\UseRawInputEncoding
\documentclass[12pt]{book}
\usepackage{amsmath, amscd, amssymb, amsthm}
\usepackage[subnum]{cases}
\usepackage{changepage}
\usepackage{marginnote}
\usepackage{hyperref}
\usepackage{cleveref}
\usepackage[all]{xy}
\usepackage{tabularx}
\usepackage{ltablex}
\usepackage[T1]{fontenc}
\usepackage{setspace}
\usepackage{fancybox}
\usepackage{ragged2e}
\usepackage[margin=1.1in]{geometry}
\hypersetup{colorlinks,linkcolor={black},citecolor={black}}
\usepackage{mathpazo}
\usepackage{enumitem}
\setlist[itemize]{leftmargin=+0.1in}

\newtheorem{theorem}{Theorem}[section]

\newtheorem{proposition}[theorem]{Proposition}

\theoremstyle{definition}
\newtheorem{definition}[theorem]{Definition}
\newtheorem{remark}[theorem]{Remark}
\numberwithin{equation}{section}
\newtheorem{example}[theorem]{Example}
\newtheorem{assumption}[theorem]{Assumption}
\newtheorem{setting}[theorem]{Setting}

\usepackage{fancyhdr}
\begin{document}

\normalfont

\title{$\infty$-Categorical Approaches to Hodge-Iwasawa Theory I: Introduction and Extensions}
\author{Xin Tong}
\date{}

\maketitle

\subsection*{Abstract}
\rm In this paper, we give the introduction to the Hodge-Iwasawa Theory introduced by the author. After that we will give some well-defined extensions to the already shaped framework established in our previous work.

\newpage

\tableofcontents

\newpage

\chapter{Preliminary Discussion}

\section{Preliminary}

\subsection{Preliminary}

This paper is a study of geometric and representation theoretic aspects of the corresponding $p$-adic motives. We have the following foundational materials.\\

\noindent  1. Noncommutative Motives: \cite{1NS},  \cite{1Ta};\\
\noindent 2. Noncommutative Harmonic Analysis, Noncommutative Microlocal Analysis and Pseudodifferential Analysis: \cite{1LVO};\\
\noindent 3. $\infty$-Categorical Topological Analysis: \cite{1BK}, \cite{1BBK}, \cite{1BBBK}, \cite{KKM}, \cite{BBM}, \cite{1CS1}, \cite{1CS2};\\
\noindent 4. Topological Division Rings and Topological Vector Spaces: \cite{1Bou};\\
\noindent 5. Topological Rings and Topological Modules: In general we need to consider very general topological rings not necessarily commutative, for instance see \cite{1AGM}, \cite{1CBCKSW}, \cite{1Hu1}, \cite{1Hu2},  \cite{1SW}, \cite{1TGI}, \cite{1U}, \cite{1W};\\
\noindent 6. Functional Analytic Rings and Functional Analytic Modules: In general we need to consider very general topological rings not necessarily commutative or archimedean, for instance see \cite{1AGM}, \cite{1BGR}, \cite{1KL1}, \cite{1KL2}, \cite{1U},   \cite{1W}. Much of the corresponding discussion in \cite{1KL1} and \cite{1KL2} works for noncommutative seminormed rings, also see \cite{1He}; \\
\noindent 7. Adic Rings:  For commutative see \cite{1CBCKSW}, \cite{1Hu1}, \cite{1Hu2},  \cite{1KL1}, \cite{1KL2}, \cite{1SW},  and for noncommutative see \cite{1FK};\\ 
\noindent 8. Distinguished Deformations of Rings:  \cite{1BMS1}, \cite{1BMS2}, \cite{1BS}, \cite{1CBCKSW}, \cite{1GR}, \cite{1KL1}, \cite{1KL2}, \cite{1Sch1}; \\
\noindent 9. Commutative Algebra: \cite{BourbakiAC}, \cite{1Lurie2}, \cite{1Lurie3}, \cite{1R}, \cite{1SP};\\
\noindent 10. Schemes: \cite{1EGAI}, \cite{1EGAII}, \cite{1EGAIII1}, \cite{1EGAIII2}, \cite{1EGAIV1},\\
 \cite{1EGAIV1}, \cite{1EGAIV1}, \cite{1EGAIV1}, \cite{1SGAI}, \cite{1SGAII}, \cite{1SGAIII1}, \cite{1SGAIII2}, \cite{1SGAIII3}, \cite{1SGAIV1}, \cite{1SGAIV2}, \cite{1SGAIV3}, \cite{1SGAIV.5}, \cite{1SGAV}, \cite{1SGAVI}, \cite{1SGAVII1},\\
  \cite{1SGAVII2}, \cite{1SP};\\
\noindent 11. Huber Spaces: \cite{1CBCKSW}, \cite{1Hu1}, \cite{1Hu2}, \cite{1KL1}, \cite{1KL2}, \cite{1Sch1}, \cite{1Sch2}, \cite{1Sch3}, \cite{1SW};\\
\noindent 12. Analytic Spaces: \cite{1BBBK}, \cite{1BK}, \cite{1BBK}, \cite{1CBCKSW}, \cite{1CS1}, \cite{1CS2}, \cite{1FS}, \cite{1Hu1}, \cite{1Hu2},  \cite{1KL1}, \cite{1KL2}, \cite{1Sch1}, \cite{1Sch2}, \cite{1Sch3}, \cite{1SW};\\
\noindent 13. Algebraic Topology: \cite{1E}, \cite{1M}, \cite{1MP}, \cite{1N};\\
\noindent 14. $\infty$-Categories and Their Models: \cite{1Bergner}, \cite{1Ci}, \cite{1J}, \cite{1Lurie1};\\
\noindent 15. Higher Algebra, Higher Toposes, Higher Geometries: \cite{1Lurie1}, \cite{1Lurie2}, \cite{1Lurie3},  \cite{1TV1}, \cite{1TV2}. Actually the corresponding $\infty$-categories of Banach ring spectra as in \cite{1BBBK} and the corresponding derived $I$-complete objects are really more relevant in the arithmetic geometry in our mind such as the corresponding objects in \cite{1BMS2}, \cite{1BS}, \cite{1BS2}, \cite{1Lurie3}, \cite{1NS}, \cite{1Po} and \cite{1Ye}\footnote{These already include many different types of $I$-adic completion in the derived sense, and the corresponding localization and completion in the algebraic topology and $\infty$-categorical theory. We would believe that the corresponding noncommutative consideration as these will be robust enough to deal with the main problems in Iwasawa theory and noncommutative derived analytic geometry.}. In fact on the $\infty$-categorical level, \cite{1BMS2} and \cite{1NS} considered the $I$-completion in the $\infty$-category theory as in \cite[Chapter 7.3]{1Lurie3}. The corresponding derived Nygaard-complete, derived Hodge-complete or other derived filtred complete objects are relevant while the corresponding derived Nygaard-incomplete, derived Hodge-incomplete or other derived filtred incomplete objects may be also relevant in certain situations.;\\
\noindent 16. EC, TMF, TAF, SAV, Sp-DG: \cite{BL1}, \cite{Lan2}, \cite{McCleary1}, \cite{DFHH1}, \cite{Lurie2}, \cite{Lurie3}, \cite{Lurie4}. This is the part for elliptic cohomology, topological modular forms, the topological automorphic forms, spectral abelian varieties, spectral $p$-divisible groups. We treat this part as a tiny degeneralization of the more general theory such as the chromatic homotopy, with restriction to derived abelian varieties and derived modular varieties. We highly recommend the reader to read Lurie's materials on elliptic cohomologies;\\
\noindent 17. Motives, Les Theories de Cohomologie d'apr\`es Weil, La Theorie de Cohomologie d'apr\`es Weil, Standard Conjectures: \cite{Motive1}-\cite{Motive33};\\
\noindent 18. Homotopty, Model Categories: \cite{Model1}, \cite{Model2}, \cite{Model3}, \cite{Model4}, \cite{Model5}, \cite{Model6}, \cite{Model7}, \cite{Model8}, \cite{Model9}, \cite{Model10}, \cite{Model11}, \cite{Model12}, \cite{Model13}, \cite{Model14}, \cite{Model15}, \cite{Model16}, \cite{Model17}, \cite{Model18}, \cite{Model19}, \cite{Model20}, \cite{Model21};\\
\noindent 19. K, THH, TCH, TP, TAQ, ArK: \cite{KTheory1}-\cite{KTheory17};\\  
\noindent 20. Motivic Homotopy, $\mathbb{A}^1$-Homotopy and $\mathbb{B}^1$-Homotopy: \cite{MotivicHomotopy1}, \cite{MotivicHomotopy2}, \cite{MotivicHomotopy3}, \cite{MotivicHomotopy4}, \cite{MotivicHomotopy5}, \cite{MotivicHomotopy6}, \cite{MotivicHomotopy7}, \cite{MotivicHomotopy8}, \cite{MotivicHomotopy9}, \cite{MotivicHomotopy10}, \cite{MotivicHomotopy11}.\\

\newpage
\subsection{Scholze's $v$-Spaces and Six Formalism}

Scholze's $v$-stacks happen over the category of perfectoid space over $\mathbb{F}_p$. We give the introduction following closely \cite{1Sch3} and \cite{1SW}. Our presentation is also following closely \cite{1Sch3}.

\begin{setting}
$v$-sheaves and $v$-stacks carry the topology which is called the $v$-topology which is finer than analytic topology, \'etale topology, pro-\'etale topology.
\end{setting}

\begin{example}
The two important situations are the following. First is the situation where the $v$-stack carries a basis of topological neighbourhoods consisting of perfectoids. The second situation is the key moduli of vector bundles over FF curves in \cite{1FS}.
\end{example}

\begin{definition}\mbox{\textbf{(Scholze \cite{1Sch3}, Analytic Prestacks)}}
Consider the following Grothendieck sites:
\begin{align}
\mathrm{Perfectoid}_{\mathbb{F}_p,\text{\'etale}},\mathrm{Perfectoid}_{\mathbb{F}_p,\text{pro\'etale}}, \mathrm{Perfectoid}_{\mathbb{F}_p,v}.
\end{align}
We just define a $(2,1)$-presheaf $\mathcal{F}$ over these sites to be any functor from 
\begin{align}
\mathrm{Perfectoid}_{\mathbb{F}_p,\text{\'etale}},\mathrm{Perfectoid}_{\mathbb{F}_p,\text{pro\'etale}}, \mathrm{Perfectoid}_{\mathbb{F}_p,v}.
\end{align}
to the groupoids. 

\end{definition}

\begin{definition}\mbox{\textbf{(Scholze \cite{1Sch3}, Analytic Stacks)}}
Consider the following Grothendieck sites:
\begin{align}
\mathrm{Perfectoid}_{\mathbb{F}_p,\text{\'etale}},\mathrm{Perfectoid}_{\mathbb{F}_p,\text{pro\'etale}}, \mathrm{Perfectoid}_{\mathbb{F}_p,v}.
\end{align}
We just define a $(2,1)$-sheaf $\mathcal{F}$ over these sites to be any functor from 
\begin{align}
\mathrm{Perfectoid}_{\mathbb{F}_p,\text{\'etale}},\mathrm{Perfectoid}_{\mathbb{F}_p,\text{pro\'etale}}, \mathrm{Perfectoid}_{\mathbb{F}_p,v}.
\end{align}
to the groupoids, which is further a stack in groupoids.

\end{definition}

\

\indent We then have the morphisms and sites of analytic stacks in certain situations.

\
\begin{definition}\mbox{\textbf{(Scholze \cite[Definition 1.20]{1Sch3}, Morphisms of Analytic Stacks)}}
The \'etale and quasi-pro-\'etale morphisms between small $v$-stacks are defined by using perfectoid coverings, and defining the corresponding \'etaleness and quasi-pro-\'etaleness after taking base changes along such perfectoid coverings.

\end{definition}

\begin{definition}\mbox{\textbf{(Scholze \cite[Definition 26.1]{1Sch3}, Sites of Analytic Stacks)}}
Let $X$ be a small $v$-stack, we have the $v$-site $X_v$\footnote{In order to study the cohomology of any small $v$-stacks, the approach taken by Scholze in \cite{1Sch3} and \cite{1FS} is not defining the particular \'etale or quasipro\'etale derived categories, instead the definition is rather not straightforward by looking at the desired subcategory of $v$-derived categories. The solidified derived category $D_\blacksquare$ in \cite{1FS} is also constructed in the same way.}. This then gives the $\infty$-category $D_{I}(X_v,\lambda)$ of derived $I$-complete objects in $\infty$-category of all $\lambda$-sheaves $D(X_v,\lambda)$ where $\lambda$ is a derived $I$-commutative algebra object.

\end{definition}

\

\begin{theorem} \mbox{\textbf{(Scholze \cite[Definition 26.1, Before the Remark 26.3]{1Sch3})}}
This $\infty$-category $D_{I}(X_v,\lambda)$ as well as the corresponding sub $\infty$-categories $D_{I,\text{qspro\'et}}(X_v,\lambda)$ and $D_{I,\text{\'et}}(X_v,\lambda)$ of the corresponding derived $I$-complete objects over the quasi-pro\'etale and \'etale sites of small $v$-stacks admit six formalism through: $\otimes, \mathrm{Hom}, f^!,f_!,f^*,f_*$ as in \cite[Definition 26.1, Before the Remark 26.3]{1Sch3}.\	
\end{theorem}

\indent Then as in \cite{1FS} for the corresponding solid $\infty$-category we can say the parallel things:

\begin{theorem} \mbox{\textbf{(Fargues-Scholze \cite[Chapter VII.2]{1FS})}}
This $\infty$-category $D_{\blacksquare}(X_v,\lambda)$ admits four formalism as in \cite[Chapter VII.2]{1FS}. Here the assumption on $\lambda$ will be definitely weaker solid condensed ring as assumed in \cite[Beginning of Chapter VII.2]{1FS}.	
\end{theorem}

\begin{remark}
In our situation, certainly what we have will be some derived $I$-completed Iwasawa modules instead of the usual Iwasawa Theoretic $I$-adic versions. This is definitely weaker, but things become robust and more well-defined. For instance if $\lambda$ is usual Iwasawa algebra, we have the following.	
\end{remark}

\begin{theorem} \mbox{\textbf{(Scholze \cite[Definition 26.1, Before the Remark 26.3]{1Sch3})}}
This $\infty$-category $D_{I}(X_v,\lambda)$ as well as the corresponding sub $\infty$-categories $D_{I,\text{qspro\'et}}(X_v,\lambda)$ and $D_{I,\text{\'et}}(X_v,\lambda)$ of the corresponding derived $I$-complete objects over the quasi-pro\'etale and \'etale sites of small $v$-stacks admit six formalism through: $\otimes, \mathrm{Hom}, f^!,f_!,f^*,f_*$ as in \cite[Definition 26.1, Before the Remark 26.3]{1Sch3}. Here we assume $\lambda$ be an Iwasawa algebra attached to some $\ell$-adic Lie group.	
\end{theorem}

\indent Then as in \cite{1FS} for the corresponding solid $\infty$-category we can say the parallel things:

\begin{theorem} \mbox{\textbf{(Fargues-Scholze \cite[Chapter VII.2]{1FS})}}
This $\infty$-category $D_{\blacksquare}(X_v,\lambda)$ admits four formalism as in \cite[Chapter VII.2]{1FS}. Here we assume $\lambda$ be an Iwasawa algebra attached to some $\ell$-adic Lie group.	
\end{theorem}

\newpage

\subsection{Introduction to the Dissertation \cite{XT}}

\subsection{The Motivation}

\

\indent Noncommutative Tamagawa number conjectures and noncommutative Iwasawa main conjectures are main topics in the field of arithmetic geometry. They give us the explicit approaches to study the motives after we take the realizations\footnote{One might want to consider motives over $\mathbb{Z}$ after some reasonable $p$-adic cohomology theories are established such as in \cite{1Ked3} and \cite{1Sch5}. Doing so might be relevant in some $p$-adic cohomological approach to Riemann Hypothesis over $\mathbb{Z}$ after  \cite{1Ked2} which is parallel to \cite{1De1}, \cite{1De2}, \cite{1Wei}.}. The associated Galois representations for instance could be used to define the corresponding $L$-functions, as well as local factors. Roughly speaking, the conjectures mention the following style  isomorphism up to certain factor\footnote{$\pi_1(*)$ is the profinite fundamental group of certain arithmetic scheme or analytic space.}:
\begin{align}
\mathrm{Determinant}_A(R\Gamma(\pi_1(*),p\mathrm{LieDef}_A(V)))\overset{}{\longrightarrow} \mathrm{Determinant}_A(0)
\end{align}
when taking some possibly trivial Iwasawa deformation along some $p$-adic Lie group tower and one could achieve the zero homotopy from the determinant of the cohomology in certain categories\footnote{As those categories in \cite{1De3}, \cite{1FK}, \cite{1MT}, \cite{1Wit3}.}, which should be closely related to special values of $L$-function and generalized functional equation.\\

\indent Our study was significantly inspired by the pictures of \cite{1BF1}, \cite{1BF2}, \cite{1FK}, \cite{1Na1}, \cite{1Wit1}, \cite{1Wit2}, \cite{1Wit3}\footnote{The commutative picture was already in \cite{1Ka1}, \cite{1Ka2}, \cite{1PR1}, \cite{1PR2}.}. On the other hand, work of Kedlaya-Liu \cite{1KL1}, \cite{1KL2} and work of Kedlaya-Pottharst \cite{1KP} also significantly inspired us to consider at least the geometrization or generalization of the work of \cite{1BF1}, \cite{1BF2}, \cite{1FK}, \cite{1Na1}, \cite{1Wit1}, \cite{1Wit2}, \cite{1Wit3}. In fact we have other important $p$-adic cohomology theories, especially in the integral setting we have the prismatic cohomology from Bhatt-Scholze \cite{1BS}. Therefore one might want to ask if we could have a well-posed prismatic Iwasawa deformation theory. To be more precise, whenever we have a well-defined $p$-adic cohomology theory, we should be able to consider some interesting Iwasawa deformation theory, possibly carrying some Banach coefficients as in \cite{1Na1}. It is actually very straightforward to first consider the integral picture, and it is very confusing to actually first understand the picture given certain Fr\'echet-Stein algebra in \cite{1ST} since the classical Iwasawa main conjecture happens over integral Iwasawa algebra $\Lambda$ as in \cite{1Iw}. \\

\newpage
\subsection{The Results of \cite{XT}}

\indent In fact our study is just as described in the consideration above, namely we actually study some non-\'etale objects which will happen over certain period sheaves over analytic spaces carrying Frobenius actions. And then we can take the corresponding Banach deformation as well as the Iwasawa deformation. Note that the previous deformation will be very crucial such as in \cite{1KPX} and \cite{1Na1}. Meanwhile we can regard them as some certain sheaves over Fargues-Fontaine stacks in some deformed way after \cite{1FF},  \cite{1KL1} and \cite{1KL2}. And we need to work with noncommutative rings as in \cite{1BF1}, \cite{1BF2}, \cite{1FK}, \cite{1Wit1}, \cite{1Wit2}, \cite{1Wit3}, \cite{1Z}. These are reflected in the following of \cite{XT}:\\

\begin{itemize}
\item 1. Hodge-Iwasawa Deformation: {corollary7.6}, {proposition2.3.10}, 
{propsition3.4.25}, {propsition3.4.26} and {proposition3.5.51}.\\ 
\item 2. Multidimensional Frobenius Modules: {theorem4.1.5} and {theorem4.1.6}.\\
\item 3. Hodge-Iwasawa Cohomology: {proposition6.3}, {proposition6.4.26}, 
{proposition6.4.27}, {proposition6.4.28}, {proposition4.42}, {proposition6.4.45}, {proposition4.45}, {proposition6.4.51}. \\
\item 4. Noncommutative Hodge-Iwasawa Deformation: {theorem3.10},  {theorem3.11}, {corollary3.12},  {theorem7.3.18}, {theorem7.3.21}, {theorem4.12}, {proposition5.13}, {theorem8.1.1}, {theorem8.1.2}, {theorem8.1.3}, {theorem8.1.4}, {theorem8.1.5}, {theorem8.1.6}, {theorem8.1.7}, {theorem9.1.1}, {theorem9.1.2}.\\
\item 5. Derived Noncommutative Hodge-Iwasawa Deformation: {theorem10.1.1}, {theorem10.1.2}, {theorem10.1.3}, {theorem10.1.4}, {theorem10.5.10}.\\
\end{itemize}

\begin{remark}
Our goal is to find certain derived Iwasawa $\infty$-categories as in \cite{1De3}, \cite{1FK}, \cite{1MT}, \cite{1Wit3} in analytic geometry over certain deformed sheaves or rings. Note that this is very complicated and particularly far beyond Iwasawa deformation of regular motives. We believe Clausen-Scholze derived category $D(\mathrm{Mod}_{\Pi,\mathrm{condensed}})$ of condensed modules over centain condensed ring $\Pi$ in \cite{1CS1} could be more robust consideration. Here in our mind $\Pi$ should be some deformed period ring such as deformation $\widetilde{\Pi}^I_{R}\otimes^{\mathbb{L}\mathrm{Solidified}}\mathcal{A}$ of the Robba ring in \cite{1KL1} and \cite{1KL2}. Also the derived category $D_\mathrm{Solidified}(X)\subset D(X_v)$ of Fargues-Scholze in \cite{1FS} carrying condensed coefficients\footnote{As well as those $D_{\text{\'et}}(X)\subset D(X_v)$ and $D_{\text{qpro\'et}}(X)\subset D(X_v)$ in \cite{1Sch3}, with application in mind to seminormal rigid analytic spaces being regarded as small $v$-stacks.} and the work \cite{1BBK} could be more robust consideration as well. What we have achieved is literally some generalization of Kedlaya-Liu abelian categories of pseudocoherent sheaves in \cite{1KL2}. 
\end{remark}

\

\indent We now discuss some examples of the spaces and local charts in our study. As mention above, they at least will be some spaces over or attached to some period rings or sheaves, which sometimes are the corresponding local charts of some stacks such as the Fargues-Fontaine spaces or the spaces before taking the quotients by equivariance coming from the motivic structures.\\


\begin{itemize}
\item 1\text{($\infty$-Categorical Rings and $\infty$-Categorical Spaces)}. The rigid analytic affinoids and spaces in \cite{1Ta} are some key examples, which is actually some initial goal in our main consideration mentioned above inspired by \cite{1BF1}, \cite{1BF2}, \cite{1FK}, \cite{1Wit1}, \cite{1Wit2}, \cite{1Wit3}, \cite{1Z};\\
\item 2\text{($\infty$-Categorical Rings and $\infty$-Categorical Spaces)}. The pseudorigid analytic affinoids and spaces over $\mathbb{Z}_p$ are also interesting to study, which is a geometric version of the arithmetic family in \cite{1Bel1} and \cite{1Bel2};\\
\item 3\text{($\infty$-Categorical Rings and $\infty$-Categorical Spaces)}. We will consider the noncommutative $I$-adically complete rings as in \cite{1FK} in the Iwasawa consideration, as well as simplicial commutative rings which are derived $I$-adically complete rings over any interesting period rings, such as the prisms in \cite{1BS}, for instance the $\mathbb{A}_\mathrm{inf}(\mathcal{O}_{\mathbb{C}_p^\flat}), W(\mathbb{F}_p)[[u]], \mathbb{Z}_p[[q-1]]$ and Robba rings in \cite{1KL1} and \cite{1KL2}. One can consider for instance the topological rings over these period rings carrying the topology induced from the period rings. Tate adic Banach rings in \cite{1KL1} and \cite{1KL2} produce certain topological adic rings satisfying the open mapping property as in \cite{1CBCKSW}.
\end{itemize}

\

\begin{remark}
Both the foundations in Bambozzi-Ben-Bassat-Kremnizer and Clausen-Scholze \cite{1BBBK}, \cite{1CS1} and \cite{1CS2} have noncommutative categories of noncommutative associative analytic rings and noncommutative associative Banach rings. Therefore this may allow one to study certain $\infty$-stacks in $\infty$-groupoids fibered over these categories under some descent consideration for instance after \cite{1KR1} and \cite{1KR2}. 
\end{remark}

\begin{itemize}

\item In chapter 2 and 3 of \cite{XT}, we study the Frobenius modules over Robba rings carrying rigid analytic coefficients and Fr\'echet-Stein coefficients, in both equal-characteristic situation and mixed-characteristic situation. We call the theory Hodge-Iwasawa since the study of the Frobenius modules over the Robba rings and sheaves are significant in Galois deformation theory, deformation of representations of fundamental groups of analytic spaces and our generalizations of the picture in \cite{1BF1}, \cite{1BF2}, \cite{1FK}, \cite{1KP}, \cite{1Na1}, \cite{1Wit1}, \cite{1Wit2}, \cite{1Wit3}. We show the equivalence between categories of finite projective or pseudocoherent $(\varphi,\Gamma)$-modules over Robba rings with rigid analytic coefficients and Fr\'echet-Stein coefficients, which also are compared to sheaves over schematic and adic Fargues-Fontaine curves in some deformed sense. Especially when we are working over analytic fields, the picture is already interesting and significant enough in the Galois representation theory and Galois deformation theory.\\

\item In chapter 4 and 5 of \cite{XT}, we study multidimensional Robba rings and multidimensional $(\varphi,\Gamma)$-modules. This point of view of taking multidimensional analogs of the Robba rings and multidimensional $(\varphi,\Gamma)$-modules is actually motivated from some programs in making progress of the local Tamagawa number conjecture of Nakamura in \cite{1Na1} literally proposed in \cite{1PZ}. Following \cite{1CKZ} and \cite{1PZ}, we define the corresponding multidimensional Robba rings and multidimensional $(\varphi,\Gamma)$-modules by taking analytic function rings over $p$-adic rigid affinoids in rigid geometry. And we define the multidimensional $(\varphi,\Gamma)$-cohomologies, multidimensional $(\psi,\Gamma)$-cohomologies and multidimensional $\psi$-cohomologies. We carefully study the complexes of the multidimensional $(\varphi,\Gamma)$-cohomologies, multidimensional $(\psi,\Gamma)$-cohomologies and multidimensional $\psi$-cohomologies, and show that they are actually living in the derived category of the bounded perfect complexes. Chapter 4 mainly focuses on imperfect Robba rings, while in chapter 5 we define perfection of Robba rings in several variables and study the comparison of multidimensional $(\varphi,\Gamma)$-modules in certain situations carefully, which is literally following \cite{1KL1}, \cite{1KL2}, \cite{1KPX} and \cite{1Ked1}. \\

\item In chapter 6 of \cite{XT}, we apply the main results in chapter 3 to study the cohomologies and categories of relative $(\varphi,\Gamma)$-modules over Robba sheaves over certain analytic spaces. We mainly discuss three applications which are crucial in our project of Iwasawa deformation of motivic structures over some higher dimensional spaces. The first is the study of abelian property of the categories of relative $(\varphi,\Gamma)$-modules over Robba sheaves over rigid analytic spaces after \cite{1KL2}, which is important whenever one would like to construct some $K$-theoretic objects to formulate Iwasawa main conjectures as in \cite{1Wit1}, \cite{1Wit2}, \cite{1Wit3}. The second is the study of families of Riemann-Hilbert correspondences after \cite{1LZ} which is crucial in further application to the arithmetic geometry along our consideration. The last one is the consideration of the equivariant version of the Iwasawa main conjecture of Nakamura \cite{1Na2}.\\

\item In chapter 7, 8 and 9 of \cite{XT}, we consider generalization of our work of the generalization of the work of Kedlaya-Liu presented in chapter 2, 3 and 6, which is also motivated by our consideration in  chapter 4 and 5 after \cite{1CKZ}, \cite{1Na1}, \cite{1PZ} and \cite{1Z} in order to make further progress. We consider the deformation in possibly noncommutative general Banach rings such as perfectoid rings, preperfectoid rings, general quotient of the noncommutative free Tate rings and so on. Certainly one thing we have to deal with is the sheafiness of the deformed rings, which will produce some difficulty to apply Kedlaya-Liu's descent in \cite{1KL1}, \cite{1KL2} for vector bundles and stably-pseudocoherent sheaves. Even in the noncommutative coefficient situation we have not worked out a theory on the noncommutative analytic toposes, which implies there is no geometric method for us to study and apply. So finding new ideas is very important. In fact, we on the representation theoretic level have the result due to Kedlaya to have the descent for vector bundles. And one can in the commutative situation use Clausen-Scholze space \cite{1CS2} to achieve the similar result by embedding Huber spaces to Clausen-Scholze spaces in \cite{1CS2} and apply Clausen-Scholze descent. Also we could consider the derived analytic spaces from Bambozzi-Kremnizer in \cite{1BK}. In the noncommutative situation we generalize results in \cite{1KL1}, \cite{1KL2} and \cite{1CBCKSW} to deform the structure sheaves directly in analytic topology, \'etale topology and pro-\'etale but not $v$-topology, which allows us to compare certain stably-pseudocoherent sheaves and modules carrying Banach deformed coefficients even if they are noncommutative, which certainly provides possibility to make further progress in the study of the noncommutative situation in chapter 4 and 5. \\

\item In chapter 10 of \cite{XT}, we initiate the project on some topics on the geometric and representation theoretic aspects of period rings. In this first paper, we consider more general base spaces. To be more precise we discuss more general perfectoid rings. Distinguished deformation of rings is a generalization notion of the Fontaine-Wintenberger idempotent correspondence. For instance in \cite{1BS}, for any quasiregular semiperfectoid ring $A$ one can canonically associate a prism $(P_A,I_{A})$. This is a very general correspondence generalized from the notions in \cite{1BMS1} and \cite{1GR}. Therefore in our consideration we will consider the adic Banach rings in \cite{1KL1}, \cite{1KL2} which are not necessarily analytic in the sense of Kedlaya's AWS lecture notes in \cite{1CBCKSW}. Assuming certain perfectoidness after \cite{1KL1} and \cite{1KL2} we study the deformed Robba rings associated. Also after \cite{1KL1} and \cite{1KL2} we studied derived deformation of the Robba rings and the descent of finite projective module spectra over them, which we will believe has some application to conjectural derived eigenvarieties and derived Galois deformation for instance in \cite{1GV}. \\

\item In chapter 11 and 12 of \cite{XT}, we consider more widely the process of taking topological and functional analytic completions, in the derived sense coming from \cite{1B1}, \cite{1BBBK}, \cite{1BK},  \cite{1BMS2}, \cite{1BS}, \cite{1CS1}. Maybe this derived consideration will let us see the hidden $\infty$-categorical structures of the representations of the Iwasawa algebras in our Iwasawa-Prismatic theory, since we are taking the deformations in some derived sense. Besides the application in mind to the prismatic or derived de Rham Iwasawa theory, chapter 11 and chapter 12 initiate also the study of some relative $p$-adic motive and Hodge theory over general derived $I$-adic spaces after \cite{1B1}, \cite{1B2}, \cite{1BMS2}, \cite{1BS}, \cite{1CBCKSW}, \cite{1DLLZ1}, \cite{1DLLZ2}, \cite{1Dr1}, \cite{1GL}, \cite{1G1}, \cite{1Hu2}, \cite{1III1}, \cite{1III2}, \cite{1KL1}, \cite{1KL2}, \cite{1NS}, \cite{1O}, \cite{1Ro}, \cite{1Sch2} such as the pseudorigid analytic spaces and more general spaces carrying some derived $I$-adic topology, as well as the prelog simplicial commutative rings in \cite{1B1} carrying some derived $I$-adic topology as well.

\end{itemize}

\newpage
\subsection{The Picture and The Future Consideration}

Let us try to discuss a little bit how one should think about our main motivation. We start from the following two settings of possible geometrizations of Iwasawa theory.\\

\begin{itemize}
\item 1. \textbf{$\infty$-Categorical Iwasawa-Prismatic Theory}\footnote{This should be slightly more relevant in the geometrization of integral Iwasawa theory, but one can invert $p$ as well.}: After Bhatt-Lurie, Bhatt-Scholze and Drinfeld \cite{1BS}, \cite{1Dr}, \cite{1Sch4};\\
\item 2. \textbf{$\infty$-Categorical Hodge-Iwasawa Theory}\footnote{This should be slightly more relevant in the geometrization of rational Iwasawa theory, but one can not invert $p$ as well.}: After Kedlaya-Liu and Kedlaya-Pottharst \cite{1KL1}, \cite{1KL2}, \cite{1KP}.
\end{itemize}

\begin{itemize}

\item \textbf{($\infty$-Categorical Iwasawa-Prismatic Theory)}
Just like in \cite{1KP}, \cite{1Wit1}, \cite{1Wit2}, \cite{1Wit3}, and by using prismatic cohomology theory in \cite{1BS},   \cite{1Sch4} we can take a quasisyntomic formal ring $R$\footnote{As in \cite{1Sch4}, one can take a $p$-adic fornal scheme which is quasisyntomic. Certainly this is already a type of spaces which is interesting enough including a point situation and many situations where motivic comparisons could happen as in \cite{1BMS2}, \cite{1BS}.} over $\mathbb{Z}_p$ and have two sites of $\mathrm{Spf}R$. The first site is the corresponding prismatic site $(X_{\mathrm{prim}},\mathcal{O}_{X_{\mathrm{prim}}})$. The second site is the quasisyntomic site $(X_\mathrm{qsyn},\mathcal{O}_{X_{\mathrm{qsyn}}})$. Recall for the second one, for any quasiregular semiperfectoid affinoid $A$ in $X_\mathrm{qsyn}$ we have that one can canonically associate a prism $(P_A,I_A)$ to $A$ where we have $\mathcal{O}_{X_{\mathrm{qsyn}}}(A):=P_A$. Now by Bhatt-Scholze \cite{1Sch4} we have more sheaves from the structure sheaves here namely we have:
\begin{align}
\mathcal{O}_{X_{\mathrm{prim}}}[1/I_{\mathcal{O}_{X_{\mathrm{prim}}}}]^\wedge_p, \mathcal{O}_{X_{\mathrm{prim}}}[1/I_{\mathcal{O}_{X_{\mathrm{prim}}}}]^\wedge_p[1/p], \mathcal{O}_{X_{\mathrm{qsyn}}}[1/I_{\mathcal{O}_{X_{\mathrm{qsyn}}}}]^\wedge_p, \mathcal{O}_{X_{\mathrm{qsyn}}}[1/I_{\mathcal{O}_{X_{\mathrm{qsyn}}}}]^\wedge_p[1/p]. 
\end{align}
Carrying some integral Iwasawa algebra $\mathbb{Z}_p[[G]]$ for some compact $p$-adic Lie group, and after taking the derived completion\footnote{At this moment we are assuming that the derived completion is possible in our current setting. One may also consider the solidification of Clausen-Scholze. We want to mention that in the noncommutative setting there are many ways to do the completion in the derived sense, which is already subtle in the commutative setting.} we have:
\begin{align}
\mathcal{O}_{X_{\mathrm{prim}}}[1/I_{\mathcal{O}_{X_{\mathrm{prim}}}}]^\wedge_p{\otimes}^{\mathbb{L}\mathrm{com}}\mathbb{Z}_p[[G]],\\
 \mathcal{O}_{X_{\mathrm{prim}}}[1/I_{\mathcal{O}_{X_{\mathrm{prim}}}}]^\wedge_p{\otimes}^{\mathbb{L}\mathrm{com}}\mathbb{Z}_p[[G]][1/p],\\
  \mathcal{O}_{X_{\mathrm{qsyn}}}[1/I_{\mathcal{O}_{X_{\mathrm{qsyn}}}}]^\wedge_p{\otimes}^{\mathbb{L}\mathrm{com}}\mathbb{Z}_p[[G]], \\
  \mathcal{O}_{X_{\mathrm{qsyn}}}[1/I_{\mathcal{O}_{X_{\mathrm{qsyn}}}}]^\wedge_p{\otimes}^{\mathbb{L}\mathrm{com}}\mathbb{Z}_p[[G]][1/p]. 
\end{align}
Then one might want to ask if one can use such style deformation to establish the parallel story in \cite{1BF1}, \cite{1BF2}, \cite{1FK}, \cite{1KP}, \cite{1Wit1}, \cite{1Wit2}, \cite{1Wit3}. For instance taking the $J\subset \mathbb{Z}_p[[G]]$-adic quotient we have the Koszul complexes parametrized by such $J$:
\begin{align}
\mathrm{Kos}_J\mathcal{O}_{X_{\mathrm{prim}}}[1/I_{\mathcal{O}_{X_{\mathrm{prim}}}}]^\wedge_p{\otimes}\mathbb{Z}_p[[G]],\\
\mathrm{Kos}_J  \mathcal{O}_{X_{\mathrm{qsyn}}}[1/I_{\mathcal{O}_{X_{\mathrm{qsyn}}}}]^\wedge_p{\otimes}\mathbb{Z}_p[[G]]. \\ 
\end{align}
Then one may define the $\infty$-category of pro-systems of the Iwasawa complexes over these $\mathbb{E}_1$-rings, and consider the associated Waldhausen categories as in \cite{1Wit1}, \cite{1Wit2}, \cite{1Wit3}.

\item \textbf{($\infty$-Categorical Iwasawa-Prismatic Theory)}
Within the same framework we take $R$ to be $\mathcal{O}_K$ for some $p$-adic local field $K$. Bhatt-Scholze \cite{1BS}, \cite{1Sch4} showed that we have the category of Galois representations of $\mathbb{Z}_p$-coefficients of $\mathrm{Gal}_K$ is equivalent to the category of prismatic $F$-crystals over 
\begin{center}
$(X_\mathrm{prim},\mathcal{O}_{X_{\mathrm{prim}}}[1/I_{\mathcal{O}_{X_{\mathrm{prim}}}}]^\wedge_p)$, 
\end{center}
while the category of Galois representations of $\mathbb{Q}_p$-coefficients of $\mathrm{Gal}_K$ is equivalent to the category of prismatic $F$-crystals over $(X_\mathrm{prim},\mathcal{O}_{X_{\mathrm{prim}}}[1/I_{\mathcal{O}_{X_{\mathrm{prim}}}}]^\wedge_p[1/p])$. Then one could ask if we could consider some Iwasawa deformation through the some $p$-adic Lie quotient of $\mathrm{Gal}_K$ to establish the parallel story in \cite{1KP}, \cite{1Wit1}, \cite{1Wit2}, \cite{1Wit3}. Namely for any such $F$-crystal $M$ with associated representation $V$ over $\mathbb{Z}_p$ or $\mathbb{Q}_p$ we take the Iwasawa deformation $p\mathrm{DfLie}(M)$ by some Iwasawa $F$-crystal\footnote{How one should define this crystal will be determined by how one forms the completed tensor product in the Iwasawa deformation.} through the quotient from $G_K$ to some compact $p$-adic Lie group $G$, then we can ask if the following:
\begin{align}
R\Gamma(X_\mathrm{prim},p\mathrm{DfLie}(M)), R\Gamma(X_\mathrm{qsyn},p\mathrm{DfLie}(M))
\end{align}
recover the classical Iwasawa theory by using the Galois cohomology 
\begin{center}
$R\Gamma(\mathrm{Gal}_K,p\mathrm{DfLie}(V))$
\end{center}
of $p\mathrm{DfLie}(V)$, as well as the \'etale cohomology $R\Gamma(\mathrm{Spec}K,p\mathrm{DfLie}(\widetilde{V}))$ of the local system $p\mathrm{DfLie}(\widetilde{V})$ attached to $p\mathrm{DfLie}(V)$.
\end{itemize}

\

\indent Beyond the somehow \'etale situations in the above picture, one could consider the corresponding category of prismatic crystals, which will be beyond the Galois representation theoretic consideration. Also one could regard these objects as certain sheaves over the prismatic stacks in \cite{1Dr}.\\ 

\indent In our study, we have the following picture. We will consider picture beyond \'etale situation, and we will study the Frobenius sheaves and Frobenius modules in the very general situation. And we will have the chance to regard the sheaves and modules with Frobenius actions as certain sheaves over Fargues-Fontaine stacks after \cite{1FF}, \cite{1KL1} and \cite{1KL2}. Actually we conjecture that the quasisyntomic descent and étale comparison results of Bhatt-Scholze \cite{1BS} will imply equivalence in some accurate sense beyond the vague similarity.\\

\chapter{Hodge-Iwasawa Theory}

\section{Motivation I}

\subsection{\text{Motivation I: Dememorization and Memorization}}
\begin{itemize}

\item<1-> Consider the cyclotomic tower $\{\mathbb{Q}_p(\zeta_{p^n})\}_n$ of $\mathbb{Q}_p$.

\item<2-> The infinite level of this tower is kind of special after the corresponding completion.

\item<3-> Over $\mathbb{Q}_p$, we could consider $\mathrm{Spa}(\mathbb{Q}_p,\mathfrak{o}_{\mathbb{Q}_p})_{\text{pro\'et}}$ due to Scholze \cite{Sch}, although the infinite level of the towers above participates in this topology but the corresponding pro-\'etale site forgets the corresponding cyclotomic tower while it is defined by using pro-systems of \'etale morphisms.

\item<4-> Work of Pottharst \cite{P1}, Kedlaya-Pottharst-Xiao \cite{KPX}, Kedlaya-Pottharst \cite{KP} implies one may see the corresponding cyclotomic tower back by considering the corresponding cyclotomic deformation as below.

\item<5-> One has the so-called $\psi$-cohomology originally dated back to Fontaine (see \cite[II.1.3]{CC}) attached to a $(\varphi,\Gamma)$-module $M$ (you could regarded this as a Galois representation):
\begin{align}
H_\psi(M)	
\end{align}
by using the operator $\psi$.

\item<6-> And we have the corresponding $(\varphi,\Gamma)$-module after Herr, but we consider the cyclotomic deformation as in Kedlaya-Pottharst-Xiao over the Robba ring $\mathcal{R}^\infty_{\mathbb{Q}_p}(\Gamma)$:
\begin{align}
H_{\varphi,\Gamma}(\mathbf{CycDef}(M)).	
\end{align}

\end{itemize}

\noindent{\text{Motivation I: Dememorization and Memorization}}
\begin{itemize}
\item<1-> This is defined by taking the corresponding external tensor product of $M$ with the corresponding module coming from the quotient $\Gamma$. This dates back to Pottharst on his analytic Iwasawa cohomology \cite{P1}.

\item<2-> Work of Kedlaya-Pottharst \cite{KP} observes that we can have the following sheaf version of the construction:
\begin{align}
H_\text{pro-\'etale}(\mathrm{Spa}(\mathbb{Q}_p,\mathfrak{o}_{\mathbb{Q}_p}),\mathbf{CycDef}(\widetilde{M})),	
\end{align}
which is defined by taking the corresponding external product of Kedlaya-Liu's sheaf $\widetilde{M}$ \cite{KL1} with the one defined by using the quotient $\Gamma$. 
\item<3-> The point is that we have the following comparison:
\begin{align}
H_{\psi}(M) \overset{\sim}{\rightarrow}	 H_{\varphi,\Gamma}(\mathbf{CycDef}(M))\overset{\sim}{\rightarrow}H_\text{pro-\'etale}(\mathrm{Spa}(\mathbb{Q}_p,\mathfrak{o}_{\mathbb{Q}_p}),\mathbf{CycDef}(\widetilde{M})).
\end{align}

\item<4-> Suppose $M(V)$ comes from a Galois representation $V$ of $G_{\mathbb{Q}_p}$ we even have the following comparison after Perrin-Riou \cite{PR}:
\begin{align}
&{\mathcal{O}_{\mathrm{Sp}\mathcal{R}^\infty_{\mathbb{Q}_p}(\Gamma)}}\widehat{\otimes}_\Lambda H_\mathrm{IW}(G_{\mathbb{Q}_p},V)\overset{\sim}{\rightarrow} \\
&\overset{\sim}{\rightarrow} H_{\psi}((M(V)) \overset{\sim}{\rightarrow}H_{\varphi,\Gamma}(\mathbf{CycDef}((M(V)))\overset{\sim}{\rightarrow}H_\text{pro-\'etale}(\mathrm{Spa}(\mathbb{Q}_p,\mathfrak{o}_{\mathbb{Q}_p}),\mathbf{CycDef}(\widetilde{(M(V)})).
\end{align}
\item<5-> Natural questions come:\\
I. How about the Lubin-Tate Iwasawa theory in Berger-Fourquaux-Schneider-Venjakob's work, observed by Kedlaya-Pottharst \cite{KP}, \cite{BF}, \cite{SV}.\\
II. How about higher dimensional toric towers and more general towers of rigid analytic spaces for instance.

\end{itemize}

\noindent{\text{Motivation I: Dememorization and Memorization}}
\begin{itemize}
\item<1-> These need us to generalize the corresponding framework to higher dimensional situation and more general deformed version. The problem is challenging, since we have some rigidized objects combined together. 
\item<2-> Rational coefficients are very complicated comparing to algebraic geometry, since sometimes we do not have the integral lattices over the \'etale sites. This is already a problem in the context of Kedlaya-Liu \cite{KL1}, \cite{KL2}. 
\item<3->It is not surprising much for us to consider generalizing the frame work of non-\'etale objects since even in the usual situations over a point work of Nakamura \cite{Nakamura1}, Kedlaya-Pottharst-Xiao \cite{KPX} and Kedlaya-Liu \cite{KL1}, \cite{KL2} implies that all kinds of families of Galois representations will be more conveniently studied by using $B$-pairs and $(\varphi,\Gamma)$-modules. 
\item<4-> If we only have some abelian group $G$ the corresponding deformation happens along the algebra $\mathbb{Q}_p[G]$ which gives rise to Galois representation of $\mathrm{Gal}(\overline{\mathbb{Q}}_p/\mathbb{Q}_p)$ with coefficient in $\mathbb{Q}_p[G]$ along the quotient $\mathrm{Gal}(\overline{\mathbb{Q}}_p/\mathbb{Q}_p)\rightarrow G$. One can then regard this as a sheave of module over some sheaf with deformed coefficient in $\mathbb{Q}_p[G]$. 
\item<5-> Note that we can also consider some deformation over an affinoid algebra in the rigid analytic geometry, which amounts to the $p$-adic families of special values. This is not available at once in archimedean functional analysis.	
\end{itemize}

\newpage


\section{\text{Motivation II: Higher Dimensional modeling of the
 Weil Conjectures}}
\begin{itemize}
\item<1-> The corresponding equivariant consideration could be obviously generalized to the relative $p$-adic Hodge theory which is aimed at the study of the \'etale local systems over rigid analytic spaces.
\item<2-> This amounts to higher dimensional modeling of the generalized Weil conjecture after the work due to many people, to name a few Deligne \cite{De1}, \cite{De2} ($\ell$-adic \'etale sheaves), Kedlaya \cite{Ked1} ($p$-adic differential modules), Abe and Caro \cite{AC} ($p$-adic arithmetic $D$-modules) and so on.

\item<3-> The invariance comes from the quotient of \'etale fundamental groups of rigid analytic spaces, or the corresponding profinite fundamental groups of rigid analytic spaces.

\item<4-> (\text{Example}) One can consider the corresponding Fr\'echet-Stein algebras associated to the group $\mathbb{Z}_p\ltimes \mathbb{Z}_p^n$ which is Galois group (quotient of the corresponding profinite fundamental group) of a local chart of smooth proper rigid analytic spaces. Note that the top of this local chart in the smooth proper setting naturally participates in some nice topology.

\item<5-> (\text{Example}) One can consider the local systems in more general sense, for instance the locally constant sheaf $\underline{A}$	
attached to a topological ring, for instance an affinoid algebra in the rigid analytic geometry after Tate. This is somewhat special in the $p$-adic setting due to the fact the corresponding Hodge structures could achieve variation in $p$-adic rigid family.

\end{itemize}

\newpage

\section{What can be learnt from noncommutative Iwasawa Theory}

\subsection{\text{Integral \'Etale Noncommutative Iwasawa Theory}}

\begin{itemize}
\item<1-> Following some idea in the noncommutative Tamagawa Number conjecture after Fukaya-Kato \cite{FK} and the noncommutative Iwasawa theory over a scheme over finite field after Witte we would like to consider the following picture after Witte \cite{Wit1}.\\

\item<2-> Let $T$ be an adic ring in the sense of Fukaya-Kato \cite{FK}, which is a compact ring with two sided ideal $I$ such that we have each $T/I^n$ is finite for $n\geq 0$ and taking the inverse limit we recover the ring $T$ itself. This ring could be noncommutative, for instance the Iwasawa algebra attached to some $p$-adic Lie group.\\
\item<3-> (\text{Definition, after Witte \cite[Definition 5.4.1]{Wit1}}) Consider a rigid analytic space or a scheme $X/\mathbb{Q}_p$ separated and of finite type we consider the category $\mathbb{D}_\mathrm{perf}(X_\sharp,T)$ ($\sharp=\text{\'et},\text{pro\'et}$) which is the category of the inverse limit of perfect complexes of abelian sheaves of left modules over quotients of $T$ by open two-sided ideals of $T$ which are $DG$-flat, parametrized by open two-sided ideals of $T$.\\

\item<4-> (\text{Theorem, Witte \cite[Proposition 6.1.5]{Wit1}}) Let $p$ be a unit in $T$. The category defined above could be endowed with the structure of Waldhausen category\footnote{Strictly Speaking, these are the complicial biWaldhausen ones.} and the total direct image functor induces a well defined functor in the situation where $X$ is a scheme and $\sharp=\text{\'et}$:
\[
\xymatrix@R+0pc@C+2pc{
\mathbb{D}_\mathrm{perf}(X_\sharp,T)\ar[r]^{R\Gamma(X_\sharp,.)}\ar[r]\ar[r] &\mathbb{D}_\mathrm{perf}(T) 
}
\]	
which induces the corresponding map on the $K$-theory space:
\[
\xymatrix@R+0pc@C+2pc{
\mathbb{K}\mathbb{D}_\mathrm{perf}(X_\sharp,T)\ar[r]^{\mathbb{K}R\Gamma(X_\sharp,.)}\ar[r]\ar[r] &\mathbb{K}\mathbb{D}_\mathrm{perf}(T). 
}
\]
Then this map is homotopic to zero in some canonical way.

\end{itemize}

\newpage

\section{Introduction to the Interactions among Motives}

\subsection{\text{Equivariant relative $p$-adic Hodge Theory}}

\begin{itemize}
\item<1->Things discussed so far have motivated the corresponding equivariant relative $p$-adic Hodge Theory in the following sense. Witte \cite{Wit1} considered general framework of Grothendieck abelian categories, for instance one can consider the following categories:\\

\item<1-> 1. The category of all the abelian sheaves over the \'etale or pro-\'etale sites of schemes of finite type over a field $k$ after Grothendieck, Scholze, Bhatt \cite{SGA4}, \cite{BS1} and etc;\\
\item<2-> 2. The category of all the abelian sheaves over the \'etale or pro-\'etale sites of adic spaces of finite type over a field $k$ after Huber, Scholze, Kedlaya-Liu \cite{Hu}, \cite{Sch}, \cite{KL1}, \cite{KL2};\\	
\item<3-> 3. The ind-category of the abelian category of the pseudocoherent Frobenius $\varphi$-sheaves over a rigid analytic space over a complete discrete valued field with perfect residue field $k$ after Kedlaya-Liu \cite{KL1}, \cite{KL2}.\\

\item<4-> 4. The category of abelian sheaves over the syntonic site  by covering of quasiregular semiperfect algebras, as in the work of Bhatt-Morrow-Scholze \cite{BMS}.  \\

\item<5-> One can naturally consider the corresponding $P$-objects throughout the categories listed above, where $P$ is noetherian for instance. For instance one can consider the third category and consider the corresponding local systems over $\underline{A}$ where $A$ is an affinoid algebra in rigid analytic geometry after Tate \cite{Ta1}, which are the $A$-objects in the corresponding category of all the abelian sheaves.

\end{itemize}

\noindent{\text{Equivariant relative $p$-adic Hodge Theory}}

\begin{itemize}

\item<1-> The corresponding $P$-objects are interesting, but in general are not that easy to study, especially we consider for instance those ring defined over $\mathbb{Q}_p$, let it alone if one would like to consider the categories of the complexes of such objects. 

\item<2-> We choose to consider the corresponding embedding of such objects into the categories of Frobenius sheaves with coefficients in $P$ after Kedlaya-Liu \cite{KL1}, \cite{KL2}. Again we expect everything will be more convenient to handle in the category of $(\varphi,\Gamma)$-modules.
	
\item<3-> Working over $R$ now a uniform Banach algebra with further structure of an adic ring over $\mathbb{F}_p$. And we assume that $R$ is perfect.
Let $\text{Robba}^\text{extended}_{I,R}$ be the Robba sheaves defined by Kedlaya-Liu \cite{KL1}, \cite{KL2}, with respect to some interval $I\subset (0,\infty)$, which are Fr\'echet completions of the ring of Witt vector of $R$ with respect to the Gauss norms induced from the norm on $R$. 

\item<4-> Taking suitable interval one can define the corresponding Robba rings $\text{Robba}^\text{extended}_{r,R}$, $\text{Robba}^\text{extended}_{\infty,R}$ and the corresponding full Robba ring $\text{Robba}^\text{extended}_{R}$.

\item<5-> We work in the category of Banach and ind-Fr\'echet spaces, which are commutative. Our generalization comes from those Banach reduced affinoid algebras $A$.

\end{itemize}

\noindent{\text{Equivariant relative $p$-adic Hodge Theory}}
\begin{itemize}

\item<1-> The $p$-adic functional analysis produces us some manageable structures within our study of relative $p$-adic Hodge theory, generalizing the original $p$-adic functional analytic framework of Kedlaya-Liu \cite{KL1}, \cite{KL2}.

\item<2-> Starting from Kedlaya-Liu's period rings, 
\begin{align}
&\text{Robba}^\text{extended}_{\infty,R},\text{Robba}^\text{extended}_{I,R},\text{Robba}^\text{extended}_{r,R},\text{Robba}^\text{extended}_{R},	\text{Robba}^\text{extended}_{{\mathrm{int},r},R},\\
&\text{Robba}^\text{extended}_{{\mathrm{int}},R},\text{Robba}^\text{extended}_{{\mathrm{bd},r},R},\text{Robba}^\text{extended}_{{\mathrm{bd}},R}
\end{align}
we can form the corresponding $A$-relative of the period rings:
\begin{align}
&\text{Robba}^\text{extended}_{\infty,R,A},\text{Robba}^\text{extended}_{I,R,A},\text{Robba}^\text{extended}_{r,R,A},\text{Robba}^\text{extended}_{R,A},	\text{Robba}^\text{extended}_{{\mathrm{int},r},R,A},\\
&\text{Robba}^\text{extended}_{\mathrm{int},R,A},\text{Robba}^\text{extended}_{{\mathrm{bd},r},R,A},\text{Robba}^\text{extended}_{{\mathrm{bd}},R,A}.	
\end{align} 
\item<3-> (\text{Remark}) There should be also many interesting contexts, for instance consider a finitely generated abelian group $G$, one can consider the group rings: 
\begin{align}
\text{Robba}^\text{extended}_{I,R}[G].	
\end{align}
\item<4-> And then consider the completion living inside the corresponding infinite direct sum Banach modules 
\begin{align}
\bigoplus\text{Robba}^\text{extended}_{I,R},	
\end{align}
over the corresponding period rings:
\begin{align}
\overline{\text{Robba}^\text{extended}_{I,R}[G]}.	
\end{align}
Then we take suitable intersection and union one can have possibly some interesting period rings $\overline{\text{Robba}^\text{extended}_{r,R}[G]}$ and $\overline{\text{Robba}^\text{extended}_{R}[G]}$.
\end{itemize}

\noindent{\text{Equivariant relative $p$-adic Hodge Theory}}
\begin{itemize}
\item<1-> The equivariant period rings in the situations we mentioned above carry relative Frobenius action $\varphi$ induced from the Witt vectors. 

\item<2-> They carry the corresponding Banach or (ind-)Fr\'echet spaces structures. So we can generalize the corresponding Kedlaya-Liu's construction to the following situations (here let $G$ be finite):

\item<3-> We can then consider the corresponding completed Frobenius modules over the rings in the equivariant setting. To be more precise over:
\begin{align}
\overline{\text{Robba}^\text{extended}_{R}[G]},\Omega_{\mathrm{int},R,A},\Omega_{R,A},\text{Robba}^\text{extended}_{R,A},\text{Robba}^\text{extended}_{\mathrm{bd},R,A}	
\end{align}
one considers the Frobenius modules finite locally free.

\item<4->  With the corresponding finite locally free models over
\begin{align}
\overline{\text{Robba}^\text{extended}_{r,R}[G]},\text{Robba}^\text{extended}_{r,R,A},\text{Robba}^\text{extended}_{{\mathrm{bd},r},R,A},	
\end{align} 
again carrying the corresponding semilinear Frobenius structures, where $r$ could be $\infty$.

\item<5-> One also consider families of Frobenius modules over 
\begin{align}
\overline{\text{Robba}^\text{extended}_{I,R}[G]},\text{Robba}^\text{extended}_{I,R,A},	
\end{align} 
in glueing fashion with obvious cocycle condition with respect to three intervals $I\subset J\subset K$. These are called the corresponding Frobenius bundles.
	
\end{itemize}

\newpage

\section{The Key Deformation}

\subsection{\text{Deformation of Schemes}}	

\begin{itemize}
\item<1-> One can consider the corresponding schemes attached to the above commutative rings, for instance
\begin{align}
\mathrm{Spec}\text{Robba}^\text{extended}_{r,R,A},\mathrm{Spec}\overline{\text{Robba}^\text{extended}_{r,R}[G]}.	
\end{align}
And consider the corresponding categories:
\begin{align}
\mathrm{Mod}(\mathcal{O}_{\mathrm{Spec}\text{Robba}^\text{extended}_{r,R,A}}),\mathrm{Mod}(\mathcal{O}_{\mathrm{Spec}\overline{\text{Robba}^\text{extended}_{r,R}[G]}}).	
\end{align}
 
\item<2-> These are very straightforward and even crucial especially when we consider 
\begin{align}
\varphi-\mathrm{Mod}(\mathcal{O}_{\mathrm{Spec}\text{Robba}^\text{extended}_{\infty,R,A}}),\varphi-\mathrm{Mod}(\mathcal{O}_{\mathrm{Spec}\overline{\text{Robba}^\text{extended}_{\infty,R}[G]}}),	
\end{align}
in some Frobenius equivariant way.

\item<3-> But on the other hand it is also very convenient to encode the Frobenius action inside the spaces themselves, which leads to Fargues-Fontaine Schemes as those in the work of Kedlaya-Liu \cite{KL1}, \cite{KL2}, \cite{FF}. 

\end{itemize}

\noindent{\text{Deformation of Schemes}}
\begin{itemize}
\item<1-> Roughly one takes the corresponding $\varphi=p^n$ equivariant elements in the full Robba ring, and putting them to be a commutative graded ring $\bigoplus P_{R,A,n}$, and then glueing them through the Proj construction by glueing subschemes taking the form of $\mathrm{Spec}P_{R,A}[1/f]_0$.
\item<2-> Roughly one takes the corresponding $\varphi=p^n$ equivariant elements in the full Robba ring $\overline{\text{Robba}^\text{extended}_{R}[G]}$, and putting them to be a commutative graded ring $\bigoplus P_{R,G,n}$, and then glueing them through the Proj construction by glueing subschemes taking the form of $\mathrm{Spec}P_{R,G}[1/f]_0$.

\item<3-> Therefore we have the natural functor:
\[
\xymatrix@R+0pc@C+0pc{
\mathrm{Mod}\mathcal{O}_{\mathrm{Proj}_{R,A}}\ar[r]\ar[r]\ar[r] &\mathrm{Mod}\mathcal{O}_{\mathrm{Spec}\text{Robba}^\text{extended}_{\infty,R,A}} , 
}
\]
defined by using the corresponding pullbacks.

\item<4-> \text{(Theorem, Tong \cite[Theorem 1.3]{T})} We have the following categories are equivalent (generalizing the work of \text{Kedlaya-Liu} \cite{KL1}, \text{Kedlaya-Pottharst} \cite{KP}): \\
I. The category of all the quasicoherent finite locally free sheaves over $\mathrm{Proj}\bigoplus P_{R,A,n}$;\\
II. The category of all the Frobenius modules of the global sections of all the $\varphi$-equivariant quasicoherent finite locally free sheaves over $\mathrm{Spec}\text{Robba}^\text{extended}_{\infty,R,A}$;\\
III.  The category of all the Frobenius modules over $\text{Robba}^\text{extended}_{R,A}$;\\
IV. The category of all the Frobenius bundles over $\text{Robba}^\text{extended}_{R,A}$.

\item<5-> For the rings for general $G$, we expect one should also be able to establish some results parallel to this once the structures are more literally investigated.	We are also interested in the noncommutative coefficients as in Z\"ahringer's thesis \cite{Z}, but we need to use noncommutative topos.
\end{itemize}

\noindent{\text{Deformation of Schemes}}
\begin{itemize}
\item<1-> \text{(Theorem, Tong \cite[Proposition 3.16, Corollary 3.17]{T})} We have the following categories are equivalent (generalizing the work of \text{Kedlaya-Liu} \cite{KL1}, \text{Kedlaya-Pottharst} \cite{KP}): \\
I. The category of pro-systems of all the quasicoherent finite locally free sheaves over $\mathrm{Proj}\bigoplus P_{R,A_\infty,n}$;\\
II. The category of pro-systems of all the Frobenius modules coming from the global sections of all the $\varphi$-equivariant quasicoherent finite locally free sheaves over $\mathrm{Spec}\text{Robba}^\text{extended}_{\infty,R,A_\infty}$;\\
III.  The category of pro-systems of all the Frobenius modules over $\text{Robba}^\text{extended}_{R,A_\infty}$;\\
IV. The category of pro-systems of all the Frobenius bundles over $\text{Robba}^\text{extended}_{R,A_\infty}$.	\\
Here $A_\infty$ is a Fr\'echet-Stein algebra attached to a compact $p$-adic Lie group such that the algebra is limit of (commutative) reduced affinoid algebras. And the finiteness is put on the infinite level of ind-scheme, actually one can also just put on each level. 

\end{itemize}

\noindent{\text{Deformation of Schemes}}

\begin{itemize}

\item<1-> \text{(Outline)} Following Kedlaya-Liu \cite{KL1}:\\
1. Construct the glueing process over the scheme $\mathrm{Spec}\text{Robba}^\text{extended}_{\infty,R,A}$;\\
2. The functors could be read off from the corresponding diagram above, namely one glues the resulting sheaves over each $\mathrm{Spec}\text{Robba}^\text{extended}_{\infty,R,A}[1/f]$ for each suitable element $f$ in the graded ring, then takes the corresponding global section;\\
3. Then from the last category back to the quasicoherent sheaves over the Fargues-Fontaine scheme we need to solve some Frobenius algebraic equation by $p$-adic analytic method to show that taking Frobenius invariance over each affine subspace is exact, where one uses Kedlaya-Liu's approach which could be dated back to Kedlaya's approach to slope filtration over extended Robba rings \cite{Ked2}.

\item<2-> Let us look back the functor:
\[
\xymatrix@R+0pc@C+0pc{
\mathrm{Mod}\mathcal{O}_{\mathrm{Proj}P_{R,A}}\ar[r]\ar[r]\ar[r] &\varphi-\mathrm{Mod}\mathcal{O}_{\mathrm{Spec}\text{Robba}^\text{extended}_{\infty,R,A}} \ar[r]\ar[r]\ar[r] &\varphi-\mathrm{Mod}\mathcal{O}_{\mathrm{Spec}\text{Robba}^\text{extended}_{R,A}} , 
}
\]
obviously one might want to generalize the picture above, which was also considered by Kedlaya-Liu in their original work \cite{KL1}, \cite{KL2}.
\item<3-> \text{(Theorem, Tong \cite[Theorem 1.4]{T})} We have the following categories are equivalent (generalizing the work of \text{Kedlaya-Liu} \cite{KL1}, \cite{KL2}, \text{Kedlaya-Pottharst} \cite{KP}): \\
I. The category of all the pseudocoherent sheaves over $\mathrm{Proj}\bigoplus P_{R,A,n}$;\\
II.  The category of all the pseudocoherent $\varphi$-equivariant modules over $\text{Robba}^\text{extended}_{R,A}$.
\end{itemize}

\newpage

\section{$K$-Theoretic Consideration}

\subsection{\text{The $K$-theory of Algebraic Relative Hodge-Iwasawa Modules}}
\begin{itemize}
\item<1-> Based on the study we did above, it should be very natural to consider more general pseudocoherent complexes in some higher categorical sense. Note that pseudocoherent objects were naturally emerging in SGA \cite{SGAVI} from some K-theoretic point of view. Also more importantly Hodge-Iwasawa theory to some extent will behave better if we forget the derived category, when we would like to study the K-theoretic aspects.

\item<2-> (\text{Definition}) Let $Ch\mathrm{Mod}\mathcal{O}_{\mathrm{Proj}P_{R}}$ denote the category of all the complexes of objects in $\mathrm{Mod}\mathcal{O}_{\mathrm{Proj}_{R}}$. 
\item<3-> (\text{Definition}) We now use the notations:
\begin{align}
D_{\mathrm{perf}}\mathrm{Proj}P_{R},D_{\mathrm{pseudo}}\mathrm{Proj}P_{R}
\end{align}
to denote the category of all the perfect and pseudocoherent complexes. 

\item<4-> (\text{Definition}) One also has the following subcategories:
\begin{align}
&D^{\mathrm{dg-flat}}_{\mathrm{perf}}\mathrm{Proj}P_{R},\\
&D^{\mathrm{str}}_{\mathrm{perf}}\mathrm{Proj}P_{R}.
\end{align}

\item<5-> (\text{Proposition, after Thomason-Trobaugh \cite{TT}}) These categories admit Waldhausen structure.

\item<6-> (\text{Question}) In the situation where $R=\widetilde{R}_\psi$ attached to the cyclotomic tower,  we would like to know if $D_{\mathrm{perf}}\mathrm{Proj}P_{R}$ and $D^{\mathrm{str}}_{\mathrm{perf}}\mathrm{Proj}P_{R}$ admit Waldhausen exact functors to $D_{\mathrm{perf}}(\mathbb{Q}_p)$ or $D^{\mathrm{str}}_{\mathrm{perf}}(\mathbb{Q}_p)$, which induce maps on the associated K-theory spaces.

\end{itemize}

\newpage

\section{Analytic $\infty$-Categorical Functional Analytic Hodge-Iwasawa Modules}

\subsection{$\infty$-Categorical Analytic Stacks and Descents I}

\noindent We now make the corresponding discussion after our previous work \cite{T2} on the homotopical functional analysis after many projects \cite{1BBBK}, \cite{1BBK}, \cite{BBM}, \cite{1BK} , \cite{1CS1}, \cite{1CS2}, \cite{KKM}. We choose to work over the Bambozzi-Kremnizer space \cite{1BK} attached to the corresponding Banach rings in our work after \cite{1BBBK}, \cite{1BBK}, \cite{BBM}, \cite{1BK}, \cite{KKM}. Note that what is happening is that attached to any Banach ring over $\mathbb{Q}_p$, say $B$, we attach a $(\infty,1)-$stack $\mathcal{X}(B)$ fibered over (in the sense of $\infty$-groupoid, and up to taking the corresponding opposite categories) after \cite{1BBBK}, \cite{1BBK}, \cite{BBM}, \cite{1BK}, \cite{KKM}:
\begin{align}
\mathrm{sComm}\mathrm{Simp}\mathrm{Ind}\mathrm{Ban}_{\mathbb{Q}_p},	
\end{align}
with 
\begin{align}
\mathrm{sComm}\mathrm{Simp}\mathrm{Ind}^m\mathrm{Ban}_{\mathbb{Q}_p}.	
\end{align}
associated with a $(\infty,1)$-ring object $\mathcal{O}_{\mathcal{X}(B)}$, such that we have the corresponding under the basic derived rational localization $\infty$-Grothendieck site
\begin{center}
 $(\mathcal{X}(B), \mathcal{O}_{\mathcal{X}(B),\mathrm{drl}})$ 
\end{center}
carrying the homotopical epimorphisms as the corresponding topology.

\begin{itemize}
\item<1-> By using this framework (certainly one can also consider \cite{1CS1} and \cite{1CS2} as the foundations, as in \cite{LBV}), we have the $\infty$-stack after Kedlaya-Liu \cite{KL1}, \cite{KL2}. Here in the following let $A$ be any Banach ring over $\mathbb{Q}_p$.

\item<2-> Generalizing Kedlaya-Liu's construction in \cite{KL1}, \cite{KL2} of the adic Fargues-Fontaine space we have a quotient (by using powers of the Frobenius operator) $X_{R,A}$ of the space 
\begin{align}
Y_{R,A}:=\bigcup_{0<s<r}\mathcal{X}(\text{Robba}^\text{extended}_{R,[s,r],A}).	
\end{align}

\item<3-> This is a locally ringed space $(X_{R,A},\mathcal{O}_{X_{R,A}})$, so one can consider the stable $\infty$-category $\mathrm{Ind}\mathrm{Banach}(\mathcal{O}_{X_{R,A}}) $ which is the $\infty$-category of all the $\mathcal{O}_{X_{R,A}}$-sheaves of inductive Banach modules over $X_{R,A}$. We have the parallel categories for $Y_{R,A}$, namely $\varphi\mathrm{Ind}\mathrm{Banach}(\mathcal{O}_{X_{R,A}})$ and so on. Here we will consider presheaves.
 
\item<4-> This is a locally ringed space $(X_{R,A},\mathcal{O}_{X_{R,A}})$, so one can consider the stable $\infty$-category $\mathrm{Ind}^m\mathrm{Banach}(\mathcal{O}_{X_{R,A}}) $ which is the $\infty$-category of all the $\mathcal{O}_{X_{R,A}}$-sheaves of inductive monomorphic Banach modules over $X_{R,A}$. We have the parallel categories for $Y_{R,A}$, namely $\varphi\mathrm{Ind}^m\mathrm{Banach}(\mathcal{O}_{X_{R,A}})$ and so on. Here we will consider presheaves.
   
\item<5-> In this context one can consider the $K$-theory as in the scheme situation by using the ideas and constructions from Blumberg-Gepner-Tabuada \cite{BGT}. Moreover we can study the Hodge Theory.

\item<6-> We expect that one can study among these big categories to find interesting relationships, since this should give us the right understanding of the $p$-adic Hodge theory. The corresponding pseudocoherent version comparison could be expected to be deduced as in Kedlaya-Liu's work \cite{KL1}, \cite{KL2}.

\begin{assumption}\label{assumtionpresheaves}
All the functors of modules or algebras below are presheaves.	
\end{assumption}

\item (\text{Proposition}) There is an equivalence between the $\infty$-categories of inductive Banach quasicoherent presheaves:
\[
\xymatrix@R+0pc@C+0pc{
\mathrm{Ind}\mathrm{Banach}(\mathcal{O}_{X_{R,A}})\ar[r]^{\mathrm{equi}}\ar[r]\ar[r] &\varphi\mathrm{Ind}\mathrm{Banach}(\mathcal{O}_{Y_{R,A}}).  
}
\]
\item (\text{Proposition}) There is an equivalence between the $\infty$-categories of monomorphic inductive Banach quasicoherent presheaves:
\[
\xymatrix@R+0pc@C+0pc{
\mathrm{Ind}^m\mathrm{Banach}(\mathcal{O}_{X_{R,A}})\ar[r]^{\mathrm{equi}}\ar[r]\ar[r] &\varphi\mathrm{Ind}^m\mathrm{Banach}(\mathcal{O}_{Y_{R,A}}).  
}
\]
\end{itemize}

\begin{itemize}

\item (\text{Proposition}) There is an equivalence between the $\infty$-categories of inductive Banach quasicoherent presheaves:
\[
\xymatrix@R+0pc@C+0pc{
\mathrm{Ind}\mathrm{Banach}(\mathcal{O}_{X_{R,A}})\ar[r]^{\mathrm{equi}}\ar[r]\ar[r] &\varphi\mathrm{Ind}\mathrm{Banach}(\mathcal{O}_{Y_{R,A}}).  
}
\]
\item (\text{Proposition}) There is an equivalence between the $\infty$-categories of monomorphic inductive Banach quasicoherent presheaves:
\[
\xymatrix@R+0pc@C+0pc{
\mathrm{Ind}^m\mathrm{Banach}(\mathcal{O}_{X_{R,A}})\ar[r]^{\mathrm{equi}}\ar[r]\ar[r] &\varphi\mathrm{Ind}^m\mathrm{Banach}(\mathcal{O}_{Y_{R,A}}).  
}
\]
\item (\text{Proposition}) There is an equivalence between the $\infty$-categories of inductive Banach quasicoherent commutative algebra $E_\infty$ objects:
\[
\xymatrix@R+0pc@C+0pc{
\mathrm{sComm}_\mathrm{simplicial}\mathrm{Ind}\mathrm{Banach}(\mathcal{O}_{X_{R,A}})\ar[r]^{\mathrm{equi}}\ar[r]\ar[r] &\mathrm{sComm}_\mathrm{simplicial}\varphi\mathrm{Ind}\mathrm{Banach}(\mathcal{O}_{Y_{R,A}}).  
}
\]
\item (\text{Proposition}) There is an equivalence between the $\infty$-categories of monomorphic inductive Banach quasicoherent commutative algebra $E_\infty$ objects:
\[
\xymatrix@R+0pc@C+0pc{
\mathrm{sComm}_\mathrm{simplicial}\mathrm{Ind}^m\mathrm{Banach}(\mathcal{O}_{X_{R,A}})\ar[r]^{\mathrm{equi}}\ar[r]\ar[r] &\mathrm{sComm}_\mathrm{simplicial}\varphi\mathrm{Ind}^m\mathrm{Banach}(\mathcal{O}_{Y_{R,A}}).  
}
\]

\item Then parallel as in \cite{LBV} we have the equivalence of the de Rham complex after \cite[Definition 5.9, Section 5.2.1]{KKM}:
\[
\xymatrix@R+0pc@C+0pc{
\mathrm{deRham}_{\mathrm{sComm}_\mathrm{simplicial}\mathrm{Ind}\mathrm{Banach}(\mathcal{O}_{X_{R,A}})\ar[r]^{\mathrm{equi}}}(-)\ar[r]\ar[r] &\mathrm{deRham}_{\mathrm{sComm}_\mathrm{simplicial}\varphi\mathrm{Ind}\mathrm{Banach}(\mathcal{O}_{Y_{R,A}})}(-), 
}
\]
\[
\xymatrix@R+0pc@C+0pc{
\mathrm{deRham}_{\mathrm{sComm}_\mathrm{simplicial}\mathrm{Ind}^m\mathrm{Banach}(\mathcal{O}_{X_{R,A}})\ar[r]^{\mathrm{equi}}}(-)\ar[r]\ar[r] &\mathrm{deRham}_{\mathrm{sComm}_\mathrm{simplicial}\varphi\mathrm{Ind}^m\mathrm{Banach}(\mathcal{O}_{Y_{R,A}})}(-). 
}
\]

\item Then we have the following equivalence of $K$-group $(\infty,1)$-spectrum from \cite{BGT}:
\[
\xymatrix@R+0pc@C+0pc{
\mathrm{K}^\mathrm{BGT}_{\mathrm{sComm}_\mathrm{simplicial}\mathrm{Ind}\mathrm{Banach}(\mathcal{O}_{X_{R,A}})\ar[r]^{\mathrm{equi}}}(-)\ar[r]\ar[r] &\mathrm{K}^\mathrm{BGT}_{\mathrm{sComm}_\mathrm{simplicial}\varphi\mathrm{Ind}\mathrm{Banach}(\mathcal{O}_{Y_{R,A}})}(-), 
}
\]
\[
\xymatrix@R+0pc@C+0pc{
\mathrm{K}^\mathrm{BGT}_{\mathrm{sComm}_\mathrm{simplicial}\mathrm{Ind}^m\mathrm{Banach}(\mathcal{O}_{X_{R,A}})\ar[r]^{\mathrm{equi}}}(-)\ar[r]\ar[r] &\mathrm{K}^\mathrm{BGT}_{\mathrm{sComm}_\mathrm{simplicial}\varphi\mathrm{Ind}^m\mathrm{Banach}(\mathcal{O}_{Y_{R,A}})}(-). 
}
\]
\end{itemize}

\noindent Now let $R=\mathbb{Q}_p(p^{1/p^\infty})^{\wedge\flat}$ and $R_k=\mathbb{Q}_p(p^{1/p^\infty})^{\wedge}\left<T_1^{\pm 1/p^{\infty}},...,T_k^{\pm 1/p^{\infty}}\right>^\flat$ we have the following Galois theoretic results with naturality along $f:\mathrm{Spa}(\mathbb{Q}_p(p^{1/p^\infty})^{\wedge}\left<T_1^{\pm 1/p^\infty},...,T_k^{\pm 1/p^\infty}\right>^\flat)\rightarrow \mathrm{Spa}(\mathbb{Q}_p(p^{1/p^\infty})^{\wedge\flat})$:

\begin{itemize}
\item (\text{Proposition}) There is an equivalence between the $\infty$-categories of inductive Banach quasicoherent presheaves:
\[
\xymatrix@R+6pc@C+0pc{
\mathrm{Ind}\mathrm{Banach}(\mathcal{O}_{X_{\mathbb{Q}_p(p^{1/p^\infty})^{\wedge}\left<T_1^{\pm 1/p^\infty},...,T_k^{\pm 1/p^\infty}\right>^\flat,A}})\ar[d]\ar[d]\ar[d]\ar[d] \ar[r]^{\mathrm{equi}}\ar[r]\ar[r] &\varphi\mathrm{Ind}\mathrm{Banach}(\mathcal{O}_{Y_{\mathbb{Q}_p(p^{1/p^\infty})^{\wedge}\left<T_1^{\pm 1/p^\infty},...,T_k^{\pm 1/p^\infty}\right>^\flat,A}}) \ar[d]\ar[d]\ar[d]\ar[d].\\
\mathrm{Ind}\mathrm{Banach}(\mathcal{O}_{X_{\mathbb{Q}_p(p^{1/p^\infty})^{\wedge\flat},A}})\ar[r]^{\mathrm{equi}}\ar[r]\ar[r] &\varphi\mathrm{Ind}\mathrm{Banach}(\mathcal{O}_{Y_{\mathbb{Q}_p(p^{1/p^\infty})^{\wedge\flat},A}}).\\ 
}
\]
\item (\text{Proposition}) There is an equivalence between the $\infty$-categories of monomorphic inductive Banach quasicoherent presheaves:
\[
\xymatrix@R+6pc@C+0pc{
\mathrm{Ind}^m\mathrm{Banach}(\mathcal{O}_{X_{R_k,A}})\ar[r]^{\mathrm{equi}}\ar[d]\ar[d]\ar[d]\ar[d]\ar[r]\ar[r] &\varphi\mathrm{Ind}^m\mathrm{Banach}(\mathcal{O}_{Y_{R_k,A}})\ar[d]\ar[d]\ar[d]\ar[d]\\
\mathrm{Ind}^m\mathrm{Banach}(\mathcal{O}_{X_{\mathbb{Q}_p(p^{1/p^\infty})^{\wedge\flat},A}})\ar[r]^{\mathrm{equi}}\ar[r]\ar[r] &\varphi\mathrm{Ind}^m\mathrm{Banach}(\mathcal{O}_{Y_{\mathbb{Q}_p(p^{1/p^\infty})^{\wedge\flat},A}}).\\  
}
\]
\item (\text{Proposition}) There is an equivalence between the $\infty$-categories of inductive Banach quasicoherent commutative algebra $E_\infty$ objects:
\[\displayindent=-0.1in
\xymatrix@R+6pc@C+0pc{
\mathrm{sComm}_\mathrm{simplicial}\mathrm{Ind}\mathrm{Banach}(\mathcal{O}_{X_{R_k,A}})\ar[d]\ar[d]\ar[d]\ar[d]\ar[r]^{\mathrm{equi}}\ar[r]\ar[r] &\mathrm{sComm}_\mathrm{simplicial}\varphi\mathrm{Ind}\mathrm{Banach}(\mathcal{O}_{Y_{R_k,A}})\ar[d]\ar[d]\ar[d]\ar[d]\\
\mathrm{sComm}_\mathrm{simplicial}\mathrm{Ind}\mathrm{Banach}(\mathcal{O}_{X_{\mathbb{Q}_p(p^{1/p^\infty})^{\wedge\flat},A}})\ar[r]^{\mathrm{equi}}\ar[r]\ar[r] &\mathrm{sComm}_\mathrm{simplicial}\varphi\mathrm{Ind}\mathrm{Banach}(\mathcal{O}_{Y_{\mathbb{Q}_p(p^{1/p^\infty})^{\wedge\flat},A}})  
}
\]
\item (\text{Proposition}) There is an equivalence between the $\infty$-categories of monomorphic inductive Banach quasicoherent commutative algebra $E_\infty$ objects:
\[\displayindent=-0.1in
\xymatrix@R+6pc@C+0pc{
\mathrm{sComm}_\mathrm{simplicial}\mathrm{Ind}^m\mathrm{Banach}(\mathcal{O}_{X_{R_k,A}})\ar[d]\ar[d]\ar[d]\ar[d]\ar[r]^{\mathrm{equi}}\ar[r]\ar[r] &\mathrm{sComm}_\mathrm{simplicial}\varphi\mathrm{Ind}^m\mathrm{Banach}(\mathcal{O}_{Y_{R_k,A}})\ar[d]\ar[d]\ar[d]\ar[d]\\
 \mathrm{sComm}_\mathrm{simplicial}\mathrm{Ind}^m\mathrm{Banach}(\mathcal{O}_{X_{\mathbb{Q}_p(p^{1/p^\infty})^{\wedge\flat},A}})\ar[r]^{\mathrm{equi}}\ar[r]\ar[r] &\mathrm{sComm}_\mathrm{simplicial}\varphi\mathrm{Ind}^m\mathrm{Banach}(\mathcal{O}_{Y_{\mathbb{Q}_p(p^{1/p^\infty})^{\wedge\flat},A}}) 
}
\]

\item Then parallel as in \cite{LBV} we have the equivalence of the de Rham complex after \cite[Definition 5.9, Section 5.2.1]{KKM}:
\[\displayindent=-0.2in
\xymatrix@R+6pc@C+0pc{
\mathrm{deRham}_{\mathrm{sComm}_\mathrm{simplicial}\mathrm{Ind}\mathrm{Banach}(\mathcal{O}_{X_{R_k,A}})\ar[r]^{\mathrm{equi}}}(-)\ar[d]\ar[d]\ar[d]\ar[d]\ar[r]\ar[r] &\mathrm{deRham}_{\mathrm{sComm}_\mathrm{simplicial}\varphi\mathrm{Ind}\mathrm{Banach}(\mathcal{O}_{Y_{R_k,A}})}(-)\ar[d]\ar[d]\ar[d]\ar[d]\\
\mathrm{deRham}_{\mathrm{sComm}_\mathrm{simplicial}\mathrm{Ind}\mathrm{Banach}(\mathcal{O}_{X_{\mathbb{Q}_p(p^{1/p^\infty})^{\wedge\flat},A}})\ar[r]^{\mathrm{equi}}}(-)\ar[r]\ar[r] &\mathrm{deRham}_{\mathrm{sComm}_\mathrm{simplicial}\varphi\mathrm{Ind}\mathrm{Banach}(\mathcal{O}_{Y_{\mathbb{Q}_p(p^{1/p^\infty})^{\wedge\flat},A}})}(-) 
}
\]
\[\displayindent=-0.2in
\xymatrix@R+6pc@C+0pc{
\mathrm{deRham}_{\mathrm{sComm}_\mathrm{simplicial}\mathrm{Ind}^m\mathrm{Banach}(\mathcal{O}_{X_{R_k,A}})\ar[r]^{\mathrm{equi}}}(-)\ar[d]\ar[d]\ar[d]\ar[d]\ar[r]\ar[r] &\mathrm{deRham}_{\mathrm{sComm}_\mathrm{simplicial}\varphi\mathrm{Ind}^m\mathrm{Banach}(\mathcal{O}_{Y_{R_k,A}})}(-)\ar[d]\ar[d]\ar[d]\ar[d]\\
\mathrm{deRham}_{\mathrm{sComm}_\mathrm{simplicial}\mathrm{Ind}^m\mathrm{Banach}(\mathcal{O}_{X_{\mathbb{Q}_p(p^{1/p^\infty})^{\wedge\flat},A}})\ar[r]^{\mathrm{equi}}}(-)\ar[r]\ar[r] &\mathrm{deRham}_{\mathrm{sComm}_\mathrm{simplicial}\varphi\mathrm{Ind}^m\mathrm{Banach}(\mathcal{O}_{Y_{\mathbb{Q}_p(p^{1/p^\infty})^{\wedge\flat},A}})}(-) 
}
\]

\item Then we have the following equivalence of $K$-group $(\infty,1)$-spectrum from \cite{BGT}:
\[
\xymatrix@R+6pc@C+0pc{
\mathrm{K}^\mathrm{BGT}_{\mathrm{sComm}_\mathrm{simplicial}\mathrm{Ind}\mathrm{Banach}(\mathcal{O}_{X_{R_k,A}})\ar[r]^{\mathrm{equi}}}(-)\ar[d]\ar[d]\ar[d]\ar[d]\ar[r]\ar[r] &\mathrm{K}^\mathrm{BGT}_{\mathrm{sComm}_\mathrm{simplicial}\varphi\mathrm{Ind}\mathrm{Banach}(\mathcal{O}_{Y_{R_k,A}})}(-)\ar[d]\ar[d]\ar[d]\ar[d]\\
\mathrm{K}^\mathrm{BGT}_{\mathrm{sComm}_\mathrm{simplicial}\mathrm{Ind}\mathrm{Banach}(\mathcal{O}_{X_{\mathbb{Q}_p(p^{1/p^\infty})^{\wedge\flat},A}})\ar[r]^{\mathrm{equi}}}(-)\ar[r]\ar[r] &\mathrm{K}^\mathrm{BGT}_{\mathrm{sComm}_\mathrm{simplicial}\varphi\mathrm{Ind}\mathrm{Banach}(\mathcal{O}_{Y_{\mathbb{Q}_p(p^{1/p^\infty})^{\wedge\flat},A}})}(-) 
}
\]
\[
\xymatrix@R+6pc@C+0pc{
\mathrm{K}^\mathrm{BGT}_{\mathrm{sComm}_\mathrm{simplicial}\mathrm{Ind}^m\mathrm{Banach}(\mathcal{O}_{X_{R_k,A}})\ar[r]^{\mathrm{equi}}}(-)\ar[d]\ar[d]\ar[d]\ar[d]\ar[r]\ar[r] &\mathrm{K}^\mathrm{BGT}_{\mathrm{sComm}_\mathrm{simplicial}\varphi\mathrm{Ind}^m\mathrm{Banach}(\mathcal{O}_{Y_{R_k,A}})}(-)\ar[d]\ar[d]\ar[d]\ar[d]\\
\mathrm{K}^\mathrm{BGT}_{\mathrm{sComm}_\mathrm{simplicial}\mathrm{Ind}^m\mathrm{Banach}(\mathcal{O}_{X_{\mathbb{Q}_p(p^{1/p^\infty})^{\wedge\flat},A}})\ar[r]^{\mathrm{equi}}}(-)\ar[r]\ar[r] &\mathrm{K}^\mathrm{BGT}_{\mathrm{sComm}_\mathrm{simplicial}\varphi\mathrm{Ind}^m\mathrm{Banach}(\mathcal{O}_{Y_{\mathbb{Q}_p(p^{1/p^\infty})^{\wedge\flat},A}})}(-) 
}
\]

\end{itemize}

\
\indent Then we consider further equivariance by considering the arithmetic profinite fundamental groups $\Gamma_{\mathbb{Q}_p}$ and $\mathrm{Gal}(\overline{\mathbb{Q}_p\left<T_1^{\pm 1},...,T_k^{\pm 1}\right>}/R_k)$ through the following diagram:

\[
\xymatrix@R+0pc@C+0pc{
\mathbb{Z}_p^k=\mathrm{Gal}(R_k/{\mathbb{Q}_p(p^{1/p^\infty})^\wedge\left<T_1^{\pm 1},...,T_k^{\pm 1}\right>}) \ar[r]\ar[r] \ar[r]\ar[r] &\Gamma_k:=\mathrm{Gal}(R_k/{\mathbb{Q}_p\left<T_1^{\pm 1},...,T_k^{\pm 1}\right>}) \ar[r] \ar[r]\ar[r] &\Gamma_{\mathbb{Q}_p}.
}
\]

\begin{itemize}
\item (\text{Proposition}) There is an equivalence between the $\infty$-categories of inductive Banach quasicoherent presheaves:
\[
\xymatrix@R+6pc@C+0pc{
\mathrm{Ind}\mathrm{Banach}_{\Gamma_k}(\mathcal{O}_{X_{\mathbb{Q}_p(p^{1/p^\infty})^{\wedge}\left<T_1^{\pm 1/p^\infty},...,T_k^{\pm 1/p^\infty}\right>^\flat,A}})\ar[d]\ar[d]\ar[d]\ar[d] \ar[r]^{\mathrm{equi}}\ar[r]\ar[r] &\varphi\mathrm{Ind}\mathrm{Banach}_{\Gamma_k}(\mathcal{O}_{Y_{\mathbb{Q}_p(p^{1/p^\infty})^{\wedge}\left<T_1^{\pm 1/p^\infty},...,T_k^{\pm 1/p^\infty}\right>^\flat,A}}) \ar[d]\ar[d]\ar[d]\ar[d].\\
\mathrm{Ind}\mathrm{Banach}(\mathcal{O}_{X_{\mathbb{Q}_p(p^{1/p^\infty})^{\wedge\flat},A}})\ar[r]^{\mathrm{equi}}\ar[r]\ar[r] &\varphi\mathrm{Ind}\mathrm{Banach}(\mathcal{O}_{Y_{\mathbb{Q}_p(p^{1/p^\infty})^{\wedge\flat},A}}).\\ 
}
\]
\item (\text{Proposition}) There is an equivalence between the $\infty$-categories of monomorphic inductive Banach quasicoherent presheaves:
\[
\xymatrix@R+6pc@C+0pc{
\mathrm{Ind}^m\mathrm{Banach}_{\Gamma_k}(\mathcal{O}_{X_{R_k,A}})\ar[r]^{\mathrm{equi}}\ar[d]\ar[d]\ar[d]\ar[d]\ar[r]\ar[r] &\varphi\mathrm{Ind}^m\mathrm{Banach}_{\Gamma_k}(\mathcal{O}_{Y_{R_k,A}})\ar[d]\ar[d]\ar[d]\ar[d]\\
\mathrm{Ind}^m\mathrm{Banach}_{\Gamma_0}(\mathcal{O}_{X_{\mathbb{Q}_p(p^{1/p^\infty})^{\wedge\flat},A}})\ar[r]^{\mathrm{equi}}\ar[r]\ar[r] &\varphi\mathrm{Ind}^m\mathrm{Banach}_{\Gamma_0}(\mathcal{O}_{Y_{\mathbb{Q}_p(p^{1/p^\infty})^{\wedge\flat},A}}).\\  
}
\]
\item (\text{Proposition}) There is an equivalence between the $\infty$-categories of inductive Banach quasicoherent commutative algebra $E_\infty$ objects:
\[
\xymatrix@R+6pc@C+0pc{
\mathrm{sComm}_\mathrm{simplicial}\mathrm{Ind}\mathrm{Banach}_{\Gamma_k}(\mathcal{O}_{X_{R_k,A}})\ar[d]\ar[d]\ar[d]\ar[d]\ar[r]^{\mathrm{equi}}\ar[r]\ar[r] &\mathrm{sComm}_\mathrm{simplicial}\varphi\mathrm{Ind}\mathrm{Banach}_{\Gamma_k}(\mathcal{O}_{Y_{R_k,A}})\ar[d]\ar[d]\ar[d]\ar[d]\\
\mathrm{sComm}_\mathrm{simplicial}\mathrm{Ind}\mathrm{Banach}_{\Gamma_0}(\mathcal{O}_{X_{\mathbb{Q}_p(p^{1/p^\infty})^{\wedge\flat},A}})\ar[r]^{\mathrm{equi}}\ar[r]\ar[r] &\mathrm{sComm}_\mathrm{simplicial}\varphi\mathrm{Ind}\mathrm{Banach}_{\Gamma_0}(\mathcal{O}_{Y_{\mathbb{Q}_p(p^{1/p^\infty})^{\wedge\flat},A}})  
}
\]
\item (\text{Proposition}) There is an equivalence between the $\infty$-categories of monomorphic inductive Banach quasicoherent commutative algebra $E_\infty$ objects:
\[
\xymatrix@R+6pc@C+0pc{
\mathrm{sComm}_\mathrm{simplicial}\mathrm{Ind}^m\mathrm{Banach}_{\Gamma_k}(\mathcal{O}_{X_{R_k,A}})\ar[d]\ar[d]\ar[d]\ar[d]\ar[r]^{\mathrm{equi}}\ar[r]\ar[r] &\mathrm{sComm}_\mathrm{simplicial}\varphi\mathrm{Ind}^m\mathrm{Banach}_{\Gamma_k}(\mathcal{O}_{Y_{R_k,A}})\ar[d]\ar[d]\ar[d]\ar[d]\\
 \mathrm{sComm}_\mathrm{simplicial}\mathrm{Ind}^m\mathrm{Banach}_{\Gamma_0}(\mathcal{O}_{X_{\mathbb{Q}_p(p^{1/p^\infty})^{\wedge\flat},A}})\ar[r]^{\mathrm{equi}}\ar[r]\ar[r] &\mathrm{sComm}_\mathrm{simplicial}\varphi\mathrm{Ind}^m\mathrm{Banach}_{\Gamma_0}(\mathcal{O}_{Y_{\mathbb{Q}_p(p^{1/p^\infty})^{\wedge\flat},A}}) 
}
\]

\item Then parallel as in \cite{LBV} we have the equivalence of the de Rham complex after \cite[Definition 5.9, Section 5.2.1]{KKM}:
\[\displayindent=-0.2in
\xymatrix@R+6pc@C+0pc{
\mathrm{deRham}_{\mathrm{sComm}_\mathrm{simplicial}\mathrm{Ind}\mathrm{Banach}_{\Gamma_k}(\mathcal{O}_{X_{R_k,A}})\ar[r]^{\mathrm{equi}}}(-)\ar[d]\ar[d]\ar[d]\ar[d]\ar[r]\ar[r] &\mathrm{deRham}_{\mathrm{sComm}_\mathrm{simplicial}\varphi\mathrm{Ind}\mathrm{Banach}_{\Gamma_k}(\mathcal{O}_{Y_{R_k,A}})}(-)\ar[d]\ar[d]\ar[d]\ar[d]\\
\mathrm{deRham}_{\mathrm{sComm}_\mathrm{simplicial}\mathrm{Ind}\mathrm{Banach}_{\Gamma_0}(\mathcal{O}_{X_{\mathbb{Q}_p(p^{1/p^\infty})^{\wedge\flat},A}})\ar[r]^{\mathrm{equi}}}(-)\ar[r]\ar[r] &\mathrm{deRham}_{\mathrm{sComm}_\mathrm{simplicial}\varphi\mathrm{Ind}\mathrm{Banach}_{\Gamma_0}(\mathcal{O}_{Y_{\mathbb{Q}_p(p^{1/p^\infty})^{\wedge\flat},A}})}(-) 
}
\]
\[\displayindent=-0.4in
\xymatrix@R+6pc@C+0pc{
\mathrm{deRham}_{\mathrm{sComm}_\mathrm{simplicial}\mathrm{Ind}^m\mathrm{Banach}_{\Gamma_k}(\mathcal{O}_{X_{R_k,A}})\ar[r]^{\mathrm{equi}}}(-)\ar[d]\ar[d]\ar[d]\ar[d]\ar[r]\ar[r] &\mathrm{deRham}_{\mathrm{sComm}_\mathrm{simplicial}\varphi\mathrm{Ind}^m\mathrm{Banach}_{\Gamma_k}(\mathcal{O}_{Y_{R_k,A}})}(-)\ar[d]\ar[d]\ar[d]\ar[d]\\
\mathrm{deRham}_{\mathrm{sComm}_\mathrm{simplicial}\mathrm{Ind}^m\mathrm{Banach}_{\Gamma_0}(\mathcal{O}_{X_{\mathbb{Q}_p(p^{1/p^\infty})^{\wedge\flat},A}})\ar[r]^{\mathrm{equi}}}(-)\ar[r]\ar[r] &\mathrm{deRham}_{\mathrm{sComm}_\mathrm{simplicial}\varphi\mathrm{Ind}^m\mathrm{Banach}_{\Gamma_0}(\mathcal{O}_{Y_{\mathbb{Q}_p(p^{1/p^\infty})^{\wedge\flat},A}})}(-) 
}
\]

\item Then we have the following equivalence of $K$-group $(\infty,1)$-spectrum from \cite{BGT}:
\[
\xymatrix@R+6pc@C+0pc{
\mathrm{K}^\mathrm{BGT}_{\mathrm{sComm}_\mathrm{simplicial}\mathrm{Ind}\mathrm{Banach}_{\Gamma_k}(\mathcal{O}_{X_{R_k,A}})\ar[r]^{\mathrm{equi}}}(-)\ar[d]\ar[d]\ar[d]\ar[d]\ar[r]\ar[r] &\mathrm{K}^\mathrm{BGT}_{\mathrm{sComm}_\mathrm{simplicial}\varphi\mathrm{Ind}\mathrm{Banach}_{\Gamma_k}(\mathcal{O}_{Y_{R_k,A}})}(-)\ar[d]\ar[d]\ar[d]\ar[d]\\
\mathrm{K}^\mathrm{BGT}_{\mathrm{sComm}_\mathrm{simplicial}\mathrm{Ind}\mathrm{Banach}_{\Gamma_0}(\mathcal{O}_{X_{\mathbb{Q}_p(p^{1/p^\infty})^{\wedge\flat},A}})\ar[r]^{\mathrm{equi}}}(-)\ar[r]\ar[r] &\mathrm{K}^\mathrm{BGT}_{\mathrm{sComm}_\mathrm{simplicial}\varphi\mathrm{Ind}\mathrm{Banach}_{\Gamma_0}(\mathcal{O}_{Y_{\mathbb{Q}_p(p^{1/p^\infty})^{\wedge\flat},A}})}(-) 
}
\]
\[
\xymatrix@R+6pc@C+0pc{
\mathrm{K}^\mathrm{BGT}_{\mathrm{sComm}_\mathrm{simplicial}\mathrm{Ind}^m\mathrm{Banach}_{\Gamma_k}(\mathcal{O}_{X_{R_k,A}})\ar[r]^{\mathrm{equi}}}(-)\ar[d]\ar[d]\ar[d]\ar[d]\ar[r]\ar[r] &\mathrm{K}^\mathrm{BGT}_{\mathrm{sComm}_\mathrm{simplicial}\varphi\mathrm{Ind}^m\mathrm{Banach}_{\Gamma_k}(\mathcal{O}_{Y_{R_k,A}})}(-)\ar[d]\ar[d]\ar[d]\ar[d]\\
\mathrm{K}^\mathrm{BGT}_{\mathrm{sComm}_\mathrm{simplicial}\mathrm{Ind}^m\mathrm{Banach}_{\Gamma_0}(\mathcal{O}_{X_{\mathbb{Q}_p(p^{1/p^\infty})^{\wedge\flat},A}})\ar[r]^{\mathrm{equi}}}(-)\ar[r]\ar[r] &\mathrm{K}^\mathrm{BGT}_{\mathrm{sComm}_\mathrm{simplicial}\varphi\mathrm{Ind}^m\mathrm{Banach}_{\Gamma_0}(\mathcal{O}_{Y_{\mathbb{Q}_p(p^{1/p^\infty})^{\wedge\flat},A}})}(-). 
}
\]

\end{itemize}

\

Furthermore we have the corresponding pro-\'etale version without the corresponding fundamental group equivariances.

\begin{itemize}
\item (\text{Proposition}) There is an equivalence between the $\infty$-categories of inductive Banach quasicoherent presheaves:
\[
\xymatrix@R+6pc@C+0pc{
\mathrm{Ind}\mathrm{Banach}(\mathcal{O}_{X_{\mathbb{Q}_p\left<T_1^{\pm 1},...,T_k^{\pm 1}\right>,\text{pro\'etale},A}})\ar[d]\ar[d]\ar[d]\ar[d] \ar[r]^{\mathrm{equi}}\ar[r]\ar[r] &\varphi\mathrm{Ind}\mathrm{Banach}(\mathcal{O}_{Y_{\mathbb{Q}_p\left<T_1^{\pm 1},...,T_k^{\pm 1}\right>,\text{pro\'etale},A}}) \ar[d]\ar[d]\ar[d]\ar[d].\\
\mathrm{Ind}\mathrm{Banach}(\mathcal{O}_{X_{\mathbb{Q}_p,\text{pro\'etale},A}})\ar[r]^{\mathrm{equi}}\ar[r]\ar[r] &\varphi\mathrm{Ind}\mathrm{Banach}(\mathcal{O}_{Y_{\mathbb{Q}_p,\text{pro\'etale},A}}).\\ 
}
\]
\item (\text{Proposition}) There is an equivalence between the $\infty$-categories of monomorphic inductive Banach quasicoherent presheaves:
\[
\xymatrix@R+6pc@C+0pc{
\mathrm{Ind}^m\mathrm{Banach}(\mathcal{O}_{X_{\mathbb{Q}_p\left<T_1^{\pm 1},...,T_k^{\pm 1}\right>,\text{pro\'etale},A}})\ar[r]^{\mathrm{equi}}\ar[d]\ar[d]\ar[d]\ar[d]\ar[r]\ar[r] &\varphi\mathrm{Ind}^m\mathrm{Banach}(\mathcal{O}_{Y_{\mathbb{Q}_p\left<T_1^{\pm 1},...,T_k^{\pm 1}\right>,\text{pro\'etale},A}})\ar[d]\ar[d]\ar[d]\ar[d]\\
\mathrm{Ind}^m\mathrm{Banach}(\mathcal{O}_{X_{\mathbb{Q}_p,\text{pro\'etale},A}})\ar[r]^{\mathrm{equi}}\ar[r]\ar[r] &\varphi\mathrm{Ind}^m\mathrm{Banach}(\mathcal{O}_{Y_{\mathbb{Q}_p,\text{pro\'etale},A}}).\\  
}
\]
\item (\text{Proposition}) There is an equivalence between the $\infty$-categories of inductive Banach quasicoherent commutative algebra $E_\infty$ objects:
\[\displayindent=-0.4in
\xymatrix@R+6pc@C+0pc{
\mathrm{sComm}_\mathrm{simplicial}\mathrm{Ind}\mathrm{Banach}(\mathcal{O}_{X_{\mathbb{Q}_p\left<T_1^{\pm 1},...,T_k^{\pm 1}\right>,\text{pro\'etale},A}})\ar[d]\ar[d]\ar[d]\ar[d]\ar[r]^{\mathrm{equi}}\ar[r]\ar[r] &\mathrm{sComm}_\mathrm{simplicial}\varphi\mathrm{Ind}\mathrm{Banach}(\mathcal{O}_{Y_{\mathbb{Q}_p\left<T_1^{\pm 1},...,T_k^{\pm 1}\right>,\text{pro\'etale},A}})\ar[d]\ar[d]\ar[d]\ar[d]\\
\mathrm{sComm}_\mathrm{simplicial}\mathrm{Ind}\mathrm{Banach}(\mathcal{O}_{X_{\mathbb{Q}_p,\text{pro\'etale},A}})\ar[r]^{\mathrm{equi}}\ar[r]\ar[r] &\mathrm{sComm}_\mathrm{simplicial}\varphi\mathrm{Ind}\mathrm{Banach}(\mathcal{O}_{Y_{\mathbb{Q}_p,\text{pro\'etale},A}})  
}
\]
\item (\text{Proposition}) There is an equivalence between the $\infty$-categories of monomorphic inductive Banach quasicoherent commutative algebra $E_\infty$ objects:
\[\displayindent=-0.7in
\xymatrix@R+6pc@C+0pc{
\mathrm{sComm}_\mathrm{simplicial}\mathrm{Ind}^m\mathrm{Banach}(\mathcal{O}_{X_{\mathbb{Q}_p\left<T_1^{\pm 1},...,T_k^{\pm 1}\right>,\text{pro\'etale},A}})\ar[d]\ar[d]\ar[d]\ar[d]\ar[r]^{\mathrm{equi}}\ar[r]\ar[r] &\mathrm{sComm}_\mathrm{simplicial}\varphi\mathrm{Ind}^m\mathrm{Banach}(\mathcal{O}_{Y_{\mathbb{Q}_p\left<T_1^{\pm 1},...,T_k^{\pm 1}\right>,\text{pro\'etale},A}})\ar[d]\ar[d]\ar[d]\ar[d]\\
 \mathrm{sComm}_\mathrm{simplicial}\mathrm{Ind}^m\mathrm{Banach}(\mathcal{O}_{X_{\mathbb{Q}_p,\text{pro\'etale},A}})\ar[r]^{\mathrm{equi}}\ar[r]\ar[r] &\mathrm{sComm}_\mathrm{simplicial}\varphi\mathrm{Ind}^m\mathrm{Banach}(\mathcal{O}_{Y_{\mathbb{Q}_p,\text{pro\'etale},A}}) 
}
\]

\item Then parallel as in \cite{LBV} we have the equivalence of the de Rham complex after \cite[Definition 5.9, Section 5.2.1]{KKM}:
\[\displayindent=-0.9in
\xymatrix@R+6pc@C+0pc{
\mathrm{deRham}_{\mathrm{sComm}_\mathrm{simplicial}\mathrm{Ind}\mathrm{Banach}(\mathcal{O}_{X_{\mathbb{Q}_p\left<T_1^{\pm 1},...,T_k^{\pm 1}\right>,\text{pro\'etale},A}})\ar[r]^{\mathrm{equi}}}(-)\ar[d]\ar[d]\ar[d]\ar[d]\ar[r]\ar[r] &\mathrm{deRham}_{\mathrm{sComm}_\mathrm{simplicial}\varphi\mathrm{Ind}\mathrm{Banach}(\mathcal{O}_{Y_{\mathbb{Q}_p\left<T_1^{\pm 1},...,T_k^{\pm 1}\right>,\text{pro\'etale},A}})}(-)\ar[d]\ar[d]\ar[d]\ar[d]\\
\mathrm{deRham}_{\mathrm{sComm}_\mathrm{simplicial}\mathrm{Ind}\mathrm{Banach}(\mathcal{O}_{X_{\mathbb{Q}_p,\text{pro\'etale},A}})\ar[r]^{\mathrm{equi}}}(-)\ar[r]\ar[r] &\mathrm{deRham}_{\mathrm{sComm}_\mathrm{simplicial}\varphi\mathrm{Ind}\mathrm{Banach}(\mathcal{O}_{Y_{\mathbb{Q}_p,\text{pro\'etale},A}})}(-) 
}
\]
\[\displayindent=-1.0in
\xymatrix@R+6pc@C+0pc{
\mathrm{deRham}_{\mathrm{sComm}_\mathrm{simplicial}\mathrm{Ind}^m\mathrm{Banach}(\mathcal{O}_{X_{\mathbb{Q}_p\left<T_1^{\pm 1},...,T_k^{\pm 1}\right>,\text{pro\'etale},A}})\ar[r]^{\mathrm{equi}}}(-)\ar[d]\ar[d]\ar[d]\ar[d]\ar[r]\ar[r] &\mathrm{deRham}_{\mathrm{sComm}_\mathrm{simplicial}\varphi\mathrm{Ind}^m\mathrm{Banach}(\mathcal{O}_{Y_{\mathbb{Q}_p\left<T_1^{\pm 1},...,T_k^{\pm 1}\right>,\text{pro\'etale},A}})}(-)\ar[d]\ar[d]\ar[d]\ar[d]\\
\mathrm{deRham}_{\mathrm{sComm}_\mathrm{simplicial}\mathrm{Ind}^m\mathrm{Banach}(\mathcal{O}_{X_{\mathbb{Q}_p,\text{pro\'etale},A}})\ar[r]^{\mathrm{equi}}}(-)\ar[r]\ar[r] &\mathrm{deRham}_{\mathrm{sComm}_\mathrm{simplicial}\varphi\mathrm{Ind}^m\mathrm{Banach}(\mathcal{O}_{Y_{\mathbb{Q}_p,\text{pro\'etale},A}})}(-) 
}
\]

\item Then we have the following equivalence of $K$-group $(\infty,1)$-spectrum from \cite{BGT}:
\[\displayindent=-0.4in
\xymatrix@R+6pc@C+0pc{
\mathrm{K}^\mathrm{BGT}_{\mathrm{sComm}_\mathrm{simplicial}\mathrm{Ind}\mathrm{Banach}(\mathcal{O}_{X_{\mathbb{Q}_p\left<T_1^{\pm 1},...,T_k^{\pm 1}\right>,\text{pro\'etale},A}})\ar[r]^{\mathrm{equi}}}(-)\ar[d]\ar[d]\ar[d]\ar[d]\ar[r]\ar[r] &\mathrm{K}^\mathrm{BGT}_{\mathrm{sComm}_\mathrm{simplicial}\varphi\mathrm{Ind}\mathrm{Banach}(\mathcal{O}_{Y_{\mathbb{Q}_p\left<T_1^{\pm 1},...,T_k^{\pm 1}\right>,\text{pro\'etale},A}})}(-)\ar[d]\ar[d]\ar[d]\ar[d]\\
\mathrm{K}^\mathrm{BGT}_{\mathrm{sComm}_\mathrm{simplicial}\mathrm{Ind}\mathrm{Banach}(\mathcal{O}_{X_{\mathbb{Q}_p,\text{pro\'etale},A}})\ar[r]^{\mathrm{equi}}}(-)\ar[r]\ar[r] &\mathrm{K}^\mathrm{BGT}_{\mathrm{sComm}_\mathrm{simplicial}\varphi\mathrm{Ind}\mathrm{Banach}(\mathcal{O}_{Y_{\mathbb{Q}_p,\text{pro\'etale},A}})}(-) 
}
\]
\[\displayindent=-0.4in
\xymatrix@R+6pc@C+0pc{
\mathrm{K}^\mathrm{BGT}_{\mathrm{sComm}_\mathrm{simplicial}\mathrm{Ind}^m\mathrm{Banach}(\mathcal{O}_{X_{\mathbb{Q}_p\left<T_1^{\pm 1},...,T_k^{\pm 1}\right>,\text{pro\'etale},A}})\ar[r]^{\mathrm{equi}}}(-)\ar[d]\ar[d]\ar[d]\ar[d]\ar[r]\ar[r] &\mathrm{K}^\mathrm{BGT}_{\mathrm{sComm}_\mathrm{simplicial}\varphi\mathrm{Ind}^m\mathrm{Banach}(\mathcal{O}_{Y_{\mathbb{Q}_p\left<T_1^{\pm 1},...,T_k^{\pm 1}\right>,\text{pro\'etale},A}})}(-)\ar[d]\ar[d]\ar[d]\ar[d]\\
\mathrm{K}^\mathrm{BGT}_{\mathrm{sComm}_\mathrm{simplicial}\mathrm{Ind}^m\mathrm{Banach}(\mathcal{O}_{X_{\mathbb{Q}_p,\text{pro\'etale},A}})\ar[r]^{\mathrm{equi}}}(-)\ar[r]\ar[r] &\mathrm{K}^\mathrm{BGT}_{\mathrm{sComm}_\mathrm{simplicial}\varphi\mathrm{Ind}^m\mathrm{Banach}(\mathcal{O}_{Y_{\mathbb{Q}_p,\text{pro\'etale},A}})}(-). 
}
\]
\\

\end{itemize}

\

\indent Now we consider \cite{1CS1} and \cite[Proposition 13.8, Theorem 14.9, Remark 14.10]{1CS2}\footnote{Note that we are motivated as well from \cite{LBV}.}, and study the corresponding solid perfect complexes, solid quasicoherent sheaves and solid vector bundles. Here we are going to use different formalism, therefore we will have different categories and functors. We use the notation $\circledcirc$ to denote any element of $\{\text{solid perfect complexes}, \text{solid quasicoherent sheaves}, \text{solid vector bundles}\}$ from \cite{1CS2} with the corresponding descent results of \cite[Proposition 13.8, Theorem 14.9, Remark 14.10]{1CS2}. Then we have the following:

\begin{itemize}
\item Generalizing Kedlaya-Liu's construction in \cite{KL1}, \cite{KL2} of the adic Fargues-Fontaine space we have a quotient (by using powers of the Frobenius operator) $X_{R,A}$ of the space by using \cite{1CS2}:
\begin{align}
Y_{R,A}:=\bigcup_{0<s<r}\mathcal{X}^\mathrm{CS}(\text{Robba}^\text{extended}_{R,[s,r],A}).	
\end{align}

\item This is a locally ringed space $(X_{R,A},\mathcal{O}_{X_{R,A}})$, so one can consider the stable $\infty$-category $\mathrm{Module}_{\text{CS},\mathrm{quasicoherent}}(\mathcal{O}_{X_{R,A}}) $ which is the $\infty$-category of all the $\mathcal{O}_{X_{R,A}}$-sheaves of solid modules over $X_{R,A}$. We have the parallel categories for $Y_{R,A}$, namely $\varphi\mathrm{Module}_{\text{CS},\mathrm{quasicoherent}}(\mathcal{O}_{X_{R,A}})$ and so on. Here we will consider sheaves.

\begin{assumption}\label{assumtionpresheaves}
All the functors of modules or algebras below are Clausen-Scholze sheaves \cite[Proposition 13.8, Theorem 14.9, Remark 14.10]{1CS2}.	
\end{assumption}

\item (\text{Proposition}) There is an equivalence between the $\infty$-categories of inductive solid sheaves:
\[
\xymatrix@R+0pc@C+0pc{
\mathrm{Module}_\circledcirc(\mathcal{O}_{X_{R,A}})\ar[r]^{\mathrm{equi}}\ar[r]\ar[r] &\varphi\mathrm{Module}_\circledcirc(\mathcal{O}_{Y_{R,A}}).  
}
\]
\end{itemize}

\begin{itemize}

\item (\text{Proposition}) There is an equivalence between the $\infty$-categories of inductive Banach quasicoherent commutative algebra $E_\infty$ objects:
\[
\xymatrix@R+0pc@C+0pc{
\mathrm{sComm}_\mathrm{simplicial}\mathrm{Module}_{\text{solidquasicoherentsheaves}}(\mathcal{O}_{X_{R,A}})\ar[r]^{\mathrm{equi}}\ar[r]\ar[r] &\mathrm{sComm}_\mathrm{simplicial}\varphi\mathrm{Module}_{\text{solidquasicoherentsheaves}}(\mathcal{O}_{Y_{R,A}}).  
}
\]

\item Then as in \cite{LBV} we have the equivalence of the de Rham complex after \cite[Definition 5.9, Section 5.2.1]{KKM}\footnote{Here $\circledcirc=\text{solidquasicoherentsheaves}$.}:
\[
\xymatrix@R+0pc@C+0pc{
\mathrm{deRham}_{\mathrm{sComm}_\mathrm{simplicial}\mathrm{Module}_\circledcirc(\mathcal{O}_{X_{R,A}})\ar[r]^{\mathrm{equi}}}(-)\ar[r]\ar[r] &\mathrm{deRham}_{\mathrm{sComm}_\mathrm{simplicial}\varphi\mathrm{Module}_\circledcirc(\mathcal{O}_{Y_{R,A}})}(-). 
}
\]

\item Then we have the following equivalence of $K$-group $(\infty,1)$-spectrum from \cite{BGT}\footnote{Here $\circledcirc=\text{solidquasicoherentsheaves}$.}:
\[
\xymatrix@R+0pc@C+0pc{
\mathrm{K}^\mathrm{BGT}_{\mathrm{sComm}_\mathrm{simplicial}\mathrm{Module}_\circledcirc(\mathcal{O}_{X_{R,A}})\ar[r]^{\mathrm{equi}}}(-)\ar[r]\ar[r] &\mathrm{K}^\mathrm{BGT}_{\mathrm{sComm}_\mathrm{simplicial}\varphi\mathrm{Module}_\circledcirc(\mathcal{O}_{Y_{R,A}})}(-). 
}
\]
\end{itemize}

\noindent Now let $R=\mathbb{Q}_p(p^{1/p^\infty})^{\wedge\flat}$ and $R_k=\mathbb{Q}_p(p^{1/p^\infty})^{\wedge}\left<T_1^{\pm 1/p^{\infty}},...,T_k^{\pm 1/p^{\infty}}\right>^\flat$ we have the following Galois theoretic results with naturality along $f:\mathrm{Spa}(\mathbb{Q}_p(p^{1/p^\infty})^{\wedge}\left<T_1^{\pm 1/p^\infty},...,T_k^{\pm 1/p^\infty}\right>^\flat)\rightarrow \mathrm{Spa}(\mathbb{Q}_p(p^{1/p^\infty})^{\wedge\flat})$:

\begin{itemize}
\item (\text{Proposition}) There is an equivalence between the $\infty$-categories of inductive Banach quasicoherent presheaves\footnote{Here $\circledcirc=\text{solidquasicoherentsheaves}$.}:
\[
\xymatrix@R+6pc@C+0pc{
\mathrm{Modules}_\circledcirc(\mathcal{O}_{X_{\mathbb{Q}_p(p^{1/p^\infty})^{\wedge}\left<T_1^{\pm 1/p^\infty},...,T_k^{\pm 1/p^\infty}\right>^\flat,A}})\ar[d]\ar[d]\ar[d]\ar[d] \ar[r]^{\mathrm{equi}}\ar[r]\ar[r] &\varphi\mathrm{Modules}_\circledcirc(\mathcal{O}_{Y_{\mathbb{Q}_p(p^{1/p^\infty})^{\wedge}\left<T_1^{\pm 1/p^\infty},...,T_k^{\pm 1/p^\infty}\right>^\flat,A}}) \ar[d]\ar[d]\ar[d]\ar[d].\\
\mathrm{Modules}_\circledcirc(\mathcal{O}_{X_{\mathbb{Q}_p(p^{1/p^\infty})^{\wedge\flat},A}})\ar[r]^{\mathrm{equi}}\ar[r]\ar[r] &\varphi\mathrm{Modules}_\circledcirc(\mathcal{O}_{Y_{\mathbb{Q}_p(p^{1/p^\infty})^{\wedge\flat},A}}).\\ 
}
\]
\item (\text{Proposition}) There is an equivalence between the $\infty$-categories of inductive Banach quasicoherent commutative algebra $E_\infty$ objects\footnote{Here $\circledcirc=\text{solidquasicoherentsheaves}$.}:
\[
\xymatrix@R+6pc@C+0pc{
\mathrm{sComm}_\mathrm{simplicial}\mathrm{Modules}_\circledcirc(\mathcal{O}_{X_{R_k,A}})\ar[d]\ar[d]\ar[d]\ar[d]\ar[r]^{\mathrm{equi}}\ar[r]\ar[r] &\mathrm{sComm}_\mathrm{simplicial}\varphi\mathrm{Modules}_\circledcirc(\mathcal{O}_{Y_{R_k,A}})\ar[d]\ar[d]\ar[d]\ar[d]\\
\mathrm{sComm}_\mathrm{simplicial}\mathrm{Modules}_\circledcirc(\mathcal{O}_{X_{\mathbb{Q}_p(p^{1/p^\infty})^{\wedge\flat},A}})\ar[r]^{\mathrm{equi}}\ar[r]\ar[r] &\mathrm{sComm}_\mathrm{simplicial}\varphi\mathrm{Modules}_\circledcirc(\mathcal{O}_{Y_{\mathbb{Q}_p(p^{1/p^\infty})^{\wedge\flat},A}}).  
}
\]
\item Then as in \cite{LBV} we have the equivalence of the de Rham complex after \cite[Definition 5.9, Section 5.2.1]{KKM}\footnote{Here $\circledcirc=\text{solidquasicoherentsheaves}$.}:
\[
\xymatrix@R+6pc@C+0pc{
\mathrm{deRham}_{\mathrm{sComm}_\mathrm{simplicial}\mathrm{Modules}_\circledcirc(\mathcal{O}_{X_{R_k,A}})\ar[r]^{\mathrm{equi}}}(-)\ar[d]\ar[d]\ar[d]\ar[d]\ar[r]\ar[r] &\mathrm{deRham}_{\mathrm{sComm}_\mathrm{simplicial}\varphi\mathrm{Modules}_\circledcirc(\mathcal{O}_{Y_{R_k,A}})}(-)\ar[d]\ar[d]\ar[d]\ar[d]\\
\mathrm{deRham}_{\mathrm{sComm}_\mathrm{simplicial}\mathrm{Modules}_\circledcirc(\mathcal{O}_{X_{\mathbb{Q}_p(p^{1/p^\infty})^{\wedge\flat},A}})\ar[r]^{\mathrm{equi}}}(-)\ar[r]\ar[r] &\mathrm{deRham}_{\mathrm{sComm}_\mathrm{simplicial}\varphi\mathrm{Modules}_\circledcirc(\mathcal{O}_{Y_{\mathbb{Q}_p(p^{1/p^\infty})^{\wedge\flat},A}})}(-). 
}
\]

\item Then we have the following equivalence of $K$-group $(\infty,1)$-spectrum from \cite{BGT}\footnote{Here $\circledcirc=\text{solidquasicoherentsheaves}$.}:
\[
\xymatrix@R+6pc@C+0pc{
\mathrm{K}^\mathrm{BGT}_{\mathrm{sComm}_\mathrm{simplicial}\mathrm{Modules}_\circledcirc(\mathcal{O}_{X_{R_k,A}})\ar[r]^{\mathrm{equi}}}(-)\ar[d]\ar[d]\ar[d]\ar[d]\ar[r]\ar[r] &\mathrm{K}^\mathrm{BGT}_{\mathrm{sComm}_\mathrm{simplicial}\varphi\mathrm{Modules}_\circledcirc(\mathcal{O}_{Y_{R_k,A}})}(-)\ar[d]\ar[d]\ar[d]\ar[d]\\
\mathrm{K}^\mathrm{BGT}_{\mathrm{sComm}_\mathrm{simplicial}\mathrm{Modules}_\circledcirc(\mathcal{O}_{X_{\mathbb{Q}_p(p^{1/p^\infty})^{\wedge\flat},A}})\ar[r]^{\mathrm{equi}}}(-)\ar[r]\ar[r] &\mathrm{K}^\mathrm{BGT}_{\mathrm{sComm}_\mathrm{simplicial}\varphi\mathrm{Modules}_\circledcirc(\mathcal{O}_{Y_{\mathbb{Q}_p(p^{1/p^\infty})^{\wedge\flat},A}})}(-).}
\]

\end{itemize}

\
\indent Then we consider further equivariance by considering the arithmetic profinite fundamental groups $\Gamma_{\mathbb{Q}_p}$ and $\mathrm{Gal}(\overline{\mathbb{Q}_p\left<T_1^{\pm 1},...,T_k^{\pm 1}\right>}/R_k)$ through the following diagram:

\[
\xymatrix@R+0pc@C+0pc{
\mathbb{Z}_p^k=\mathrm{Gal}(R_k/{\mathbb{Q}_p(p^{1/p^\infty})^\wedge\left<T_1^{\pm 1},...,T_k^{\pm 1}\right>}) \ar[r]\ar[r] \ar[r]\ar[r] &\Gamma_k:=\mathrm{Gal}(R_k/{\mathbb{Q}_p\left<T_1^{\pm 1},...,T_k^{\pm 1}\right>}) \ar[r] \ar[r]\ar[r] &\Gamma_{\mathbb{Q}_p}.
}
\]

\begin{itemize}
\item (\text{Proposition}) There is an equivalence between the $\infty$-categories of inductive Banach quasicoherent presheaves\footnote{Here $\circledcirc=\text{solidquasicoherentsheaves}$.}:
\[
\xymatrix@R+6pc@C+0pc{
\mathrm{Modules}_{\circledcirc,\Gamma_k}(\mathcal{O}_{X_{\mathbb{Q}_p(p^{1/p^\infty})^{\wedge}\left<T_1^{\pm 1/p^\infty},...,T_k^{\pm 1/p^\infty}\right>^\flat,A}})\ar[d]\ar[d]\ar[d]\ar[d] \ar[r]^{\mathrm{equi}}\ar[r]\ar[r] &\varphi\mathrm{Modules}_{\circledcirc,\Gamma_k}(\mathcal{O}_{Y_{\mathbb{Q}_p(p^{1/p^\infty})^{\wedge}\left<T_1^{\pm 1/p^\infty},...,T_k^{\pm 1/p^\infty}\right>^\flat,A}}) \ar[d]\ar[d]\ar[d]\ar[d].\\
\mathrm{Modules}_{\circledcirc,\Gamma_0}(\mathcal{O}_{X_{\mathbb{Q}_p(p^{1/p^\infty})^{\wedge\flat},A}})\ar[r]^{\mathrm{equi}}\ar[r]\ar[r] &\varphi\mathrm{Modules}_{\circledcirc,\Gamma_0}(\mathcal{O}_{Y_{\mathbb{Q}_p(p^{1/p^\infty})^{\wedge\flat},A}}).\\ 
}
\]

\item (\text{Proposition}) There is an equivalence between the $\infty$-categories of inductive Banach quasicoherent commutative algebra $E_\infty$ objects\footnote{Here $\circledcirc=\text{solidquasicoherentsheaves}$.}:
\[
\xymatrix@R+6pc@C+0pc{
\mathrm{sComm}_\mathrm{simplicial}\mathrm{Modules}_{\circledcirc,\Gamma_k}(\mathcal{O}_{X_{R_k,A}})\ar[d]\ar[d]\ar[d]\ar[d]\ar[r]^{\mathrm{equi}}\ar[r]\ar[r] &\mathrm{sComm}_\mathrm{simplicial}\varphi\mathrm{Modules}_{\circledcirc,\Gamma_k}(\mathcal{O}_{Y_{R_k,A}})\ar[d]\ar[d]\ar[d]\ar[d]\\
\mathrm{sComm}_\mathrm{simplicial}\mathrm{Modules}_{\circledcirc,\Gamma_0}(\mathcal{O}_{X_{\mathbb{Q}_p(p^{1/p^\infty})^{\wedge\flat},A}})\ar[r]^{\mathrm{equi}}\ar[r]\ar[r] &\mathrm{sComm}_\mathrm{simplicial}\varphi\mathrm{Modules}_{\circledcirc,\Gamma_0}(\mathcal{O}_{Y_{\mathbb{Q}_p(p^{1/p^\infty})^{\wedge\flat},A}}).  
}
\]

\item Then as in \cite{LBV} we have the equivalence of the de Rham complex after \cite[Definition 5.9, Section 5.2.1]{KKM}\footnote{Here $\circledcirc=\text{solidquasicoherentsheaves}$.}:
\[\displayindent=-0.4in
\xymatrix@R+6pc@C+0pc{
\mathrm{deRham}_{\mathrm{sComm}_\mathrm{simplicial}\mathrm{Modules}_{\circledcirc,\Gamma_k}(\mathcal{O}_{X_{R_k,A}})\ar[r]^{\mathrm{equi}}}(-)\ar[d]\ar[d]\ar[d]\ar[d]\ar[r]\ar[r] &\mathrm{deRham}_{\mathrm{sComm}_\mathrm{simplicial}\varphi\mathrm{Modules}_{\circledcirc,\Gamma_k}(\mathcal{O}_{Y_{R_k,A}})}(-)\ar[d]\ar[d]\ar[d]\ar[d]\\
\mathrm{deRham}_{\mathrm{sComm}_\mathrm{simplicial}\mathrm{Modules}_{\circledcirc,\Gamma_0}(\mathcal{O}_{X_{\mathbb{Q}_p(p^{1/p^\infty})^{\wedge\flat},A}})\ar[r]^{\mathrm{equi}}}(-)\ar[r]\ar[r] &\mathrm{deRham}_{\mathrm{sComm}_\mathrm{simplicial}\varphi\mathrm{Modules}_{\circledcirc,\Gamma_0}(\mathcal{O}_{Y_{\mathbb{Q}_p(p^{1/p^\infty})^{\wedge\flat},A}})}(-). 
}
\]

\item Then we have the following equivalence of $K$-group $(\infty,1)$-spectrum from \cite{BGT}\footnote{Here $\circledcirc=\text{solidquasicoherentsheaves}$.}:
\[
\xymatrix@R+6pc@C+0pc{
\mathrm{K}^\mathrm{BGT}_{\mathrm{sComm}_\mathrm{simplicial}\mathrm{Modules}_{\circledcirc,\Gamma_k}(\mathcal{O}_{X_{R_k,A}})\ar[r]^{\mathrm{equi}}}(-)\ar[d]\ar[d]\ar[d]\ar[d]\ar[r]\ar[r] &\mathrm{K}^\mathrm{BGT}_{\mathrm{sComm}_\mathrm{simplicial}\varphi\mathrm{Modules}_{\circledcirc,\Gamma_k}(\mathcal{O}_{Y_{R_k,A}})}(-)\ar[d]\ar[d]\ar[d]\ar[d]\\
\mathrm{K}^\mathrm{BGT}_{\mathrm{sComm}_\mathrm{simplicial}\mathrm{Modules}_{\circledcirc,\Gamma_0}(\mathcal{O}_{X_{\mathbb{Q}_p(p^{1/p^\infty})^{\wedge\flat},A}})\ar[r]^{\mathrm{equi}}}(-)\ar[r]\ar[r] &\mathrm{K}^\mathrm{BGT}_{\mathrm{sComm}_\mathrm{simplicial}\varphi\mathrm{Modules}_{\circledcirc,\Gamma_0}(\mathcal{O}_{Y_{\mathbb{Q}_p(p^{1/p^\infty})^{\wedge\flat},A}})}(-).
}
\]

\end{itemize}

\

Furthermore we have the corresponding pro-\'etale version without the corresponding fundamental group equivariances.

\begin{itemize}
\item (\text{Proposition}) There is an equivalence between the $\infty$-categories of inductive Banach quasicoherent presheaves\footnote{Here $\circledcirc=\text{solidquasicoherentsheaves}$.}:
\[
\xymatrix@R+6pc@C+0pc{
\mathrm{Modules}_\circledcirc(\mathcal{O}_{X_{\mathbb{Q}_p\left<T_1^{\pm 1},...,T_k^{\pm 1}\right>,\text{pro\'etale},A}})\ar[d]\ar[d]\ar[d]\ar[d] \ar[r]^{\mathrm{equi}}\ar[r]\ar[r] &\varphi\mathrm{Modules}_\circledcirc(\mathcal{O}_{Y_{\mathbb{Q}_p\left<T_1^{\pm 1},...,T_k^{\pm 1}\right>,\text{pro\'etale},A}}) \ar[d]\ar[d]\ar[d]\ar[d].\\
\mathrm{Modules}_\circledcirc(\mathcal{O}_{X_{\mathbb{Q}_p,\text{pro\'etale},A}})\ar[r]^{\mathrm{equi}}\ar[r]\ar[r] &\varphi\mathrm{Modules}_\circledcirc(\mathcal{O}_{Y_{\mathbb{Q}_p,\text{pro\'etale},A}}).\\ 
}
\]

\item (\text{Proposition}) There is an equivalence between the $\infty$-categories of inductive Banach quasicoherent commutative algebra $E_\infty$ objects\footnote{Here $\circledcirc=\text{solidquasicoherentsheaves}$.}:
\[\displayindent=-0.4in
\xymatrix@R+6pc@C+0pc{
\mathrm{sComm}_\mathrm{simplicial}\mathrm{Modules}_\circledcirc(\mathcal{O}_{X_{\mathbb{Q}_p\left<T_1^{\pm 1},...,T_k^{\pm 1}\right>,\text{pro\'etale},A}})\ar[d]\ar[d]\ar[d]\ar[d]\ar[r]^{\mathrm{equi}}\ar[r]\ar[r] &\mathrm{sComm}_\mathrm{simplicial}\varphi\mathrm{Modules}_\circledcirc(\mathcal{O}_{Y_{\mathbb{Q}_p\left<T_1^{\pm 1},...,T_k^{\pm 1}\right>,\text{pro\'etale},A}})\ar[d]\ar[d]\ar[d]\ar[d]\\
\mathrm{sComm}_\mathrm{simplicial}\mathrm{Modules}_\circledcirc(\mathcal{O}_{X_{\mathbb{Q}_p,\text{pro\'etale},A}})\ar[r]^{\mathrm{equi}}\ar[r]\ar[r] &\mathrm{sComm}_\mathrm{simplicial}\varphi\mathrm{Modules}_\circledcirc(\mathcal{O}_{Y_{\mathbb{Q}_p,\text{pro\'etale},A}}).  
}
\]

\item Then as in \cite{LBV} we have the equivalence of the de Rham complex after \cite[Definition 5.9, Section 5.2.1]{KKM}\footnote{Here $\circledcirc=\text{solidquasicoherentsheaves}$.}:
\[\displayindent=-0.9in
\xymatrix@R+6pc@C+0pc{
\mathrm{deRham}_{\mathrm{sComm}_\mathrm{simplicial}\mathrm{Modules}_\circledcirc(\mathcal{O}_{X_{\mathbb{Q}_p\left<T_1^{\pm 1},...,T_k^{\pm 1}\right>,\text{pro\'etale},A}})\ar[r]^{\mathrm{equi}}}(-)\ar[d]\ar[d]\ar[d]\ar[d]\ar[r]\ar[r] &\mathrm{deRham}_{\mathrm{sComm}_\mathrm{simplicial}\varphi\mathrm{Modules}_\circledcirc(\mathcal{O}_{Y_{\mathbb{Q}_p\left<T_1^{\pm 1},...,T_k^{\pm 1}\right>,\text{pro\'etale},A}})}(-)\ar[d]\ar[d]\ar[d]\ar[d]\\
\mathrm{deRham}_{\mathrm{sComm}_\mathrm{simplicial}\mathrm{Modules}_\circledcirc(\mathcal{O}_{X_{\mathbb{Q}_p,\text{pro\'etale},A}})\ar[r]^{\mathrm{equi}}}(-)\ar[r]\ar[r] &\mathrm{deRham}_{\mathrm{sComm}_\mathrm{simplicial}\varphi\mathrm{Modules}_\circledcirc(\mathcal{O}_{Y_{\mathbb{Q}_p,\text{pro\'etale},A}})}(-). 
}
\]

\item Then we have the following equivalence of $K$-group $(\infty,1)$-spectrum from \cite{BGT}\footnote{Here $\circledcirc=\text{solidquasicoherentsheaves}$.}:
\[\displayindent=-0.4in
\xymatrix@R+6pc@C+0pc{
\mathrm{K}^\mathrm{BGT}_{\mathrm{sComm}_\mathrm{simplicial}\mathrm{Modules}_\circledcirc(\mathcal{O}_{X_{\mathbb{Q}_p\left<T_1^{\pm 1},...,T_k^{\pm 1}\right>,\text{pro\'etale},A}})\ar[r]^{\mathrm{equi}}}(-)\ar[d]\ar[d]\ar[d]\ar[d]\ar[r]\ar[r] &\mathrm{K}^\mathrm{BGT}_{\mathrm{sComm}_\mathrm{simplicial}\varphi\mathrm{Modules}_\circledcirc(\mathcal{O}_{Y_{\mathbb{Q}_p\left<T_1^{\pm 1},...,T_k^{\pm 1}\right>,\text{pro\'etale},A}})}(-)\ar[d]\ar[d]\ar[d]\ar[d]\\
\mathrm{K}^\mathrm{BGT}_{\mathrm{sComm}_\mathrm{simplicial}\mathrm{Modules}_\circledcirc(\mathcal{O}_{X_{\mathbb{Q}_p,\text{pro\'etale},A}})\ar[r]^{\mathrm{equi}}}(-)\ar[r]\ar[r] &\mathrm{K}^\mathrm{BGT}_{\mathrm{sComm}_\mathrm{simplicial}\varphi\mathrm{Modules}_\circledcirc(\mathcal{O}_{Y_{\mathbb{Q}_p,\text{pro\'etale},A}})}(-). 
}
\]
	
\end{itemize}

\newpage
\subsection{$\infty$-Categorical Analytic Stacks and Descents II}

\indent Then by the corresponding \v{C}ech $\infty$-descent in \cite[Section 9.3]{KKM} and \cite{BBM} we have the following objects by directly taking the corresponding \v{C}ech $\infty$-descent. In the following the right had of each row in each diagram will be the corresponding quasicoherent Robba bundles over the Robba ring carrying the corresponding action from the Frobenius or the fundamental groups, defined by directly applying \cite[Section 9.3]{KKM} and \cite{BBM}. We then have the following global section functors:

\begin{itemize}

\item (\text{Proposition}) There is a functor (global section) between the $\infty$-categories of inductive Banach quasicoherent presheaves:
\[
\xymatrix@R+0pc@C+0pc{
\mathrm{Ind}\mathrm{Banach}(\mathcal{O}_{X_{R,A}})\ar[r]^{\mathrm{global}}\ar[r]\ar[r] &\varphi\mathrm{Ind}\mathrm{Banach}(\{\mathrm{Robba}^\mathrm{extended}_{{R,A,I}}\}_I).  
}
\]
\item (\text{Proposition}) There is a functor (global section) between the $\infty$-categories of monomorphic inductive Banach quasicoherent presheaves:
\[
\xymatrix@R+0pc@C+0pc{
\mathrm{Ind}^m\mathrm{Banach}(\mathcal{O}_{X_{R,A}})\ar[r]^{\mathrm{global}}\ar[r]\ar[r] &\varphi\mathrm{Ind}^m\mathrm{Banach}(\{\mathrm{Robba}^\mathrm{extended}_{{R,A,I}}\}_I).  
}
\]

\item (\text{Proposition}) There is a functor (global section) between the $\infty$-categories of inductive Banach quasicoherent presheaves:
\[
\xymatrix@R+0pc@C+0pc{
\mathrm{Ind}\mathrm{Banach}(\mathcal{O}_{X_{R,A}})\ar[r]^{\mathrm{global}}\ar[r]\ar[r] &\varphi\mathrm{Ind}\mathrm{Banach}(\{\mathrm{Robba}^\mathrm{extended}_{{R,A,I}}\}_I).  
}
\]
\item (\text{Proposition}) There is a functor (global section) between the $\infty$-categories of monomorphic inductive Banach quasicoherent presheaves:
\[
\xymatrix@R+0pc@C+0pc{
\mathrm{Ind}^m\mathrm{Banach}(\mathcal{O}_{X_{R,A}})\ar[r]^{\mathrm{global}}\ar[r]\ar[r] &\varphi\mathrm{Ind}^m\mathrm{Banach}(\{\mathrm{Robba}^\mathrm{extended}_{{R,A,I}}\}_I).  
}
\]
\item (\text{Proposition}) There is a functor (global section) between the $\infty$-categories of inductive Banach quasicoherent commutative algebra $E_\infty$ objects:
\[
\xymatrix@R+0pc@C+0pc{
\mathrm{sComm}_\mathrm{simplicial}\mathrm{Ind}\mathrm{Banach}(\mathcal{O}_{X_{R,A}})\ar[r]^{\mathrm{global}}\ar[r]\ar[r] &\mathrm{sComm}_\mathrm{simplicial}\varphi\mathrm{Ind}\mathrm{Banach}(\{\mathrm{Robba}^\mathrm{extended}_{{R,A,I}}\}_I).  
}
\]
\item (\text{Proposition}) There is a functor (global section) between the $\infty$-categories of monomorphic inductive Banach quasicoherent commutative algebra $E_\infty$ objects:
\[
\xymatrix@R+0pc@C+0pc{
\mathrm{sComm}_\mathrm{simplicial}\mathrm{Ind}^m\mathrm{Banach}(\mathcal{O}_{X_{R,A}})\ar[r]^{\mathrm{global}}\ar[r]\ar[r] &\mathrm{sComm}_\mathrm{simplicial}\varphi\mathrm{Ind}^m\mathrm{Banach}(\{\mathrm{Robba}^\mathrm{extended}_{{R,A,I}}\}_I).  
}
\]

\item Then parallel as in \cite{LBV} we have a functor (global section) of the de Rham complex after \cite[Definition 5.9, Section 5.2.1]{KKM}:
\[
\xymatrix@R+0pc@C+0pc{
\mathrm{deRham}_{\mathrm{sComm}_\mathrm{simplicial}\mathrm{Ind}\mathrm{Banach}(\mathcal{O}_{X_{R,A}})\ar[r]^{\mathrm{global}}}(-)\ar[r]\ar[r] &\mathrm{deRham}_{\mathrm{sComm}_\mathrm{simplicial}\varphi\mathrm{Ind}\mathrm{Banach}(\{\mathrm{Robba}^\mathrm{extended}_{{R,A,I}}\}_I)}(-), 
}
\]
\[
\xymatrix@R+0pc@C+0pc{
\mathrm{deRham}_{\mathrm{sComm}_\mathrm{simplicial}\mathrm{Ind}^m\mathrm{Banach}(\mathcal{O}_{X_{R,A}})\ar[r]^{\mathrm{global}}}(-)\ar[r]\ar[r] &\mathrm{deRham}_{\mathrm{sComm}_\mathrm{simplicial}\varphi\mathrm{Ind}^m\mathrm{Banach}(\{\mathrm{Robba}^\mathrm{extended}_{{R,A,I}}\}_I)}(-). 
}
\]

\item Then we have the following a functor (global section) of $K$-group $(\infty,1)$-spectrum from \cite{BGT}:
\[
\xymatrix@R+0pc@C+0pc{
\mathrm{K}^\mathrm{BGT}_{\mathrm{sComm}_\mathrm{simplicial}\mathrm{Ind}\mathrm{Banach}(\mathcal{O}_{X_{R,A}})\ar[r]^{\mathrm{global}}}(-)\ar[r]\ar[r] &\mathrm{K}^\mathrm{BGT}_{\mathrm{sComm}_\mathrm{simplicial}\varphi\mathrm{Ind}\mathrm{Banach}(\{\mathrm{Robba}^\mathrm{extended}_{{R,A,I}}\}_I)}(-), 
}
\]
\[
\xymatrix@R+0pc@C+0pc{
\mathrm{K}^\mathrm{BGT}_{\mathrm{sComm}_\mathrm{simplicial}\mathrm{Ind}^m\mathrm{Banach}(\mathcal{O}_{X_{R,A}})\ar[r]^{\mathrm{global}}}(-)\ar[r]\ar[r] &\mathrm{K}^\mathrm{BGT}_{\mathrm{sComm}_\mathrm{simplicial}\varphi\mathrm{Ind}^m\mathrm{Banach}(\{\mathrm{Robba}^\mathrm{extended}_{{R,A,I}}\}_I)}(-). 
}
\]
\end{itemize}

\noindent Now let $R=\mathbb{Q}_p(p^{1/p^\infty})^{\wedge\flat}$ and $R_k=\mathbb{Q}_p(p^{1/p^\infty})^{\wedge}\left<T_1^{\pm 1/p^{\infty}},...,T_k^{\pm 1/p^{\infty}}\right>^\flat$ we have the following Galois theoretic results with naturality along $f:\mathrm{Spa}(\mathbb{Q}_p(p^{1/p^\infty})^{\wedge}\left<T_1^{\pm 1/p^\infty},...,T_k^{\pm 1/p^\infty}\right>^\flat)\rightarrow \mathrm{Spa}(\mathbb{Q}_p(p^{1/p^\infty})^{\wedge\flat})$:

\begin{itemize}
\item (\text{Proposition}) There is a functor (global section) between the $\infty$-categories of inductive Banach quasicoherent presheaves:
\[
\xymatrix@R+6pc@C+0pc{
\mathrm{Ind}\mathrm{Banach}(\mathcal{O}_{X_{\mathbb{Q}_p(p^{1/p^\infty})^{\wedge}\left<T_1^{\pm 1/p^\infty},...,T_k^{\pm 1/p^\infty}\right>^\flat,A}})\ar[d]\ar[d]\ar[d]\ar[d] \ar[r]^{\mathrm{global}}\ar[r]\ar[r] &\varphi\mathrm{Ind}\mathrm{Banach}(\{\mathrm{Robba}^\mathrm{extended}_{{R_k,A,I}}\}_I) \ar[d]\ar[d]\ar[d]\ar[d].\\
\mathrm{Ind}\mathrm{Banach}(\mathcal{O}_{X_{\mathbb{Q}_p(p^{1/p^\infty})^{\wedge\flat},A}})\ar[r]^{\mathrm{global}}\ar[r]\ar[r] &\varphi\mathrm{Ind}\mathrm{Banach}(\{\mathrm{Robba}^\mathrm{extended}_{{R_0,A,I}}\}_I).\\ 
}
\]
\item (\text{Proposition}) There is a functor (global section) between the $\infty$-categories of monomorphic inductive Banach quasicoherent presheaves:
\[
\xymatrix@R+6pc@C+0pc{
\mathrm{Ind}^m\mathrm{Banach}(\mathcal{O}_{X_{R_k,A}})\ar[r]^{\mathrm{global}}\ar[d]\ar[d]\ar[d]\ar[d]\ar[r]\ar[r] &\varphi\mathrm{Ind}^m\mathrm{Banach}(\{\mathrm{Robba}^\mathrm{extended}_{{R_k,A,I}}\}_I)\ar[d]\ar[d]\ar[d]\ar[d]\\
\mathrm{Ind}^m\mathrm{Banach}(\mathcal{O}_{X_{\mathbb{Q}_p(p^{1/p^\infty})^{\wedge\flat},A}})\ar[r]^{\mathrm{global}}\ar[r]\ar[r] &\varphi\mathrm{Ind}^m\mathrm{Banach}(\{\mathrm{Robba}^\mathrm{extended}_{{R_0,A,I}}\}_I).\\  
}
\]
\item (\text{Proposition}) There is a functor (global section) between the $\infty$-categories of inductive Banach quasicoherent commutative algebra $E_\infty$ objects:
\[
\xymatrix@R+6pc@C+0pc{
\mathrm{sComm}_\mathrm{simplicial}\mathrm{Ind}\mathrm{Banach}(\mathcal{O}_{X_{R_k,A}})\ar[d]\ar[d]\ar[d]\ar[d]\ar[r]^{\mathrm{global}}\ar[r]\ar[r] &\mathrm{sComm}_\mathrm{simplicial}\varphi\mathrm{Ind}\mathrm{Banach}(\{\mathrm{Robba}^\mathrm{extended}_{{R_k,A,I}}\}_I)\ar[d]\ar[d]\ar[d]\ar[d]\\
\mathrm{sComm}_\mathrm{simplicial}\mathrm{Ind}\mathrm{Banach}(\mathcal{O}_{X_{\mathbb{Q}_p(p^{1/p^\infty})^{\wedge\flat},A}})\ar[r]^{\mathrm{global}}\ar[r]\ar[r] &\mathrm{sComm}_\mathrm{simplicial}\varphi\mathrm{Ind}\mathrm{Banach}(\{\mathrm{Robba}^\mathrm{extended}_{{R_0,A,I}}\}_I).  
}
\]
\item (\text{Proposition}) There is a functor (global section) between the $\infty$-categories of monomorphic inductive Banach quasicoherent commutative algebra $E_\infty$ objects:
\[
\xymatrix@R+6pc@C+0pc{
\mathrm{sComm}_\mathrm{simplicial}\mathrm{Ind}^m\mathrm{Banach}(\mathcal{O}_{X_{R_k,A}})\ar[d]\ar[d]\ar[d]\ar[d]\ar[r]^{\mathrm{global}}\ar[r]\ar[r] &\mathrm{sComm}_\mathrm{simplicial}\varphi\mathrm{Ind}^m\mathrm{Banach}(\{\mathrm{Robba}^\mathrm{extended}_{{R_k,A,I}}\}_I)\ar[d]\ar[d]\ar[d]\ar[d]\\
 \mathrm{sComm}_\mathrm{simplicial}\mathrm{Ind}^m\mathrm{Banach}(\mathcal{O}_{X_{\mathbb{Q}_p(p^{1/p^\infty})^{\wedge\flat},A}})\ar[r]^{\mathrm{global}}\ar[r]\ar[r] &\mathrm{sComm}_\mathrm{simplicial}\varphi\mathrm{Ind}^m\mathrm{Banach}(\{\mathrm{Robba}^\mathrm{extended}_{{R_0,A,I}}\}_I). 
}
\]

\item Then parallel as in \cite{LBV} we have a functor (global section) of the de Rham complex after \cite[Definition 5.9, Section 5.2.1]{KKM}:
\[\displayindent=-0.4in
\xymatrix@R+6pc@C+0pc{
\mathrm{deRham}_{\mathrm{sComm}_\mathrm{simplicial}\mathrm{Ind}\mathrm{Banach}(\mathcal{O}_{X_{R_k,A}})\ar[r]^{\mathrm{global}}}(-)\ar[d]\ar[d]\ar[d]\ar[d]\ar[r]\ar[r] &\mathrm{deRham}_{\mathrm{sComm}_\mathrm{simplicial}\varphi\mathrm{Ind}\mathrm{Banach}(\{\mathrm{Robba}^\mathrm{extended}_{{R_k,A,I}}\}_I)}(-)\ar[d]\ar[d]\ar[d]\ar[d]\\
\mathrm{deRham}_{\mathrm{sComm}_\mathrm{simplicial}\mathrm{Ind}\mathrm{Banach}(\mathcal{O}_{X_{\mathbb{Q}_p(p^{1/p^\infty})^{\wedge\flat},A}})\ar[r]^{\mathrm{global}}}(-)\ar[r]\ar[r] &\mathrm{deRham}_{\mathrm{sComm}_\mathrm{simplicial}\varphi\mathrm{Ind}\mathrm{Banach}(\{\mathrm{Robba}^\mathrm{extended}_{{R_0,A,I}}\}_I)}(-), 
}
\]
\[\displayindent=-0.4in
\xymatrix@R+6pc@C+0pc{
\mathrm{deRham}_{\mathrm{sComm}_\mathrm{simplicial}\mathrm{Ind}^m\mathrm{Banach}(\mathcal{O}_{X_{R_k,A}})\ar[r]^{\mathrm{global}}}(-)\ar[d]\ar[d]\ar[d]\ar[d]\ar[r]\ar[r] &\mathrm{deRham}_{\mathrm{sComm}_\mathrm{simplicial}\varphi\mathrm{Ind}^m\mathrm{Banach}(\{\mathrm{Robba}^\mathrm{extended}_{{R_k,A,I}}\}_I)}(-)\ar[d]\ar[d]\ar[d]\ar[d]\\
\mathrm{deRham}_{\mathrm{sComm}_\mathrm{simplicial}\mathrm{Ind}^m\mathrm{Banach}(\mathcal{O}_{X_{\mathbb{Q}_p(p^{1/p^\infty})^{\wedge\flat},A}})\ar[r]^{\mathrm{global}}}(-)\ar[r]\ar[r] &\mathrm{deRham}_{\mathrm{sComm}_\mathrm{simplicial}\varphi\mathrm{Ind}^m\mathrm{Banach}(\{\mathrm{Robba}^\mathrm{extended}_{{R_0,A,I}}\}_I)}(-). 
}
\]

\item Then we have the following a functor (global section) of $K$-group $(\infty,1)$-spectrum from \cite{BGT}:
\[
\xymatrix@R+6pc@C+0pc{
\mathrm{K}^\mathrm{BGT}_{\mathrm{sComm}_\mathrm{simplicial}\mathrm{Ind}\mathrm{Banach}(\mathcal{O}_{X_{R_k,A}})\ar[r]^{\mathrm{global}}}(-)\ar[d]\ar[d]\ar[d]\ar[d]\ar[r]\ar[r] &\mathrm{K}^\mathrm{BGT}_{\mathrm{sComm}_\mathrm{simplicial}\varphi\mathrm{Ind}\mathrm{Banach}(\{\mathrm{Robba}^\mathrm{extended}_{{R_k,A,I}}\}_I)}(-)\ar[d]\ar[d]\ar[d]\ar[d]\\
\mathrm{K}^\mathrm{BGT}_{\mathrm{sComm}_\mathrm{simplicial}\mathrm{Ind}\mathrm{Banach}(\mathcal{O}_{X_{\mathbb{Q}_p(p^{1/p^\infty})^{\wedge\flat},A}})\ar[r]^{\mathrm{global}}}(-)\ar[r]\ar[r] &\mathrm{K}^\mathrm{BGT}_{\mathrm{sComm}_\mathrm{simplicial}\varphi\mathrm{Ind}\mathrm{Banach}(\{\mathrm{Robba}^\mathrm{extended}_{{R_0,A,I}}\}_I)}(-), 
}
\]
\[
\xymatrix@R+6pc@C+0pc{
\mathrm{K}^\mathrm{BGT}_{\mathrm{sComm}_\mathrm{simplicial}\mathrm{Ind}^m\mathrm{Banach}(\mathcal{O}_{X_{R_k,A}})\ar[r]^{\mathrm{global}}}(-)\ar[d]\ar[d]\ar[d]\ar[d]\ar[r]\ar[r] &\mathrm{K}^\mathrm{BGT}_{\mathrm{sComm}_\mathrm{simplicial}\varphi\mathrm{Ind}^m\mathrm{Banach}(\{\mathrm{Robba}^\mathrm{extended}_{{R_k,A,I}}\}_I)}(-)\ar[d]\ar[d]\ar[d]\ar[d]\\
\mathrm{K}^\mathrm{BGT}_{\mathrm{sComm}_\mathrm{simplicial}\mathrm{Ind}^m\mathrm{Banach}(\mathcal{O}_{X_{\mathbb{Q}_p(p^{1/p^\infty})^{\wedge\flat},A}})\ar[r]^{\mathrm{global}}}(-)\ar[r]\ar[r] &\mathrm{K}^\mathrm{BGT}_{\mathrm{sComm}_\mathrm{simplicial}\varphi\mathrm{Ind}^m\mathrm{Banach}(\{\mathrm{Robba}^\mathrm{extended}_{{R_0,A,I}}\}_I)}(-). 
}
\]

\end{itemize}

\
\indent Then we consider further equivariance by considering the arithmetic profinite fundamental groups $\Gamma_{\mathbb{Q}_p}$ and $\mathrm{Gal}(\overline{\mathbb{Q}_p\left<T_1^{\pm 1},...,T_k^{\pm 1}\right>}/R_k)$ through the following diagram:

\[
\xymatrix@R+0pc@C+0pc{
\mathbb{Z}_p^k=\mathrm{Gal}(R_k/{\mathbb{Q}_p(p^{1/p^\infty})^\wedge\left<T_1^{\pm 1},...,T_k^{\pm 1}\right>}) \ar[r]\ar[r] \ar[r]\ar[r] &\Gamma_k:=\mathrm{Gal}(R_k/{\mathbb{Q}_p\left<T_1^{\pm 1},...,T_k^{\pm 1}\right>}) \ar[r] \ar[r]\ar[r] &\Gamma_{\mathbb{Q}_p}.
}
\]

\begin{itemize}
\item (\text{Proposition}) There is a functor (global section) between the $\infty$-categories of inductive Banach quasicoherent presheaves:
\[
\xymatrix@R+6pc@C+0pc{
\mathrm{Ind}\mathrm{Banach}_{\Gamma_k}(\mathcal{O}_{X_{\mathbb{Q}_p(p^{1/p^\infty})^{\wedge}\left<T_1^{\pm 1/p^\infty},...,T_k^{\pm 1/p^\infty}\right>^\flat,A}})\ar[d]\ar[d]\ar[d]\ar[d] \ar[r]^{\mathrm{global}}\ar[r]\ar[r] &\varphi\mathrm{Ind}\mathrm{Banach}_{\Gamma_k}(\{\mathrm{Robba}^\mathrm{extended}_{{R_k,A,I}}\}_I) \ar[d]\ar[d]\ar[d]\ar[d].\\
\mathrm{Ind}\mathrm{Banach}(\mathcal{O}_{X_{\mathbb{Q}_p(p^{1/p^\infty})^{\wedge\flat},A}})\ar[r]^{\mathrm{global}}\ar[r]\ar[r] &\varphi\mathrm{Ind}\mathrm{Banach}(\{\mathrm{Robba}^\mathrm{extended}_{{R_0,A,I}}\}_I).\\ 
}
\]
\item (\text{Proposition}) There is a functor (global section) between the $\infty$-categories of monomorphic inductive Banach quasicoherent presheaves:
\[
\xymatrix@R+6pc@C+0pc{
\mathrm{Ind}^m\mathrm{Banach}_{\Gamma_k}(\mathcal{O}_{X_{R_k,A}})\ar[r]^{\mathrm{global}}\ar[d]\ar[d]\ar[d]\ar[d]\ar[r]\ar[r] &\varphi\mathrm{Ind}^m\mathrm{Banach}_{\Gamma_k}(\{\mathrm{Robba}^\mathrm{extended}_{{R_k,A,I}}\}_I)\ar[d]\ar[d]\ar[d]\ar[d]\\
\mathrm{Ind}^m\mathrm{Banach}_{\Gamma_0}(\mathcal{O}_{X_{\mathbb{Q}_p(p^{1/p^\infty})^{\wedge\flat},A}})\ar[r]^{\mathrm{global}}\ar[r]\ar[r] &\varphi\mathrm{Ind}^m\mathrm{Banach}_{\Gamma_0}(\{\mathrm{Robba}^\mathrm{extended}_{{R_0,A,I}}\}_I).\\  
}
\]
\item (\text{Proposition}) There is a functor (global section) between the $\infty$-categories of inductive Banach quasicoherent commutative algebra $E_\infty$ objects:
\[\displayindent=-0.4in
\xymatrix@R+6pc@C+0pc{
\mathrm{sComm}_\mathrm{simplicial}\mathrm{Ind}\mathrm{Banach}_{\Gamma_k}(\mathcal{O}_{X_{R_k,A}})\ar[d]\ar[d]\ar[d]\ar[d]\ar[r]^{\mathrm{global}}\ar[r]\ar[r] &\mathrm{sComm}_\mathrm{simplicial}\varphi\mathrm{Ind}\mathrm{Banach}_{\Gamma_k}(\{\mathrm{Robba}^\mathrm{extended}_{{R_k,A,I}}\}_I)\ar[d]\ar[d]\ar[d]\ar[d]\\
\mathrm{sComm}_\mathrm{simplicial}\mathrm{Ind}\mathrm{Banach}_{\Gamma_0}(\mathcal{O}_{X_{\mathbb{Q}_p(p^{1/p^\infty})^{\wedge\flat},A}})\ar[r]^{\mathrm{global}}\ar[r]\ar[r] &\mathrm{sComm}_\mathrm{simplicial}\varphi\mathrm{Ind}\mathrm{Banach}_{\Gamma_0}(\{\mathrm{Robba}^\mathrm{extended}_{{R_0,A,I}}\}_I).  
}
\]
\item (\text{Proposition}) There is a functor (global section) between the $\infty$-categories of monomorphic inductive Banach quasicoherent commutative algebra $E_\infty$ objects:
\[\displayindent=-0.4in
\xymatrix@R+6pc@C+0pc{
\mathrm{sComm}_\mathrm{simplicial}\mathrm{Ind}^m\mathrm{Banach}_{\Gamma_k}(\mathcal{O}_{X_{R_k,A}})\ar[d]\ar[d]\ar[d]\ar[d]\ar[r]^{\mathrm{equi}}\ar[r]\ar[r] &\mathrm{sComm}_\mathrm{simplicial}\varphi\mathrm{Ind}^m\mathrm{Banach}_{\Gamma_k}(\{\mathrm{Robba}^\mathrm{extended}_{{R_k,A,I}}\}_I)\ar[d]\ar[d]\ar[d]\ar[d]\\
 \mathrm{sComm}_\mathrm{simplicial}\mathrm{Ind}^m\mathrm{Banach}_{\Gamma_0}(\mathcal{O}_{X_{\mathbb{Q}_p(p^{1/p^\infty})^{\wedge\flat},A}})\ar[r]^{\mathrm{equi}}\ar[r]\ar[r] &\mathrm{sComm}_\mathrm{simplicial}\varphi\mathrm{Ind}^m\mathrm{Banach}_{\Gamma_0}(\{\mathrm{Robba}^\mathrm{extended}_{{R_0,A,I}}\}_I). 
}
\]

\item Then parallel as in \cite{LBV} we have a functor (global section) of the de Rham complex after \cite[Definition 5.9, Section 5.2.1]{KKM}:
\[\displayindent=-0.4in
\xymatrix@R+6pc@C+0pc{
\mathrm{deRham}_{\mathrm{sComm}_\mathrm{simplicial}\mathrm{Ind}\mathrm{Banach}_{\Gamma_k}(\mathcal{O}_{X_{R_k,A}})\ar[r]^{\mathrm{global}}}(-)\ar[d]\ar[d]\ar[d]\ar[d]\ar[r]\ar[r] &\mathrm{deRham}_{\mathrm{sComm}_\mathrm{simplicial}\varphi\mathrm{Ind}\mathrm{Banach}_{\Gamma_k}(\{\mathrm{Robba}^\mathrm{extended}_{{R_k,A,I}}\}_I)}(-)\ar[d]\ar[d]\ar[d]\ar[d]\\
\mathrm{deRham}_{\mathrm{sComm}_\mathrm{simplicial}\mathrm{Ind}\mathrm{Banach}_{\Gamma_0}(\mathcal{O}_{X_{\mathbb{Q}_p(p^{1/p^\infty})^{\wedge\flat},A}})\ar[r]^{\mathrm{global}}}(-)\ar[r]\ar[r] &\mathrm{deRham}_{\mathrm{sComm}_\mathrm{simplicial}\varphi\mathrm{Ind}\mathrm{Banach}_{\Gamma_0}(\{\mathrm{Robba}^\mathrm{extended}_{{R_0,A,I}}\}_I)}(-), 
}
\]
\[\displayindent=-0.4in
\xymatrix@R+6pc@C+0pc{
\mathrm{deRham}_{\mathrm{sComm}_\mathrm{simplicial}\mathrm{Ind}^m\mathrm{Banach}_{\Gamma_k}(\mathcal{O}_{X_{R_k,A}})\ar[r]^{\mathrm{global}}}(-)\ar[d]\ar[d]\ar[d]\ar[d]\ar[r]\ar[r] &\mathrm{deRham}_{\mathrm{sComm}_\mathrm{simplicial}\varphi\mathrm{Ind}^m\mathrm{Banach}_{\Gamma_k}(\{\mathrm{Robba}^\mathrm{extended}_{{R_k,A,I}}\}_I)}(-)\ar[d]\ar[d]\ar[d]\ar[d]\\
\mathrm{deRham}_{\mathrm{sComm}_\mathrm{simplicial}\mathrm{Ind}^m\mathrm{Banach}_{\Gamma_0}(\mathcal{O}_{X_{\mathbb{Q}_p(p^{1/p^\infty})^{\wedge\flat},A}})\ar[r]^{\mathrm{global}}}(-)\ar[r]\ar[r] &\mathrm{deRham}_{\mathrm{sComm}_\mathrm{simplicial}\varphi\mathrm{Ind}^m\mathrm{Banach}_{\Gamma_0}(\{\mathrm{Robba}^\mathrm{extended}_{{R_0,A,I}}\}_I)}(-). 
}
\]

\item Then we have the following a functor (global section) of $K$-group $(\infty,1)$-spectrum from \cite{BGT}:
\[
\xymatrix@R+6pc@C+0pc{
\mathrm{K}^\mathrm{BGT}_{\mathrm{sComm}_\mathrm{simplicial}\mathrm{Ind}\mathrm{Banach}_{\Gamma_k}(\mathcal{O}_{X_{R_k,A}})\ar[r]^{\mathrm{global}}}(-)\ar[d]\ar[d]\ar[d]\ar[d]\ar[r]\ar[r] &\mathrm{K}^\mathrm{BGT}_{\mathrm{sComm}_\mathrm{simplicial}\varphi\mathrm{Ind}\mathrm{Banach}_{\Gamma_k}(\{\mathrm{Robba}^\mathrm{extended}_{{R_k,A,I}}\}_I)}(-)\ar[d]\ar[d]\ar[d]\ar[d]\\
\mathrm{K}^\mathrm{BGT}_{\mathrm{sComm}_\mathrm{simplicial}\mathrm{Ind}\mathrm{Banach}_{\Gamma_0}(\mathcal{O}_{X_{\mathbb{Q}_p(p^{1/p^\infty})^{\wedge\flat},A}})\ar[r]^{\mathrm{global}}}(-)\ar[r]\ar[r] &\mathrm{K}^\mathrm{BGT}_{\mathrm{sComm}_\mathrm{simplicial}\varphi\mathrm{Ind}\mathrm{Banach}_{\Gamma_0}(\{\mathrm{Robba}^\mathrm{extended}_{{R_0,A,I}}\}_I)}(-), 
}
\]
\[
\xymatrix@R+6pc@C+0pc{
\mathrm{K}^\mathrm{BGT}_{\mathrm{sComm}_\mathrm{simplicial}\mathrm{Ind}^m\mathrm{Banach}_{\Gamma_k}(\mathcal{O}_{X_{R_k,A}})\ar[r]^{\mathrm{global}}}(-)\ar[d]\ar[d]\ar[d]\ar[d]\ar[r]\ar[r] &\mathrm{K}^\mathrm{BGT}_{\mathrm{sComm}_\mathrm{simplicial}\varphi\mathrm{Ind}^m\mathrm{Banach}_{\Gamma_k}(\{\mathrm{Robba}^\mathrm{extended}_{{R_k,A,I}}\}_I)}(-)\ar[d]\ar[d]\ar[d]\ar[d]\\
\mathrm{K}^\mathrm{BGT}_{\mathrm{sComm}_\mathrm{simplicial}\mathrm{Ind}^m\mathrm{Banach}_{\Gamma_0}(\mathcal{O}_{X_{\mathbb{Q}_p(p^{1/p^\infty})^{\wedge\flat},A}})\ar[r]^{\mathrm{global}}}(-)\ar[r]\ar[r] &\mathrm{K}^\mathrm{BGT}_{\mathrm{sComm}_\mathrm{simplicial}\varphi\mathrm{Ind}^m\mathrm{Banach}_{\Gamma_0}(\{\mathrm{Robba}^\mathrm{extended}_{{R_0,A,I}}\}_I)}(-). 
}
\]

\end{itemize}

\

\

\begin{remark}
\noindent We can certainly consider the quasicoherent sheaves in \cite[Lemma 7.11, Remark 7.12]{1BBK}, therefore all the quasicoherent presheaves and modules will be those satisfying \cite[Lemma 7.11, Remark 7.12]{1BBK} if one would like to consider the the quasicoherent sheaves. That being all as this said, we would believe that the big quasicoherent presheaves are automatically quasicoherent sheaves (namely satisfying the corresponding \v{C}ech $\infty$-descent as in \cite[Section 9.3]{KKM} and \cite[Lemma 7.11, Remark 7.12]{1BBK}) and the corresponding global section functors are automatically equivalence of $\infty$-categories. 
\end{remark}

\

\indent In Clausen-Scholze formalism we have the following:

\begin{itemize}

\item (\text{Proposition}) There is a functor (global section) between the $\infty$-categories of inductive Banach quasicoherent sheaves:
\[
\xymatrix@R+0pc@C+0pc{
\mathrm{Module}_\circledcirc(\mathcal{O}_{X_{R,A}})\ar[r]^{\mathrm{global}}\ar[r]\ar[r] &\varphi\mathrm{Module}_\circledcirc(\{\mathrm{Robba}^\mathrm{extended}_{{R,A,I}}\}_I).  
}
\]

\item (\text{Proposition}) There is a functor (global section) between the $\infty$-categories of inductive Banach quasicoherent sheaves:
\[
\xymatrix@R+0pc@C+0pc{
\mathrm{Module}_\circledcirc(\mathcal{O}_{X_{R,A}})\ar[r]^{\mathrm{global}}\ar[r]\ar[r] &\varphi\mathrm{Module}_\circledcirc(\{\mathrm{Robba}^\mathrm{extended}_{{R,A,I}}\}_I).  
}
\]

\item (\text{Proposition}) There is a functor (global section) between the $\infty$-categories of inductive Banach quasicoherent commutative algebra $E_\infty$ objects\footnote{Here $\circledcirc=\text{solidquasicoherentsheaves}$.}:
\[
\xymatrix@R+0pc@C+0pc{
\mathrm{sComm}_\mathrm{simplicial}\mathrm{Module}_\circledcirc(\mathcal{O}_{X_{R,A}})\ar[r]^{\mathrm{global}}\ar[r]\ar[r] &\mathrm{sComm}_\mathrm{simplicial}\varphi\mathrm{Module}_\circledcirc(\{\mathrm{Robba}^\mathrm{extended}_{{R,A,I}}\}_I).  
}
\]

\item Then as in \cite{LBV} we have a functor (global section) of the de Rham complex after \cite[Definition 5.9, Section 5.2.1]{KKM}\footnote{Here $\circledcirc=\text{solidquasicoherentsheaves}$.}:
\[
\xymatrix@R+0pc@C+0pc{
\mathrm{deRham}_{\mathrm{sComm}_\mathrm{simplicial}\mathrm{Module}_\circledcirc(\mathcal{O}_{X_{R,A}})\ar[r]^{\mathrm{global}}}(-)\ar[r]\ar[r] &\mathrm{deRham}_{\mathrm{sComm}_\mathrm{simplicial}\varphi\mathrm{Module}_\circledcirc(\{\mathrm{Robba}^\mathrm{extended}_{{R,A,I}}\}_I)}(-). 
}
\]

\item Then we have the following a functor (global section) of $K$-group $(\infty,1)$-spectrum from \cite{BGT}\footnote{Here $\circledcirc=\text{solidquasicoherentsheaves}$.}:
\[
\xymatrix@R+0pc@C+0pc{
\mathrm{K}^\mathrm{BGT}_{\mathrm{sComm}_\mathrm{simplicial}\mathrm{Module}_\circledcirc(\mathcal{O}_{X_{R,A}})\ar[r]^{\mathrm{global}}}(-)\ar[r]\ar[r] &\mathrm{K}^\mathrm{BGT}_{\mathrm{sComm}_\mathrm{simplicial}\varphi\mathrm{Module}_\circledcirc(\{\mathrm{Robba}^\mathrm{extended}_{{R,A,I}}\}_I)}(-). 
}
\]
\end{itemize}

\noindent Now let $R=\mathbb{Q}_p(p^{1/p^\infty})^{\wedge\flat}$ and $R_k=\mathbb{Q}_p(p^{1/p^\infty})^{\wedge}\left<T_1^{\pm 1/p^{\infty}},...,T_k^{\pm 1/p^{\infty}}\right>^\flat$ we have the following Galois theoretic results with naturality along $f:\mathrm{Spa}(\mathbb{Q}_p(p^{1/p^\infty})^{\wedge}\left<T_1^{\pm 1/p^\infty},...,T_k^{\pm 1/p^\infty}\right>^\flat)\rightarrow \mathrm{Spa}(\mathbb{Q}_p(p^{1/p^\infty})^{\wedge\flat})$:

\begin{itemize}
\item (\text{Proposition}) There is a functor (global section) between the $\infty$-categories of inductive Banach quasicoherent sheaves\footnote{Here $\circledcirc=\text{solidquasicoherentsheaves}$.}:
\[
\xymatrix@R+6pc@C+0pc{
\mathrm{Module}_\circledcirc(\mathcal{O}_{X_{\mathbb{Q}_p(p^{1/p^\infty})^{\wedge}\left<T_1^{\pm 1/p^\infty},...,T_k^{\pm 1/p^\infty}\right>^\flat,A}})\ar[d]\ar[d]\ar[d]\ar[d] \ar[r]^{\mathrm{global}}\ar[r]\ar[r] &\varphi\mathrm{Module}_\circledcirc(\{\mathrm{Robba}^\mathrm{extended}_{{R_k,A,I}}\}_I) \ar[d]\ar[d]\ar[d]\ar[d].\\
\mathrm{Module}_\circledcirc(\mathcal{O}_{X_{\mathbb{Q}_p(p^{1/p^\infty})^{\wedge\flat},A}})\ar[r]^{\mathrm{global}}\ar[r]\ar[r] &\varphi\mathrm{Module}_\circledcirc(\{\mathrm{Robba}^\mathrm{extended}_{{R_0,A,I}}\}_I).\\ 
}
\]

\item (\text{Proposition}) There is a functor (global section) between the $\infty$-categories of inductive Banach quasicoherent commutative algebra $E_\infty$ objects\footnote{Here $\circledcirc=\text{solidquasicoherentsheaves}$.}:
\[
\xymatrix@R+6pc@C+0pc{
\mathrm{sComm}_\mathrm{simplicial}\mathrm{Module}_\circledcirc(\mathcal{O}_{X_{R_k,A}})\ar[d]\ar[d]\ar[d]\ar[d]\ar[r]^{\mathrm{global}}\ar[r]\ar[r] &\mathrm{sComm}_\mathrm{simplicial}\varphi\mathrm{Module}_\circledcirc(\{\mathrm{Robba}^\mathrm{extended}_{{R_k,A,I}}\}_I)\ar[d]\ar[d]\ar[d]\ar[d]\\
\mathrm{sComm}_\mathrm{simplicial}\mathrm{Module}_\circledcirc(\mathcal{O}_{X_{\mathbb{Q}_p(p^{1/p^\infty})^{\wedge\flat},A}})\ar[r]^{\mathrm{global}}\ar[r]\ar[r] &\mathrm{sComm}_\mathrm{simplicial}\varphi\mathrm{Module}_\circledcirc(\{\mathrm{Robba}^\mathrm{extended}_{{R_0,A,I}}\}_I).  
}
\]

\item Then as in \cite{LBV} we have a functor (global section) of the de Rham complex after \cite[Definition 5.9, Section 5.2.1]{KKM}\footnote{Here $\circledcirc=\text{solidquasicoherentsheaves}$.}:
\[
\xymatrix@R+6pc@C+0pc{
\mathrm{deRham}_{\mathrm{sComm}_\mathrm{simplicial}\mathrm{Module}_\circledcirc(\mathcal{O}_{X_{R_k,A}})\ar[r]^{\mathrm{global}}}(-)\ar[d]\ar[d]\ar[d]\ar[d]\ar[r]\ar[r] &\mathrm{deRham}_{\mathrm{sComm}_\mathrm{simplicial}\varphi\mathrm{Module}_\circledcirc(\{\mathrm{Robba}^\mathrm{extended}_{{R_k,A,I}}\}_I)}(-)\ar[d]\ar[d]\ar[d]\ar[d]\\
\mathrm{deRham}_{\mathrm{sComm}_\mathrm{simplicial}\mathrm{Module}_\circledcirc(\mathcal{O}_{X_{\mathbb{Q}_p(p^{1/p^\infty})^{\wedge\flat},A}})\ar[r]^{\mathrm{global}}}(-)\ar[r]\ar[r] &\mathrm{deRham}_{\mathrm{sComm}_\mathrm{simplicial}\varphi\mathrm{Module}_\circledcirc(\{\mathrm{Robba}^\mathrm{extended}_{{R_0,A,I}}\}_I)}(-). 
}
\]

\item Then we have the following a functor (global section) of $K$-group $(\infty,1)$-spectrum from \cite{BGT}\footnote{Here $\circledcirc=\text{solidquasicoherentsheaves}$.}:
\[
\xymatrix@R+6pc@C+0pc{
\mathrm{K}^\mathrm{BGT}_{\mathrm{sComm}_\mathrm{simplicial}\mathrm{Module}_\circledcirc(\mathcal{O}_{X_{R_k,A}})\ar[r]^{\mathrm{global}}}(-)\ar[d]\ar[d]\ar[d]\ar[d]\ar[r]\ar[r] &\mathrm{K}^\mathrm{BGT}_{\mathrm{sComm}_\mathrm{simplicial}\varphi\mathrm{Module}_\circledcirc(\{\mathrm{Robba}^\mathrm{extended}_{{R_k,A,I}}\}_I)}(-)\ar[d]\ar[d]\ar[d]\ar[d]\\
\mathrm{K}^\mathrm{BGT}_{\mathrm{sComm}_\mathrm{simplicial}\mathrm{Module}_\circledcirc(\mathcal{O}_{X_{\mathbb{Q}_p(p^{1/p^\infty})^{\wedge\flat},A}})\ar[r]^{\mathrm{global}}}(-)\ar[r]\ar[r] &\mathrm{K}^\mathrm{BGT}_{\mathrm{sComm}_\mathrm{simplicial}\varphi\mathrm{Module}_\circledcirc(\{\mathrm{Robba}^\mathrm{extended}_{{R_0,A,I}}\}_I)}(-). 
}
\]

\end{itemize}

\
\indent Then we consider further equivariance by considering the arithmetic profinite fundamental groups $\Gamma_{\mathbb{Q}_p}$ and $\mathrm{Gal}(\overline{\mathbb{Q}_p\left<T_1^{\pm 1},...,T_k^{\pm 1}\right>}/R_k)$ through the following diagram:

\[
\xymatrix@R+0pc@C+0pc{
\mathbb{Z}_p^k=\mathrm{Gal}(R_k/{\mathbb{Q}_p(p^{1/p^\infty})^\wedge\left<T_1^{\pm 1},...,T_k^{\pm 1}\right>}) \ar[r]\ar[r] \ar[r]\ar[r] &\Gamma_k:=\mathrm{Gal}(R_k/{\mathbb{Q}_p\left<T_1^{\pm 1},...,T_k^{\pm 1}\right>}) \ar[r] \ar[r]\ar[r] &\Gamma_{\mathbb{Q}_p}.
}
\]

\begin{itemize}
\item (\text{Proposition}) There is a functor (global section) between the $\infty$-categories of inductive Banach quasicoherent sheaves\footnote{Here $\circledcirc=\text{solidquasicoherentsheaves}$.}:
\[
\xymatrix@R+6pc@C+0pc{
{\mathrm{Module}_\circledcirc}_{\Gamma_k}(\mathcal{O}_{X_{\mathbb{Q}_p(p^{1/p^\infty})^{\wedge}\left<T_1^{\pm 1/p^\infty},...,T_k^{\pm 1/p^\infty}\right>^\flat,A}})\ar[d]\ar[d]\ar[d]\ar[d] \ar[r]^{\mathrm{global}}\ar[r]\ar[r] &\varphi{\mathrm{Module}_\circledcirc}_{\Gamma_k}(\{\mathrm{Robba}^\mathrm{extended}_{{R_k,A,I}}\}_I) \ar[d]\ar[d]\ar[d]\ar[d].\\
{\mathrm{Module}_\circledcirc}(\mathcal{O}_{X_{\mathbb{Q}_p(p^{1/p^\infty})^{\wedge\flat},A}})\ar[r]^{\mathrm{global}}\ar[r]\ar[r] &\varphi{\mathrm{Module}_\circledcirc}(\{\mathrm{Robba}^\mathrm{extended}_{{R_0,A,I}}\}_I).\\ 
}
\]

\item (\text{Proposition}) There is a functor (global section) between the $\infty$-categories of inductive Banach quasicoherent commutative algebra $E_\infty$ objects\footnote{Here $\circledcirc=\text{solidquasicoherentsheaves}$.}:
\[
\xymatrix@R+6pc@C+0pc{
\mathrm{sComm}_\mathrm{simplicial}{\mathrm{Module}_\circledcirc}_{\Gamma_k}(\mathcal{O}_{X_{R_k,A}})\ar[d]\ar[d]\ar[d]\ar[d]\ar[r]^{\mathrm{global}}\ar[r]\ar[r] &\mathrm{sComm}_\mathrm{simplicial}\varphi{\mathrm{Module}_\circledcirc}(\{\mathrm{Robba}^\mathrm{extended}_{{R_k,A,I}}\}_I)\ar[d]\ar[d]\ar[d]\ar[d]\\
\mathrm{sComm}_\mathrm{simplicial}{\mathrm{Module}_\circledcirc}_{\Gamma_0}(\mathcal{O}_{X_{\mathbb{Q}_p(p^{1/p^\infty})^{\wedge\flat},A}})\ar[r]^{\mathrm{global}}\ar[r]\ar[r] &\mathrm{sComm}_\mathrm{simplicial}\varphi{\mathrm{Modules}_\circledcirc}_{\Gamma_0}(\{\mathrm{Robba}^\mathrm{extended}_{{R_0,A,I}}\}_I).  
}
\]

\item Then as in \cite{LBV} we have a functor (global section) of the de Rham complex after \cite[Definition 5.9, Section 5.2.1]{KKM}\footnote{Here $\circledcirc=\text{solidquasicoherentsheaves}$.}:
\[\displayindent=-0.4in
\xymatrix@R+6pc@C+0pc{
\mathrm{deRham}_{\mathrm{sComm}_\mathrm{simplicial}{\mathrm{Modules}_\circledcirc}_{\Gamma_k}(\mathcal{O}_{X_{R_k,A}})\ar[r]^{\mathrm{global}}}(-)\ar[d]\ar[d]\ar[d]\ar[d]\ar[r]\ar[r] &\mathrm{deRham}_{\mathrm{sComm}_\mathrm{simplicial}\varphi{\mathrm{Modules}_\circledcirc}_{\Gamma_k}(\{\mathrm{Robba}^\mathrm{extended}_{{R_k,A,I}}\}_I)}(-)\ar[d]\ar[d]\ar[d]\ar[d]\\
\mathrm{deRham}_{\mathrm{sComm}_\mathrm{simplicial}{\mathrm{Modules}_\circledcirc}_{\Gamma_0}(\mathcal{O}_{X_{\mathbb{Q}_p(p^{1/p^\infty})^{\wedge\flat},A}})\ar[r]^{\mathrm{global}}}(-)\ar[r]\ar[r] &\mathrm{deRham}_{\mathrm{sComm}_\mathrm{simplicial}\varphi{\mathrm{Modules}_\circledcirc}_{\Gamma_0}(\{\mathrm{Robba}^\mathrm{extended}_{{R_0,A,I}}\}_I)}(-). 
}
\]

\item Then we have the following a functor (global section) of $K$-group $(\infty,1)$-spectrum from \cite{BGT}\footnote{Here $\circledcirc=\text{solidquasicoherentsheaves}$.}:
\[
\xymatrix@R+6pc@C+0pc{
\mathrm{K}^\mathrm{BGT}_{\mathrm{sComm}_\mathrm{simplicial}{\mathrm{Modules}_\circledcirc}_{\Gamma_k}(\mathcal{O}_{X_{R_k,A}})\ar[r]^{\mathrm{global}}}(-)\ar[d]\ar[d]\ar[d]\ar[d]\ar[r]\ar[r] &\mathrm{K}^\mathrm{BGT}_{\mathrm{sComm}_\mathrm{simplicial}\varphi{\mathrm{Modules}_\circledcirc}_{\Gamma_k}(\{\mathrm{Robba}^\mathrm{extended}_{{R_k,A,I}}\}_I)}(-)\ar[d]\ar[d]\ar[d]\ar[d]\\
\mathrm{K}^\mathrm{BGT}_{\mathrm{sComm}_\mathrm{simplicial}{\mathrm{Modules}_\circledcirc}_{\Gamma_0}(\mathcal{O}_{X_{\mathbb{Q}_p(p^{1/p^\infty})^{\wedge\flat},A}})\ar[r]^{\mathrm{global}}}(-)\ar[r]\ar[r] &\mathrm{K}^\mathrm{BGT}_{\mathrm{sComm}_\mathrm{simplicial}\varphi{\mathrm{Modules}_\circledcirc}_{\Gamma_0}(\{\mathrm{Robba}^\mathrm{extended}_{{R_0,A,I}}\}_I)}(-). 
}
\]
	
\end{itemize}

\begin{proposition}
All the global functors from \cite[Proposition 13.8, Theorem 14.9, Remark 14.10]{1CS2} achieve the equivalences on both sides.	
\end{proposition}

\newpage
\subsection{$\infty$-Categorical Analytic Stacks and Descents III}

\indent In the following the right had of each row in each diagram will be the corresponding quasicoherent Robba bundles over the Robba ring carrying the corresponding action from the Frobenius or the fundamental groups, defined by directly applying \cite[Section 9.3]{KKM} and \cite{BBM}. We now let $\mathcal{A}$ be any commutative algebra objects in the corresponding $\infty$-toposes over ind-Banach commutative algebra objects over $\mathbb{Q}_p$ or the corresponding born\'e commutative algebra objects over $\mathbb{Q}_p$ carrying the Grothendieck topology defined by essentially the corresponding monomorphism homotopy in the opposite category. Then we promote the construction to the corresponding $\infty$-stack over the same $\infty$-categories of affinoids. Let $\mathcal{A}$ vary in the category of all the Banach algebras over $\mathbb{Q}_p$ we have the following.

\begin{itemize}

\item (\text{Proposition}) There is a functor (global section) between the $\infty$-prestacks of inductive Banach quasicoherent presheaves:
\[
\xymatrix@R+0pc@C+0pc{
\mathrm{Ind}\mathrm{Banach}(\mathcal{O}_{X_{R,-}})\ar[r]^{\mathrm{global}}\ar[r]\ar[r] &\varphi\mathrm{Ind}\mathrm{Banach}(\{\mathrm{Robba}^\mathrm{extended}_{{R,-,I}}\}_I).  
}
\]
\item (\text{Proposition}) There is a functor (global section) between the $\infty$-prestacks of monomorphic inductive Banach quasicoherent presheaves:
\[
\xymatrix@R+0pc@C+0pc{
\mathrm{Ind}^m\mathrm{Banach}(\mathcal{O}_{X_{R,-}})\ar[r]^{\mathrm{global}}\ar[r]\ar[r] &\varphi\mathrm{Ind}^m\mathrm{Banach}(\{\mathrm{Robba}^\mathrm{extended}_{{R,-,I}}\}_I).  
}
\]

\item (\text{Proposition}) There is a functor (global section) between the $\infty$-prestacks of inductive Banach quasicoherent presheaves:
\[
\xymatrix@R+0pc@C+0pc{
\mathrm{Ind}\mathrm{Banach}(\mathcal{O}_{X_{R,-}})\ar[r]^{\mathrm{global}}\ar[r]\ar[r] &\varphi\mathrm{Ind}\mathrm{Banach}(\{\mathrm{Robba}^\mathrm{extended}_{{R,-,I}}\}_I).  
}
\]
\item (\text{Proposition}) There is a functor (global section) between the $\infty$-stacks of monomorphic inductive Banach quasicoherent presheaves:
\[
\xymatrix@R+0pc@C+0pc{
\mathrm{Ind}^m\mathrm{Banach}(\mathcal{O}_{X_{R,-}})\ar[r]^{\mathrm{global}}\ar[r]\ar[r] &\varphi\mathrm{Ind}^m\mathrm{Banach}(\{\mathrm{Robba}^\mathrm{extended}_{{R,-,I}}\}_I).  
}
\]
\item (\text{Proposition}) There is a functor (global section) between the $\infty$-prestacks of inductive Banach quasicoherent commutative algebra $E_\infty$ objects:
\[
\xymatrix@R+0pc@C+0pc{
\mathrm{sComm}_\mathrm{simplicial}\mathrm{Ind}\mathrm{Banach}(\mathcal{O}_{X_{R,-}})\ar[r]^{\mathrm{global}}\ar[r]\ar[r] &\mathrm{sComm}_\mathrm{simplicial}\varphi\mathrm{Ind}\mathrm{Banach}(\{\mathrm{Robba}^\mathrm{extended}_{{R,-,I}}\}_I).  
}
\]
\item (\text{Proposition}) There is a functor (global section) between the $\infty$-prestacks of monomorphic inductive Banach quasicoherent commutative algebra $E_\infty$ objects:
\[
\xymatrix@R+0pc@C+0pc{
\mathrm{sComm}_\mathrm{simplicial}\mathrm{Ind}^m\mathrm{Banach}(\mathcal{O}_{X_{R,-}})\ar[r]^{\mathrm{global}}\ar[r]\ar[r] &\mathrm{sComm}_\mathrm{simplicial}\varphi\mathrm{Ind}^m\mathrm{Banach}(\{\mathrm{Robba}^\mathrm{extended}_{{R,-,I}}\}_I).  
}
\]

\item Then parallel as in \cite{LBV} we have a functor (global section) of the de Rham complex after \cite[Definition 5.9, Section 5.2.1]{KKM}:
\[
\xymatrix@R+0pc@C+0pc{
\mathrm{deRham}_{\mathrm{sComm}_\mathrm{simplicial}\mathrm{Ind}\mathrm{Banach}(\mathcal{O}_{X_{R,-}})\ar[r]^{\mathrm{global}}}(-)\ar[r]\ar[r] &\mathrm{deRham}_{\mathrm{sComm}_\mathrm{simplicial}\varphi\mathrm{Ind}\mathrm{Banach}(\{\mathrm{Robba}^\mathrm{extended}_{{R,-,I}}\}_I)}(-), 
}
\]
\[
\xymatrix@R+0pc@C+0pc{
\mathrm{deRham}_{\mathrm{sComm}_\mathrm{simplicial}\mathrm{Ind}^m\mathrm{Banach}(\mathcal{O}_{X_{R,-}})\ar[r]^{\mathrm{global}}}(-)\ar[r]\ar[r] &\mathrm{deRham}_{\mathrm{sComm}_\mathrm{simplicial}\varphi\mathrm{Ind}^m\mathrm{Banach}(\{\mathrm{Robba}^\mathrm{extended}_{{R,-,I}}\}_I)}(-). 
}
\]

\item Then we have the following a functor (global section) of $K$-group $(\infty,1)$-spectrum from \cite{BGT}:
\[
\xymatrix@R+0pc@C+0pc{
\mathrm{K}^\mathrm{BGT}_{\mathrm{sComm}_\mathrm{simplicial}\mathrm{Ind}\mathrm{Banach}(\mathcal{O}_{X_{R,-}})\ar[r]^{\mathrm{global}}}(-)\ar[r]\ar[r] &\mathrm{K}^\mathrm{BGT}_{\mathrm{sComm}_\mathrm{simplicial}\varphi\mathrm{Ind}\mathrm{Banach}(\{\mathrm{Robba}^\mathrm{extended}_{{R,-,I}}\}_I)}(-), 
}
\]
\[
\xymatrix@R+0pc@C+0pc{
\mathrm{K}^\mathrm{BGT}_{\mathrm{sComm}_\mathrm{simplicial}\mathrm{Ind}^m\mathrm{Banach}(\mathcal{O}_{X_{R,-}})\ar[r]^{\mathrm{global}}}(-)\ar[r]\ar[r] &\mathrm{K}^\mathrm{BGT}_{\mathrm{sComm}_\mathrm{simplicial}\varphi\mathrm{Ind}^m\mathrm{Banach}(\{\mathrm{Robba}^\mathrm{extended}_{{R,-,I}}\}_I)}(-). 
}
\]
\end{itemize}

\noindent Now let $R=\mathbb{Q}_p(p^{1/p^\infty})^{\wedge\flat}$ and $R_k=\mathbb{Q}_p(p^{1/p^\infty})^{\wedge}\left<T_1^{\pm 1/p^{\infty}},...,T_k^{\pm 1/p^{\infty}}\right>^\flat$ we have the following Galois theoretic results with naturality along $f:\mathrm{Spa}(\mathbb{Q}_p(p^{1/p^\infty})^{\wedge}\left<T_1^{\pm 1/p^\infty},...,T_k^{\pm 1/p^\infty}\right>^\flat)\rightarrow \mathrm{Spa}(\mathbb{Q}_p(p^{1/p^\infty})^{\wedge\flat})$:

\begin{itemize}
\item (\text{Proposition}) There is a functor (global section) between the $\infty$-prestacks of inductive Banach quasicoherent presheaves:
\[
\xymatrix@R+6pc@C+0pc{
\mathrm{Ind}\mathrm{Banach}(\mathcal{O}_{X_{\mathbb{Q}_p(p^{1/p^\infty})^{\wedge}\left<T_1^{\pm 1/p^\infty},...,T_k^{\pm 1/p^\infty}\right>^\flat,-}})\ar[d]\ar[d]\ar[d]\ar[d] \ar[r]^{\mathrm{global}}\ar[r]\ar[r] &\varphi\mathrm{Ind}\mathrm{Banach}(\{\mathrm{Robba}^\mathrm{extended}_{{R_k,-,I}}\}_I) \ar[d]\ar[d]\ar[d]\ar[d].\\
\mathrm{Ind}\mathrm{Banach}(\mathcal{O}_{X_{\mathbb{Q}_p(p^{1/p^\infty})^{\wedge\flat},-}})\ar[r]^{\mathrm{global}}\ar[r]\ar[r] &\varphi\mathrm{Ind}\mathrm{Banach}(\{\mathrm{Robba}^\mathrm{extended}_{{R_0,-,I}}\}_I).\\ 
}
\]
\item (\text{Proposition}) There is a functor (global section) between the $\infty$-prestacks of monomorphic inductive Banach quasicoherent presheaves:
\[
\xymatrix@R+6pc@C+0pc{
\mathrm{Ind}^m\mathrm{Banach}(\mathcal{O}_{X_{R_k,-}})\ar[r]^{\mathrm{global}}\ar[d]\ar[d]\ar[d]\ar[d]\ar[r]\ar[r] &\varphi\mathrm{Ind}^m\mathrm{Banach}(\{\mathrm{Robba}^\mathrm{extended}_{{R_k,-,I}}\}_I)\ar[d]\ar[d]\ar[d]\ar[d]\\
\mathrm{Ind}^m\mathrm{Banach}(\mathcal{O}_{X_{\mathbb{Q}_p(p^{1/p^\infty})^{\wedge\flat},-}})\ar[r]^{\mathrm{global}}\ar[r]\ar[r] &\varphi\mathrm{Ind}^m\mathrm{Banach}(\{\mathrm{Robba}^\mathrm{extended}_{{R_0,-,I}}\}_I).\\  
}
\]
\item (\text{Proposition}) There is a functor (global section) between the $\infty$-prestacks of inductive Banach quasicoherent commutative algebra $E_\infty$ objects:
\[
\xymatrix@R+6pc@C+0pc{
\mathrm{sComm}_\mathrm{simplicial}\mathrm{Ind}\mathrm{Banach}(\mathcal{O}_{X_{R_k,-}})\ar[d]\ar[d]\ar[d]\ar[d]\ar[r]^{\mathrm{global}}\ar[r]\ar[r] &\mathrm{sComm}_\mathrm{simplicial}\varphi\mathrm{Ind}\mathrm{Banach}(\{\mathrm{Robba}^\mathrm{extended}_{{R_k,-,I}}\}_I)\ar[d]\ar[d]\ar[d]\ar[d]\\
\mathrm{sComm}_\mathrm{simplicial}\mathrm{Ind}\mathrm{Banach}(\mathcal{O}_{X_{\mathbb{Q}_p(p^{1/p^\infty})^{\wedge\flat},-}})\ar[r]^{\mathrm{global}}\ar[r]\ar[r] &\mathrm{sComm}_\mathrm{simplicial}\varphi\mathrm{Ind}\mathrm{Banach}(\{\mathrm{Robba}^\mathrm{extended}_{{R_0,-,I}}\}_I).  
}
\]
\item (\text{Proposition}) There is a functor (global section) between the $\infty$-prestacks of monomorphic inductive Banach quasicoherent commutative algebra $E_\infty$ objects:
\[
\xymatrix@R+6pc@C+0pc{
\mathrm{sComm}_\mathrm{simplicial}\mathrm{Ind}^m\mathrm{Banach}(\mathcal{O}_{X_{R_k,-}})\ar[d]\ar[d]\ar[d]\ar[d]\ar[r]^{\mathrm{global}}\ar[r]\ar[r] &\mathrm{sComm}_\mathrm{simplicial}\varphi\mathrm{Ind}^m\mathrm{Banach}(\{\mathrm{Robba}^\mathrm{extended}_{{R_k,-,I}}\}_I)\ar[d]\ar[d]\ar[d]\ar[d]\\
 \mathrm{sComm}_\mathrm{simplicial}\mathrm{Ind}^m\mathrm{Banach}(\mathcal{O}_{X_{\mathbb{Q}_p(p^{1/p^\infty})^{\wedge\flat},-}})\ar[r]^{\mathrm{global}}\ar[r]\ar[r] &\mathrm{sComm}_\mathrm{simplicial}\varphi\mathrm{Ind}^m\mathrm{Banach}(\{\mathrm{Robba}^\mathrm{extended}_{{R_0,-,I}}\}_I).
}
\]

\item Then parallel as in \cite{LBV} we have a functor (global section) of the de Rham complex after \cite[Definition 5.9, Section 5.2.1]{KKM}:
\[\displayindent=-0.4in
\xymatrix@R+6pc@C+0pc{
\mathrm{deRham}_{\mathrm{sComm}_\mathrm{simplicial}\mathrm{Ind}\mathrm{Banach}(\mathcal{O}_{X_{R_k,-}})\ar[r]^{\mathrm{global}}}(-)\ar[d]\ar[d]\ar[d]\ar[d]\ar[r]\ar[r] &\mathrm{deRham}_{\mathrm{sComm}_\mathrm{simplicial}\varphi\mathrm{Ind}\mathrm{Banach}(\{\mathrm{Robba}^\mathrm{extended}_{{R_k,-,I}}\}_I)}(-)\ar[d]\ar[d]\ar[d]\ar[d]\\
\mathrm{deRham}_{\mathrm{sComm}_\mathrm{simplicial}\mathrm{Ind}\mathrm{Banach}(\mathcal{O}_{X_{\mathbb{Q}_p(p^{1/p^\infty})^{\wedge\flat},-}})\ar[r]^{\mathrm{global}}}(-)\ar[r]\ar[r] &\mathrm{deRham}_{\mathrm{sComm}_\mathrm{simplicial}\varphi\mathrm{Ind}\mathrm{Banach}(\{\mathrm{Robba}^\mathrm{extended}_{{R_0,-,I}}\}_I)}(-), 
}
\]
\[\displayindent=-0.4in
\xymatrix@R+6pc@C+0pc{
\mathrm{deRham}_{\mathrm{sComm}_\mathrm{simplicial}\mathrm{Ind}^m\mathrm{Banach}(\mathcal{O}_{X_{R_k,-}})\ar[r]^{\mathrm{global}}}(-)\ar[d]\ar[d]\ar[d]\ar[d]\ar[r]\ar[r] &\mathrm{deRham}_{\mathrm{sComm}_\mathrm{simplicial}\varphi\mathrm{Ind}^m\mathrm{Banach}(\{\mathrm{Robba}^\mathrm{extended}_{{R_k,-,I}}\}_I)}(-)\ar[d]\ar[d]\ar[d]\ar[d]\\
\mathrm{deRham}_{\mathrm{sComm}_\mathrm{simplicial}\mathrm{Ind}^m\mathrm{Banach}(\mathcal{O}_{X_{\mathbb{Q}_p(p^{1/p^\infty})^{\wedge\flat},-}})\ar[r]^{\mathrm{global}}}(-)\ar[r]\ar[r] &\mathrm{deRham}_{\mathrm{sComm}_\mathrm{simplicial}\varphi\mathrm{Ind}^m\mathrm{Banach}(\{\mathrm{Robba}^\mathrm{extended}_{{R_0,-,I}}\}_I)}(-). 
}
\]

\item Then we have the following a functor (global section) of $K$-group $(\infty,1)$-spectrum from \cite{BGT}:
\[
\xymatrix@R+6pc@C+0pc{
\mathrm{K}^\mathrm{BGT}_{\mathrm{sComm}_\mathrm{simplicial}\mathrm{Ind}\mathrm{Banach}(\mathcal{O}_{X_{R_k,-}})\ar[r]^{\mathrm{global}}}(-)\ar[d]\ar[d]\ar[d]\ar[d]\ar[r]\ar[r] &\mathrm{K}^\mathrm{BGT}_{\mathrm{sComm}_\mathrm{simplicial}\varphi\mathrm{Ind}\mathrm{Banach}(\{\mathrm{Robba}^\mathrm{extended}_{{R_k,-,I}}\}_I)}(-)\ar[d]\ar[d]\ar[d]\ar[d]\\
\mathrm{K}^\mathrm{BGT}_{\mathrm{sComm}_\mathrm{simplicial}\mathrm{Ind}\mathrm{Banach}(\mathcal{O}_{X_{\mathbb{Q}_p(p^{1/p^\infty})^{\wedge\flat},-}})\ar[r]^{\mathrm{global}}}(-)\ar[r]\ar[r] &\mathrm{K}^\mathrm{BGT}_{\mathrm{sComm}_\mathrm{simplicial}\varphi\mathrm{Ind}\mathrm{Banach}(\{\mathrm{Robba}^\mathrm{extended}_{{R_0,-,I}}\}_I)}(-), 
}
\]
\[
\xymatrix@R+6pc@C+0pc{
\mathrm{K}^\mathrm{BGT}_{\mathrm{sComm}_\mathrm{simplicial}\mathrm{Ind}^m\mathrm{Banach}(\mathcal{O}_{X_{R_k,-}})\ar[r]^{\mathrm{global}}}(-)\ar[d]\ar[d]\ar[d]\ar[d]\ar[r]\ar[r] &\mathrm{K}^\mathrm{BGT}_{\mathrm{sComm}_\mathrm{simplicial}\varphi\mathrm{Ind}^m\mathrm{Banach}(\{\mathrm{Robba}^\mathrm{extended}_{{R_k,-,I}}\}_I)}(-)\ar[d]\ar[d]\ar[d]\ar[d]\\
\mathrm{K}^\mathrm{BGT}_{\mathrm{sComm}_\mathrm{simplicial}\mathrm{Ind}^m\mathrm{Banach}(\mathcal{O}_{X_{\mathbb{Q}_p(p^{1/p^\infty})^{\wedge\flat},-}})\ar[r]^{\mathrm{global}}}(-)\ar[r]\ar[r] &\mathrm{K}^\mathrm{BGT}_{\mathrm{sComm}_\mathrm{simplicial}\varphi\mathrm{Ind}^m\mathrm{Banach}(\{\mathrm{Robba}^\mathrm{extended}_{{R_0,-,I}}\}_I)}(-). 
}
\]

\end{itemize}

\
\indent Then we consider further equivariance by considering the arithmetic profinite fundamental groups $\Gamma_{\mathbb{Q}_p}$ and $\mathrm{Gal}(\overline{{Q}_p\left<T_1^{\pm 1},...,T_k^{\pm 1}\right>}/R_k)$ through the following diagram:

\[
\xymatrix@R+0pc@C+0pc{
\mathbb{Z}_p^k=\mathrm{Gal}(R_k/{\mathbb{Q}_p(p^{1/p^\infty})^\wedge\left<T_1^{\pm 1},...,T_k^{\pm 1}\right>}) \ar[r]\ar[r] \ar[r]\ar[r] &\Gamma_k:=\mathrm{Gal}(R_k/{\mathbb{Q}_p\left<T_1^{\pm 1},...,T_k^{\pm 1}\right>}) \ar[r] \ar[r]\ar[r] &\Gamma_{\mathbb{Q}_p}.
}
\]

\begin{itemize}
\item (\text{Proposition}) There is a functor (global section) between the $\infty$-prestacks of inductive Banach quasicoherent presheaves:
\[
\xymatrix@R+6pc@C+0pc{
\mathrm{Ind}\mathrm{Banach}_{\Gamma_k}(\mathcal{O}_{X_{\mathbb{Q}_p(p^{1/p^\infty})^{\wedge}\left<T_1^{\pm 1/p^\infty},...,T_k^{\pm 1/p^\infty}\right>^\flat,-}})\ar[d]\ar[d]\ar[d]\ar[d] \ar[r]^{\mathrm{global}}\ar[r]\ar[r] &\varphi\mathrm{Ind}\mathrm{Banach}_{\Gamma_k}(\{\mathrm{Robba}^\mathrm{extended}_{{R_k,-,I}}\}_I) \ar[d]\ar[d]\ar[d]\ar[d].\\
\mathrm{Ind}\mathrm{Banach}(\mathcal{O}_{X_{\mathbb{Q}_p(p^{1/p^\infty})^{\wedge\flat},-}})\ar[r]^{\mathrm{global}}\ar[r]\ar[r] &\varphi\mathrm{Ind}\mathrm{Banach}(\{\mathrm{Robba}^\mathrm{extended}_{{R_0,-,I}}\}_I).\\ 
}
\]
\item (\text{Proposition}) There is a functor (global section) between the $\infty$-prestacks of monomorphic inductive Banach quasicoherent presheaves:
\[
\xymatrix@R+6pc@C+0pc{
\mathrm{Ind}^m\mathrm{Banach}_{\Gamma_k}(\mathcal{O}_{X_{R_k,-}})\ar[r]^{\mathrm{global}}\ar[d]\ar[d]\ar[d]\ar[d]\ar[r]\ar[r] &\varphi\mathrm{Ind}^m\mathrm{Banach}_{\Gamma_k}(\{\mathrm{Robba}^\mathrm{extended}_{{R_k,-,I}}\}_I)\ar[d]\ar[d]\ar[d]\ar[d]\\
\mathrm{Ind}^m\mathrm{Banach}_{\Gamma_0}(\mathcal{O}_{X_{\mathbb{Q}_p(p^{1/p^\infty})^{\wedge\flat},-}})\ar[r]^{\mathrm{global}}\ar[r]\ar[r] &\varphi\mathrm{Ind}^m\mathrm{Banach}_{\Gamma_0}(\{\mathrm{Robba}^\mathrm{extended}_{{R_0,-,I}}\}_I).\\  
}
\]
\item (\text{Proposition}) There is a functor (global section) between the $\infty$-stacks of inductive Banach quasicoherent commutative algebra $E_\infty$ objects:
\[
\xymatrix@R+6pc@C+0pc{
\mathrm{sComm}_\mathrm{simplicial}\mathrm{Ind}\mathrm{Banach}_{\Gamma_k}(\mathcal{O}_{X_{R_k,-}})\ar[d]\ar[d]\ar[d]\ar[d]\ar[r]^{\mathrm{global}}\ar[r]\ar[r] &\mathrm{sComm}_\mathrm{simplicial}\varphi\mathrm{Ind}\mathrm{Banach}_{\Gamma_k}(\{\mathrm{Robba}^\mathrm{extended}_{{R_k,-,I}}\}_I)\ar[d]\ar[d]\ar[d]\ar[d]\\
\mathrm{sComm}_\mathrm{simplicial}\mathrm{Ind}\mathrm{Banach}_{\Gamma_0}(\mathcal{O}_{X_{\mathbb{Q}_p(p^{1/p^\infty})^{\wedge\flat},-}})\ar[r]^{\mathrm{global}}\ar[r]\ar[r] &\mathrm{sComm}_\mathrm{simplicial}\varphi\mathrm{Ind}\mathrm{Banach}_{\Gamma_0}(\{\mathrm{Robba}^\mathrm{extended}_{{R_0,-,I}}\}_I).  
}
\]
\item (\text{Proposition}) There is a functor (global section) between the $\infty$-prestacks of monomorphic inductive Banach quasicoherent commutative algebra $E_\infty$ objects:
\[\displayindent=-0.4in
\xymatrix@R+6pc@C+0pc{
\mathrm{sComm}_\mathrm{simplicial}\mathrm{Ind}^m\mathrm{Banach}_{\Gamma_k}(\mathcal{O}_{X_{R_k,-}})\ar[d]\ar[d]\ar[d]\ar[d]\ar[r]^{\mathrm{global}}\ar[r]\ar[r] &\mathrm{sComm}_\mathrm{simplicial}\varphi\mathrm{Ind}^m\mathrm{Banach}_{\Gamma_k}(\{\mathrm{Robba}^\mathrm{extended}_{{R_k,-,I}}\}_I)\ar[d]\ar[d]\ar[d]\ar[d]\\
 \mathrm{sComm}_\mathrm{simplicial}\mathrm{Ind}^m\mathrm{Banach}_{\Gamma_0}(\mathcal{O}_{X_{\mathbb{Q}_p(p^{1/p^\infty})^{\wedge\flat},-}})\ar[r]^{\mathrm{global}}\ar[r]\ar[r] &\mathrm{sComm}_\mathrm{simplicial}\varphi\mathrm{Ind}^m\mathrm{Banach}_{\Gamma_0}(\{\mathrm{Robba}^\mathrm{extended}_{{R_0,-,I}}\}_I). 
}
\]

\item Then parallel as in \cite{LBV} we have a functor (global section) of the de Rham complex after \cite[Definition 5.9, Section 5.2.1]{KKM}:
\[\displayindent=-0.4in
\xymatrix@R+6pc@C+0pc{
\mathrm{deRham}_{\mathrm{sComm}_\mathrm{simplicial}\mathrm{Ind}\mathrm{Banach}_{\Gamma_k}(\mathcal{O}_{X_{R_k,-}})\ar[r]^{\mathrm{global}}}(-)\ar[d]\ar[d]\ar[d]\ar[d]\ar[r]\ar[r] &\mathrm{deRham}_{\mathrm{sComm}_\mathrm{simplicial}\varphi\mathrm{Ind}\mathrm{Banach}_{\Gamma_k}(\{\mathrm{Robba}^\mathrm{extended}_{{R_k,-,I}}\}_I)}(-)\ar[d]\ar[d]\ar[d]\ar[d]\\
\mathrm{deRham}_{\mathrm{sComm}_\mathrm{simplicial}\mathrm{Ind}\mathrm{Banach}_{\Gamma_0}(\mathcal{O}_{X_{\mathbb{Q}_p(p^{1/p^\infty})^{\wedge\flat},-}})\ar[r]^{\mathrm{global}}}(-)\ar[r]\ar[r] &\mathrm{deRham}_{\mathrm{sComm}_\mathrm{simplicial}\varphi\mathrm{Ind}\mathrm{Banach}_{\Gamma_0}(\{\mathrm{Robba}^\mathrm{extended}_{{R_0,-,I}}\}_I)}(-), 
}
\]
\[\displayindent=-0.4in
\xymatrix@R+6pc@C+0pc{
\mathrm{deRham}_{\mathrm{sComm}_\mathrm{simplicial}\mathrm{Ind}^m\mathrm{Banach}_{\Gamma_k}(\mathcal{O}_{X_{R_k,-}})\ar[r]^{\mathrm{global}}}(-)\ar[d]\ar[d]\ar[d]\ar[d]\ar[r]\ar[r] &\mathrm{deRham}_{\mathrm{sComm}_\mathrm{simplicial}\varphi\mathrm{Ind}^m\mathrm{Banach}_{\Gamma_k}(\{\mathrm{Robba}^\mathrm{extended}_{{R_k,-,I}}\}_I)}(-)\ar[d]\ar[d]\ar[d]\ar[d]\\
\mathrm{deRham}_{\mathrm{sComm}_\mathrm{simplicial}\mathrm{Ind}^m\mathrm{Banach}_{\Gamma_0}(\mathcal{O}_{X_{\mathbb{Q}_p(p^{1/p^\infty})^{\wedge\flat},-}})\ar[r]^{\mathrm{global}}}(-)\ar[r]\ar[r] &\mathrm{deRham}_{\mathrm{sComm}_\mathrm{simplicial}\varphi\mathrm{Ind}^m\mathrm{Banach}_{\Gamma_0}(\{\mathrm{Robba}^\mathrm{extended}_{{R_0,-,I}}\}_I)}(-). 
}
\]

\item Then we have the following a functor (global section) of $K$-group $(\infty,1)$-spectrum from \cite{BGT}:
\[
\xymatrix@R+6pc@C+0pc{
\mathrm{K}^\mathrm{BGT}_{\mathrm{sComm}_\mathrm{simplicial}\mathrm{Ind}\mathrm{Banach}_{\Gamma_k}(\mathcal{O}_{X_{R_k,-}})\ar[r]^{\mathrm{global}}}(-)\ar[d]\ar[d]\ar[d]\ar[d]\ar[r]\ar[r] &\mathrm{K}^\mathrm{BGT}_{\mathrm{sComm}_\mathrm{simplicial}\varphi\mathrm{Ind}\mathrm{Banach}_{\Gamma_k}(\{\mathrm{Robba}^\mathrm{extended}_{{R_k,-,I}}\}_I)}(-)\ar[d]\ar[d]\ar[d]\ar[d]\\
\mathrm{K}^\mathrm{BGT}_{\mathrm{sComm}_\mathrm{simplicial}\mathrm{Ind}\mathrm{Banach}_{\Gamma_0}(\mathcal{O}_{X_{\mathbb{Q}_p(p^{1/p^\infty})^{\wedge\flat},-}})\ar[r]^{\mathrm{global}}}(-)\ar[r]\ar[r] &\mathrm{K}^\mathrm{BGT}_{\mathrm{sComm}_\mathrm{simplicial}\varphi\mathrm{Ind}\mathrm{Banach}_{\Gamma_0}(\{\mathrm{Robba}^\mathrm{extended}_{{R_0,-,I}}\}_I)}(-), 
}
\]
\[
\xymatrix@R+6pc@C+0pc{
\mathrm{K}^\mathrm{BGT}_{\mathrm{sComm}_\mathrm{simplicial}\mathrm{Ind}^m\mathrm{Banach}_{\Gamma_k}(\mathcal{O}_{X_{R_k,-}})\ar[r]^{\mathrm{global}}}(-)\ar[d]\ar[d]\ar[d]\ar[d]\ar[r]\ar[r] &\mathrm{K}^\mathrm{BGT}_{\mathrm{sComm}_\mathrm{simplicial}\varphi\mathrm{Ind}^m\mathrm{Banach}_{\Gamma_k}(\{\mathrm{Robba}^\mathrm{extended}_{{R_k,-,I}}\}_I)}(-)\ar[d]\ar[d]\ar[d]\ar[d]\\
\mathrm{K}^\mathrm{BGT}_{\mathrm{sComm}_\mathrm{simplicial}\mathrm{Ind}^m\mathrm{Banach}_{\Gamma_0}(\mathcal{O}_{X_{\mathbb{Q}_p(p^{1/p^\infty})^{\wedge\flat},-}})\ar[r]^{\mathrm{global}}}(-)\ar[r]\ar[r] &\mathrm{K}^\mathrm{BGT}_{\mathrm{sComm}_\mathrm{simplicial}\varphi\mathrm{Ind}^m\mathrm{Banach}_{\Gamma_0}(\{\mathrm{Robba}^\mathrm{extended}_{{R_0,-,I}}\}_I)}(-). 
}
\]

\end{itemize}

\

\begin{remark}
\noindent We can certainly consider the quasicoherent sheaves in \cite[Lemma 7.11, Remark 7.12]{1BBK}, therefore all the quasicoherent presheaves and modules will be those satisfying \cite[Lemma 7.11, Remark 7.12]{1BBK} if one would like to consider the the quasicoherent sheaves. That being all as this said, we would believe that the big quasicoherent presheaves are automatically quasicoherent sheaves (namely satisfying the corresponding \v{C}ech $\infty$-descent as in \cite[Section 9.3]{KKM} and \cite[Lemma 7.11, Remark 7.12]{1BBK}) and the corresponding global section functors are automatically equivalence of $\infty$-categories. \\
\end{remark}

\

\indent In Clausen-Scholze formalism we have the following:

\begin{itemize}

\item (\text{Proposition}) There is a functor (global section) between the $\infty$-prestacks of inductive Banach quasicoherent sheaves:
\[
\xymatrix@R+0pc@C+0pc{
{\mathrm{Modules}_\circledcirc}(\mathcal{O}_{X_{R,-}})\ar[r]^{\mathrm{global}}\ar[r]\ar[r] &\varphi{\mathrm{Modules}_\circledcirc}(\{\mathrm{Robba}^\mathrm{extended}_{{R,-,I}}\}_I).  
}
\]

\item (\text{Proposition}) There is a functor (global section) between the $\infty$-prestacks of inductive Banach quasicoherent sheaves:
\[
\xymatrix@R+0pc@C+0pc{
{\mathrm{Modules}_\circledcirc}(\mathcal{O}_{X_{R,-}})\ar[r]^{\mathrm{global}}\ar[r]\ar[r] &\varphi{\mathrm{Modules}_\circledcirc}(\{\mathrm{Robba}^\mathrm{extended}_{{R,-,I}}\}_I).  
}
\]

\item (\text{Proposition}) There is a functor (global section) between the $\infty$-prestacks of inductive Banach quasicoherent commutative algebra $E_\infty$ objects\footnote{Here $\circledcirc=\text{solidquasicoherentsheaves}$.}:
\[
\xymatrix@R+0pc@C+0pc{
\mathrm{sComm}_\mathrm{simplicial}{\mathrm{Modules}_\circledcirc}(\mathcal{O}_{X_{R,-}})\ar[r]^{\mathrm{global}}\ar[r]\ar[r] &\mathrm{sComm}_\mathrm{simplicial}\varphi{\mathrm{Modules}_\circledcirc}(\{\mathrm{Robba}^\mathrm{extended}_{{R,-,I}}\}_I).  
}
\]

\item Then as in \cite{LBV} we have a functor (global section) of the de Rham complex after \cite[Definition 5.9, Section 5.2.1]{KKM}\footnote{Here $\circledcirc=\text{solidquasicoherentsheaves}$.}:
\[
\xymatrix@R+0pc@C+0pc{
\mathrm{deRham}_{\mathrm{sComm}_\mathrm{simplicial}{\mathrm{Modules}_\circledcirc}(\mathcal{O}_{X_{R,-}})\ar[r]^{\mathrm{global}}}(-)\ar[r]\ar[r] &\mathrm{deRham}_{\mathrm{sComm}_\mathrm{simplicial}\varphi{\mathrm{Modules}_\circledcirc}(\{\mathrm{Robba}^\mathrm{extended}_{{R,-,I}}\}_I)}(-). 
}
\]

\item Then we have the following a functor (global section) of $K$-group $(\infty,1)$-spectrum from \cite{BGT}\footnote{Here $\circledcirc=\text{solidquasicoherentsheaves}$.}:
\[
\xymatrix@R+0pc@C+0pc{
\mathrm{K}^\mathrm{BGT}_{\mathrm{sComm}_\mathrm{simplicial}{\mathrm{Modules}_\circledcirc}(\mathcal{O}_{X_{R,-}})\ar[r]^{\mathrm{global}}}(-)\ar[r]\ar[r] &\mathrm{K}^\mathrm{BGT}_{\mathrm{sComm}_\mathrm{simplicial}\varphi{\mathrm{Modules}_\circledcirc}(\{\mathrm{Robba}^\mathrm{extended}_{{R,-,I}}\}_I)}(-). 
}
\]
\end{itemize}

\noindent Now let $R=\mathbb{Q}_p(p^{1/p^\infty})^{\wedge\flat}$ and $R_k=\mathbb{Q}_p(p^{1/p^\infty})^{\wedge}\left<T_1^{\pm 1/p^{\infty}},...,T_k^{\pm 1/p^{\infty}}\right>^\flat$ we have the following Galois theoretic results with naturality along $f:\mathrm{Spa}(\mathbb{Q}_p(p^{1/p^\infty})^{\wedge}\left<T_1^{\pm 1/p^\infty},...,T_k^{\pm 1/p^\infty}\right>^\flat)\rightarrow \mathrm{Spa}(\mathbb{Q}_p(p^{1/p^\infty})^{\wedge\flat})$:

\begin{itemize}
\item (\text{Proposition}) There is a functor (global section) between the $\infty$-prestacks of inductive Banach quasicoherent sheaves\footnote{Here $\circledcirc=\text{solidquasicoherentsheaves}$.}:
\[
\xymatrix@R+6pc@C+0pc{
{\mathrm{Modules}_\circledcirc}(\mathcal{O}_{X_{\mathbb{Q}_p(p^{1/p^\infty})^{\wedge}\left<T_1^{\pm 1/p^\infty},...,T_k^{\pm 1/p^\infty}\right>^\flat,-}})\ar[d]\ar[d]\ar[d]\ar[d] \ar[r]^{\mathrm{global}}\ar[r]\ar[r] &\varphi{\mathrm{Modules}_\circledcirc}(\{\mathrm{Robba}^\mathrm{extended}_{{R_k,-,I}}\}_I) \ar[d]\ar[d]\ar[d]\ar[d].\\
{\mathrm{Modules}_\circledcirc}(\mathcal{O}_{X_{\mathbb{Q}_p(p^{1/p^\infty})^{\wedge\flat},-}})\ar[r]^{\mathrm{global}}\ar[r]\ar[r] &\varphi{\mathrm{Modules}_\circledcirc}(\{\mathrm{Robba}^\mathrm{extended}_{{R_0,-,I}}\}_I).\\ 
}
\]
\item (\text{Proposition}) There is a functor (global section) between the $\infty$-prestacks of inductive Banach quasicoherent commutative algebra $E_\infty$ objects\footnote{Here $\circledcirc=\text{solidquasicoherentsheaves}$.}:
\[
\xymatrix@R+6pc@C+0pc{
\mathrm{sComm}_\mathrm{simplicial}{\mathrm{Modules}_\circledcirc}(\mathcal{O}_{X_{R_k,-}})\ar[d]\ar[d]\ar[d]\ar[d]\ar[r]^{\mathrm{global}}\ar[r]\ar[r] &\mathrm{sComm}_\mathrm{simplicial}\varphi{\mathrm{Modules}_\circledcirc}(\{\mathrm{Robba}^\mathrm{extended}_{{R_k,-,I}}\}_I)\ar[d]\ar[d]\ar[d]\ar[d]\\
\mathrm{sComm}_\mathrm{simplicial}{\mathrm{Modules}_\circledcirc}(\mathcal{O}_{X_{\mathbb{Q}_p(p^{1/p^\infty})^{\wedge\flat},-}})\ar[r]^{\mathrm{global}}\ar[r]\ar[r] &\mathrm{sComm}_\mathrm{simplicial}\varphi{\mathrm{Modules}_\circledcirc}(\{\mathrm{Robba}^\mathrm{extended}_{{R_0,-,I}}\}_I).  
}
\]

\item Then as in \cite{LBV} we have a functor (global section) of the de Rham complex after \cite[Definition 5.9, Section 5.2.1]{KKM}\footnote{Here $\circledcirc=\text{solidquasicoherentsheaves}$.}:
\[\displayindent=-0.4in
\xymatrix@R+6pc@C+0pc{
\mathrm{deRham}_{\mathrm{sComm}_\mathrm{simplicial}{\mathrm{Modules}_\circledcirc}(\mathcal{O}_{X_{R_k,-}})\ar[r]^{\mathrm{global}}}(-)\ar[d]\ar[d]\ar[d]\ar[d]\ar[r]\ar[r] &\mathrm{deRham}_{\mathrm{sComm}_\mathrm{simplicial}\varphi{\mathrm{Modules}_\circledcirc}(\{\mathrm{Robba}^\mathrm{extended}_{{R_k,-,I}}\}_I)}(-)\ar[d]\ar[d]\ar[d]\ar[d]\\
\mathrm{deRham}_{\mathrm{sComm}_\mathrm{simplicial}{\mathrm{Modules}_\circledcirc}(\mathcal{O}_{X_{\mathbb{Q}_p(p^{1/p^\infty})^{\wedge\flat},-}})\ar[r]^{\mathrm{global}}}(-)\ar[r]\ar[r] &\mathrm{deRham}_{\mathrm{sComm}_\mathrm{simplicial}\varphi{\mathrm{Modules}_\circledcirc}(\{\mathrm{Robba}^\mathrm{extended}_{{R_0,-,I}}\}_I)}(-). 
}
\]

\item Then we have the following a functor (global section) of $K$-group $(\infty,1)$-spectrum from \cite{BGT}\footnote{Here $\circledcirc=\text{solidquasicoherentsheaves}$.}:
\[
\xymatrix@R+6pc@C+0pc{
\mathrm{K}^\mathrm{BGT}_{\mathrm{sComm}_\mathrm{simplicial}{\mathrm{Modules}_\circledcirc}(\mathcal{O}_{X_{R_k,-}})\ar[r]^{\mathrm{global}}}(-)\ar[d]\ar[d]\ar[d]\ar[d]\ar[r]\ar[r] &\mathrm{K}^\mathrm{BGT}_{\mathrm{sComm}_\mathrm{simplicial}\varphi{\mathrm{Modules}_\circledcirc}(\{\mathrm{Robba}^\mathrm{extended}_{{R_k,-,I}}\}_I)}(-)\ar[d]\ar[d]\ar[d]\ar[d]\\
\mathrm{K}^\mathrm{BGT}_{\mathrm{sComm}_\mathrm{simplicial}{\mathrm{Modules}_\circledcirc}(\mathcal{O}_{X_{\mathbb{Q}_p(p^{1/p^\infty})^{\wedge\flat},-}})\ar[r]^{\mathrm{global}}}(-)\ar[r]\ar[r] &\mathrm{K}^\mathrm{BGT}_{\mathrm{sComm}_\mathrm{simplicial}\varphi{\mathrm{Modules}_\circledcirc}(\{\mathrm{Robba}^\mathrm{extended}_{{R_0,-,I}}\}_I)}(-). 
}
\]

\end{itemize}

\
\indent Then we consider further equivariance by considering the arithmetic profinite fundamental groups $\Gamma_{\mathbb{Q}_p}$ and $\mathrm{Gal}(\overline{{Q}_p\left<T_1^{\pm 1},...,T_k^{\pm 1}\right>}/R_k)$ through the following diagram:

\[
\xymatrix@R+0pc@C+0pc{
\mathbb{Z}_p^k=\mathrm{Gal}(R_k/{\mathbb{Q}_p(p^{1/p^\infty})^\wedge\left<T_1^{\pm 1},...,T_k^{\pm 1}\right>}) \ar[r]\ar[r] \ar[r]\ar[r] &\Gamma_k:=\mathrm{Gal}(R_k/{\mathbb{Q}_p\left<T_1^{\pm 1},...,T_k^{\pm 1}\right>}) \ar[r] \ar[r]\ar[r] &\Gamma_{\mathbb{Q}_p}.
}
\]

\begin{itemize}
\item (\text{Proposition}) There is a functor (global section) between the $\infty$-prestacks of inductive Banach quasicoherent sheaves\footnote{Here $\circledcirc=\text{solidquasicoherentsheaves}$.}:
\[
\xymatrix@R+6pc@C+0pc{
{\mathrm{Modules}_\circledcirc}_{\Gamma_k}(\mathcal{O}_{X_{\mathbb{Q}_p(p^{1/p^\infty})^{\wedge}\left<T_1^{\pm 1/p^\infty},...,T_k^{\pm 1/p^\infty}\right>^\flat,-}})\ar[d]\ar[d]\ar[d]\ar[d] \ar[r]^{\mathrm{global}}\ar[r]\ar[r] &\varphi{\mathrm{Modules}_\circledcirc}_{\Gamma_k}(\{\mathrm{Robba}^\mathrm{extended}_{{R_k,-,I}}\}_I) \ar[d]\ar[d]\ar[d]\ar[d].\\
{\mathrm{Modules}_\circledcirc}(\mathcal{O}_{X_{\mathbb{Q}_p(p^{1/p^\infty})^{\wedge\flat},-}})\ar[r]^{\mathrm{global}}\ar[r]\ar[r] &\varphi{\mathrm{Modules}_\circledcirc}(\{\mathrm{Robba}^\mathrm{extended}_{{R_0,-,I}}\}_I).\\ 
}
\]

\item (\text{Proposition}) There is a functor (global section) between the $\infty$-stacks of inductive Banach quasicoherent commutative algebra $E_\infty$ objects\footnote{Here $\circledcirc=\text{solidquasicoherentsheaves}$.}:
\[
\xymatrix@R+6pc@C+0pc{
\mathrm{sComm}_\mathrm{simplicial}{\mathrm{Modules}_\circledcirc}_{\Gamma_k}(\mathcal{O}_{X_{R_k,-}})\ar[d]\ar[d]\ar[d]\ar[d]\ar[r]^{\mathrm{global}}\ar[r]\ar[r] &\mathrm{sComm}_\mathrm{simplicial}\varphi{\mathrm{Modules}_\circledcirc}_{\Gamma_k}(\{\mathrm{Robba}^\mathrm{extended}_{{R_k,-,I}}\}_I)\ar[d]\ar[d]\ar[d]\ar[d]\\
\mathrm{sComm}_\mathrm{simplicial}{\mathrm{Modules}_\circledcirc}_{\Gamma_0}(\mathcal{O}_{X_{\mathbb{Q}_p(p^{1/p^\infty})^{\wedge\flat},-}})\ar[r]^{\mathrm{global}}\ar[r]\ar[r] &\mathrm{sComm}_\mathrm{simplicial}\varphi{\mathrm{Modules}_\circledcirc}_{\Gamma_0}(\{\mathrm{Robba}^\mathrm{extended}_{{R_0,-,I}}\}_I).  
}
\]

\item Then as in \cite{LBV} we have a functor (global section) of the de Rham complex after \cite[Definition 5.9, Section 5.2.1]{KKM}\footnote{Here $\circledcirc=\text{solidquasicoherentsheaves}$.}:
\[\displayindent=-0.4in
\xymatrix@R+6pc@C+0pc{
\mathrm{deRham}_{\mathrm{sComm}_\mathrm{simplicial}{\mathrm{Modules}_\circledcirc}_{\Gamma_k}(\mathcal{O}_{X_{R_k,-}})\ar[r]^{\mathrm{global}}}(-)\ar[d]\ar[d]\ar[d]\ar[d]\ar[r]\ar[r] &\mathrm{deRham}_{\mathrm{sComm}_\mathrm{simplicial}\varphi{\mathrm{Modules}_\circledcirc}_{\Gamma_k}(\{\mathrm{Robba}^\mathrm{extended}_{{R_k,-,I}}\}_I)}(-)\ar[d]\ar[d]\ar[d]\ar[d]\\
\mathrm{deRham}_{\mathrm{sComm}_\mathrm{simplicial}{\mathrm{Modules}_\circledcirc}_{\Gamma_0}(\mathcal{O}_{X_{\mathbb{Q}_p(p^{1/p^\infty})^{\wedge\flat},-}})\ar[r]^{\mathrm{global}}}(-)\ar[r]\ar[r] &\mathrm{deRham}_{\mathrm{sComm}_\mathrm{simplicial}\varphi{\mathrm{Modules}_\circledcirc}_{\Gamma_0}(\{\mathrm{Robba}^\mathrm{extended}_{{R_0,-,I}}\}_I)}(-). 
}
\]

\item Then we have the following a functor (global section) of $K$-group $(\infty,1)$-spectrum from \cite{BGT}\footnote{Here $\circledcirc=\text{solidquasicoherentsheaves}$.}:
\[
\xymatrix@R+6pc@C+0pc{
\mathrm{K}^\mathrm{BGT}_{\mathrm{sComm}_\mathrm{simplicial}{\mathrm{Modules}_\circledcirc}_{\Gamma_k}(\mathcal{O}_{X_{R_k,-}})\ar[r]^{\mathrm{global}}}(-)\ar[d]\ar[d]\ar[d]\ar[d]\ar[r]\ar[r] &\mathrm{K}^\mathrm{BGT}_{\mathrm{sComm}_\mathrm{simplicial}\varphi{\mathrm{Modules}_\circledcirc}_{\Gamma_k}(\{\mathrm{Robba}^\mathrm{extended}_{{R_k,-,I}}\}_I)}(-)\ar[d]\ar[d]\ar[d]\ar[d]\\
\mathrm{K}^\mathrm{BGT}_{\mathrm{sComm}_\mathrm{simplicial}{\mathrm{Modules}_\circledcirc}_{\Gamma_0}(\mathcal{O}_{X_{\mathbb{Q}_p(p^{1/p^\infty})^{\wedge\flat},-}})\ar[r]^{\mathrm{global}}}(-)\ar[r]\ar[r] &\mathrm{K}^\mathrm{BGT}_{\mathrm{sComm}_\mathrm{simplicial}\varphi{\mathrm{Modules}_\circledcirc}_{\Gamma_0}(\{\mathrm{Robba}^\mathrm{extended}_{{R_0,-,I}}\}_I)}(-). 
}
\]

\end{itemize}

\begin{proposition}
All the global functors from \cite[Proposition 13.8, Theorem 14.9, Remark 14.10]{1CS2} achieve the equivalences on both sides.	
\end{proposition}

\newpage
\subsection{$\infty$-Categorical Analytic Stacks and Descents IV}

\indent In the following the right had of each row in each diagram will be the corresponding quasicoherent Robba bundles over the Robba ring carrying the corresponding action from the Frobenius or the fundamental groups, defined by directly applying \cite[Section 9.3]{KKM} and \cite{BBM}. We now let $\mathcal{A}$ be any commutative algebra objects in the corresponding $\infty$-toposes over ind-Banach commutative algebra objects over $\mathbb{Q}_p$ or the corresponding born\'e commutative algebra objects over $\mathbb{Q}_p$ carrying the Grothendieck topology defined by essentially the corresponding monomorphism homotopy in the opposite category. Then we promote the construction to the corresponding $\infty$-stack over the same $\infty$-categories of affinoids. We now take the corresponding colimit through all the $(\infty,1)$-categories. Therefore all the corresponding $(\infty,1)$-functors into $(\infty,1)$-categories or $(\infty,1)$-groupoids are from the homotopy closure of $\mathbb{Q}_p\left<C_1,...,C_\ell\right>$ $\ell=1,2,...$ in $\mathrm{sComm}\mathrm{Ind}\mathrm{Banach}_{\mathbb{Q}_p}$ or $\mathbb{Q}_p\left<C_1,...,C_\ell\right>$ $\ell=1,2,...$ in $\mathrm{sComm}\mathrm{Ind}^m\mathrm{Banach}_{\mathbb{Q}_p}$ as in \cite[Section 4.2]{BBM}:
\begin{align}
\mathrm{Ind}^{\mathbb{Q}_p\left<C_1,...,C_\ell\right>,\ell=1,2,...}\mathrm{sComm}\mathrm{Ind}\mathrm{Banach}_{\mathbb{Q}_p},\\
\mathrm{Ind}^{\mathbb{Q}_p\left<C_1,...,C_\ell\right>,\ell=1,2,...}\mathrm{sComm}\mathrm{Ind}\mathrm{Banach}_{\mathbb{Q}_p}	.
\end{align}

\begin{itemize}

\item (\text{Proposition}) There is a functor (global section) between the $\infty$-prestacks of inductive Banach quasicoherent presheaves:
\[
\xymatrix@R+0pc@C+0pc{
\mathrm{Ind}\mathrm{Banach}(\mathcal{O}_{X_{R,-}})\ar[r]^{\mathrm{global}}\ar[r]\ar[r] &\varphi\mathrm{Ind}\mathrm{Banach}(\{\mathrm{Robba}^\mathrm{extended}_{{R,-,I}}\}_I).  
}
\]
The definition is given by the following:
\[
\xymatrix@R+0pc@C+0pc{
\mathrm{homotopycolimit}_i(\mathrm{Ind}\mathrm{Banach}(\mathcal{O}_{X_{R,-}})\ar[r]^{\mathrm{global}}\ar[r]\ar[r] &\varphi\mathrm{Ind}\mathrm{Banach}(\{\mathrm{Robba}^\mathrm{extended}_{{R,-,I}}\}_I))(\mathcal{O}_i),  
}
\]
each $\mathcal{O}_i$ is just as $\mathbb{Q}_p\left<C_1,...,C_\ell\right>,\ell=1,2,...$.
\item (\text{Proposition}) There is a functor (global section) between the $\infty$-prestacks of monomorphic inductive Banach quasicoherent presheaves:
\[
\xymatrix@R+0pc@C+0pc{
\mathrm{Ind}^m\mathrm{Banach}(\mathcal{O}_{X_{R,-}})\ar[r]^{\mathrm{global}}\ar[r]\ar[r] &\varphi\mathrm{Ind}^m\mathrm{Banach}(\{\mathrm{Robba}^\mathrm{extended}_{{R,-,I}}\}_I).  
}
\]
The definition is given by the following:
\[
\xymatrix@R+0pc@C+0pc{
\mathrm{homotopycolimit}_i(\mathrm{Ind}^m\mathrm{Banach}(\mathcal{O}_{X_{R,-}})\ar[r]^{\mathrm{global}}\ar[r]\ar[r] &\varphi\mathrm{Ind}^m\mathrm{Banach}(\{\mathrm{Robba}^\mathrm{extended}_{{R,-,I}}\}_I))(\mathcal{O}_i),  
}
\]
each $\mathcal{O}_i$ is just as $\mathbb{Q}_p\left<C_1,...,C_\ell\right>,\ell=1,2,...$.

\item (\text{Proposition}) There is a functor (global section) between the $\infty$-prestacks of inductive Banach quasicoherent presheaves:
\[
\xymatrix@R+0pc@C+0pc{
\mathrm{Ind}\mathrm{Banach}(\mathcal{O}_{X_{R,-}})\ar[r]^{\mathrm{global}}\ar[r]\ar[r] &\varphi\mathrm{Ind}\mathrm{Banach}(\{\mathrm{Robba}^\mathrm{extended}_{{R,-,I}}\}_I).  
}
\]
The definition is given by the following:
\[
\xymatrix@R+0pc@C+0pc{
\mathrm{homotopycolimit}_i(\mathrm{Ind}\mathrm{Banach}(\mathcal{O}_{X_{R,-}})\ar[r]^{\mathrm{global}}\ar[r]\ar[r] &\varphi\mathrm{Ind}\mathrm{Banach}(\{\mathrm{Robba}^\mathrm{extended}_{{R,-,I}}\}_I))(\mathcal{O}_i),  
}
\]
each $\mathcal{O}_i$ is just as $\mathbb{Q}_p\left<C_1,...,C_\ell\right>,\ell=1,2,...$.
\item (\text{Proposition}) There is a functor (global section) between the $\infty$-stacks of monomorphic inductive Banach quasicoherent presheaves:
\[
\xymatrix@R+0pc@C+0pc{
\mathrm{Ind}^m\mathrm{Banach}(\mathcal{O}_{X_{R,-}})\ar[r]^{\mathrm{global}}\ar[r]\ar[r] &\varphi\mathrm{Ind}^m\mathrm{Banach}(\{\mathrm{Robba}^\mathrm{extended}_{{R,-,I}}\}_I).  
}
\]
The definition is given by the following:
\[
\xymatrix@R+0pc@C+0pc{
\mathrm{homotopycolimit}_i(\mathrm{Ind}^m\mathrm{Banach}(\mathcal{O}_{X_{R,-}})\ar[r]^{\mathrm{global}}\ar[r]\ar[r] &\varphi\mathrm{Ind}^m\mathrm{Banach}(\{\mathrm{Robba}^\mathrm{extended}_{{R,-,I}}\}_I))(\mathcal{O}_i),  
}
\]
each $\mathcal{O}_i$ is just as $\mathbb{Q}_p\left<C_1,...,C_\ell\right>,\ell=1,2,...$.
\item (\text{Proposition}) There is a functor (global section) between the $\infty$-prestacks of inductive Banach quasicoherent commutative algebra $E_\infty$ objects:
\[
\xymatrix@R+0pc@C+0pc{
\mathrm{sComm}_\mathrm{simplicial}\mathrm{Ind}\mathrm{Banach}(\mathcal{O}_{X_{R,-}})\ar[r]^{\mathrm{global}}\ar[r]\ar[r] &\mathrm{sComm}_\mathrm{simplicial}\varphi\mathrm{Ind}\mathrm{Banach}(\{\mathrm{Robba}^\mathrm{extended}_{{R,-,I}}\}_I).  
}
\]
The definition is given by the following:
\[\displayindent=-0.4in
\xymatrix@R+0pc@C+0pc{
\mathrm{homotopycolimit}_i(\mathrm{sComm}_\mathrm{simplicial}\mathrm{Ind}\mathrm{Banach}(\mathcal{O}_{X_{R,-}})\ar[r]^{\mathrm{global}}\ar[r]\ar[r] &\mathrm{sComm}_\mathrm{simplicial}\varphi\mathrm{Ind}\mathrm{Banach}(\{\mathrm{Robba}^\mathrm{extended}_{{R,-,I}}\}_I))(\mathcal{O}_i),  
}
\]
each $\mathcal{O}_i$ is just as $\mathbb{Q}_p\left<C_1,...,C_\ell\right>,\ell=1,2,...$.
\item (\text{Proposition}) There is a functor (global section) between the $\infty$-prestacks of monomorphic inductive Banach quasicoherent commutative algebra $E_\infty$ objects:
\[
\xymatrix@R+0pc@C+0pc{
\mathrm{sComm}_\mathrm{simplicial}\mathrm{Ind}^m\mathrm{Banach}(\mathcal{O}_{X_{R,-}})\ar[r]^{\mathrm{global}}\ar[r]\ar[r] &\mathrm{sComm}_\mathrm{simplicial}\varphi\mathrm{Ind}^m\mathrm{Banach}(\{\mathrm{Robba}^\mathrm{extended}_{{R,-,I}}\}_I).  
}
\]
The definition is given by the following:
\[
\xymatrix@R+0pc@C+0pc{
\mathrm{homotopycolimit}_i(\mathrm{sComm}_\mathrm{simplicial}\mathrm{Ind}^m\mathrm{Banach}(\mathcal{O}_{X_{R,-}})\ar[r]^{\mathrm{global}}\ar[r]\ar[r] &\mathrm{sComm}_\mathrm{simplicial}\varphi\mathrm{Ind}^m\mathrm{Banach}(\{\mathrm{Robba}^\mathrm{extended}_{{R,-,I}}\}_I))(\mathcal{O}_i),  
}
\]
each $\mathcal{O}_i$ is just as $\mathbb{Q}_p\left<C_1,...,C_\ell\right>,\ell=1,2,...$.

\item Then parallel as in \cite{LBV} we have a functor (global section ) of the de Rham complex after \cite[Definition 5.9, Section 5.2.1]{KKM}:
\[
\xymatrix@R+0pc@C+0pc{
\mathrm{deRham}_{\mathrm{sComm}_\mathrm{simplicial}\mathrm{Ind}\mathrm{Banach}(\mathcal{O}_{X_{R,-}})\ar[r]^{\mathrm{global}}}(-)\ar[r]\ar[r] &\mathrm{deRham}_{\mathrm{sComm}_\mathrm{simplicial}\varphi\mathrm{Ind}\mathrm{Banach}(\{\mathrm{Robba}^\mathrm{extended}_{{R,-,I}}\}_I)}(-), 
}
\]
\[
\xymatrix@R+0pc@C+0pc{
\mathrm{deRham}_{\mathrm{sComm}_\mathrm{simplicial}\mathrm{Ind}^m\mathrm{Banach}(\mathcal{O}_{X_{R,-}})\ar[r]^{\mathrm{global}}}(-)\ar[r]\ar[r] &\mathrm{deRham}_{\mathrm{sComm}_\mathrm{simplicial}\varphi\mathrm{Ind}^m\mathrm{Banach}(\{\mathrm{Robba}^\mathrm{extended}_{{R,-,I}}\}_I)}(-). 
}
\]
The definition is given by the following:
\[
\xymatrix@R+0pc@C+0pc{
\mathrm{homotopycolimit}_i\\
(\mathrm{deRham}_{\mathrm{sComm}_\mathrm{simplicial}\mathrm{Ind}\mathrm{Banach}(\mathcal{O}_{X_{R,-}})\ar[r]^{\mathrm{global}}}(-)\ar[r]\ar[r] &\mathrm{deRham}_{\mathrm{sComm}_\mathrm{simplicial}\varphi\mathrm{Ind}\mathrm{Banach}(\{\mathrm{Robba}^\mathrm{extended}_{{R,-,I}}\}_I)}(-))(\mathcal{O}_i),  
}
\]
\[
\xymatrix@R+0pc@C+0pc{
\mathrm{homotopycolimit}_i\\
(\mathrm{deRham}_{\mathrm{sComm}_\mathrm{simplicial}\mathrm{Ind}^m\mathrm{Banach}(\mathcal{O}_{X_{R,-}})\ar[r]^{\mathrm{global}}}(-)\ar[r]\ar[r] &\mathrm{deRham}_{\mathrm{sComm}_\mathrm{simplicial}\varphi\mathrm{Ind}^m\mathrm{Banach}(\{\mathrm{Robba}^\mathrm{extended}_{{R,-,I}}\}_I)}(-))(\mathcal{O}_i),  
}
\]
each $\mathcal{O}_i$ is just as $\mathbb{Q}_p\left<C_1,...,C_\ell\right>,\ell=1,2,...$.\item Then we have the following a functor (global section) of $K$-group $(\infty,1)$-spectrum from \cite{BGT}:
\[
\xymatrix@R+0pc@C+0pc{
\mathrm{K}^\mathrm{BGT}_{\mathrm{sComm}_\mathrm{simplicial}\mathrm{Ind}\mathrm{Banach}(\mathcal{O}_{X_{R,-}})\ar[r]^{\mathrm{global}}}(-)\ar[r]\ar[r] &\mathrm{K}^\mathrm{BGT}_{\mathrm{sComm}_\mathrm{simplicial}\varphi\mathrm{Ind}\mathrm{Banach}(\{\mathrm{Robba}^\mathrm{extended}_{{R,-,I}}\}_I)}(-), 
}
\]
\[
\xymatrix@R+0pc@C+0pc{
\mathrm{K}^\mathrm{BGT}_{\mathrm{sComm}_\mathrm{simplicial}\mathrm{Ind}^m\mathrm{Banach}(\mathcal{O}_{X_{R,-}})\ar[r]^{\mathrm{global}}}(-)\ar[r]\ar[r] &\mathrm{K}^\mathrm{BGT}_{\mathrm{sComm}_\mathrm{simplicial}\varphi\mathrm{Ind}^m\mathrm{Banach}(\{\mathrm{Robba}^\mathrm{extended}_{{R,-,I}}\}_I)}(-). 
}
\]
The definition is given by the following:
\[\displayindent=-0.4in
\xymatrix@R+0pc@C+0pc{
\mathrm{homotopycolimit}_i(\mathrm{K}^\mathrm{BGT}_{\mathrm{sComm}_\mathrm{simplicial}\mathrm{Ind}\mathrm{Banach}(\mathcal{O}_{X_{R,-}})\ar[r]^{\mathrm{global}}}(-)\ar[r]\ar[r] &\mathrm{K}^\mathrm{BGT}_{\mathrm{sComm}_\mathrm{simplicial}\varphi\mathrm{Ind}\mathrm{Banach}(\{\mathrm{Robba}^\mathrm{extended}_{{R,-,I}}\}_I)}(-))(\mathcal{O}_i),  
}
\]
\[\displayindent=-0.4in
\xymatrix@R+0pc@C+0pc{
\mathrm{homotopycolimit}_i(\mathrm{K}^\mathrm{BGT}_{\mathrm{sComm}_\mathrm{simplicial}\mathrm{Ind}^m\mathrm{Banach}(\mathcal{O}_{X_{R,-}})\ar[r]^{\mathrm{global}}}(-)\ar[r]\ar[r] &\mathrm{K}^\mathrm{BGT}_{\mathrm{sComm}_\mathrm{simplicial}\varphi\mathrm{Ind}^m\mathrm{Banach}(\{\mathrm{Robba}^\mathrm{extended}_{{R,-,I}}\}_I)}(-))(\mathcal{O}_i),  
}
\]
each $\mathcal{O}_i$ is just as $\mathbb{Q}_p\left<C_1,...,C_\ell\right>,\ell=1,2,...$.
\end{itemize}

\noindent Now let $R=\mathbb{Q}_p(p^{1/p^\infty})^{\wedge\flat}$ and $R_k=\mathbb{Q}_p(p^{1/p^\infty})^{\wedge}\left<T_1^{\pm 1/p^{\infty}},...,T_k^{\pm 1/p^{\infty}}\right>^\flat$ we have the following Galois theoretic results with naturality along $f:\mathrm{Spa}(\mathbb{Q}_p(p^{1/p^\infty})^{\wedge}\left<T_1^{\pm 1/p^\infty},...,T_k^{\pm 1/p^\infty}\right>^\flat)\rightarrow \mathrm{Spa}(\mathbb{Q}_p(p^{1/p^\infty})^{\wedge\flat})$:

\begin{itemize}
\item (\text{Proposition}) There is a functor (global section) between the $\infty$-prestacks of inductive Banach quasicoherent presheaves:
\[
\xymatrix@R+6pc@C+0pc{
\mathrm{Ind}\mathrm{Banach}(\mathcal{O}_{X_{\mathbb{Q}_p(p^{1/p^\infty})^{\wedge}\left<T_1^{\pm 1/p^\infty},...,T_k^{\pm 1/p^\infty}\right>^\flat,-}})\ar[d]\ar[d]\ar[d]\ar[d] \ar[r]^{\mathrm{global}}\ar[r]\ar[r] &\varphi\mathrm{Ind}\mathrm{Banach}(\{\mathrm{Robba}^\mathrm{extended}_{{R_k,-,I}}\}_I) \ar[d]\ar[d]\ar[d]\ar[d].\\
\mathrm{Ind}\mathrm{Banach}(\mathcal{O}_{X_{\mathbb{Q}_p(p^{1/p^\infty})^{\wedge\flat},-}})\ar[r]^{\mathrm{global}}\ar[r]\ar[r] &\varphi\mathrm{Ind}\mathrm{Banach}(\{\mathrm{Robba}^\mathrm{extended}_{{R_0,-,I}}\}_I).\\ 
}
\]
The definition is given by the following:
\[
\xymatrix@R+0pc@C+0pc{
\mathrm{homotopycolimit}_i(\square)(\mathcal{O}_i),  
}
\]
each $\mathcal{O}_i$ is just as $\mathbb{Q}_p\left<C_1,...,C_\ell\right>,\ell=1,2,...$ and $\square$ is the relative diagram of $\infty$-functors.
\item (\text{Proposition}) There is a functor (global section) between the $\infty$-prestacks of monomorphic inductive Banach quasicoherent presheaves:
\[
\xymatrix@R+6pc@C+0pc{
\mathrm{Ind}^m\mathrm{Banach}(\mathcal{O}_{X_{R_k,-}})\ar[r]^{\mathrm{global}}\ar[d]\ar[d]\ar[d]\ar[d]\ar[r]\ar[r] &\varphi\mathrm{Ind}^m\mathrm{Banach}(\{\mathrm{Robba}^\mathrm{extended}_{{R_k,-,I}}\}_I)\ar[d]\ar[d]\ar[d]\ar[d]\\
\mathrm{Ind}^m\mathrm{Banach}(\mathcal{O}_{X_{\mathbb{Q}_p(p^{1/p^\infty})^{\wedge\flat},-}})\ar[r]^{\mathrm{global}}\ar[r]\ar[r] &\varphi\mathrm{Ind}^m\mathrm{Banach}(\{\mathrm{Robba}^\mathrm{extended}_{{R_0,-,I}}\}_I).\\  
}
\]
The definition is given by the following:
\[
\xymatrix@R+0pc@C+0pc{
\mathrm{homotopycolimit}_i(\square)(\mathcal{O}_i),  
}
\]
each $\mathcal{O}_i$ is just as $\mathbb{Q}_p\left<C_1,...,C_\ell\right>,\ell=1,2,...$ and $\square$ is the relative diagram of $\infty$-functors.

\item (\text{Proposition}) There is a functor (global section) between the $\infty$-prestacks of inductive Banach quasicoherent commutative algebra $E_\infty$ objects:
\[
\xymatrix@R+6pc@C+0pc{
\mathrm{sComm}_\mathrm{simplicial}\mathrm{Ind}\mathrm{Banach}(\mathcal{O}_{X_{R_k,-}})\ar[d]\ar[d]\ar[d]\ar[d]\ar[r]^{\mathrm{global}}\ar[r]\ar[r] &\mathrm{sComm}_\mathrm{simplicial}\varphi\mathrm{Ind}\mathrm{Banach}(\{\mathrm{Robba}^\mathrm{extended}_{{R_k,-,I}}\}_I)\ar[d]\ar[d]\ar[d]\ar[d]\\
\mathrm{sComm}_\mathrm{simplicial}\mathrm{Ind}\mathrm{Banach}(\mathcal{O}_{X_{\mathbb{Q}_p(p^{1/p^\infty})^{\wedge\flat},-}})\ar[r]^{\mathrm{global}}\ar[r]\ar[r] &\mathrm{sComm}_\mathrm{simplicial}\varphi\mathrm{Ind}\mathrm{Banach}(\{\mathrm{Robba}^\mathrm{extended}_{{R_0,-,I}}\}_I).  
}
\]
The definition is given by the following:
\[
\xymatrix@R+0pc@C+0pc{
\mathrm{homotopycolimit}_i(\square)(\mathcal{O}_i),  
}
\]
each $\mathcal{O}_i$ is just as $\mathbb{Q}_p\left<C_1,...,C_\ell\right>,\ell=1,2,...$ and $\square$ is the relative diagram of $\infty$-functors.

\item (\text{Proposition}) There is a functor (global section) between the $\infty$-prestacks of monomorphic inductive Banach quasicoherent commutative algebra $E_\infty$ objects:
\[
\xymatrix@R+6pc@C+0pc{
\mathrm{sComm}_\mathrm{simplicial}\mathrm{Ind}^m\mathrm{Banach}(\mathcal{O}_{X_{R_k,-}})\ar[d]\ar[d]\ar[d]\ar[d]\ar[r]^{\mathrm{global}}\ar[r]\ar[r] &\mathrm{sComm}_\mathrm{simplicial}\varphi\mathrm{Ind}^m\mathrm{Banach}(\{\mathrm{Robba}^\mathrm{extended}_{{R_k,-,I}}\}_I)\ar[d]\ar[d]\ar[d]\ar[d]\\
 \mathrm{sComm}_\mathrm{simplicial}\mathrm{Ind}^m\mathrm{Banach}(\mathcal{O}_{X_{\mathbb{Q}_p(p^{1/p^\infty})^{\wedge\flat},-}})\ar[r]^{\mathrm{global}}\ar[r]\ar[r] &\mathrm{sComm}_\mathrm{simplicial}\varphi\mathrm{Ind}^m\mathrm{Banach}(\{\mathrm{Robba}^\mathrm{extended}_{{R_0,-,I}}\}_I).
}
\]
The definition is given by the following:
\[
\xymatrix@R+0pc@C+0pc{
\mathrm{homotopycolimit}_i(\square)(\mathcal{O}_i),  
}
\]
each $\mathcal{O}_i$ is just as $\mathbb{Q}_p\left<C_1,...,C_\ell\right>,\ell=1,2,...$ and $\square$ is the relative diagram of $\infty$-functors.

\item Then parallel as in \cite{LBV} we have a functor (global section) of the de Rham complex after \cite[Definition 5.9, Section 5.2.1]{KKM}:
\[\displayindent=-0.4in
\xymatrix@R+6pc@C+0pc{
\mathrm{deRham}_{\mathrm{sComm}_\mathrm{simplicial}\mathrm{Ind}\mathrm{Banach}(\mathcal{O}_{X_{R_k,-}})\ar[r]^{\mathrm{global}}}(-)\ar[d]\ar[d]\ar[d]\ar[d]\ar[r]\ar[r] &\mathrm{deRham}_{\mathrm{sComm}_\mathrm{simplicial}\varphi\mathrm{Ind}\mathrm{Banach}(\{\mathrm{Robba}^\mathrm{extended}_{{R_k,-,I}}\}_I)}(-)\ar[d]\ar[d]\ar[d]\ar[d]\\
\mathrm{deRham}_{\mathrm{sComm}_\mathrm{simplicial}\mathrm{Ind}\mathrm{Banach}(\mathcal{O}_{X_{\mathbb{Q}_p(p^{1/p^\infty})^{\wedge\flat},-}})\ar[r]^{\mathrm{global}}}(-)\ar[r]\ar[r] &\mathrm{deRham}_{\mathrm{sComm}_\mathrm{simplicial}\varphi\mathrm{Ind}\mathrm{Banach}(\{\mathrm{Robba}^\mathrm{extended}_{{R_0,-,I}}\}_I)}(-), 
}
\]
\[\displayindent=-0.4in
\xymatrix@R+6pc@C+0pc{
\mathrm{deRham}_{\mathrm{sComm}_\mathrm{simplicial}\mathrm{Ind}^m\mathrm{Banach}(\mathcal{O}_{X_{R_k,-}})\ar[r]^{\mathrm{global}}}(-)\ar[d]\ar[d]\ar[d]\ar[d]\ar[r]\ar[r] &\mathrm{deRham}_{\mathrm{sComm}_\mathrm{simplicial}\varphi\mathrm{Ind}^m\mathrm{Banach}(\{\mathrm{Robba}^\mathrm{extended}_{{R_k,-,I}}\}_I)}(-)\ar[d]\ar[d]\ar[d]\ar[d]\\
\mathrm{deRham}_{\mathrm{sComm}_\mathrm{simplicial}\mathrm{Ind}^m\mathrm{Banach}(\mathcal{O}_{X_{\mathbb{Q}_p(p^{1/p^\infty})^{\wedge\flat},-}})\ar[r]^{\mathrm{global}}}(-)\ar[r]\ar[r] &\mathrm{deRham}_{\mathrm{sComm}_\mathrm{simplicial}\varphi\mathrm{Ind}^m\mathrm{Banach}(\{\mathrm{Robba}^\mathrm{extended}_{{R_0,-,I}}\}_I)}(-). 
}
\]

\item Then we have the following a functor (global section) of $K$-group $(\infty,1)$-spectrum from \cite{BGT}:
\[
\xymatrix@R+6pc@C+0pc{
\mathrm{K}^\mathrm{BGT}_{\mathrm{sComm}_\mathrm{simplicial}\mathrm{Ind}\mathrm{Banach}(\mathcal{O}_{X_{R_k,-}})\ar[r]^{\mathrm{global}}}(-)\ar[d]\ar[d]\ar[d]\ar[d]\ar[r]\ar[r] &\mathrm{K}^\mathrm{BGT}_{\mathrm{sComm}_\mathrm{simplicial}\varphi\mathrm{Ind}\mathrm{Banach}(\{\mathrm{Robba}^\mathrm{extended}_{{R_k,-,I}}\}_I)}(-)\ar[d]\ar[d]\ar[d]\ar[d]\\
\mathrm{K}^\mathrm{BGT}_{\mathrm{sComm}_\mathrm{simplicial}\mathrm{Ind}\mathrm{Banach}(\mathcal{O}_{X_{\mathbb{Q}_p(p^{1/p^\infty})^{\wedge\flat},-}})\ar[r]^{\mathrm{global}}}(-)\ar[r]\ar[r] &\mathrm{K}^\mathrm{BGT}_{\mathrm{sComm}_\mathrm{simplicial}\varphi\mathrm{Ind}\mathrm{Banach}(\{\mathrm{Robba}^\mathrm{extended}_{{R_0,-,I}}\}_I)}(-), 
}
\]
\[
\xymatrix@R+6pc@C+0pc{
\mathrm{K}^\mathrm{BGT}_{\mathrm{sComm}_\mathrm{simplicial}\mathrm{Ind}^m\mathrm{Banach}(\mathcal{O}_{X_{R_k,-}})\ar[r]^{\mathrm{global}}}(-)\ar[d]\ar[d]\ar[d]\ar[d]\ar[r]\ar[r] &\mathrm{K}^\mathrm{BGT}_{\mathrm{sComm}_\mathrm{simplicial}\varphi\mathrm{Ind}^m\mathrm{Banach}(\{\mathrm{Robba}^\mathrm{extended}_{{R_k,-,I}}\}_I)}(-)\ar[d]\ar[d]\ar[d]\ar[d]\\
\mathrm{K}^\mathrm{BGT}_{\mathrm{sComm}_\mathrm{simplicial}\mathrm{Ind}^m\mathrm{Banach}(\mathcal{O}_{X_{\mathbb{Q}_p(p^{1/p^\infty})^{\wedge\flat},-}})\ar[r]^{\mathrm{global}}}(-)\ar[r]\ar[r] &\mathrm{K}^\mathrm{BGT}_{\mathrm{sComm}_\mathrm{simplicial}\varphi\mathrm{Ind}^m\mathrm{Banach}(\{\mathrm{Robba}^\mathrm{extended}_{{R_0,-,I}}\}_I)}(-). 
}
\]
The definition is given by the following:
\[
\xymatrix@R+0pc@C+0pc{
\mathrm{homotopycolimit}_i(\square)(\mathcal{O}_i),  
}
\]
each $\mathcal{O}_i$ is just as $\mathbb{Q}_p\left<C_1,...,C_\ell\right>,\ell=1,2,...$ and $\square$ is the relative diagram of $\infty$-functors.

\end{itemize}

\
\indent Then we consider further equivariance by considering the arithmetic profinite fundamental groups $\Gamma_{\mathbb{Q}_p}$ and $\mathrm{Gal}(\overline{{Q}_p\left<T_1^{\pm 1},...,T_k^{\pm 1}\right>}/R_k)$ through the following diagram:

\[
\xymatrix@R+0pc@C+0pc{
\mathbb{Z}_p^k=\mathrm{Gal}(R_k/{\mathbb{Q}_p(p^{1/p^\infty})^\wedge\left<T_1^{\pm 1},...,T_k^{\pm 1}\right>}) \ar[r]\ar[r] \ar[r]\ar[r] &\Gamma_k:=\mathrm{Gal}(R_k/{\mathbb{Q}_p\left<T_1^{\pm 1},...,T_k^{\pm 1}\right>}) \ar[r] \ar[r]\ar[r] &\Gamma_{\mathbb{Q}_p}.
}
\]

\begin{itemize}
\item (\text{Proposition}) There is a functor (global section) between the $\infty$-prestacks of inductive Banach quasicoherent presheaves:
\[
\xymatrix@R+6pc@C+0pc{
\mathrm{Ind}\mathrm{Banach}_{\Gamma_k}(\mathcal{O}_{X_{\mathbb{Q}_p(p^{1/p^\infty})^{\wedge}\left<T_1^{\pm 1/p^\infty},...,T_k^{\pm 1/p^\infty}\right>^\flat,-}})\ar[d]\ar[d]\ar[d]\ar[d] \ar[r]^{\mathrm{global}}\ar[r]\ar[r] &\varphi\mathrm{Ind}\mathrm{Banach}_{\Gamma_k}(\{\mathrm{Robba}^\mathrm{extended}_{{R_k,-,I}}\}_I) \ar[d]\ar[d]\ar[d]\ar[d].\\
\mathrm{Ind}\mathrm{Banach}(\mathcal{O}_{X_{\mathbb{Q}_p(p^{1/p^\infty})^{\wedge\flat},-}})\ar[r]^{\mathrm{global}}\ar[r]\ar[r] &\varphi\mathrm{Ind}\mathrm{Banach}(\{\mathrm{Robba}^\mathrm{extended}_{{R_0,-,I}}\}_I).\\ 
}
\]
The definition is given by the following:
\[
\xymatrix@R+0pc@C+0pc{
\mathrm{homotopycolimit}_i(\square)(\mathcal{O}_i),  
}
\]
each $\mathcal{O}_i$ is just as $\mathbb{Q}_p\left<C_1,...,C_\ell\right>,\ell=1,2,...$ and $\square$ is the relative diagram of $\infty$-functors.

\item (\text{Proposition}) There is a functor (global section) between the $\infty$-prestacks of monomorphic inductive Banach quasicoherent presheaves:
\[
\xymatrix@R+6pc@C+0pc{
\mathrm{Ind}^m\mathrm{Banach}_{\Gamma_k}(\mathcal{O}_{X_{R_k,-}})\ar[r]^{\mathrm{global}}\ar[d]\ar[d]\ar[d]\ar[d]\ar[r]\ar[r] &\varphi\mathrm{Ind}^m\mathrm{Banach}_{\Gamma_k}(\{\mathrm{Robba}^\mathrm{extended}_{{R_k,-,I}}\}_I)\ar[d]\ar[d]\ar[d]\ar[d]\\
\mathrm{Ind}^m\mathrm{Banach}_{\Gamma_0}(\mathcal{O}_{X_{\mathbb{Q}_p(p^{1/p^\infty})^{\wedge\flat},-}})\ar[r]^{\mathrm{global}}\ar[r]\ar[r] &\varphi\mathrm{Ind}^m\mathrm{Banach}_{\Gamma_0}(\{\mathrm{Robba}^\mathrm{extended}_{{R_0,-,I}}\}_I).\\  
}
\]
The definition is given by the following:
\[
\xymatrix@R+0pc@C+0pc{
\mathrm{homotopycolimit}_i(\square)(\mathcal{O}_i),  
}
\]
each $\mathcal{O}_i$ is just as $\mathbb{Q}_p\left<C_1,...,C_\ell\right>,\ell=1,2,...$ and $\square$ is the relative diagram of $\infty$-functors.

\item (\text{Proposition}) There is a functor (global section) between the $\infty$-stacks of inductive Banach quasicoherent commutative algebra $E_\infty$ objects:
\[
\xymatrix@R+6pc@C+0pc{
\mathrm{sComm}_\mathrm{simplicial}\mathrm{Ind}\mathrm{Banach}_{\Gamma_k}(\mathcal{O}_{X_{R_k,-}})\ar[d]\ar[d]\ar[d]\ar[d]\ar[r]^{\mathrm{global}}\ar[r]\ar[r] &\mathrm{sComm}_\mathrm{simplicial}\varphi\mathrm{Ind}\mathrm{Banach}_{\Gamma_k}(\{\mathrm{Robba}^\mathrm{extended}_{{R_k,-,I}}\}_I)\ar[d]\ar[d]\ar[d]\ar[d]\\
\mathrm{sComm}_\mathrm{simplicial}\mathrm{Ind}\mathrm{Banach}_{\Gamma_0}(\mathcal{O}_{X_{\mathbb{Q}_p(p^{1/p^\infty})^{\wedge\flat},-}})\ar[r]^{\mathrm{global}}\ar[r]\ar[r] &\mathrm{sComm}_\mathrm{simplicial}\varphi\mathrm{Ind}\mathrm{Banach}_{\Gamma_0}(\{\mathrm{Robba}^\mathrm{extended}_{{R_0,-,I}}\}_I).  
}
\]
The definition is given by the following:
\[
\xymatrix@R+0pc@C+0pc{
\mathrm{homotopycolimit}_i(\square)(\mathcal{O}_i),  
}
\]
each $\mathcal{O}_i$ is just as $\mathbb{Q}_p\left<C_1,...,C_\ell\right>,\ell=1,2,...$ and $\square$ is the relative diagram of $\infty$-functors.

\item (\text{Proposition}) There is a functor (global section) between the $\infty$-prestacks of monomorphic inductive Banach quasicoherent commutative algebra $E_\infty$ objects:
\[\displayindent=-0.4in
\xymatrix@R+6pc@C+0pc{
\mathrm{sComm}_\mathrm{simplicial}\mathrm{Ind}^m\mathrm{Banach}_{\Gamma_k}(\mathcal{O}_{X_{R_k,-}})\ar[d]\ar[d]\ar[d]\ar[d]\ar[r]^{\mathrm{global}}\ar[r]\ar[r] &\mathrm{sComm}_\mathrm{simplicial}\varphi\mathrm{Ind}^m\mathrm{Banach}_{\Gamma_k}(\{\mathrm{Robba}^\mathrm{extended}_{{R_k,-,I}}\}_I)\ar[d]\ar[d]\ar[d]\ar[d]\\
 \mathrm{sComm}_\mathrm{simplicial}\mathrm{Ind}^m\mathrm{Banach}_{\Gamma_0}(\mathcal{O}_{X_{\mathbb{Q}_p(p^{1/p^\infty})^{\wedge\flat},-}})\ar[r]^{\mathrm{global}}\ar[r]\ar[r] &\mathrm{sComm}_\mathrm{simplicial}\varphi\mathrm{Ind}^m\mathrm{Banach}_{\Gamma_0}(\{\mathrm{Robba}^\mathrm{extended}_{{R_0,-,I}}\}_I). 
}
\]
The definition is given by the following:
\[
\xymatrix@R+0pc@C+0pc{
\mathrm{homotopycolimit}_i(\square)(\mathcal{O}_i),  
}
\]
each $\mathcal{O}_i$ is just as $\mathbb{Q}_p\left<C_1,...,C_\ell\right>,\ell=1,2,...$ and $\square$ is the relative diagram of $\infty$-functors.

\item Then parallel as in \cite{LBV} we have a functor (global section) of the de Rham complex after \cite[Definition 5.9, Section 5.2.1]{KKM}:
\[\displayindent=-0.4in
\xymatrix@R+6pc@C+0pc{
\mathrm{deRham}_{\mathrm{sComm}_\mathrm{simplicial}\mathrm{Ind}\mathrm{Banach}_{\Gamma_k}(\mathcal{O}_{X_{R_k,-}})\ar[r]^{\mathrm{global}}}(-)\ar[d]\ar[d]\ar[d]\ar[d]\ar[r]\ar[r] &\mathrm{deRham}_{\mathrm{sComm}_\mathrm{simplicial}\varphi\mathrm{Ind}\mathrm{Banach}_{\Gamma_k}(\{\mathrm{Robba}^\mathrm{extended}_{{R_k,-,I}}\}_I)}(-)\ar[d]\ar[d]\ar[d]\ar[d]\\
\mathrm{deRham}_{\mathrm{sComm}_\mathrm{simplicial}\mathrm{Ind}\mathrm{Banach}_{\Gamma_0}(\mathcal{O}_{X_{\mathbb{Q}_p(p^{1/p^\infty})^{\wedge\flat},-}})\ar[r]^{\mathrm{global}}}(-)\ar[r]\ar[r] &\mathrm{deRham}_{\mathrm{sComm}_\mathrm{simplicial}\varphi\mathrm{Ind}\mathrm{Banach}_{\Gamma_0}(\{\mathrm{Robba}^\mathrm{extended}_{{R_0,-,I}}\}_I)}(-), 
}
\]
\[\displayindent=-0.4in
\xymatrix@R+6pc@C+0pc{
\mathrm{deRham}_{\mathrm{sComm}_\mathrm{simplicial}\mathrm{Ind}^m\mathrm{Banach}_{\Gamma_k}(\mathcal{O}_{X_{R_k,-}})\ar[r]^{\mathrm{global}}}(-)\ar[d]\ar[d]\ar[d]\ar[d]\ar[r]\ar[r] &\mathrm{deRham}_{\mathrm{sComm}_\mathrm{simplicial}\varphi\mathrm{Ind}^m\mathrm{Banach}_{\Gamma_k}(\{\mathrm{Robba}^\mathrm{extended}_{{R_k,-,I}}\}_I)}(-)\ar[d]\ar[d]\ar[d]\ar[d]\\
\mathrm{deRham}_{\mathrm{sComm}_\mathrm{simplicial}\mathrm{Ind}^m\mathrm{Banach}_{\Gamma_0}(\mathcal{O}_{X_{\mathbb{Q}_p(p^{1/p^\infty})^{\wedge\flat},-}})\ar[r]^{\mathrm{global}}}(-)\ar[r]\ar[r] &\mathrm{deRham}_{\mathrm{sComm}_\mathrm{simplicial}\varphi\mathrm{Ind}^m\mathrm{Banach}_{\Gamma_0}(\{\mathrm{Robba}^\mathrm{extended}_{{R_0,-,I}}\}_I)}(-). 
}
\]
The definition is given by the following:
\[
\xymatrix@R+0pc@C+0pc{
\mathrm{homotopycolimit}_i(\square)(\mathcal{O}_i),  
}
\]
each $\mathcal{O}_i$ is just as $\mathbb{Q}_p\left<C_1,...,C_\ell\right>,\ell=1,2,...$ and $\square$ is the relative diagram of $\infty$-functors.

\item Then we have the following a functor (global section) of $K$-group $(\infty,1)$-spectrum from \cite{BGT}:
\[
\xymatrix@R+6pc@C+0pc{
\mathrm{K}^\mathrm{BGT}_{\mathrm{sComm}_\mathrm{simplicial}\mathrm{Ind}\mathrm{Banach}_{\Gamma_k}(\mathcal{O}_{X_{R_k,-}})\ar[r]^{\mathrm{global}}}(-)\ar[d]\ar[d]\ar[d]\ar[d]\ar[r]\ar[r] &\mathrm{K}^\mathrm{BGT}_{\mathrm{sComm}_\mathrm{simplicial}\varphi\mathrm{Ind}\mathrm{Banach}_{\Gamma_k}(\{\mathrm{Robba}^\mathrm{extended}_{{R_k,-,I}}\}_I)}(-)\ar[d]\ar[d]\ar[d]\ar[d]\\
\mathrm{K}^\mathrm{BGT}_{\mathrm{sComm}_\mathrm{simplicial}\mathrm{Ind}\mathrm{Banach}_{\Gamma_0}(\mathcal{O}_{X_{\mathbb{Q}_p(p^{1/p^\infty})^{\wedge\flat},-}})\ar[r]^{\mathrm{global}}}(-)\ar[r]\ar[r] &\mathrm{K}^\mathrm{BGT}_{\mathrm{sComm}_\mathrm{simplicial}\varphi\mathrm{Ind}\mathrm{Banach}_{\Gamma_0}(\{\mathrm{Robba}^\mathrm{extended}_{{R_0,-,I}}\}_I)}(-), 
}
\]
\[
\xymatrix@R+6pc@C+0pc{
\mathrm{K}^\mathrm{BGT}_{\mathrm{sComm}_\mathrm{simplicial}\mathrm{Ind}^m\mathrm{Banach}_{\Gamma_k}(\mathcal{O}_{X_{R_k,-}})\ar[r]^{\mathrm{global}}}(-)\ar[d]\ar[d]\ar[d]\ar[d]\ar[r]\ar[r] &\mathrm{K}^\mathrm{BGT}_{\mathrm{sComm}_\mathrm{simplicial}\varphi\mathrm{Ind}^m\mathrm{Banach}_{\Gamma_k}(\{\mathrm{Robba}^\mathrm{extended}_{{R_k,-,I}}\}_I)}(-)\ar[d]\ar[d]\ar[d]\ar[d]\\
\mathrm{K}^\mathrm{BGT}_{\mathrm{sComm}_\mathrm{simplicial}\mathrm{Ind}^m\mathrm{Banach}_{\Gamma_0}(\mathcal{O}_{X_{\mathbb{Q}_p(p^{1/p^\infty})^{\wedge\flat},-}})\ar[r]^{\mathrm{global}}}(-)\ar[r]\ar[r] &\mathrm{K}^\mathrm{BGT}_{\mathrm{sComm}_\mathrm{simplicial}\varphi\mathrm{Ind}^m\mathrm{Banach}_{\Gamma_0}(\{\mathrm{Robba}^\mathrm{extended}_{{R_0,-,I}}\}_I)}(-). 
}
\]
The definition is given by the following:
\[
\xymatrix@R+0pc@C+0pc{
\mathrm{homotopycolimit}_i(\square)(\mathcal{O}_i),  
}
\]
each $\mathcal{O}_i$ is just as $\mathbb{Q}_p\left<C_1,...,C_\ell\right>,\ell=1,2,...$ and $\square$ is the relative diagram of $\infty$-functors.

\end{itemize}

\

\begin{remark}
\noindent We can certainly consider the quasicoherent sheaves in \cite[Lemma 7.11, Remark 7.12]{1BBK}, therefore all the quasicoherent presheaves and modules will be those satisfying \cite[Lemma 7.11, Remark 7.12]{1BBK} if one would like to consider the the quasicoherent sheaves. That being all as this said, we would believe that the big quasicoherent presheaves are automatically quasicoherent sheaves (namely satisfying the corresponding \v{C}ech $\infty$-descent as in \cite[Section 9.3]{KKM} and \cite[Lemma 7.11, Remark 7.12]{1BBK}) and the corresponding global section functors are automatically equivalence of $\infty$-categories. 
\end{remark}

\

\indent In Clausen-Scholze formalism we have the following\footnote{Certainly the homotopy colimit in the rings side will be within the condensed solid animated analytic rings from \cite{1CS2}.}:

\begin{itemize}
\item (\text{Proposition}) There is a functor (global section) between the $\infty$-prestacks of inductive Banach quasicoherent sheaves:
\[
\xymatrix@R+0pc@C+0pc{
{\mathrm{Modules}_\circledcirc}(\mathcal{O}_{X_{R,-}})\ar[r]^{\mathrm{global}}\ar[r]\ar[r] &\varphi{\mathrm{Modules}_\circledcirc}(\{\mathrm{Robba}^\mathrm{extended}_{{R,-,I}}\}_I).  
}
\]
The definition is given by the following:
\[
\xymatrix@R+0pc@C+0pc{
\mathrm{homotopycolimit}_i({\mathrm{Modules}_\circledcirc}(\mathcal{O}_{X_{R,-}})\ar[r]^{\mathrm{global}}\ar[r]\ar[r] &\varphi{\mathrm{Modules}_\circledcirc}(\{\mathrm{Robba}^\mathrm{extended}_{{R,-,I}}\}_I))(\mathcal{O}_i),  
}
\]
each $\mathcal{O}_i$ is just as $\mathbb{Q}_p\left<C_1,...,C_\ell\right>,\ell=1,2,...$.

\item (\text{Proposition}) There is a functor (global section) between the $\infty$-prestacks of inductive Banach quasicoherent sheaves:
\[
\xymatrix@R+0pc@C+0pc{
{\mathrm{Modules}_\circledcirc}(\mathcal{O}_{X_{R,-}})\ar[r]^{\mathrm{global}}\ar[r]\ar[r] &\varphi{\mathrm{Modules}_\circledcirc}(\{\mathrm{Robba}^\mathrm{extended}_{{R,-,I}}\}_I).  
}
\]
The definition is given by the following:
\[
\xymatrix@R+0pc@C+0pc{
\mathrm{homotopycolimit}_i({\mathrm{Modules}_\circledcirc}(\mathcal{O}_{X_{R,-}})\ar[r]^{\mathrm{global}}\ar[r]\ar[r] &\varphi{\mathrm{Modules}_\circledcirc}(\{\mathrm{Robba}^\mathrm{extended}_{{R,-,I}}\}_I))(\mathcal{O}_i),  
}
\]
each $\mathcal{O}_i$ is just as $\mathbb{Q}_p\left<C_1,...,C_\ell\right>,\ell=1,2,...$.

\item (\text{Proposition}) There is a functor (global section) between the $\infty$-prestacks of inductive Banach quasicoherent commutative algebra $E_\infty$ objects\footnote{Here $\circledcirc=\text{solidquasicoherentsheaves}$.}:
\[
\xymatrix@R+0pc@C+0pc{
\mathrm{sComm}_\mathrm{simplicial}{\mathrm{Modules}_\circledcirc}(\mathcal{O}_{X_{R,-}})\ar[r]^{\mathrm{global}}\ar[r]\ar[r] &\mathrm{sComm}_\mathrm{simplicial}\varphi{\mathrm{Modules}_\circledcirc}(\{\mathrm{Robba}^\mathrm{extended}_{{R,-,I}}\}_I).  
}
\]
The definition is given by the following:
\[
\xymatrix@R+0pc@C+0pc{
\mathrm{homotopycolimit}_i(\mathrm{sComm}_\mathrm{simplicial}{\mathrm{Modules}_\circledcirc}(\mathcal{O}_{X_{R,-}})\ar[r]^{\mathrm{global}}\ar[r]\ar[r] &\mathrm{sComm}_\mathrm{simplicial}\varphi{\mathrm{Modules}_\circledcirc}(\{\mathrm{Robba}^\mathrm{extended}_{{R,-,I}}\}_I))(\mathcal{O}_i),  
}
\]
each $\mathcal{O}_i$ is just as $\mathbb{Q}_p\left<C_1,...,C_\ell\right>,\ell=1,2,...$.
\item Then as in \cite{LBV} we have a functor (global section ) of the de Rham complex after \cite[Definition 5.9, Section 5.2.1]{KKM}\footnote{Here $\circledcirc=\text{solidquasicoherentsheaves}$.}:
\[
\xymatrix@R+0pc@C+0pc{
\mathrm{deRham}_{\mathrm{sComm}_\mathrm{simplicial}{\mathrm{Modules}_\circledcirc}(\mathcal{O}_{X_{R,-}})\ar[r]^{\mathrm{global}}}(-)\ar[r]\ar[r] &\mathrm{deRham}_{\mathrm{sComm}_\mathrm{simplicial}\varphi{\mathrm{Modules}_\circledcirc}(\{\mathrm{Robba}^\mathrm{extended}_{{R,-,I}}\}_I)}(-), 
}
\]
The definition is given by the following:
\[
\xymatrix@R+0pc@C+0pc{
\mathrm{homotopycolimit}_i\\
(\mathrm{deRham}_{\mathrm{sComm}_\mathrm{simplicial}{\mathrm{Modules}_\circledcirc}(\mathcal{O}_{X_{R,-}})\ar[r]^{\mathrm{global}}}(-)\ar[r]\ar[r] &\mathrm{deRham}_{\mathrm{sComm}_\mathrm{simplicial}\varphi{\mathrm{Modules}_\circledcirc}(\{\mathrm{Robba}^\mathrm{extended}_{{R,-,I}}\}_I)}(-))(\mathcal{O}_i),  
}
\]
each $\mathcal{O}_i$ is just as $\mathbb{Q}_p\left<C_1,...,C_\ell\right>,\ell=1,2,...$.\item Then we have the following a functor (global section) of $K$-group $(\infty,1)$-spectrum from \cite{BGT}\footnote{Here $\circledcirc=\text{solidquasicoherentsheaves}$.}:
\[
\xymatrix@R+0pc@C+0pc{
\mathrm{K}^\mathrm{BGT}_{\mathrm{sComm}_\mathrm{simplicial}{\mathrm{Modules}_\circledcirc}(\mathcal{O}_{X_{R,-}})\ar[r]^{\mathrm{global}}}(-)\ar[r]\ar[r] &\mathrm{K}^\mathrm{BGT}_{\mathrm{sComm}_\mathrm{simplicial}\varphi{\mathrm{Modules}_\circledcirc}(\{\mathrm{Robba}^\mathrm{extended}_{{R,-,I}}\}_I)}(-). 
}
\]
The definition is given by the following:
\[\displayindent=-0.4in
\xymatrix@R+0pc@C+0pc{
\mathrm{homotopycolimit}_i(\mathrm{K}^\mathrm{BGT}_{\mathrm{sComm}_\mathrm{simplicial}{\mathrm{Modules}_\circledcirc}(\mathcal{O}_{X_{R,-}})\ar[r]^{\mathrm{global}}}(-)\ar[r]\ar[r] &\mathrm{K}^\mathrm{BGT}_{\mathrm{sComm}_\mathrm{simplicial}\varphi{\mathrm{Modules}_\circledcirc}(\{\mathrm{Robba}^\mathrm{extended}_{{R,-,I}}\}_I)}(-))(\mathcal{O}_i),  
}
\]
each $\mathcal{O}_i$ is just as $\mathbb{Q}_p\left<C_1,...,C_\ell\right>,\ell=1,2,...$.
\end{itemize}

\noindent Now let $R=\mathbb{Q}_p(p^{1/p^\infty})^{\wedge\flat}$ and $R_k=\mathbb{Q}_p(p^{1/p^\infty})^{\wedge}\left<T_1^{\pm 1/p^{\infty}},...,T_k^{\pm 1/p^{\infty}}\right>^\flat$ we have the following Galois theoretic results with naturality along $f:\mathrm{Spa}(\mathbb{Q}_p(p^{1/p^\infty})^{\wedge}\left<T_1^{\pm 1/p^\infty},...,T_k^{\pm 1/p^\infty}\right>^\flat)\rightarrow \mathrm{Spa}(\mathbb{Q}_p(p^{1/p^\infty})^{\wedge\flat})$:

\begin{itemize}
\item (\text{Proposition}) There is a functor (global section) between the $\infty$-prestacks of inductive Banach quasicoherent sheaves\footnote{Here $\circledcirc=\text{solidquasicoherentsheaves}$.}:
\[
\xymatrix@R+6pc@C+0pc{
{\mathrm{Modules}_\circledcirc}(\mathcal{O}_{X_{\mathbb{Q}_p(p^{1/p^\infty})^{\wedge}\left<T_1^{\pm 1/p^\infty},...,T_k^{\pm 1/p^\infty}\right>^\flat,-}})\ar[d]\ar[d]\ar[d]\ar[d] \ar[r]^{\mathrm{global}}\ar[r]\ar[r] &\varphi{\mathrm{Modules}_\circledcirc}(\{\mathrm{Robba}^\mathrm{extended}_{{R_k,-,I}}\}_I) \ar[d]\ar[d]\ar[d]\ar[d].\\
{\mathrm{Modules}_\circledcirc}(\mathcal{O}_{X_{\mathbb{Q}_p(p^{1/p^\infty})^{\wedge\flat},-}})\ar[r]^{\mathrm{global}}\ar[r]\ar[r] &\varphi{\mathrm{Modules}_\circledcirc}(\{\mathrm{Robba}^\mathrm{extended}_{{R_0,-,I}}\}_I).\\ 
}
\]
The definition is given by the following:
\[
\xymatrix@R+0pc@C+0pc{
\mathrm{homotopycolimit}_i(\square)(\mathcal{O}_i),  
}
\]
each $\mathcal{O}_i$ is just as $\mathbb{Q}_p\left<C_1,...,C_\ell\right>,\ell=1,2,...$ and $\square$ is the relative diagram of $\infty$-functors.

\item (\text{Proposition}) There is a functor (global section) between the $\infty$-prestacks of inductive Banach quasicoherent commutative algebra $E_\infty$ objects\footnote{Here $\circledcirc=\text{solidquasicoherentsheaves}$.}:
\[
\xymatrix@R+6pc@C+0pc{
\mathrm{sComm}_\mathrm{simplicial}{\mathrm{Modules}_\circledcirc}(\mathcal{O}_{X_{R_k,-}})\ar[d]\ar[d]\ar[d]\ar[d]\ar[r]^{\mathrm{global}}\ar[r]\ar[r] &\mathrm{sComm}_\mathrm{simplicial}\varphi{\mathrm{Modules}_\circledcirc}(\{\mathrm{Robba}^\mathrm{extended}_{{R_k,-,I}}\}_I)\ar[d]\ar[d]\ar[d]\ar[d]\\
\mathrm{sComm}_\mathrm{simplicial}{\mathrm{Modules}_\circledcirc}(\mathcal{O}_{X_{\mathbb{Q}_p(p^{1/p^\infty})^{\wedge\flat},-}})\ar[r]^{\mathrm{global}}\ar[r]\ar[r] &\mathrm{sComm}_\mathrm{simplicial}\varphi{\mathrm{Modules}_\circledcirc}(\{\mathrm{Robba}^\mathrm{extended}_{{R_0,-,I}}\}_I).  
}
\]
The definition is given by the following:
\[
\xymatrix@R+0pc@C+0pc{
\mathrm{homotopycolimit}_i(\square)(\mathcal{O}_i),  
}
\]
each $\mathcal{O}_i$ is just as $\mathbb{Q}_p\left<C_1,...,C_\ell\right>,\ell=1,2,...$ and $\square$ is the relative diagram of $\infty$-functors.

\item Then as in \cite{LBV} we have a functor (global section) of the de Rham complex after \cite[Definition 5.9, Section 5.2.1]{KKM}\footnote{Here $\circledcirc=\text{solidquasicoherentsheaves}$.}:
\[\displayindent=-0.4in
\xymatrix@R+6pc@C+0pc{
\mathrm{deRham}_{\mathrm{sComm}_\mathrm{simplicial}{\mathrm{Modules}_\circledcirc}(\mathcal{O}_{X_{R_k,-}})\ar[r]^{\mathrm{global}}}(-)\ar[d]\ar[d]\ar[d]\ar[d]\ar[r]\ar[r] &\mathrm{deRham}_{\mathrm{sComm}_\mathrm{simplicial}\varphi{\mathrm{Modules}_\circledcirc}(\{\mathrm{Robba}^\mathrm{extended}_{{R_k,-,I}}\}_I)}(-)\ar[d]\ar[d]\ar[d]\ar[d]\\
\mathrm{deRham}_{\mathrm{sComm}_\mathrm{simplicial}{\mathrm{Modules}_\circledcirc}(\mathcal{O}_{X_{\mathbb{Q}_p(p^{1/p^\infty})^{\wedge\flat},-}})\ar[r]^{\mathrm{global}}}(-)\ar[r]\ar[r] &\mathrm{deRham}_{\mathrm{sComm}_\mathrm{simplicial}\varphi{\mathrm{Modules}_\circledcirc}(\{\mathrm{Robba}^\mathrm{extended}_{{R_0,-,I}}\}_I)}(-), 
}
\]

\item Then we have the following a functor (global section) of $K$-group $(\infty,1)$-spectrum from \cite{BGT}\footnote{Here $\circledcirc=\text{solidquasicoherentsheaves}$.}:
\[
\xymatrix@R+6pc@C+0pc{
\mathrm{K}^\mathrm{BGT}_{\mathrm{sComm}_\mathrm{simplicial}{\mathrm{Modules}_\circledcirc}(\mathcal{O}_{X_{R_k,-}})\ar[r]^{\mathrm{global}}}(-)\ar[d]\ar[d]\ar[d]\ar[d]\ar[r]\ar[r] &\mathrm{K}^\mathrm{BGT}_{\mathrm{sComm}_\mathrm{simplicial}\varphi{\mathrm{Modules}_\circledcirc}(\{\mathrm{Robba}^\mathrm{extended}_{{R_k,-,I}}\}_I)}(-)\ar[d]\ar[d]\ar[d]\ar[d]\\
\mathrm{K}^\mathrm{BGT}_{\mathrm{sComm}_\mathrm{simplicial}{\mathrm{Modules}_\circledcirc}(\mathcal{O}_{X_{\mathbb{Q}_p(p^{1/p^\infty})^{\wedge\flat},-}})\ar[r]^{\mathrm{global}}}(-)\ar[r]\ar[r] &\mathrm{K}^\mathrm{BGT}_{\mathrm{sComm}_\mathrm{simplicial}\varphi{\mathrm{Modules}_\circledcirc}(\{\mathrm{Robba}^\mathrm{extended}_{{R_0,-,I}}\}_I)}(-), 
}
\]

The definition is given by the following:
\[
\xymatrix@R+0pc@C+0pc{
\mathrm{homotopycolimit}_i(\square)(\mathcal{O}_i),  
}
\]
each $\mathcal{O}_i$ is just as $\mathbb{Q}_p\left<C_1,...,C_\ell\right>,\ell=1,2,...$ and $\square$ is the relative diagram of $\infty$-functors.

\end{itemize}

\
\indent Then we consider further equivariance by considering the arithmetic profinite fundamental groups $\Gamma_{\mathbb{Q}_p}$ and $\mathrm{Gal}(\overline{{Q}_p\left<T_1^{\pm 1},...,T_k^{\pm 1}\right>}/R_k)$ through the following diagram:

\[
\xymatrix@R+0pc@C+0pc{
\mathbb{Z}_p^k=\mathrm{Gal}(R_k/{\mathbb{Q}_p(p^{1/p^\infty})^\wedge\left<T_1^{\pm 1},...,T_k^{\pm 1}\right>}) \ar[r]\ar[r] \ar[r]\ar[r] &\Gamma_k:=\mathrm{Gal}(R_k/{\mathbb{Q}_p\left<T_1^{\pm 1},...,T_k^{\pm 1}\right>}) \ar[r] \ar[r]\ar[r] &\Gamma_{\mathbb{Q}_p}.
}
\]

\begin{itemize}
\item (\text{Proposition}) There is a functor (global section) between the $\infty$-prestacks of inductive Banach quasicoherent sheaves\footnote{Here $\circledcirc=\text{solidquasicoherentsheaves}$.}:
\[
\xymatrix@R+6pc@C+0pc{
{\mathrm{Modules}_\circledcirc}_{\Gamma_k}(\mathcal{O}_{X_{\mathbb{Q}_p(p^{1/p^\infty})^{\wedge}\left<T_1^{\pm 1/p^\infty},...,T_k^{\pm 1/p^\infty}\right>^\flat,-}})\ar[d]\ar[d]\ar[d]\ar[d] \ar[r]^{\mathrm{global}}\ar[r]\ar[r] &\varphi{\mathrm{Modules}_\circledcirc}_{\Gamma_k}(\{\mathrm{Robba}^\mathrm{extended}_{{R_k,-,I}}\}_I) \ar[d]\ar[d]\ar[d]\ar[d].\\
{\mathrm{Modules}_\circledcirc}(\mathcal{O}_{X_{\mathbb{Q}_p(p^{1/p^\infty})^{\wedge\flat},-}})\ar[r]^{\mathrm{global}}\ar[r]\ar[r] &\varphi{\mathrm{Modules}_\circledcirc}(\{\mathrm{Robba}^\mathrm{extended}_{{R_0,-,I}}\}_I).\\ 
}
\]
The definition is given by the following:
\[
\xymatrix@R+0pc@C+0pc{
\mathrm{homotopycolimit}_i(\square)(\mathcal{O}_i),  
}
\]
each $\mathcal{O}_i$ is just as $\mathbb{Q}_p\left<C_1,...,C_\ell\right>,\ell=1,2,...$ and $\square$ is the relative diagram of $\infty$-functors.

\item (\text{Proposition}) There is a functor (global section) between the $\infty$-stacks of inductive Banach quasicoherent commutative algebra $E_\infty$ objects\footnote{Here $\circledcirc=\text{solidquasicoherentsheaves}$.}:
\[
\xymatrix@R+6pc@C+0pc{
\mathrm{sComm}_\mathrm{simplicial}{\mathrm{Modules}_\circledcirc}_{\Gamma_k}(\mathcal{O}_{X_{R_k,-}})\ar[d]\ar[d]\ar[d]\ar[d]\ar[r]^{\mathrm{global}}\ar[r]\ar[r] &\mathrm{sComm}_\mathrm{simplicial}\varphi{\mathrm{Modules}_\circledcirc}_{\Gamma_k}(\{\mathrm{Robba}^\mathrm{extended}_{{R_k,-,I}}\}_I)\ar[d]\ar[d]\ar[d]\ar[d]\\
\mathrm{sComm}_\mathrm{simplicial}{\mathrm{Modules}_\circledcirc}_{\Gamma_0}(\mathcal{O}_{X_{\mathbb{Q}_p(p^{1/p^\infty})^{\wedge\flat},-}})\ar[r]^{\mathrm{global}}\ar[r]\ar[r] &\mathrm{sComm}_\mathrm{simplicial}\varphi{\mathrm{Modules}_\circledcirc}_{\Gamma_0}(\{\mathrm{Robba}^\mathrm{extended}_{{R_0,-,I}}\}_I).  
}
\]
The definition is given by the following:
\[
\xymatrix@R+0pc@C+0pc{
\mathrm{homotopycolimit}_i(\square)(\mathcal{O}_i),  
}
\]
each $\mathcal{O}_i$ is just as $\mathbb{Q}_p\left<C_1,...,C_\ell\right>,\ell=1,2,...$ and $\square$ is the relative diagram of $\infty$-functors.

\item Then as in \cite{LBV} we have a functor (global section) of the de Rham complex after \cite[Definition 5.9, Section 5.2.1]{KKM}\footnote{Here $\circledcirc=\text{solidquasicoherentsheaves}$.}:
\[\displayindent=-0.4in
\xymatrix@R+6pc@C+0pc{
\mathrm{deRham}_{\mathrm{sComm}_\mathrm{simplicial}{\mathrm{Modules}_\circledcirc}_{\Gamma_k}(\mathcal{O}_{X_{R_k,-}})\ar[r]^{\mathrm{global}}}(-)\ar[d]\ar[d]\ar[d]\ar[d]\ar[r]\ar[r] &\mathrm{deRham}_{\mathrm{sComm}_\mathrm{simplicial}\varphi{\mathrm{Modules}_\circledcirc}_{\Gamma_k}(\{\mathrm{Robba}^\mathrm{extended}_{{R_k,-,I}}\}_I)}(-)\ar[d]\ar[d]\ar[d]\ar[d]\\
\mathrm{deRham}_{\mathrm{sComm}_\mathrm{simplicial}{\mathrm{Modules}_\circledcirc}_{\Gamma_0}(\mathcal{O}_{X_{\mathbb{Q}_p(p^{1/p^\infty})^{\wedge\flat},-}})\ar[r]^{\mathrm{global}}}(-)\ar[r]\ar[r] &\mathrm{deRham}_{\mathrm{sComm}_\mathrm{simplicial}\varphi{\mathrm{Modules}_\circledcirc}_{\Gamma_0}(\{\mathrm{Robba}^\mathrm{extended}_{{R_0,-,I}}\}_I)}(-), 
}
\]

The definition is given by the following:
\[
\xymatrix@R+0pc@C+0pc{
\mathrm{homotopycolimit}_i(\square)(\mathcal{O}_i),  
}
\]
each $\mathcal{O}_i$ is just as $\mathbb{Q}_p\left<C_1,...,C_\ell\right>,\ell=1,2,...$ and $\square$ is the relative diagram of $\infty$-functors.

\item Then we have the following a functor (global section) of $K$-group $(\infty,1)$-spectrum from \cite{BGT}\footnote{Here $\circledcirc=\text{solidquasicoherentsheaves}$.}:
\[
\xymatrix@R+6pc@C+0pc{
\mathrm{K}^\mathrm{BGT}_{\mathrm{sComm}_\mathrm{simplicial}{\mathrm{Modules}_\circledcirc}_{\Gamma_k}(\mathcal{O}_{X_{R_k,-}})\ar[r]^{\mathrm{global}}}(-)\ar[d]\ar[d]\ar[d]\ar[d]\ar[r]\ar[r] &\mathrm{K}^\mathrm{BGT}_{\mathrm{sComm}_\mathrm{simplicial}\varphi{\mathrm{Modules}_\circledcirc}_{\Gamma_k}(\{\mathrm{Robba}^\mathrm{extended}_{{R_k,-,I}}\}_I)}(-)\ar[d]\ar[d]\ar[d]\ar[d]\\
\mathrm{K}^\mathrm{BGT}_{\mathrm{sComm}_\mathrm{simplicial}{\mathrm{Modules}_\circledcirc}_{\Gamma_0}(\mathcal{O}_{X_{\mathbb{Q}_p(p^{1/p^\infty})^{\wedge\flat},-}})\ar[r]^{\mathrm{global}}}(-)\ar[r]\ar[r] &\mathrm{K}^\mathrm{BGT}_{\mathrm{sComm}_\mathrm{simplicial}\varphi{\mathrm{Modules}_\circledcirc}_{\Gamma_0}(\{\mathrm{Robba}^\mathrm{extended}_{{R_0,-,I}}\}_I)}(-), 
}
\]

The definition is given by the following:
\[
\xymatrix@R+0pc@C+0pc{
\mathrm{homotopycolimit}_i(\square)(\mathcal{O}_i),  
}
\]
each $\mathcal{O}_i$ is just as $\mathbb{Q}_p\left<C_1,...,C_\ell\right>,\ell=1,2,...$ and $\square$ is the relative diagram of $\infty$-functors.\\

\end{itemize}

\begin{proposition}
All the global functors from \cite[Proposition 13.8, Theorem 14.9, Remark 14.10]{1CS2} achieve the equivalences on both sides.	\\
\end{proposition}

\newpage
\subsection{$\infty$-Categorical Analytic Stacks and Descents V}

Here we consider the corresponding archimedean picture, after \cite[Problem A.4, Kedlaya's Lecture]{1CBCKSW}. Recall for any algebraic variety $R$ over $\mathbb{R}$ this $X_R(\mathbb{C})$ is defined to be the corresponding quotient:
\begin{align}
X_{R}(\mathbb{C}):=R(\mathbb{C})\times P^1(\mathbb{C})/\varphi,\\
Y_R(\mathbb{C}):=R(\mathbb{C})\times P^1(\mathbb{C}).	
\end{align}
The Hodge structure is given by $\varphi$. We define the relative version by considering a further algebraic variety over $\mathbb{C}$, say $A$ as in the following:
\begin{align}
X_{R,A}(\mathbb{C}):=R(\mathbb{C})\times P^1(\mathbb{C})\times A(\mathbb{C})/\varphi,\\
Y_{R,A}(\mathbb{C}):=R(\mathbb{C})\times P^1(\mathbb{C})\times A(\mathbb{C}).	
\end{align}

Then by \cite{1BBK} and \cite{1CS2} we have the corresponding $\infty$-category of $\infty$-sheaves of simplicial ind-Banach quasicoherent modules which in our situation will be assumed to the modules in \cite{1BBK}, as well as the corresponding associated Clausen-Scholze spaces:
\begin{align}
X_{R}(\mathbb{C}):=R(\mathbb{C})\times P^1(\mathbb{C})^\blacksquare/\varphi,\\
Y_R(\mathbb{C}):=R(\mathbb{C})\times P^1(\mathbb{C})^\blacksquare.	
\end{align}
\begin{align}
X_{R,A}(\mathbb{C}):=R(\mathbb{C})\times P^1(\mathbb{C})\times A(\mathbb{C})^\blacksquare/\varphi,\\
Y_{R,A}(\mathbb{C}):=R(\mathbb{C})\times P^1(\mathbb{C})\times A(\mathbb{C})^\blacksquare,	
\end{align}
with the $\infty$-category of $\infty$-sheaves of simplicial liquid quasicoherent modules, liquid vector bundles and liquid perfect complexes, with further descent \cite[Proposition 13.8, Theorem 14.9, Remark 14.10]{1CS2}. We call the resulting global sections are the corresponding $c$-equivariant Hodge Modules. Then we have the following direct analogy:

\begin{itemize}
\item (\text{Proposition}) There is an equivalence between the $\infty$-categories of inductive Banach quasicoherent presheaves:
\[
\xymatrix@R+0pc@C+0pc{
\mathrm{Ind}\mathrm{Banach}(\mathcal{O}_{X_{R,A}})\ar[r]^{\mathrm{equi}}\ar[r]\ar[r] &\varphi\mathrm{Ind}\mathrm{Banach}(\mathcal{O}_{Y_{R,A}}).  
}
\]
\item (\text{Proposition}) There is an equivalence between the $\infty$-categories of monomorphic inductive Banach quasicoherent presheaves:
\[
\xymatrix@R+0pc@C+0pc{
\mathrm{Ind}^m\mathrm{Banach}(\mathcal{O}_{X_{R,A}})\ar[r]^{\mathrm{equi}}\ar[r]\ar[r] &\varphi\mathrm{Ind}^m\mathrm{Banach}(\mathcal{O}_{Y_{R,A}}).  
}
\]
\end{itemize}

\begin{itemize}

\item (\text{Proposition}) There is an equivalence between the $\infty$-categories of inductive Banach quasicoherent presheaves:
\[
\xymatrix@R+0pc@C+0pc{
\mathrm{Ind}\mathrm{Banach}(\mathcal{O}_{X_{R,A}})\ar[r]^{\mathrm{equi}}\ar[r]\ar[r] &\varphi\mathrm{Ind}\mathrm{Banach}(\mathcal{O}_{Y_{R,A}}).  
}
\]
\item (\text{Proposition}) There is an equivalence between the $\infty$-categories of monomorphic inductive Banach quasicoherent presheaves:
\[
\xymatrix@R+0pc@C+0pc{
\mathrm{Ind}^m\mathrm{Banach}(\mathcal{O}_{X_{R,A}})\ar[r]^{\mathrm{equi}}\ar[r]\ar[r] &\varphi\mathrm{Ind}^m\mathrm{Banach}(\mathcal{O}_{Y_{R,A}}).  
}
\]
\item (\text{Proposition}) There is an equivalence between the $\infty$-categories of inductive Banach quasicoherent commutative algebra $E_\infty$ objects:
\[
\xymatrix@R+0pc@C+0pc{
\mathrm{sComm}_\mathrm{simplicial}\mathrm{Ind}\mathrm{Banach}(\mathcal{O}_{X_{R,A}})\ar[r]^{\mathrm{equi}}\ar[r]\ar[r] &\mathrm{sComm}_\mathrm{simplicial}\varphi\mathrm{Ind}\mathrm{Banach}(\mathcal{O}_{Y_{R,A}}).  
}
\]
\item (\text{Proposition}) There is an equivalence between the $\infty$-categories of monomorphic inductive Banach quasicoherent commutative algebra $E_\infty$ objects:
\[
\xymatrix@R+0pc@C+0pc{
\mathrm{sComm}_\mathrm{simplicial}\mathrm{Ind}^m\mathrm{Banach}(\mathcal{O}_{X_{R,A}})\ar[r]^{\mathrm{equi}}\ar[r]\ar[r] &\mathrm{sComm}_\mathrm{simplicial}\varphi\mathrm{Ind}^m\mathrm{Banach}(\mathcal{O}_{Y_{R,A}}).  
}
\]

\item Then parallel as in \cite{LBV} we have the equivalence of the de Rham complex after \cite[Definition 5.9, Section 5.2.1]{KKM}:
\[
\xymatrix@R+0pc@C+0pc{
\mathrm{deRham}_{\mathrm{sComm}_\mathrm{simplicial}\mathrm{Ind}\mathrm{Banach}(\mathcal{O}_{X_{R,A}})\ar[r]^{\mathrm{equi}}}(-)\ar[r]\ar[r] &\mathrm{deRham}_{\mathrm{sComm}_\mathrm{simplicial}\varphi\mathrm{Ind}\mathrm{Banach}(\mathcal{O}_{Y_{R,A}})}(-), 
}
\]
\[
\xymatrix@R+0pc@C+0pc{
\mathrm{deRham}_{\mathrm{sComm}_\mathrm{simplicial}\mathrm{Ind}^m\mathrm{Banach}(\mathcal{O}_{X_{R,A}})\ar[r]^{\mathrm{equi}}}(-)\ar[r]\ar[r] &\mathrm{deRham}_{\mathrm{sComm}_\mathrm{simplicial}\varphi\mathrm{Ind}^m\mathrm{Banach}(\mathcal{O}_{Y_{R,A}})}(-). 
}
\]

\item Then we have the following equivalence of $K$-group $(\infty,1)$-spectrum from \cite{BGT}:
\[
\xymatrix@R+0pc@C+0pc{
\mathrm{K}^\mathrm{BGT}_{\mathrm{sComm}_\mathrm{simplicial}\mathrm{Ind}\mathrm{Banach}(\mathcal{O}_{X_{R,A}})\ar[r]^{\mathrm{equi}}}(-)\ar[r]\ar[r] &\mathrm{K}^\mathrm{BGT}_{\mathrm{sComm}_\mathrm{simplicial}\varphi\mathrm{Ind}\mathrm{Banach}(\mathcal{O}_{Y_{R,A}})}(-), 
}
\]
\[
\xymatrix@R+0pc@C+0pc{
\mathrm{K}^\mathrm{BGT}_{\mathrm{sComm}_\mathrm{simplicial}\mathrm{Ind}^m\mathrm{Banach}(\mathcal{O}_{X_{R,A}})\ar[r]^{\mathrm{equi}}}(-)\ar[r]\ar[r] &\mathrm{K}^\mathrm{BGT}_{\mathrm{sComm}_\mathrm{simplicial}\varphi\mathrm{Ind}^m\mathrm{Banach}(\mathcal{O}_{Y_{R,A}})}(-). 
}
\]
\end{itemize}

\begin{assumption}\label{assumtionpresheaves}
All the functors of modules or algebras below are Clausen-Scholze sheaves \cite[Proposition 13.8, Theorem 14.9, Remark 14.10]{1CS2}.	
\end{assumption}

\begin{itemize}
\item (\text{Proposition}) There is an equivalence between the $\infty$-categories of inductive liquid sheaves:
\[
\xymatrix@R+0pc@C+0pc{
\mathrm{Module}_\circledcirc(\mathcal{O}_{X_{R,A}})\ar[r]^{\mathrm{equi}}\ar[r]\ar[r] &\varphi\mathrm{Module}_\circledcirc(\mathcal{O}_{Y_{R,A}}).  
}
\]
\end{itemize}

\begin{itemize}

\item (\text{Proposition}) There is an equivalence between the $\infty$-categories of inductive Banach quasicoherent commutative algebra $E_\infty$ objects:
\[\displayindent=-0.4in
\xymatrix@R+0pc@C+0pc{
\mathrm{sComm}_\mathrm{simplicial}\mathrm{Module}_{\text{liquidquasicoherentsheaves}}(\mathcal{O}_{X_{R,A}})\ar[r]^{\mathrm{equi}}\ar[r]\ar[r] &\mathrm{sComm}_\mathrm{simplicial}\varphi\mathrm{Module}_{\text{liquidquasicoherentsheaves}}(\mathcal{O}_{Y_{R,A}}).  
}
\]

\item Then as in \cite{LBV} we have the equivalence of the de Rham complex after \cite[Definition 5.9, Section 5.2.1]{KKM}\footnote{Here $\circledcirc=\text{liquidquasicoherentsheaves}$.}:
\[
\xymatrix@R+0pc@C+0pc{
\mathrm{deRham}_{\mathrm{sComm}_\mathrm{simplicial}\mathrm{Module}_\circledcirc(\mathcal{O}_{X_{R,A}})\ar[r]^{\mathrm{equi}}}(-)\ar[r]\ar[r] &\mathrm{deRham}_{\mathrm{sComm}_\mathrm{simplicial}\varphi\mathrm{Module}_\circledcirc(\mathcal{O}_{Y_{R,A}})}(-). 
}
\]

\item Then we have the following equivalence of $K$-group $(\infty,1)$-spectrum from \cite{BGT}\footnote{Here $\circledcirc=\text{liquidquasicoherentsheaves}$.}:
\[
\xymatrix@R+0pc@C+0pc{
\mathrm{K}^\mathrm{BGT}_{\mathrm{sComm}_\mathrm{simplicial}\mathrm{Module}_\circledcirc(\mathcal{O}_{X_{R,A}})\ar[r]^{\mathrm{equi}}}(-)\ar[r]\ar[r] &\mathrm{K}^\mathrm{BGT}_{\mathrm{sComm}_\mathrm{simplicial}\varphi\mathrm{Module}_\circledcirc(\mathcal{O}_{Y_{R,A}})}(-). 
}
\]
\end{itemize}

\newpage

\chapter{Hodge-Iwasawa Theory: The Extensions}

\section{Introduction to the Interactions among Motives}

\subsection{\text{Equivariant relative $p$-adic Hodge Theory}}

\noindent{\text{Equivariant relative $p$-adic Hodge Theory}}

\begin{itemize}

\item<1-> The corresponding $P$-objects are interesting, but in general are not that easy to study, especially we consider for instance those ring defined over $\mathbb{Q}_p$, let it alone if one would like to consider the categories of the complexes of such objects. 

\item<2-> We choose to consider the corresponding embedding of such objects into the categories of Frobenius sheaves with coefficients in $P$ after Kedlaya-Liu \cite{KL1}, \cite{KL2}. Again we expect everything will be more convenient to handle in the category of $(\varphi,\Gamma)$-modules.
	
\item<3-> Working over $R$ now a uniform Banach algebra with further structure of an adic ring over $\mathbb{F}_p$. And we assume that $R$ is perfect.
Let $\mathrm{Robba}^{\mathrm{extended}}_{I,R}$ be the Robba sheaves defined by Kedlaya-Liu \cite{KL1}, \cite{KL2}, with respect to some interval $I\subset (0,\infty)$, which are Fr\'echet completions of the ring of Witt vector of $R$ with respect to the Gauss norms induced from the norm on $R$. Here we consider the following assumption:

\begin{assumption}
We now assume $R$ comes from the local chart of a rigid space over $\mathbb{Q}_p$. \footnote{This could be made more general, but at this moment let us be closer to classical $p$-adic Hodge Theory.} This will give us the chance to consider the following period rings from \cite{1Sch2} and \cite[Definition 8.6.5]{KL2}:
\begin{align}
B^+_{R,\mathrm{dR}},B_{R,\mathrm{dR}},\\
\mathcal{O}B^+_{R,\mathrm{dR}},\mathcal{O}B_{R,\mathrm{dR}}	
\end{align}
over:
\begin{align}
B^+_{\mathbb{Q}_p(p^{1/p^\infty})^{\wedge,\flat},\mathrm{dR}},B_{\mathbb{Q}_p(p^{1/p^\infty})^{\wedge,\flat},\mathrm{dR}},\\
\mathcal{O}B^+_{\mathbb{Q}_p(p^{1/p^\infty})^{\wedge,\flat},\mathrm{dR}},\mathcal{O}B_{\mathbb{Q}_p(p^{1/p^\infty})^{\wedge,\flat},\mathrm{dR}}	
\end{align}	
where the smaller rings contain element $t=\mathrm{log}([1+\overline{\pi}])$. As in \cite{BCM} we can take the corresponding self $q$-th power product. Then we have by taking the corresponding self $q$-th power product following:
\begin{align}
B^+_{R,\mathrm{dR},q},B_{R,\mathrm{dR},q},\\
\mathcal{O}B^+_{R,\mathrm{dR},q},\mathcal{O}B_{R,\mathrm{dR},q}	
\end{align}
over:
\begin{align}
B^+_{\mathbb{Q}_p(p^{1/p^\infty})^{\wedge,\flat},\mathrm{dR},q},B_{\mathbb{Q}_p(p^{1/p^\infty})^{\wedge,\flat},\mathrm{dR},q},\\
\mathcal{O}B^+_{\mathbb{Q}_p(p^{1/p^\infty})^{\wedge,\flat},\mathrm{dR},q},\mathcal{O}B_{\mathbb{Q}_p(p^{1/p^\infty})^{\wedge,\flat},\mathrm{dR},q}.	
\end{align}
Here we have the action from the product of arithmetic profinite fundamental groups and the product of Frobenius operators.
\end{assumption}

\item<4-> Following Carter-Kedlaya-Z\'abr\'adi and Pal-Z\'abr\'adi \cite{1CKZ} and \cite{1PZ}, taking suitable interval one can define the corresponding Robba rings $\mathrm{Robba}^{\mathrm{extended},q}_{r,R}$, $\mathrm{Robba}^{\mathrm{extended},q}_{\infty,R}$ and the corresponding full Robba ring $\mathrm{Robba}^{\mathrm{extended},q}_{R}$ by the corresponding self $q$-th power product. 

\item<5-> We work in the category of Banach and ind-Fr\'echet spaces, which are commutative. Our generalization comes from those Banach reduced affinoid algebras $A$ over $\mathbb{Q}_p$.

\end{itemize}

\noindent{\text{Equivariant relative $p$-adic Hodge Theory}}
\begin{itemize}

\item<1-> The $p$-adic functional analysis produces us some manageable structures within our study of relative $p$-adic Hodge theory, generalizing the original $p$-adic functional analytic framework of Kedlaya-Liu \cite{KL1}, \cite{KL2}.

\item<2-> Starting from Kedlaya-Liu's period rings after taking product\footnote{When we are talking about $q$-th power as in this chapter, the radius $r$ is then multiradius which are allowed to be different in different components, and the interval $I$ is then multiinterval which are allowed to be different in different components.}, 
\begin{align}
&\mathrm{Robba}^{\mathrm{extended},q}_{\infty,R},\mathrm{Robba}^{\mathrm{extended},q}_{I,R},\mathrm{Robba}^{\mathrm{extended},q}_{r,R},\mathrm{Robba}^{\mathrm{extended},q}_{R},	\mathrm{Robba}^{\mathrm{extended},q}_{{\mathrm{int},r},R},\\
&\mathrm{Robba}^{\mathrm{extended},q}_{{\mathrm{int}},R},\mathrm{Robba}^{\mathrm{extended},q}_{{\mathrm{bd},r},R},\mathrm{Robba}^{\mathrm{extended},q}_{{\mathrm{bd}},R}
\end{align}
we can form the corresponding $A$-relative of the period rings\footnote{Taking products over $\mathbb{Q}_p$.}:
\begin{align}
&\mathrm{Robba}^{\mathrm{extended},q}_{\infty,R,A},\mathrm{Robba}^{\mathrm{extended},q}_{I,R,A},\mathrm{Robba}^{\mathrm{extended},q}_{r,R,A},\mathrm{Robba}^{\mathrm{extended},q}_{R,A},	\mathrm{Robba}^{\mathrm{extended},q}_{{\mathrm{int},r},R,A},\\
&\mathrm{Robba}^{\mathrm{extended},q}_{\mathrm{int},R,A},\mathrm{Robba}^{\mathrm{extended},q}_{{\mathrm{bd},r},R,A},\mathrm{Robba}^{\mathrm{extended},q}_{{\mathrm{bd}},R,A}.	
\end{align} 
\item<3-> (\text{Remark}) There should be also many interesting contexts, for instance consider a finitely generated abelian group $G$, one can consider the group rings: 
\begin{align}
\mathrm{Robba}^{\mathrm{extended},q}_{I,R}[G].	
\end{align}
\item<4-> And then consider the completion living inside the corresponding infinite direct sum Banach modules 
\begin{align}
\bigoplus\mathrm{Robba}^{\mathrm{extended},q}_{I,R},	
\end{align}
over the corresponding period rings:
\begin{align}
\overline{\mathrm{Robba}^{\mathrm{extended},q}_{I,R}[G]}.	
\end{align}
Then we take suitable intersection and union one can have possibly some interesting period rings $\overline{\mathrm{Robba}^{\mathrm{extended},q}_{r,R}[G]}$ and $\overline{\mathrm{Robba}^{\mathrm{extended},q}_{R}[G]}$.
\end{itemize}

\noindent{\text{Equivariant relative $p$-adic Hodge Theory}}
\begin{itemize}
\item<1-> The equivariant period rings in the situations we mentioned above carry relative multi-Frobenius action $\varphi_q$ induced from the Witt vectors. 

\item<2-> They carry the corresponding Banach or (ind-)Fr\'echet spaces structures. So we can generalize the corresponding Kedlaya-Liu's construction to the following situations (here let $G$ be finite):

\item<3-> We can then consider the corresponding completed Frobenius modules over the rings in the equivariant setting. To be more precise over:
\begin{align}
\overline{\mathrm{Robba}^{\mathrm{extended},q}_{R,A}[G]},\Omega_{\mathrm{int},R,A},\Omega_{R,A},\mathrm{Robba}^{\mathrm{extended},q}_{R,A},\mathrm{Robba}^{\mathrm{extended},q}_{\mathrm{bd},R,A}	
\end{align}
one considers the Frobenius modules finite locally free.

\item<4->  With the corresponding finite locally free models over
\begin{align}
\overline{\mathrm{Robba}^{\mathrm{extended},q}_{r,R,A}[G]},\mathrm{Robba}^{\mathrm{extended},q}_{r,R,A},\mathrm{Robba}^{\mathrm{extended},q}_{{\mathrm{bd},r},R,A},	
\end{align} 
again carrying the corresponding semilinear Frobenius structures, where $r$ could be $\infty$.

\item<5-> One also consider families of Frobenius modules over 
\begin{align}
\overline{\mathrm{Robba}^{\mathrm{extended},q}_{I,R,A}[G]},\mathrm{Robba}^{\mathrm{extended},q}_{I,R,A},	
\end{align} 
in glueing fashion with obvious cocycle condition with respect to three multi-intervals $I\subset J\subset K$. These are called the corresponding Frobenius bundles.
	
\end{itemize}

\newpage

\section{Analytic $\infty$-Categorical Functional Analytic Hodge-Iwasawa Modules}

\subsection{$\infty$-Categorical Analytic Stacks and Descents I}

\noindent We now make the corresponding discussion after our previous work \cite{T2} on the homotopical functional analysis after many projects \cite{1BBBK}, \cite{1BBK}, \cite{BBM}, \cite{1BK} , \cite{1CS1}, \cite{1CS2}, \cite{KKM}. We choose to work over the Bambozzi-Kremnizer space \cite{1BK} attached to the corresponding Banach rings in our work after \cite{1BBBK}, \cite{1BBK}, \cite{BBM}, \cite{1BK}, \cite{KKM}. Note that what is happening is that attached to any Banach ring over $\mathbb{Q}_p$, say $B$, we attach a $(\infty,1)-$stack $\mathcal{X}(B)$ fibered over (in the sense of $\infty$-groupoid, and up to taking the corresponding opposite categories) after \cite{1BBBK}, \cite{1BBK}, \cite{BBM}, \cite{1BK}, \cite{KKM}:
\begin{align}
\mathrm{sComm}\mathrm{Simp}\mathrm{Ind}\mathrm{Ban}_{\mathbb{Q}_p},	
\end{align}
with 
\begin{align}
\mathrm{sComm}\mathrm{Simp}\mathrm{Ind}^m\mathrm{Ban}_{\mathbb{Q}_p}.	
\end{align}
associated with a $(\infty,1)$-ring object $\mathcal{O}_{\mathcal{X}(B)}$, such that we have the corresponding under the basic derived rational localization $\infty$-Grothendieck site
\begin{center}
 $(\mathcal{X}(B), \mathcal{O}_{\mathcal{X}(B),\mathrm{drl}})$ 
\end{center}
carrying the homotopical epimorphisms as the corresponding topology.

\begin{itemize}
\item By using this framework (certainly one can also consider \cite{1CS1} and \cite{1CS2} as the foundations, as in \cite{LBV}), we have the $\infty$-stack after Kedlaya-Liu \cite{KL1}, \cite{KL2}. Here in the following let $A$ be any Banach ring over $\mathbb{Q}_p$.

\item Generalizing Kedlaya-Liu's construction in \cite{KL1}, \cite{KL2} of the adic Fargues-Fontaine space we have a quotient (by using power $\varphi_q$ of the Frobenius operator) $X_{R,A,q}$ of the space 
\begin{align}
Y_{R,A,q}:=\bigcup_{I,\text{multi}}\mathcal{X}(\mathrm{Robba}^{\mathrm{extended},q}_{R,I,A}).	
\end{align}

\item This is a locally ringed space $(X_{R,A,q},\mathcal{O}_{X_{R,A,q}})$, so one can consider the stable $\infty$-category $\mathrm{Ind}\mathrm{Banach}(\mathcal{O}_{X_{R,A,q}}) $ which is the $\infty$-category of all the $\mathcal{O}_{X_{R,A,q}}$-sheaves of inductive Banach modules over $X_{R,A,q}$. We have the parallel categories for $Y_{R,A,q}$, namely $\varphi\mathrm{Ind}\mathrm{Banach}(\mathcal{O}_{X_{R,A,q}})$ and so on. Here we only consider presheaves.

\item This is a locally ringed space $(X_{R,A,q},\mathcal{O}_{X_{R,A,q}})$, so one can consider the stable $\infty$-category $\mathrm{Ind}^m\mathrm{Banach}(\mathcal{O}_{X_{R,A,q}}) $ which is the $\infty$-category of all the $\mathcal{O}_{X_{R,A,q}}$-sheaves of inductive monomorphic Banach modules over $X_{R,A,q}$. We have the parallel categories for $Y_{R,A,q}$, namely $\varphi\mathrm{Ind}^m\mathrm{Banach}(\mathcal{O}_{X_{R,A,q}})$ and so on. Here we only consider presheaves.

\begin{assumption}\label{assumtionpresheaves}
All the functors of modules or algebras below are presheaves.	
\end{assumption}

\item In this context one can consider the $K$-theory as in the scheme situation by using the ideas and constructions from Blumberg-Gepner-Tabuada \cite{BGT}. Moreover we can study the Hodge Theory.

\item We expect that one can study among these big categories to find interesting relationships, since this should give us the right understanding of the $p$-adic Hodge theory. The corresponding pseudocoherent version comparison could be expected to be deduced as in Kedlaya-Liu's work \cite{KL1}, \cite{KL2}.

\item (\text{Proposition}) There is an equivalence between the $\infty$-categories of inductive Banach quasicoherent presheaves:
\[
\xymatrix@R+0pc@C+0pc{
\mathrm{Ind}\mathrm{Banach}(\mathcal{O}_{X_{R,A,q}})\ar[r]^{\mathrm{equi}}\ar[r]\ar[r] &\varphi_q\mathrm{Ind}\mathrm{Banach}(\mathcal{O}_{Y_{R,A,q}}).  
}
\]
\item (\text{Proposition}) There is an equivalence between the $\infty$-categories of monomorphic inductive Banach quasicoherent presheaves:
\[
\xymatrix@R+0pc@C+0pc{
\mathrm{Ind}^m\mathrm{Banach}(\mathcal{O}_{X_{R,A,q}})\ar[r]^{\mathrm{equi}}\ar[r]\ar[r] &\varphi_q\mathrm{Ind}^m\mathrm{Banach}(\mathcal{O}_{Y_{R,A,q}}).  
}
\]
\item (\text{Proposition}) There is an equivalence between the $\infty$-categories of inductive Banach quasicoherent commutative algebra $E_\infty$ objects:
\[
\xymatrix@R+0pc@C+0pc{
\mathrm{sComm}_\mathrm{simplicial}\mathrm{Ind}\mathrm{Banach}(\mathcal{O}_{X_{R,A,q}})\ar[r]^{\mathrm{equi}}\ar[r]\ar[r] &\mathrm{sComm}_\mathrm{simplicial}\varphi_q\mathrm{Ind}\mathrm{Banach}(\mathcal{O}_{Y_{R,A,q}}).  
}
\]
\item (\text{Proposition}) There is an equivalence between the $\infty$-categories of monomorphic inductive Banach quasicoherent commutative algebra $E_\infty$ objects:
\[
\xymatrix@R+0pc@C+0pc{
\mathrm{sComm}_\mathrm{simplicial}\mathrm{Ind}^m\mathrm{Banach}(\mathcal{O}_{X_{R,A,q}})\ar[r]^{\mathrm{equi}}\ar[r]\ar[r] &\mathrm{sComm}_\mathrm{simplicial}\varphi_q\mathrm{Ind}^m\mathrm{Banach}(\mathcal{O}_{Y_{R,A,q}}).  
}
\]

\item Then parallel as in \cite{LBV} we have the equivalence of the de Rham complex after \cite[Definition 5.9, Section 5.2.1]{KKM}:
\[
\xymatrix@R+0pc@C+0pc{
\mathrm{deRham}_{\mathrm{sComm}_\mathrm{simplicial}\mathrm{Ind}\mathrm{Banach}(\mathcal{O}_{X_{R,A,q}})\ar[r]^{\mathrm{equi}}}(-)\ar[r]\ar[r] &\mathrm{deRham}_{\mathrm{sComm}_\mathrm{simplicial}\varphi_q\mathrm{Ind}\mathrm{Banach}(\mathcal{O}_{Y_{R,A,q}})}(-), 
}
\]
\[
\xymatrix@R+0pc@C+0pc{
\mathrm{deRham}_{\mathrm{sComm}_\mathrm{simplicial}\mathrm{Ind}^m\mathrm{Banach}(\mathcal{O}_{X_{R,A,q}})\ar[r]^{\mathrm{equi}}}(-)\ar[r]\ar[r] &\mathrm{deRham}_{\mathrm{sComm}_\mathrm{simplicial}\varphi_q\mathrm{Ind}^m\mathrm{Banach}(\mathcal{O}_{Y_{R,A,q}})}(-). 
}
\]

\item Then we have the following equivalence of $K$-group $(\infty,1)$-spectrum from \cite{BGT}:
\[
\xymatrix@R+0pc@C+0pc{
\mathrm{K}^\mathrm{BGT}_{\mathrm{sComm}_\mathrm{simplicial}\mathrm{Ind}\mathrm{Banach}(\mathcal{O}_{X_{R,A,q}})\ar[r]^{\mathrm{equi}}}(-)\ar[r]\ar[r] &\mathrm{K}^\mathrm{BGT}_{\mathrm{sComm}_\mathrm{simplicial}\varphi_q\mathrm{Ind}\mathrm{Banach}(\mathcal{O}_{Y_{R,A,q}})}(-), 
}
\]
\[
\xymatrix@R+0pc@C+0pc{
\mathrm{K}^\mathrm{BGT}_{\mathrm{sComm}_\mathrm{simplicial}\mathrm{Ind}^m\mathrm{Banach}(\mathcal{O}_{X_{R,A,q}})\ar[r]^{\mathrm{equi}}}(-)\ar[r]\ar[r] &\mathrm{K}^\mathrm{BGT}_{\mathrm{sComm}_\mathrm{simplicial}\varphi_q\mathrm{Ind}^m\mathrm{Banach}(\mathcal{O}_{Y_{R,A,q}})}(-). 
}
\]
\end{itemize}

\noindent Now let $R=\mathbb{Q}_p(p^{1/p^\infty})^{\wedge\flat}$ and $R_k=\mathbb{Q}_p(p^{1/p^\infty})^{\wedge}\left<T_1^{\pm 1/p^{\infty}},...,T_k^{\pm 1/p^{\infty}}\right>^\flat$ we have the following Galois theoretic results with naturality along $f:\mathrm{Spa}(\mathbb{Q}_p(p^{1/p^\infty})^{\wedge}\left<T_1^{\pm 1/p^\infty},...,T_k^{\pm 1/p^\infty}\right>^\flat)\rightarrow \mathrm{Spa}(\mathbb{Q}_p(p^{1/p^\infty})^{\wedge\flat})$:

\begin{itemize}
\item (\text{Proposition}) There is an equivalence between the $\infty$-categories of inductive Banach quasicoherent presheaves:
\[
\xymatrix@R+6pc@C+0pc{
\mathrm{Ind}\mathrm{Banach}(\mathcal{O}_{X_{\mathbb{Q}_p(p^{1/p^\infty})^{\wedge}\left<T_1^{\pm 1/p^\infty},...,T_k^{\pm 1/p^\infty}\right>^\flat,A,q}})\ar[d]\ar[d]\ar[d]\ar[d] \ar[r]^{\mathrm{equi}}\ar[r]\ar[r] &\varphi_q\mathrm{Ind}\mathrm{Banach}(\mathcal{O}_{Y_{\mathbb{Q}_p(p^{1/p^\infty})^{\wedge}\left<T_1^{\pm 1/p^\infty},...,T_k^{\pm 1/p^\infty}\right>^\flat,A,q}}) \ar[d]\ar[d]\ar[d]\ar[d].\\
\mathrm{Ind}\mathrm{Banach}(\mathcal{O}_{X_{\mathbb{Q}_p(p^{1/p^\infty})^{\wedge\flat},A,q}})\ar[r]^{\mathrm{equi}}\ar[r]\ar[r] &\varphi_q\mathrm{Ind}\mathrm{Banach}(\mathcal{O}_{Y_{\mathbb{Q}_p(p^{1/p^\infty})^{\wedge\flat},A,q}}).\\ 
}
\]
\item (\text{Proposition}) There is an equivalence between the $\infty$-categories of monomorphic inductive Banach quasicoherent presheaves:
\[
\xymatrix@R+6pc@C+0pc{
\mathrm{Ind}^m\mathrm{Banach}(\mathcal{O}_{X_{R_k,A,q}})\ar[r]^{\mathrm{equi}}\ar[d]\ar[d]\ar[d]\ar[d]\ar[r]\ar[r] &\varphi_q\mathrm{Ind}^m\mathrm{Banach}(\mathcal{O}_{Y_{R_k,A,q}})\ar[d]\ar[d]\ar[d]\ar[d]\\
\mathrm{Ind}^m\mathrm{Banach}(\mathcal{O}_{X_{\mathbb{Q}_p(p^{1/p^\infty})^{\wedge\flat},A,q}})\ar[r]^{\mathrm{equi}}\ar[r]\ar[r] &\varphi_q\mathrm{Ind}^m\mathrm{Banach}(\mathcal{O}_{Y_{\mathbb{Q}_p(p^{1/p^\infty})^{\wedge\flat},A,q}}).\\  
}
\]
\item (\text{Proposition}) There is an equivalence between the $\infty$-categories of inductive Banach quasicoherent commutative algebra $E_\infty$ objects:
\[
\xymatrix@R+6pc@C+0pc{
\mathrm{sComm}_\mathrm{simplicial}\mathrm{Ind}\mathrm{Banach}(\mathcal{O}_{X_{R_k,A,q}})\ar[d]\ar[d]\ar[d]\ar[d]\ar[r]^{\mathrm{equi}}\ar[r]\ar[r] &\mathrm{sComm}_\mathrm{simplicial}\varphi_q\mathrm{Ind}\mathrm{Banach}(\mathcal{O}_{Y_{R_k,A,q}})\ar[d]\ar[d]\ar[d]\ar[d]\\
\mathrm{sComm}_\mathrm{simplicial}\mathrm{Ind}\mathrm{Banach}(\mathcal{O}_{X_{\mathbb{Q}_p(p^{1/p^\infty})^{\wedge\flat},A,q}})\ar[r]^{\mathrm{equi}}\ar[r]\ar[r] &\mathrm{sComm}_\mathrm{simplicial}\varphi_q\mathrm{Ind}\mathrm{Banach}(\mathcal{O}_{Y_{\mathbb{Q}_p(p^{1/p^\infty})^{\wedge\flat},A,q}})  
}
\]
\item (\text{Proposition}) There is an equivalence between the $\infty$-categories of monomorphic inductive Banach quasicoherent commutative algebra $E_\infty$ objects:
\[
\xymatrix@R+6pc@C+0pc{
\mathrm{sComm}_\mathrm{simplicial}\mathrm{Ind}^m\mathrm{Banach}(\mathcal{O}_{X_{R_k,A,q}})\ar[d]\ar[d]\ar[d]\ar[d]\ar[r]^{\mathrm{equi}}\ar[r]\ar[r] &\mathrm{sComm}_\mathrm{simplicial}\varphi_q\mathrm{Ind}^m\mathrm{Banach}(\mathcal{O}_{Y_{R_k,A,q}})\ar[d]\ar[d]\ar[d]\ar[d]\\
 \mathrm{sComm}_\mathrm{simplicial}\mathrm{Ind}^m\mathrm{Banach}(\mathcal{O}_{X_{\mathbb{Q}_p(p^{1/p^\infty})^{\wedge\flat},A,q}})\ar[r]^{\mathrm{equi}}\ar[r]\ar[r] &\mathrm{sComm}_\mathrm{simplicial}\varphi_q\mathrm{Ind}^m\mathrm{Banach}(\mathcal{O}_{Y_{\mathbb{Q}_p(p^{1/p^\infty})^{\wedge\flat},A,q}}) 
}
\]

\item Then parallel as in \cite{LBV} we have the equivalence of the de Rham complex after \cite[Definition 5.9, Section 5.2.1]{KKM}:
\[\displayindent=-0.4in
\xymatrix@R+6pc@C+0pc{
\mathrm{deRham}_{\mathrm{sComm}_\mathrm{simplicial}\mathrm{Ind}\mathrm{Banach}(\mathcal{O}_{X_{R_k,A,q}})\ar[r]^{\mathrm{equi}}}(-)\ar[d]\ar[d]\ar[d]\ar[d]\ar[r]\ar[r] &\mathrm{deRham}_{\mathrm{sComm}_\mathrm{simplicial}\varphi_q\mathrm{Ind}\mathrm{Banach}(\mathcal{O}_{Y_{R_k,A,q}})}(-)\ar[d]\ar[d]\ar[d]\ar[d]\\
\mathrm{deRham}_{\mathrm{sComm}_\mathrm{simplicial}\mathrm{Ind}\mathrm{Banach}(\mathcal{O}_{X_{\mathbb{Q}_p(p^{1/p^\infty})^{\wedge\flat},A,q}})\ar[r]^{\mathrm{equi}}}(-)\ar[r]\ar[r] &\mathrm{deRham}_{\mathrm{sComm}_\mathrm{simplicial}\varphi_q\mathrm{Ind}\mathrm{Banach}(\mathcal{O}_{Y_{\mathbb{Q}_p(p^{1/p^\infty})^{\wedge\flat},A,q}})}(-) 
}
\]
\[\displayindent=-0.4in
\xymatrix@R+6pc@C+0pc{
\mathrm{deRham}_{\mathrm{sComm}_\mathrm{simplicial}\mathrm{Ind}^m\mathrm{Banach}(\mathcal{O}_{X_{R_k,A,q}})\ar[r]^{\mathrm{equi}}}(-)\ar[d]\ar[d]\ar[d]\ar[d]\ar[r]\ar[r] &\mathrm{deRham}_{\mathrm{sComm}_\mathrm{simplicial}\varphi_q\mathrm{Ind}^m\mathrm{Banach}(\mathcal{O}_{Y_{R_k,A,q}})}(-)\ar[d]\ar[d]\ar[d]\ar[d]\\
\mathrm{deRham}_{\mathrm{sComm}_\mathrm{simplicial}\mathrm{Ind}^m\mathrm{Banach}(\mathcal{O}_{X_{\mathbb{Q}_p(p^{1/p^\infty})^{\wedge\flat},A,q}})\ar[r]^{\mathrm{equi}}}(-)\ar[r]\ar[r] &\mathrm{deRham}_{\mathrm{sComm}_\mathrm{simplicial}\varphi_q\mathrm{Ind}^m\mathrm{Banach}(\mathcal{O}_{Y_{\mathbb{Q}_p(p^{1/p^\infty})^{\wedge\flat},A,q}})}(-) 
}
\]

\item Then we have the following equivalence of $K$-group $(\infty,1)$-spectrum from \cite{BGT}:
\[
\xymatrix@R+6pc@C+0pc{
\mathrm{K}^\mathrm{BGT}_{\mathrm{sComm}_\mathrm{simplicial}\mathrm{Ind}\mathrm{Banach}(\mathcal{O}_{X_{R_k,A,q}})\ar[r]^{\mathrm{equi}}}(-)\ar[d]\ar[d]\ar[d]\ar[d]\ar[r]\ar[r] &\mathrm{K}^\mathrm{BGT}_{\mathrm{sComm}_\mathrm{simplicial}\varphi_q\mathrm{Ind}\mathrm{Banach}(\mathcal{O}_{Y_{R_k,A,q}})}(-)\ar[d]\ar[d]\ar[d]\ar[d]\\
\mathrm{K}^\mathrm{BGT}_{\mathrm{sComm}_\mathrm{simplicial}\mathrm{Ind}\mathrm{Banach}(\mathcal{O}_{X_{\mathbb{Q}_p(p^{1/p^\infty})^{\wedge\flat},A,q}})\ar[r]^{\mathrm{equi}}}(-)\ar[r]\ar[r] &\mathrm{K}^\mathrm{BGT}_{\mathrm{sComm}_\mathrm{simplicial}\varphi_q\mathrm{Ind}\mathrm{Banach}(\mathcal{O}_{Y_{\mathbb{Q}_p(p^{1/p^\infty})^{\wedge\flat},A,q}})}(-) 
}
\]
\[
\xymatrix@R+6pc@C+0pc{
\mathrm{K}^\mathrm{BGT}_{\mathrm{sComm}_\mathrm{simplicial}\mathrm{Ind}^m\mathrm{Banach}(\mathcal{O}_{X_{R_k,A,q}})\ar[r]^{\mathrm{equi}}}(-)\ar[d]\ar[d]\ar[d]\ar[d]\ar[r]\ar[r] &\mathrm{K}^\mathrm{BGT}_{\mathrm{sComm}_\mathrm{simplicial}\varphi_q\mathrm{Ind}^m\mathrm{Banach}(\mathcal{O}_{Y_{R_k,A,q}})}(-)\ar[d]\ar[d]\ar[d]\ar[d]\\
\mathrm{K}^\mathrm{BGT}_{\mathrm{sComm}_\mathrm{simplicial}\mathrm{Ind}^m\mathrm{Banach}(\mathcal{O}_{X_{\mathbb{Q}_p(p^{1/p^\infty})^{\wedge\flat},A,q}})\ar[r]^{\mathrm{equi}}}(-)\ar[r]\ar[r] &\mathrm{K}^\mathrm{BGT}_{\mathrm{sComm}_\mathrm{simplicial}\varphi_q\mathrm{Ind}^m\mathrm{Banach}(\mathcal{O}_{Y_{\mathbb{Q}_p(p^{1/p^\infty})^{\wedge\flat},A,q}})}(-) 
}
\]

\end{itemize}

\indent Then we consider further equivariance by considering the arithmetic profinite fundamental group and actually its $q$-th power $\mathrm{Gal}(\overline{{Q}_p\left<T_1^{\pm 1},...,T_k^{\pm 1}\right>}/R_k)^{\times q}$ through the following diagram:\\

\[
\xymatrix@R+6pc@C+0pc{
\mathbb{Z}_p^k=\mathrm{Gal}(R_k/{\mathbb{Q}_p(p^{1/p^\infty})^\wedge\left<T_1^{\pm 1},...,T_k^{\pm 1}\right>})\ar[d]\ar[d]\ar[d]\ar[d] \ar[r]\ar[r] \ar[r]\ar[r] &\mathrm{Gal}(\overline{{Q}_p\left<T_1^{\pm 1},...,T_k^{\pm 1}\right>}/R_k) \ar[d]\ar[d]\ar[d] \ar[r]\ar[r] &\Gamma_{\mathbb{Q}_p} \ar[d]\ar[d]\ar[d]\ar[d]\\
(\mathbb{Z}_p^k=\mathrm{Gal}(R_k/{\mathbb{Q}_p(p^{1/p^\infty})^\wedge\left<T_1^{\pm 1},...,T_k^{\pm 1}\right>}))^{\times q} \ar[r]\ar[r] \ar[r]\ar[r] &\Gamma_k^{\times q}:=\mathrm{Gal}(R_k/{\mathbb{Q}_p\left<T_1^{\pm 1},...,T_k^{\pm 1}\right>})^{\times q} \ar[r] \ar[r]\ar[r] &\Gamma_{\mathbb{Q}_p}^{\times q}.
}
\]

\

We then have the correspond arithmetic profinite fundamental groups equivariant versions:

\begin{itemize}
\item (\text{Proposition}) There is an equivalence between the $\infty$-categories of inductive Banach quasicoherent presheaves:
\[
\xymatrix@R+6pc@C+0pc{
\mathrm{Ind}\mathrm{Banach}_{\Gamma_{k}^{\times q}}(\mathcal{O}_{X_{\mathbb{Q}_p(p^{1/p^\infty})^{\wedge}\left<T_1^{\pm 1/p^\infty},...,T_k^{\pm 1/p^\infty}\right>^\flat,A,q}})\ar[d]\ar[d]\ar[d]\ar[d] \ar[r]^{\mathrm{equi}}\ar[r]\ar[r] &\varphi_q\mathrm{Ind}\mathrm{Banach}_{\Gamma_{k}^{\times q}}(\mathcal{O}_{Y_{\mathbb{Q}_p(p^{1/p^\infty})^{\wedge}\left<T_1^{\pm 1/p^\infty},...,T_k^{\pm 1/p^\infty}\right>^\flat,A,q}}) \ar[d]\ar[d]\ar[d]\ar[d].\\
\mathrm{Ind}\mathrm{Banach}_{\Gamma_{0}^{\times q}}(\mathcal{O}_{X_{\mathbb{Q}_p(p^{1/p^\infty})^{\wedge\flat},A,q}})\ar[r]^{\mathrm{equi}}\ar[r]\ar[r] &\varphi_q\mathrm{Ind}\mathrm{Banach}_{\Gamma_{0}^{\times q}}(\mathcal{O}_{Y_{\mathbb{Q}_p(p^{1/p^\infty})^{\wedge\flat},A,q}}).\\ 
}
\]
\item (\text{Proposition}) There is an equivalence between the $\infty$-categories of monomorphic inductive Banach quasicoherent presheaves:
\[
\xymatrix@R+6pc@C+0pc{
\mathrm{Ind}^m\mathrm{Banach}_{\Gamma_{k}^{\times q}}(\mathcal{O}_{X_{R_k,A,q}})\ar[r]^{\mathrm{equi}}\ar[d]\ar[d]\ar[d]\ar[d]\ar[r]\ar[r] &\varphi_q\mathrm{Ind}^m\mathrm{Banach}_{\Gamma_{k}^{\times q}}(\mathcal{O}_{Y_{R_k,A,q}})\ar[d]\ar[d]\ar[d]\ar[d]\\
\mathrm{Ind}^m\mathrm{Banach}_{\Gamma_{0}^{\times q}}(\mathcal{O}_{X_{\mathbb{Q}_p(p^{1/p^\infty})^{\wedge\flat},A,q}})\ar[r]^{\mathrm{equi}}\ar[r]\ar[r] &\varphi_q\mathrm{Ind}^m\mathrm{Banach}_{\Gamma_{0}^{\times q}}(\mathcal{O}_{Y_{\mathbb{Q}_p(p^{1/p^\infty})^{\wedge\flat},A,q}}).\\  
}
\]	
\end{itemize}

\

\

\indent Now we consider \cite{1CS1} and \cite[Proposition 13.8, Theorem 14.9, Remark 14.10]{1CS2}\footnote{Note that we are motivated as well from \cite{LBV}.}, and study the corresponding solid perfect complexes, solid quasicoherent sheaves and solid vector bundles. Here we are going to use different formalism, therefore we will have different categories and functors. We use the notation $\circledcirc$ to denote any element of $\{\text{solid perfect complexes}, \text{solid quasicoherent sheaves}, \text{solid vector bundles}\}$ from \cite{1CS2} with the corresponding descent results of \cite[Proposition 13.8, Theorem 14.9, Remark 14.10]{1CS2}. Then we have the following:

\begin{itemize}
\item Generalizing Kedlaya-Liu's construction in \cite{KL1}, \cite{KL2} of the adic Fargues-Fontaine space we have a quotient (by using powers of the Frobenius operator) $X_{R,A,q}$ of the space by using \cite{1CS2}:
\begin{align}
Y_{R,A,q}:=\bigcup_{I,\mathrm{multi}}\mathcal{X}^\mathrm{CS}(\text{Robba}^\text{extended}_{R,I,A,q}).	
\end{align}

\item This is a locally ringed space $(X_{R,A,q},\mathcal{O}_{X_{R,A,q}})$, so one can consider the stable $\infty$-category $\mathrm{Module}_{\text{CS},\mathrm{quasicoherent}}(\mathcal{O}_{X_{R,A,q}}) $ which is the $\infty$-category of all the $\mathcal{O}_{X_{R,A,q}}$-sheaves of solid modules over $X_{R,A,q}$. We have the parallel categories for $Y_{R,A,q}$, namely $\varphi\mathrm{Module}_{\text{CS},\mathrm{quasicoherent}}(\mathcal{O}_{X_{R,A,q}})$ and so on. Here we will consider sheaves.

\begin{assumption}\label{assumtionpresheaves}
All the functors of modules or algebras below are Clausen-Scholze sheaves \cite[Proposition 13.8, Theorem 14.9, Remark 14.10]{1CS2}.	
\end{assumption}

\item (\text{Proposition}) There is an equivalence between the $\infty$-categories of inductive solid sheaves:
\[
\xymatrix@R+0pc@C+0pc{
\mathrm{Modules}_\circledcirc(\mathcal{O}_{X_{R,A,q}})\ar[r]^{\mathrm{equi}}\ar[r]\ar[r] &\varphi_q\mathrm{Modules}_\circledcirc(\mathcal{O}_{Y_{R,A,q}}).  
}
\]
\end{itemize}

\begin{itemize}

\item (\text{Proposition}) There is an equivalence between the $\infty$-categories of inductive Banach quasicoherent commutative algebra $E_\infty$ objects:
\[\displayindent=-1.0in
\xymatrix@R+0pc@C+0pc{
\mathrm{sComm}_\mathrm{simplicial}\mathrm{Modules}_{\text{solidquasicoherentsheaves}}(\mathcal{O}_{X_{R,A,q}})\ar[r]^{\mathrm{equi}}\ar[r]\ar[r] &\mathrm{sComm}_\mathrm{simplicial}\varphi_q\mathrm{Modules}_{\text{solidquasicoherentsheaves}}(\mathcal{O}_{Y_{R,A,q}}).  
}
\]

\item Then as in \cite{LBV} we have the equivalence of the de Rham complex after \cite[Definition 5.9, Section 5.2.1]{KKM}\footnote{Here $\circledcirc=\text{solidquasicoherentsheaves}$.}:
\[
\xymatrix@R+0pc@C+0pc{
\mathrm{deRham}_{\mathrm{sComm}_\mathrm{simplicial}\mathrm{Modules}_\circledcirc(\mathcal{O}_{X_{R,A,q}})\ar[r]^{\mathrm{equi}}}(-)\ar[r]\ar[r] &\mathrm{deRham}_{\mathrm{sComm}_\mathrm{simplicial}\varphi_q\mathrm{Modules}_\circledcirc(\mathcal{O}_{Y_{R,A,q}})}(-). 
}
\]

\item Then we have the following equivalence of $K$-group $(\infty,1)$-spectrum from \cite{BGT}\footnote{Here $\circledcirc=\text{solidquasicoherentsheaves}$.}:
\[
\xymatrix@R+0pc@C+0pc{
\mathrm{K}^\mathrm{BGT}_{\mathrm{sComm}_\mathrm{simplicial}\mathrm{Modules}_\circledcirc(\mathcal{O}_{X_{R,A,q}})\ar[r]^{\mathrm{equi}}}(-)\ar[r]\ar[r] &\mathrm{K}^\mathrm{BGT}_{\mathrm{sComm}_\mathrm{simplicial}\varphi_q\mathrm{Modules}_\circledcirc(\mathcal{O}_{Y_{R,A,q}})}(-). 
}
\]
\end{itemize}

\noindent Now let $R=\mathbb{Q}_p(p^{1/p^\infty})^{\wedge\flat}$ and $R_k=\mathbb{Q}_p(p^{1/p^\infty})^{\wedge}\left<T_1^{\pm 1/p^{\infty}},...,T_k^{\pm 1/p^{\infty}}\right>^\flat$ we have the following Galois theoretic results with naturality along $f:\mathrm{Spa}(\mathbb{Q}_p(p^{1/p^\infty})^{\wedge}\left<T_1^{\pm 1/p^\infty},...,T_k^{\pm 1/p^\infty}\right>^\flat)\rightarrow \mathrm{Spa}(\mathbb{Q}_p(p^{1/p^\infty})^{\wedge\flat})$:

\begin{itemize}
\item (\text{Proposition}) There is an equivalence between the $\infty$-categories of inductive Banach quasicoherent presheaves\footnote{Here $\circledcirc=\text{solidquasicoherentsheaves}$.}:
\[
\xymatrix@R+6pc@C+0pc{
\mathrm{Modules}_\circledcirc(\mathcal{O}_{X_{\mathbb{Q}_p(p^{1/p^\infty})^{\wedge}\left<T_1^{\pm 1/p^\infty},...,T_k^{\pm 1/p^\infty}\right>^\flat,A,q}})\ar[d]\ar[d]\ar[d]\ar[d] \ar[r]^{\mathrm{equi}}\ar[r]\ar[r] &\varphi_q\mathrm{Modules}_\circledcirc(\mathcal{O}_{Y_{\mathbb{Q}_p(p^{1/p^\infty})^{\wedge}\left<T_1^{\pm 1/p^\infty},...,T_k^{\pm 1/p^\infty}\right>^\flat,A,q}}) \ar[d]\ar[d]\ar[d]\ar[d].\\
\mathrm{Modules}_\circledcirc(\mathcal{O}_{X_{\mathbb{Q}_p(p^{1/p^\infty})^{\wedge\flat},A,q}})\ar[r]^{\mathrm{equi}}\ar[r]\ar[r] &\varphi_q\mathrm{Modules}_\circledcirc(\mathcal{O}_{Y_{\mathbb{Q}_p(p^{1/p^\infty})^{\wedge\flat},A,q}}).\\ 
}
\]
\item (\text{Proposition}) There is an equivalence between the $\infty$-categories of inductive Banach quasicoherent commutative algebra $E_\infty$ objects\footnote{Here $\circledcirc=\text{solidquasicoherentsheaves}$.}:
\[
\xymatrix@R+6pc@C+0pc{
\mathrm{sComm}_\mathrm{simplicial}\mathrm{Modules}_\circledcirc(\mathcal{O}_{X_{R_k,A,q}})\ar[d]\ar[d]\ar[d]\ar[d]\ar[r]^{\mathrm{equi}}\ar[r]\ar[r] &\mathrm{sComm}_\mathrm{simplicial}\varphi_q\mathrm{Modules}_\circledcirc(\mathcal{O}_{Y_{R_k,A,q}})\ar[d]\ar[d]\ar[d]\ar[d]\\
\mathrm{sComm}_\mathrm{simplicial}\mathrm{Modules}_\circledcirc(\mathcal{O}_{X_{\mathbb{Q}_p(p^{1/p^\infty})^{\wedge\flat},A,q}})\ar[r]^{\mathrm{equi}}\ar[r]\ar[r] &\mathrm{sComm}_\mathrm{simplicial}\varphi_q\mathrm{Modules}_\circledcirc(\mathcal{O}_{Y_{\mathbb{Q}_p(p^{1/p^\infty})^{\wedge\flat},A,q}}).  
}
\]
\item Then as in \cite{LBV} we have the equivalence of the de Rham complex after \cite[Definition 5.9, Section 5.2.1]{KKM}\footnote{Here $\circledcirc=\text{solidquasicoherentsheaves}$.}:
\[\displayindent=-0.4in
\xymatrix@R+6pc@C+0pc{
\mathrm{deRham}_{\mathrm{sComm}_\mathrm{simplicial}\mathrm{Modules}_\circledcirc(\mathcal{O}_{X_{R_k,A,q}})\ar[r]^{\mathrm{equi}}}(-)\ar[d]\ar[d]\ar[d]\ar[d]\ar[r]\ar[r] &\mathrm{deRham}_{\mathrm{sComm}_\mathrm{simplicial}\varphi_q\mathrm{Modules}_\circledcirc(\mathcal{O}_{Y_{R_k,A,q}})}(-)\ar[d]\ar[d]\ar[d]\ar[d]\\
\mathrm{deRham}_{\mathrm{sComm}_\mathrm{simplicial}\mathrm{Modules}_\circledcirc(\mathcal{O}_{X_{\mathbb{Q}_p(p^{1/p^\infty})^{\wedge\flat},A,q}})\ar[r]^{\mathrm{equi}}}(-)\ar[r]\ar[r] &\mathrm{deRham}_{\mathrm{sComm}_\mathrm{simplicial}\varphi_q\mathrm{Modules}_\circledcirc(\mathcal{O}_{Y_{\mathbb{Q}_p(p^{1/p^\infty})^{\wedge\flat},A,q}})}(-). 
}
\]

\item Then we have the following equivalence of $K$-group $(\infty,1)$-spectrum from \cite{BGT}\footnote{Here $\circledcirc=\text{solidquasicoherentsheaves}$.}:
\[
\xymatrix@R+6pc@C+0pc{
\mathrm{K}^\mathrm{BGT}_{\mathrm{sComm}_\mathrm{simplicial}\mathrm{Modules}_\circledcirc(\mathcal{O}_{X_{R_k,A,q}})\ar[r]^{\mathrm{equi}}}(-)\ar[d]\ar[d]\ar[d]\ar[d]\ar[r]\ar[r] &\mathrm{K}^\mathrm{BGT}_{\mathrm{sComm}_\mathrm{simplicial}\varphi_q\mathrm{Modules}_\circledcirc(\mathcal{O}_{Y_{R_k,A,q}})}(-)\ar[d]\ar[d]\ar[d]\ar[d]\\
\mathrm{K}^\mathrm{BGT}_{\mathrm{sComm}_\mathrm{simplicial}\mathrm{Modules}_\circledcirc(\mathcal{O}_{X_{\mathbb{Q}_p(p^{1/p^\infty})^{\wedge\flat},A,q}})\ar[r]^{\mathrm{equi}}}(-)\ar[r]\ar[r] &\mathrm{K}^\mathrm{BGT}_{\mathrm{sComm}_\mathrm{simplicial}\varphi_q\mathrm{Modules}_\circledcirc(\mathcal{O}_{Y_{\mathbb{Q}_p(p^{1/p^\infty})^{\wedge\flat},A,q}})}(-).}
\]

\end{itemize}

\
\indent Then we consider further equivariance by considering the arithmetic profinite fundamental groups $\Gamma_{\mathbb{Q}_p}$ and $\mathrm{Gal}(\overline{\mathbb{Q}_p\left<T_1^{\pm 1},...,T_k^{\pm 1}\right>}/R_k)$ through the following diagram:

\[
\xymatrix@R+0pc@C+0pc{
\mathbb{Z}_p^k=\mathrm{Gal}(R_k/{\mathbb{Q}_p(p^{1/p^\infty})^\wedge\left<T_1^{\pm 1},...,T_k^{\pm 1}\right>}) \ar[r]\ar[r] \ar[r]\ar[r] &\Gamma_k:=\mathrm{Gal}(R_k/{\mathbb{Q}_p\left<T_1^{\pm 1},...,T_k^{\pm 1}\right>}) \ar[r] \ar[r]\ar[r] &\Gamma_{\mathbb{Q}_p}.
}
\]

\begin{itemize}
\item (\text{Proposition}) There is an equivalence between the $\infty$-categories of inductive Banach quasicoherent presheaves\footnote{Here $\circledcirc=\text{solidquasicoherentsheaves}$.}:
\[
\xymatrix@R+6pc@C+0pc{
\mathrm{Modules}_{\circledcirc,\Gamma_k^{\times q}}(\mathcal{O}_{X_{\mathbb{Q}_p(p^{1/p^\infty})^{\wedge}\left<T_1^{\pm 1/p^\infty},...,T_k^{\pm 1/p^\infty}\right>^\flat,A,q}})\ar[d]\ar[d]\ar[d]\ar[d] \ar[r]^{\mathrm{equi}}\ar[r]\ar[r] &\varphi_q\mathrm{Modules}_{\circledcirc,\Gamma_k^{\times q}}(\mathcal{O}_{Y_{\mathbb{Q}_p(p^{1/p^\infty})^{\wedge}\left<T_1^{\pm 1/p^\infty},...,T_k^{\pm 1/p^\infty}\right>^\flat,A,q}}) \ar[d]\ar[d]\ar[d]\ar[d].\\
\mathrm{Modules}_{\circledcirc,\Gamma_0^{\times q}}(\mathcal{O}_{X_{\mathbb{Q}_p(p^{1/p^\infty})^{\wedge\flat},A,q}})\ar[r]^{\mathrm{equi}}\ar[r]\ar[r] &\varphi_q\mathrm{Modules}_{\circledcirc,\Gamma_0^{\times q}}(\mathcal{O}_{Y_{\mathbb{Q}_p(p^{1/p^\infty})^{\wedge\flat},A,q}}).\\ 
}
\]

\item (\text{Proposition}) There is an equivalence between the $\infty$-categories of inductive Banach quasicoherent commutative algebra $E_\infty$ objects\footnote{Here $\circledcirc=\text{solidquasicoherentsheaves}$.}:
\[
\xymatrix@R+6pc@C+0pc{
\mathrm{sComm}_\mathrm{simplicial}\mathrm{Modules}_{\circledcirc,\Gamma_k^{\times q}}(\mathcal{O}_{X_{R_k,A,q}})\ar[d]\ar[d]\ar[d]\ar[d]\ar[r]^{\mathrm{equi}}\ar[r]\ar[r] &\mathrm{sComm}_\mathrm{simplicial}\varphi_q\mathrm{Modules}_{\circledcirc,\Gamma_k^{\times q}}(\mathcal{O}_{Y_{R_k,A,q}})\ar[d]\ar[d]\ar[d]\ar[d]\\
\mathrm{sComm}_\mathrm{simplicial}\mathrm{Modules}_{\circledcirc,\Gamma_0^{\times q}}(\mathcal{O}_{X_{\mathbb{Q}_p(p^{1/p^\infty})^{\wedge\flat},A,q}})\ar[r]^{\mathrm{equi}}\ar[r]\ar[r] &\mathrm{sComm}_\mathrm{simplicial}\varphi_q\mathrm{Modules}_{\circledcirc,\Gamma_0^{\times q}}(\mathcal{O}_{Y_{\mathbb{Q}_p(p^{1/p^\infty})^{\wedge\flat},A,q}}).  
}
\]

\item Then as in \cite{LBV} we have the equivalence of the de Rham complex after \cite[Definition 5.9, Section 5.2.1]{KKM}\footnote{Here $\circledcirc=\text{solidquasicoherentsheaves}$.}:
\[\displayindent=-0.4in
\xymatrix@R+6pc@C+0pc{
\mathrm{deRham}_{\mathrm{sComm}_\mathrm{simplicial}\mathrm{Modules}_{\circledcirc,\Gamma_k^{\times q}}(\mathcal{O}_{X_{R_k,A,q}})\ar[r]^{\mathrm{equi}}}(-)\ar[d]\ar[d]\ar[d]\ar[d]\ar[r]\ar[r] &\mathrm{deRham}_{\mathrm{sComm}_\mathrm{simplicial}\varphi_q\mathrm{Modules}_{\circledcirc,\Gamma_k^{\times q}}(\mathcal{O}_{Y_{R_k,A,q}})}(-)\ar[d]\ar[d]\ar[d]\ar[d]\\
\mathrm{deRham}_{\mathrm{sComm}_\mathrm{simplicial}\mathrm{Modules}_{\circledcirc,\Gamma_0^{\times q}}(\mathcal{O}_{X_{\mathbb{Q}_p(p^{1/p^\infty})^{\wedge\flat},A,q}})\ar[r]^{\mathrm{equi}}}(-)\ar[r]\ar[r] &\mathrm{deRham}_{\mathrm{sComm}_\mathrm{simplicial}\varphi_q\mathrm{Modules}_{\circledcirc,\Gamma_0^{\times q}}(\mathcal{O}_{Y_{\mathbb{Q}_p(p^{1/p^\infty})^{\wedge\flat},A,q}})}(-). 
}
\]

\item Then we have the following equivalence of $K$-group $(\infty,1)$-spectrum from \cite{BGT}\footnote{Here $\circledcirc=\text{solidquasicoherentsheaves}$.}:
\[
\xymatrix@R+6pc@C+0pc{
\mathrm{K}^\mathrm{BGT}_{\mathrm{sComm}_\mathrm{simplicial}\mathrm{Modules}_{\circledcirc,\Gamma_k^{\times q}}(\mathcal{O}_{X_{R_k,A,q}})\ar[r]^{\mathrm{equi}}}(-)\ar[d]\ar[d]\ar[d]\ar[d]\ar[r]\ar[r] &\mathrm{K}^\mathrm{BGT}_{\mathrm{sComm}_\mathrm{simplicial}\varphi_q\mathrm{Modules}_{\circledcirc,\Gamma_k^{\times q}}(\mathcal{O}_{Y_{R_k,A,q}})}(-)\ar[d]\ar[d]\ar[d]\ar[d]\\
\mathrm{K}^\mathrm{BGT}_{\mathrm{sComm}_\mathrm{simplicial}\mathrm{Modules}_{\circledcirc,\Gamma_0^{\times q}}(\mathcal{O}_{X_{\mathbb{Q}_p(p^{1/p^\infty})^{\wedge\flat},A,q}})\ar[r]^{\mathrm{equi}}}(-)\ar[r]\ar[r] &\mathrm{K}^\mathrm{BGT}_{\mathrm{sComm}_\mathrm{simplicial}\varphi_q\mathrm{Modules}_{\circledcirc,\Gamma_0^{\times q}}(\mathcal{O}_{Y_{\mathbb{Q}_p(p^{1/p^\infty})^{\wedge\flat},A,q}})}(-).
}
\]

\end{itemize}

\newpage
\subsection{$\infty$-Categorical Analytic Stacks and Descents II}

As before, we have the following:

\begin{itemize}

\item (\text{Proposition}) There is a functor (global section) between the $\infty$-categories of inductive Banach quasicoherent presheaves:
\[
\xymatrix@R+0pc@C+0pc{
\mathrm{Ind}\mathrm{Banach}(\mathcal{O}_{X_{R,A,q}})\ar[r]^{\mathrm{global}}\ar[r]\ar[r] &\varphi_q\mathrm{Ind}\mathrm{Banach}(\{\mathrm{Robba}^\mathrm{extended}_{{R,A,q,I}}\}_I).  
}
\]
\item (\text{Proposition}) There is a functor (global section) between the $\infty$-categories of monomorphic inductive Banach quasicoherent presheaves:
\[
\xymatrix@R+0pc@C+0pc{
\mathrm{Ind}^m\mathrm{Banach}(\mathcal{O}_{X_{R,A,q}})\ar[r]^{\mathrm{global}}\ar[r]\ar[r] &\varphi_q\mathrm{Ind}^m\mathrm{Banach}(\{\mathrm{Robba}^\mathrm{extended}_{{R,A,q,I}}\}_I).  
}
\]

\item (\text{Proposition}) There is a functor (global section) between the $\infty$-categories of inductive Banach quasicoherent presheaves:
\[
\xymatrix@R+0pc@C+0pc{
\mathrm{Ind}\mathrm{Banach}(\mathcal{O}_{X_{R,A,q}})\ar[r]^{\mathrm{global}}\ar[r]\ar[r] &\varphi_q\mathrm{Ind}\mathrm{Banach}(\{\mathrm{Robba}^\mathrm{extended}_{{R,A,q,I}}\}_I).  
}
\]
\item (\text{Proposition}) There is a functor (global section) between the $\infty$-categories of monomorphic inductive Banach quasicoherent presheaves:
\[
\xymatrix@R+0pc@C+0pc{
\mathrm{Ind}^m\mathrm{Banach}(\mathcal{O}_{X_{R,A,q}})\ar[r]^{\mathrm{global}}\ar[r]\ar[r] &\varphi_q\mathrm{Ind}^m\mathrm{Banach}(\{\mathrm{Robba}^\mathrm{extended}_{{R,A,q,I}}\}_I).  
}
\]
\item (\text{Proposition}) There is a functor (global section) between the $\infty$-categories of inductive Banach quasicoherent commutative algebra $E_\infty$ objects:
\[
\xymatrix@R+0pc@C+0pc{
\mathrm{sComm}_\mathrm{simplicial}\mathrm{Ind}\mathrm{Banach}(\mathcal{O}_{X_{R,A,q}})\ar[r]^{\mathrm{global}}\ar[r]\ar[r] &\mathrm{sComm}_\mathrm{simplicial}\varphi_q\mathrm{Ind}\mathrm{Banach}(\{\mathrm{Robba}^\mathrm{extended}_{{R,A,q,I}}\}_I).  
}
\]
\item (\text{Proposition}) There is a functor (global section) between the $\infty$-categories of monomorphic inductive Banach quasicoherent commutative algebra $E_\infty$ objects:
\[
\xymatrix@R+0pc@C+0pc{
\mathrm{sComm}_\mathrm{simplicial}\mathrm{Ind}^m\mathrm{Banach}(\mathcal{O}_{X_{R,A,q}})\ar[r]^{\mathrm{global}}\ar[r]\ar[r] &\mathrm{sComm}_\mathrm{simplicial}\varphi_q\mathrm{Ind}^m\mathrm{Banach}(\{\mathrm{Robba}^\mathrm{extended}_{{R,A,q,I}}\}_I).  
}
\]

\item Then parallel as in \cite{LBV} we have a functor (global section) of the de Rham complex after \cite[Definition 5.9, Section 5.2.1]{KKM}:
\[
\xymatrix@R+0pc@C+0pc{
\mathrm{deRham}_{\mathrm{sComm}_\mathrm{simplicial}\mathrm{Ind}\mathrm{Banach}(\mathcal{O}_{X_{R,A,q}})\ar[r]^{\mathrm{global}}}(-)\ar[r]\ar[r] &\mathrm{deRham}_{\mathrm{sComm}_\mathrm{simplicial}\varphi_q\mathrm{Ind}\mathrm{Banach}(\{\mathrm{Robba}^\mathrm{extended}_{{R,A,q,I}}\}_I)}(-), 
}
\]
\[
\xymatrix@R+0pc@C+0pc{
\mathrm{deRham}_{\mathrm{sComm}_\mathrm{simplicial}\mathrm{Ind}^m\mathrm{Banach}(\mathcal{O}_{X_{R,A,q}})\ar[r]^{\mathrm{global}}}(-)\ar[r]\ar[r] &\mathrm{deRham}_{\mathrm{sComm}_\mathrm{simplicial}\varphi_q\mathrm{Ind}^m\mathrm{Banach}(\{\mathrm{Robba}^\mathrm{extended}_{{R,A,q,I}}\}_I)}(-). 
}
\]

\item Then we have the following a functor (global section) of $K$-group $(\infty,1)$-spectrum from \cite{BGT}:
\[
\xymatrix@R+0pc@C+0pc{
\mathrm{K}^\mathrm{BGT}_{\mathrm{sComm}_\mathrm{simplicial}\mathrm{Ind}\mathrm{Banach}(\mathcal{O}_{X_{R,A,q}})\ar[r]^{\mathrm{global}}}(-)\ar[r]\ar[r] &\mathrm{K}^\mathrm{BGT}_{\mathrm{sComm}_\mathrm{simplicial}\varphi_q\mathrm{Ind}\mathrm{Banach}(\{\mathrm{Robba}^\mathrm{extended}_{{R,A,q,I}}\}_I)}(-), 
}
\]
\[
\xymatrix@R+0pc@C+0pc{
\mathrm{K}^\mathrm{BGT}_{\mathrm{sComm}_\mathrm{simplicial}\mathrm{Ind}^m\mathrm{Banach}(\mathcal{O}_{X_{R,A,q}})\ar[r]^{\mathrm{global}}}(-)\ar[r]\ar[r] &\mathrm{K}^\mathrm{BGT}_{\mathrm{sComm}_\mathrm{simplicial}\varphi_q\mathrm{Ind}^m\mathrm{Banach}(\{\mathrm{Robba}^\mathrm{extended}_{{R,A,q,I}}\}_I)}(-). 
}
\]
\end{itemize}

\noindent Now let $R=\mathbb{Q}_p(p^{1/p^\infty})^{\wedge\flat}$ and $R_k=\mathbb{Q}_p(p^{1/p^\infty})^{\wedge}\left<T_1^{\pm 1/p^{\infty}},...,T_k^{\pm 1/p^{\infty}}\right>^\flat$ we have the following Galois theoretic results with naturality along $f:\mathrm{Spa}(\mathbb{Q}_p(p^{1/p^\infty})^{\wedge}\left<T_1^{\pm 1/p^\infty},...,T_k^{\pm 1/p^\infty}\right>^\flat)\rightarrow \mathrm{Spa}(\mathbb{Q}_p(p^{1/p^\infty})^{\wedge\flat})$:

\begin{itemize}
\item (\text{Proposition}) There is a functor (global section) between the $\infty$-categories of inductive Banach quasicoherent presheaves:
\[
\xymatrix@R+6pc@C+0pc{
\mathrm{Ind}\mathrm{Banach}(\mathcal{O}_{X_{\mathbb{Q}_p(p^{1/p^\infty})^{\wedge}\left<T_1^{\pm 1/p^\infty},...,T_k^{\pm 1/p^\infty}\right>^\flat,A,q}})\ar[d]\ar[d]\ar[d]\ar[d] \ar[r]^{\mathrm{global}}\ar[r]\ar[r] &\varphi_q\mathrm{Ind}\mathrm{Banach}(\{\mathrm{Robba}^\mathrm{extended}_{{R_k,A,q,I}}\}_I) \ar[d]\ar[d]\ar[d]\ar[d].\\
\mathrm{Ind}\mathrm{Banach}(\mathcal{O}_{X_{\mathbb{Q}_p(p^{1/p^\infty})^{\wedge\flat},A,q}})\ar[r]^{\mathrm{global}}\ar[r]\ar[r] &\varphi_q\mathrm{Ind}\mathrm{Banach}(\{\mathrm{Robba}^\mathrm{extended}_{{R_0,A,q,I}}\}_I).\\ 
}
\]
\item (\text{Proposition}) There is a functor (global section) between the $\infty$-categories of monomorphic inductive Banach quasicoherent presheaves:
\[
\xymatrix@R+6pc@C+0pc{
\mathrm{Ind}^m\mathrm{Banach}(\mathcal{O}_{X_{R_k,A,q}})\ar[r]^{\mathrm{global}}\ar[d]\ar[d]\ar[d]\ar[d]\ar[r]\ar[r] &\varphi_q\mathrm{Ind}^m\mathrm{Banach}(\{\mathrm{Robba}^\mathrm{extended}_{{R_k,A,q,I}}\}_I)\ar[d]\ar[d]\ar[d]\ar[d]\\
\mathrm{Ind}^m\mathrm{Banach}(\mathcal{O}_{X_{\mathbb{Q}_p(p^{1/p^\infty})^{\wedge\flat},A,q}})\ar[r]^{\mathrm{global}}\ar[r]\ar[r] &\varphi_q\mathrm{Ind}^m\mathrm{Banach}(\{\mathrm{Robba}^\mathrm{extended}_{{R_0,A,q,I}}\}_I).\\  
}
\]
\item (\text{Proposition}) There is a functor (global section) between the $\infty$-categories of inductive Banach quasicoherent commutative algebra $E_\infty$ objects:
\[
\xymatrix@R+6pc@C+0pc{
\mathrm{sComm}_\mathrm{simplicial}\mathrm{Ind}\mathrm{Banach}(\mathcal{O}_{X_{R_k,A,q}})\ar[d]\ar[d]\ar[d]\ar[d]\ar[r]^{\mathrm{global}}\ar[r]\ar[r] &\mathrm{sComm}_\mathrm{simplicial}\varphi_q\mathrm{Ind}\mathrm{Banach}(\{\mathrm{Robba}^\mathrm{extended}_{{R_k,A,q,I}}\}_I)\ar[d]\ar[d]\ar[d]\ar[d]\\
\mathrm{sComm}_\mathrm{simplicial}\mathrm{Ind}\mathrm{Banach}(\mathcal{O}_{X_{\mathbb{Q}_p(p^{1/p^\infty})^{\wedge\flat},A,q}})\ar[r]^{\mathrm{global}}\ar[r]\ar[r] &\mathrm{sComm}_\mathrm{simplicial}\varphi_q\mathrm{Ind}\mathrm{Banach}(\{\mathrm{Robba}^\mathrm{extended}_{{R_0,A,q,I}}\}_I).  
}
\]
\item (\text{Proposition}) There is a functor (global section) between the $\infty$-categories of monomorphic inductive Banach quasicoherent commutative algebra $E_\infty$ objects:
\[\displayindent=-0.4in
\xymatrix@R+6pc@C+0pc{
\mathrm{sComm}_\mathrm{simplicial}\mathrm{Ind}^m\mathrm{Banach}(\mathcal{O}_{X_{R_k,A,q}})\ar[d]\ar[d]\ar[d]\ar[d]\ar[r]^{\mathrm{global}}\ar[r]\ar[r] &\mathrm{sComm}_\mathrm{simplicial}\varphi_q\mathrm{Ind}^m\mathrm{Banach}(\{\mathrm{Robba}^\mathrm{extended}_{{R_k,A,q,I}}\}_I)\ar[d]\ar[d]\ar[d]\ar[d]\\
 \mathrm{sComm}_\mathrm{simplicial}\mathrm{Ind}^m\mathrm{Banach}(\mathcal{O}_{X_{\mathbb{Q}_p(p^{1/p^\infty})^{\wedge\flat},A,q}})\ar[r]^{\mathrm{global}}\ar[r]\ar[r] &\mathrm{sComm}_\mathrm{simplicial}\varphi_q\mathrm{Ind}^m\mathrm{Banach}(\{\mathrm{Robba}^\mathrm{extended}_{{R_0,A,q,I}}\}_I). 
}
\]

\item Then parallel as in \cite{LBV} we have a functor (global section) of the de Rham complex after \cite[Definition 5.9, Section 5.2.1]{KKM}:
\[\displayindent=-0.4in
\xymatrix@R+6pc@C+0pc{
\mathrm{deRham}_{\mathrm{sComm}_\mathrm{simplicial}\mathrm{Ind}\mathrm{Banach}(\mathcal{O}_{X_{R_k,A,q}})\ar[r]^{\mathrm{global}}}(-)\ar[d]\ar[d]\ar[d]\ar[d]\ar[r]\ar[r] &\mathrm{deRham}_{\mathrm{sComm}_\mathrm{simplicial}\varphi_q\mathrm{Ind}\mathrm{Banach}(\{\mathrm{Robba}^\mathrm{extended}_{{R_k,A,q,I}}\}_I)}(-)\ar[d]\ar[d]\ar[d]\ar[d]\\
\mathrm{deRham}_{\mathrm{sComm}_\mathrm{simplicial}\mathrm{Ind}\mathrm{Banach}(\mathcal{O}_{X_{\mathbb{Q}_p(p^{1/p^\infty})^{\wedge\flat},A,q}})\ar[r]^{\mathrm{global}}}(-)\ar[r]\ar[r] &\mathrm{deRham}_{\mathrm{sComm}_\mathrm{simplicial}\varphi_q\mathrm{Ind}\mathrm{Banach}(\{\mathrm{Robba}^\mathrm{extended}_{{R_0,A,q,I}}\}_I)}(-), 
}
\]
\[\displayindent=-0.4in
\xymatrix@R+6pc@C+0pc{
\mathrm{deRham}_{\mathrm{sComm}_\mathrm{simplicial}\mathrm{Ind}^m\mathrm{Banach}(\mathcal{O}_{X_{R_k,A,q}})\ar[r]^{\mathrm{global}}}(-)\ar[d]\ar[d]\ar[d]\ar[d]\ar[r]\ar[r] &\mathrm{deRham}_{\mathrm{sComm}_\mathrm{simplicial}\varphi_q\mathrm{Ind}^m\mathrm{Banach}(\{\mathrm{Robba}^\mathrm{extended}_{{R_k,A,q,I}}\}_I)}(-)\ar[d]\ar[d]\ar[d]\ar[d]\\
\mathrm{deRham}_{\mathrm{sComm}_\mathrm{simplicial}\mathrm{Ind}^m\mathrm{Banach}(\mathcal{O}_{X_{\mathbb{Q}_p(p^{1/p^\infty})^{\wedge\flat},A,q}})\ar[r]^{\mathrm{global}}}(-)\ar[r]\ar[r] &\mathrm{deRham}_{\mathrm{sComm}_\mathrm{simplicial}\varphi_q\mathrm{Ind}^m\mathrm{Banach}(\{\mathrm{Robba}^\mathrm{extended}_{{R_0,A,q,I}}\}_I)}(-). 
}
\]

\item Then we have the following a functor (global section) of $K$-group $(\infty,1)$-spectrum from \cite{BGT}:
\[
\xymatrix@R+6pc@C+0pc{
\mathrm{K}^\mathrm{BGT}_{\mathrm{sComm}_\mathrm{simplicial}\mathrm{Ind}\mathrm{Banach}(\mathcal{O}_{X_{R_k,A,q}})\ar[r]^{\mathrm{global}}}(-)\ar[d]\ar[d]\ar[d]\ar[d]\ar[r]\ar[r] &\mathrm{K}^\mathrm{BGT}_{\mathrm{sComm}_\mathrm{simplicial}\varphi_q\mathrm{Ind}\mathrm{Banach}(\{\mathrm{Robba}^\mathrm{extended}_{{R_k,A,q,I}}\}_I)}(-)\ar[d]\ar[d]\ar[d]\ar[d]\\
\mathrm{K}^\mathrm{BGT}_{\mathrm{sComm}_\mathrm{simplicial}\mathrm{Ind}\mathrm{Banach}(\mathcal{O}_{X_{\mathbb{Q}_p(p^{1/p^\infty})^{\wedge\flat},A,q}})\ar[r]^{\mathrm{global}}}(-)\ar[r]\ar[r] &\mathrm{K}^\mathrm{BGT}_{\mathrm{sComm}_\mathrm{simplicial}\varphi_q\mathrm{Ind}\mathrm{Banach}(\{\mathrm{Robba}^\mathrm{extended}_{{R_0,A,q,I}}\}_I)}(-), 
}
\]
\[
\xymatrix@R+6pc@C+0pc{
\mathrm{K}^\mathrm{BGT}_{\mathrm{sComm}_\mathrm{simplicial}\mathrm{Ind}^m\mathrm{Banach}(\mathcal{O}_{X_{R_k,A,q}})\ar[r]^{\mathrm{global}}}(-)\ar[d]\ar[d]\ar[d]\ar[d]\ar[r]\ar[r] &\mathrm{K}^\mathrm{BGT}_{\mathrm{sComm}_\mathrm{simplicial}\varphi_q\mathrm{Ind}^m\mathrm{Banach}(\{\mathrm{Robba}^\mathrm{extended}_{{R_k,A,q,I}}\}_I)}(-)\ar[d]\ar[d]\ar[d]\ar[d]\\
\mathrm{K}^\mathrm{BGT}_{\mathrm{sComm}_\mathrm{simplicial}\mathrm{Ind}^m\mathrm{Banach}(\mathcal{O}_{X_{\mathbb{Q}_p(p^{1/p^\infty})^{\wedge\flat},A,q}})\ar[r]^{\mathrm{global}}}(-)\ar[r]\ar[r] &\mathrm{K}^\mathrm{BGT}_{\mathrm{sComm}_\mathrm{simplicial}\varphi_q\mathrm{Ind}^m\mathrm{Banach}(\{\mathrm{Robba}^\mathrm{extended}_{{R_0,A,q,I}}\}_I)}(-). 
}
\]

\end{itemize}

\
\indent Then we consider further equivariance by considering the arithmetic profinite fundamental group and actually its $q$-th power $\mathrm{Gal}(\overline{{Q}_p\left<T_1^{\pm 1},...,T_k^{\pm 1}\right>}/R_k)^{\times q}$ through the following diagram\:

\[
\xymatrix@R+6pc@C+0pc{
\mathbb{Z}_p^k=\mathrm{Gal}(R_k/{\mathbb{Q}_p(p^{1/p^\infty})^\wedge\left<T_1^{\pm 1},...,T_k^{\pm 1}\right>})\ar[d]\ar[d]\ar[d]\ar[d] \ar[r]\ar[r] \ar[r]\ar[r] &\mathrm{Gal}(\overline{{Q}_p\left<T_1^{\pm 1},...,T_k^{\pm 1}\right>}/R_k) \ar[d]\ar[d]\ar[d] \ar[r]\ar[r] &\Gamma_{\mathbb{Q}_p} \ar[d]\ar[d]\ar[d]\ar[d]\\
(\mathbb{Z}_p^k=\mathrm{Gal}(R_k/{\mathbb{Q}_p(p^{1/p^\infty})^\wedge\left<T_1^{\pm 1},...,T_k^{\pm 1}\right>}))^{\times q} \ar[r]\ar[r] \ar[r]\ar[r] &\Gamma_k^{\times q}:=\mathrm{Gal}(R_k/{\mathbb{Q}_p\left<T_1^{\pm 1},...,T_k^{\pm 1}\right>})^{\times q} \ar[r] \ar[r]\ar[r] &\Gamma_{\mathbb{Q}_p}^{\times q}.
}
\]

\

We then have the correspond arithmetic profinite fundamental groups equivariant versions:
\begin{itemize}
\item (\text{Proposition}) There is a functor (global section) between the $\infty$-categories of inductive Banach quasicoherent presheaves:
\[
\xymatrix@R+6pc@C+0pc{
\mathrm{Ind}\mathrm{Banach}_{\Gamma_{k}^{\times q}}(\mathcal{O}_{X_{\mathbb{Q}_p(p^{1/p^\infty})^{\wedge}\left<T_1^{\pm 1/p^\infty},...,T_k^{\pm 1/p^\infty}\right>^\flat,A,q}})\ar[d]\ar[d]\ar[d]\ar[d] \ar[r]^{\mathrm{global}}\ar[r]\ar[r] &\varphi_q\mathrm{Ind}\mathrm{Banach}_{\Gamma_{k}^{\times q}}(\{\mathrm{Robba}^\mathrm{extended}_{{R_k,A,q,I}}\}_I) \ar[d]\ar[d]\ar[d]\ar[d].\\
\mathrm{Ind}\mathrm{Banach}_{\Gamma_{0}^{\times q}}(\mathcal{O}_{X_{\mathbb{Q}_p(p^{1/p^\infty})^{\wedge\flat},A,q}})\ar[r]^{\mathrm{global}}\ar[r]\ar[r] &\varphi_q\mathrm{Ind}\mathrm{Banach}_{\Gamma_{0}^{\times q}}(\{\mathrm{Robba}^\mathrm{extended}_{{R_0,A,q,I}}\}_I).\\ 
}
\]
\item (\text{Proposition}) There is a functor (global section) between the $\infty$-categories of monomorphic inductive Banach quasicoherent presheaves:
\[
\xymatrix@R+6pc@C+0pc{
\mathrm{Ind}^m\mathrm{Banach}_{\Gamma_{k}^{\times q}}(\mathcal{O}_{X_{R_k,A,q}})\ar[r]^{\mathrm{global}}\ar[d]\ar[d]\ar[d]\ar[d]\ar[r]\ar[r] &\varphi_q\mathrm{Ind}^m\mathrm{Banach}_{\Gamma_{k}^{\times q}}(\{\mathrm{Robba}^\mathrm{extended}_{{R_k,A,q,I}}\}_I)\ar[d]\ar[d]\ar[d]\ar[d]\\
\mathrm{Ind}^m\mathrm{Banach}_{\Gamma_{0}^{\times q}}(\mathcal{O}_{X_{\mathbb{Q}_p(p^{1/p^\infty})^{\wedge\flat},A,q}})\ar[r]^{\mathrm{global}}\ar[r]\ar[r] &\varphi_q\mathrm{Ind}^m\mathrm{Banach}_{\Gamma_{0}^{\times q}}(\{\mathrm{Robba}^\mathrm{extended}_{{R_0,A,q,I}}\}_I).\\  
}
\]
\item (\text{Proposition}) There is a functor (global section) between the $\infty$-categories of inductive Banach quasicoherent commutative algebra $E_\infty$ objects:
\[\displayindent=-0.4in
\xymatrix@R+6pc@C+0pc{
\mathrm{sComm}_\mathrm{simplicial}\mathrm{Ind}\mathrm{Banach}_{\Gamma_{k}^{\times q}}(\mathcal{O}_{X_{R_k,A,q}})\ar[d]\ar[d]\ar[d]\ar[d]\ar[r]^{\mathrm{global}}\ar[r]\ar[r] &\mathrm{sComm}_\mathrm{simplicial}\varphi_q\mathrm{Ind}\mathrm{Banach}_{\Gamma_{k}^{\times q}}(\{\mathrm{Robba}^\mathrm{extended}_{{R_k,A,q,I}}\}_I)\ar[d]\ar[d]\ar[d]\ar[d]\\
\mathrm{sComm}_\mathrm{simplicial}\mathrm{Ind}\mathrm{Banach}_{\Gamma_{0}^{\times q}}(\mathcal{O}_{X_{\mathbb{Q}_p(p^{1/p^\infty})^{\wedge\flat},A,q}})\ar[r]^{\mathrm{global}}\ar[r]\ar[r] &\mathrm{sComm}_\mathrm{simplicial}\varphi_q\mathrm{Ind}\mathrm{Banach}_{\Gamma_{0}^{\times q}}(\{\mathrm{Robba}^\mathrm{extended}_{{R_0,A,q,I}}\}_I).  
}
\]
\item (\text{Proposition}) There is a functor (global section) between the $\infty$-categories of monomorphic inductive Banach quasicoherent commutative algebra $E_\infty$ objects:
\[\displayindent=-0.4in
\xymatrix@R+6pc@C+0pc{
\mathrm{sComm}_\mathrm{simplicial}\mathrm{Ind}^m\mathrm{Banach}_{\Gamma_{k}^{\times q}}(\mathcal{O}_{X_{R_k,A,q}})\ar[d]\ar[d]\ar[d]\ar[d]\ar[r]^{\mathrm{global}}\ar[r]\ar[r] &\mathrm{sComm}_\mathrm{simplicial}\varphi_q\mathrm{Ind}^m\mathrm{Banach}_{\Gamma_{k}^{\times q}}(\{\mathrm{Robba}^\mathrm{extended}_{{R_k,A,q,I}}\}_I)\ar[d]\ar[d]\ar[d]\ar[d]\\
 \mathrm{sComm}_\mathrm{simplicial}\mathrm{Ind}^m\mathrm{Banach}_{\Gamma_{0}^{\times q}}(\mathcal{O}_{X_{\mathbb{Q}_p(p^{1/p^\infty})^{\wedge\flat},A,q}})\ar[r]^{\mathrm{global}}\ar[r]\ar[r] &\mathrm{sComm}_\mathrm{simplicial}\varphi_q\mathrm{Ind}^m\mathrm{Banach}_{\Gamma_{0}^{\times q}}(\{\mathrm{Robba}^\mathrm{extended}_{{R_0,A,q,I}}\}_I). 
}
\]

\item Then parallel as in \cite{LBV} we have a functor (global section) of the de Rham complex after \cite[Definition 5.9, Section 5.2.1]{KKM}:
\[\displayindent=-0.5in
\xymatrix@R+6pc@C+0pc{
\mathrm{deRham}_{\mathrm{sComm}_\mathrm{simplicial}\mathrm{Ind}\mathrm{Banach}_{\Gamma_{k}^{\times q}}(\mathcal{O}_{X_{R_k,A,q}})\ar[r]^{\mathrm{global}}}(-)\ar[d]\ar[d]\ar[d]\ar[d]\ar[r]\ar[r] &\mathrm{deRham}_{\mathrm{sComm}_\mathrm{simplicial}\varphi_q\mathrm{Ind}\mathrm{Banach}_{\Gamma_{k}^{\times q}}(\{\mathrm{Robba}^\mathrm{extended}_{{R_k,A,q,I}}\}_I)}(-)\ar[d]\ar[d]\ar[d]\ar[d]\\
\mathrm{deRham}_{\mathrm{sComm}_\mathrm{simplicial}\mathrm{Ind}\mathrm{Banach}_{\Gamma_{0}^{\times q}}(\mathcal{O}_{X_{\mathbb{Q}_p(p^{1/p^\infty})^{\wedge\flat},A,q}})\ar[r]^{\mathrm{global}}}(-)\ar[r]\ar[r] &\mathrm{deRham}_{\mathrm{sComm}_\mathrm{simplicial}\varphi_q\mathrm{Ind}\mathrm{Banach}_{\Gamma_{0}^{\times q}}(\{\mathrm{Robba}^\mathrm{extended}_{{R_0,A,q,I}}\}_I)}(-), 
}
\]
\[\displayindent=-0.7in
\xymatrix@R+6pc@C+0pc{
\mathrm{deRham}_{\mathrm{sComm}_\mathrm{simplicial}\mathrm{Ind}^m\mathrm{Banach}_{\Gamma_{k}^{\times q}}(\mathcal{O}_{X_{R_k,A,q}})\ar[r]^{\mathrm{global}}}(-)\ar[d]\ar[d]\ar[d]\ar[d]\ar[r]\ar[r] &\mathrm{deRham}_{\mathrm{sComm}_\mathrm{simplicial}\varphi_q\mathrm{Ind}^m\mathrm{Banach}_{\Gamma_{k}^{\times q}}(\{\mathrm{Robba}^\mathrm{extended}_{{R_k,A,q,I}}\}_I)}(-)\ar[d]\ar[d]\ar[d]\ar[d]\\
\mathrm{deRham}_{\mathrm{sComm}_\mathrm{simplicial}\mathrm{Ind}^m\mathrm{Banach}_{\Gamma_{0}^{\times q}}(\mathcal{O}_{X_{\mathbb{Q}_p(p^{1/p^\infty})^{\wedge\flat},A,q}})\ar[r]^{\mathrm{global}}}(-)\ar[r]\ar[r] &\mathrm{deRham}_{\mathrm{sComm}_\mathrm{simplicial}\varphi_q\mathrm{Ind}^m\mathrm{Banach}_{\Gamma_{0}^{\times q}}(\{\mathrm{Robba}^\mathrm{extended}_{{R_0,A,q,I}}\}_I)}(-). 
}
\]

\item Then we have the following a functor (global section) of $K$-group $(\infty,1)$-spectrum from \cite{BGT}:
\[
\xymatrix@R+6pc@C+0pc{
\mathrm{K}^\mathrm{BGT}_{\mathrm{sComm}_\mathrm{simplicial}\mathrm{Ind}\mathrm{Banach}_{\Gamma_{k}^{\times q}}(\mathcal{O}_{X_{R_k,A,q}})\ar[r]^{\mathrm{global}}}(-)\ar[d]\ar[d]\ar[d]\ar[d]\ar[r]\ar[r] &\mathrm{K}^\mathrm{BGT}_{\mathrm{sComm}_\mathrm{simplicial}\varphi_q\mathrm{Ind}\mathrm{Banach}_{\Gamma_{k}^{\times q}}(\{\mathrm{Robba}^\mathrm{extended}_{{R_k,A,q,I}}\}_I)}(-)\ar[d]\ar[d]\ar[d]\ar[d]\\
\mathrm{K}^\mathrm{BGT}_{\mathrm{sComm}_\mathrm{simplicial}\mathrm{Ind}\mathrm{Banach}_{\Gamma_{0}^{\times q}}(\mathcal{O}_{X_{\mathbb{Q}_p(p^{1/p^\infty})^{\wedge\flat},A,q}})\ar[r]^{\mathrm{global}}}(-)\ar[r]\ar[r] &\mathrm{K}^\mathrm{BGT}_{\mathrm{sComm}_\mathrm{simplicial}\varphi_q\mathrm{Ind}\mathrm{Banach}_{\Gamma_{0}^{\times q}}(\{\mathrm{Robba}^\mathrm{extended}_{{R_0,A,q,I}}\}_I)}(-), 
}
\]
\[
\xymatrix@R+6pc@C+0pc{
\mathrm{K}^\mathrm{BGT}_{\mathrm{sComm}_\mathrm{simplicial}\mathrm{Ind}^m\mathrm{Banach}_{\Gamma_{k}^{\times q}}(\mathcal{O}_{X_{R_k,A,q}})\ar[r]^{\mathrm{global}}}(-)\ar[d]\ar[d]\ar[d]\ar[d]\ar[r]\ar[r] &\mathrm{K}^\mathrm{BGT}_{\mathrm{sComm}_\mathrm{simplicial}\varphi_q\mathrm{Ind}^m\mathrm{Banach}_{\Gamma_{k}^{\times q}}(\{\mathrm{Robba}^\mathrm{extended}_{{R_k,A,q,I}}\}_I)}(-)\ar[d]\ar[d]\ar[d]\ar[d]\\
\mathrm{K}^\mathrm{BGT}_{\mathrm{sComm}_\mathrm{simplicial}\mathrm{Ind}^m\mathrm{Banach}_{\Gamma_{0}^{\times q}}(\mathcal{O}_{X_{\mathbb{Q}_p(p^{1/p^\infty})^{\wedge\flat},A,q}})\ar[r]^{\mathrm{global}}}(-)\ar[r]\ar[r] &\mathrm{K}^\mathrm{BGT}_{\mathrm{sComm}_\mathrm{simplicial}\varphi_q\mathrm{Ind}^m\mathrm{Banach}_{\Gamma_{0}^{\times q}}(\{\mathrm{Robba}^\mathrm{extended}_{{R_0,A,q,I}}\}_I)}(-). 
}
\]

\end{itemize}

\

\

\begin{remark}
\noindent We can certainly consider the quasicoherent sheaves in \cite[Lemma 7.11, Remark 7.12]{1BBK}, therefore all the quasicoherent presheaves and modules will be those satisfying \cite[Lemma 7.11, Remark 7.12]{1BBK} if one would like to consider the the quasicoherent sheaves. That being all as this said, we would believe that the big quasicoherent presheaves are automatically quasicoherent sheaves (namely satisfying the corresponding \v{C}ech $\infty$-descent as in \cite[Section 9.3]{KKM} and \cite[Lemma 7.11, Remark 7.12]{1BBK}) and the corresponding global section functors are automatically equivalence of $\infty$-categories. 
\end{remark}

\

\indent In Clausen-Scholze formalism we have the following:

\begin{itemize}

\item (\text{Proposition}) There is a functor (global section) between the $\infty$-categories of inductive Banach quasicoherent sheaves:
\[
\xymatrix@R+0pc@C+0pc{
\mathrm{Modules}_\circledcirc(\mathcal{O}_{X_{R,A,q}})\ar[r]^{\mathrm{global}}\ar[r]\ar[r] &\varphi_q\mathrm{Modules}_\circledcirc(\{\mathrm{Robba}^{\mathrm{extended},q}_{{R,A,I}}\}_I).  
}
\]

\item (\text{Proposition}) There is a functor (global section) between the $\infty$-categories of inductive Banach quasicoherent sheaves:
\[
\xymatrix@R+0pc@C+0pc{
\mathrm{Modules}_\circledcirc(\mathcal{O}_{X_{R,A,q}})\ar[r]^{\mathrm{global}}\ar[r]\ar[r] &\varphi_q\mathrm{Modules}_\circledcirc(\{\mathrm{Robba}^{\mathrm{extended},q}_{{R,A,I}}\}_I).  
}
\]

\item (\text{Proposition}) There is a functor (global section) between the $\infty$-categories of inductive Banach quasicoherent commutative algebra $E_\infty$ objects\footnote{Here $\circledcirc=\text{solidquasicoherentsheaves}$.}:
\[
\xymatrix@R+0pc@C+0pc{
\mathrm{sComm}_\mathrm{simplicial}\mathrm{Modules}_\circledcirc(\mathcal{O}_{X_{R,A,q}})\ar[r]^{\mathrm{global}}\ar[r]\ar[r] &\mathrm{sComm}_\mathrm{simplicial}\varphi_q\mathrm{Modules}_\circledcirc(\{\mathrm{Robba}^{\mathrm{extended},q}_{{R,A,I}}\}_I).  
}
\]

\item Then as in \cite{LBV} we have a functor (global section) of the de Rham complex after \cite[Definition 5.9, Section 5.2.1]{KKM}\footnote{Here $\circledcirc=\text{solidquasicoherentsheaves}$.}:
\[
\xymatrix@R+0pc@C+0pc{
\mathrm{deRham}_{\mathrm{sComm}_\mathrm{simplicial}\mathrm{Modules}_\circledcirc(\mathcal{O}_{X_{R,A,q}})\ar[r]^{\mathrm{global}}}(-)\ar[r]\ar[r] &\mathrm{deRham}_{\mathrm{sComm}_\mathrm{simplicial}\varphi_q\mathrm{Modules}_\circledcirc(\{\mathrm{Robba}^{\mathrm{extended},q}_{{R,A,I}}\}_I)}(-). 
}
\]

\item Then we have the following a functor (global section) of $K$-group $(\infty,1)$-spectrum from \cite{BGT}\footnote{Here $\circledcirc=\text{solidquasicoherentsheaves}$.}:
\[
\xymatrix@R+0pc@C+0pc{
\mathrm{K}^\mathrm{BGT}_{\mathrm{sComm}_\mathrm{simplicial}\mathrm{Module}_\circledcirc(\mathcal{O}_{X_{R,A,q}})\ar[r]^{\mathrm{global}}}(-)\ar[r]\ar[r] &\mathrm{K}^\mathrm{BGT}_{\mathrm{sComm}_\mathrm{simplicial}\varphi_q\mathrm{Module}_\circledcirc(\{\mathrm{Robba}^{\mathrm{extended},q}_{{R,A,I}}\}_I)}(-). 
}
\]
\end{itemize}

\noindent Now let $R=\mathbb{Q}_p(p^{1/p^\infty})^{\wedge\flat}$ and $R_k=\mathbb{Q}_p(p^{1/p^\infty})^{\wedge}\left<T_1^{\pm 1/p^{\infty}},...,T_k^{\pm 1/p^{\infty}}\right>^\flat$ we have the following Galois theoretic results with naturality along $f:\mathrm{Spa}(\mathbb{Q}_p(p^{1/p^\infty})^{\wedge}\left<T_1^{\pm 1/p^\infty},...,T_k^{\pm 1/p^\infty}\right>^\flat)\rightarrow \mathrm{Spa}(\mathbb{Q}_p(p^{1/p^\infty})^{\wedge\flat})$:

\begin{itemize}
\item (\text{Proposition}) There is a functor (global section) between the $\infty$-categories of inductive Banach quasicoherent sheaves\footnote{Here $\circledcirc=\text{solidquasicoherentsheaves}$.}:
\[
\xymatrix@R+6pc@C+0pc{
\mathrm{Module}_\circledcirc(\mathcal{O}_{X_{\mathbb{Q}_p(p^{1/p^\infty})^{\wedge}\left<T_1^{\pm 1/p^\infty},...,T_k^{\pm 1/p^\infty}\right>^\flat,A,q}})\ar[d]\ar[d]\ar[d]\ar[d] \ar[r]^{\mathrm{global}}\ar[r]\ar[r] &\varphi_q\mathrm{Module}_\circledcirc(\{\mathrm{Robba}^{\mathrm{extended},q}_{{R_k,A,I}}\}_I) \ar[d]\ar[d]\ar[d]\ar[d].\\
\mathrm{Module}_\circledcirc(\mathcal{O}_{X_{\mathbb{Q}_p(p^{1/p^\infty})^{\wedge\flat},A,q}})\ar[r]^{\mathrm{global}}\ar[r]\ar[r] &\varphi_q\mathrm{Module}_\circledcirc(\{\mathrm{Robba}^{\mathrm{extended},q}_{{R_0,A,I}}\}_I).\\ 
}
\]

\item (\text{Proposition}) There is a functor (global section) between the $\infty$-categories of inductive Banach quasicoherent commutative algebra $E_\infty$ objects\footnote{Here $\circledcirc=\text{solidquasicoherentsheaves}$.}:
\[
\xymatrix@R+6pc@C+0pc{
\mathrm{sComm}_\mathrm{simplicial}\mathrm{Module}_\circledcirc(\mathcal{O}_{X_{R_k,A,q}})\ar[d]\ar[d]\ar[d]\ar[d]\ar[r]^{\mathrm{global}}\ar[r]\ar[r] &\mathrm{sComm}_\mathrm{simplicial}\varphi_q\mathrm{Module}_\circledcirc(\{\mathrm{Robba}^{\mathrm{extended},q}_{{R_k,A,I}}\}_I)\ar[d]\ar[d]\ar[d]\ar[d]\\
\mathrm{sComm}_\mathrm{simplicial}\mathrm{Module}_\circledcirc(\mathcal{O}_{X_{\mathbb{Q}_p(p^{1/p^\infty})^{\wedge\flat},A,q}})\ar[r]^{\mathrm{global}}\ar[r]\ar[r] &\mathrm{sComm}_\mathrm{simplicial}\varphi_q\mathrm{Module}_\circledcirc(\{\mathrm{Robba}^{\mathrm{extended},q}_{{R_0,A,I}}\}_I).  
}
\]

\item Then as in \cite{LBV} we have a functor (global section) of the de Rham complex after \cite[Definition 5.9, Section 5.2.1]{KKM}\footnote{Here $\circledcirc=\text{solidquasicoherentsheaves}$.}:
\[\displayindent=-0.4in
\xymatrix@R+6pc@C+0pc{
\mathrm{deRham}_{\mathrm{sComm}_\mathrm{simplicial}\mathrm{Module}_\circledcirc(\mathcal{O}_{X_{R_k,A,q}})\ar[r]^{\mathrm{global}}}(-)\ar[d]\ar[d]\ar[d]\ar[d]\ar[r]\ar[r] &\mathrm{deRham}_{\mathrm{sComm}_\mathrm{simplicial}\varphi_q\mathrm{Module}_\circledcirc(\{\mathrm{Robba}^{\mathrm{extended},q}_{{R_k,A,I}}\}_I)}(-)\ar[d]\ar[d]\ar[d]\ar[d]\\
\mathrm{deRham}_{\mathrm{sComm}_\mathrm{simplicial}\mathrm{Module}_\circledcirc(\mathcal{O}_{X_{\mathbb{Q}_p(p^{1/p^\infty})^{\wedge\flat},A,q}})\ar[r]^{\mathrm{global}}}(-)\ar[r]\ar[r] &\mathrm{deRham}_{\mathrm{sComm}_\mathrm{simplicial}\varphi_q\mathrm{Module}_\circledcirc(\{\mathrm{Robba}^{\mathrm{extended},q}_{{R_0,A,I}}\}_I)}(-). 
}
\]

\item Then we have the following a functor (global section) of $K$-group $(\infty,1)$-spectrum from \cite{BGT}\footnote{Here $\circledcirc=\text{solidquasicoherentsheaves}$.}:
\[
\xymatrix@R+6pc@C+0pc{
\mathrm{K}^\mathrm{BGT}_{\mathrm{sComm}_\mathrm{simplicial}\mathrm{Module}_\circledcirc(\mathcal{O}_{X_{R_k,A,q}})\ar[r]^{\mathrm{global}}}(-)\ar[d]\ar[d]\ar[d]\ar[d]\ar[r]\ar[r] &\mathrm{K}^\mathrm{BGT}_{\mathrm{sComm}_\mathrm{simplicial}\varphi_q\mathrm{Module}_\circledcirc(\{\mathrm{Robba}^{\mathrm{extended},q}_{{R_k,A,I}}\}_I)}(-)\ar[d]\ar[d]\ar[d]\ar[d]\\
\mathrm{K}^\mathrm{BGT}_{\mathrm{sComm}_\mathrm{simplicial}\mathrm{Module}_\circledcirc(\mathcal{O}_{X_{\mathbb{Q}_p(p^{1/p^\infty})^{\wedge\flat},A,q}})\ar[r]^{\mathrm{global}}}(-)\ar[r]\ar[r] &\mathrm{K}^\mathrm{BGT}_{\mathrm{sComm}_\mathrm{simplicial}\varphi_q\mathrm{Module}_\circledcirc(\{\mathrm{Robba}^{\mathrm{extended},q}_{{R_0,A,I}}\}_I)}(-). 
}
\]

\end{itemize}

\
\indent Then we consider further equivariance by considering the arithmetic profinite fundamental groups $\Gamma_{\mathbb{Q}_p}$ and $\mathrm{Gal}(\overline{\mathbb{Q}_p\left<T_1^{\pm 1},...,T_k^{\pm 1}\right>}/R_k)$ through the following diagram:

\[
\xymatrix@R+0pc@C+0pc{
\mathbb{Z}_p^k=\mathrm{Gal}(R_k/{\mathbb{Q}_p(p^{1/p^\infty})^\wedge\left<T_1^{\pm 1},...,T_k^{\pm 1}\right>}) \ar[r]\ar[r] \ar[r]\ar[r] &\Gamma_k^{\times q}:=\mathrm{Gal}(R_k/{\mathbb{Q}_p\left<T_1^{\pm 1},...,T_k^{\pm 1}\right>}) \ar[r] \ar[r]\ar[r] &\Gamma_{\mathbb{Q}_p}.
}
\]

\begin{itemize}
\item (\text{Proposition}) There is a functor (global section) between the $\infty$-categories of inductive Banach quasicoherent sheaves\footnote{Here $\circledcirc=\text{solidquasicoherentsheaves}$.}:
\[
\xymatrix@R+6pc@C+0pc{
{\mathrm{Module}_\circledcirc}_{\Gamma_k^{\times q}}(\mathcal{O}_{X_{\mathbb{Q}_p(p^{1/p^\infty})^{\wedge}\left<T_1^{\pm 1/p^\infty},...,T_k^{\pm 1/p^\infty}\right>^\flat,A,q}})\ar[d]\ar[d]\ar[d]\ar[d] \ar[r]^{\mathrm{global}}\ar[r]\ar[r] &\varphi_q{\mathrm{Module}_\circledcirc}_{\Gamma_k^{\times q}}(\{\mathrm{Robba}^{\mathrm{extended},q}_{{R_k,A,I}}\}_I) \ar[d]\ar[d]\ar[d]\ar[d].\\
{\mathrm{Module}_\circledcirc}(\mathcal{O}_{X_{\mathbb{Q}_p(p^{1/p^\infty})^{\wedge\flat},A,q}})\ar[r]^{\mathrm{global}}\ar[r]\ar[r] &\varphi_q{\mathrm{Module}_\circledcirc}(\{\mathrm{Robba}^{\mathrm{extended},q}_{{R_0,A,I}}\}_I).\\ 
}
\]

\item (\text{Proposition}) There is a functor (global section) between the $\infty$-categories of inductive Banach quasicoherent commutative algebra $E_\infty$ objects\footnote{Here $\circledcirc=\text{solidquasicoherentsheaves}$.}:
\[\displayindent=-0.4in
\xymatrix@R+6pc@C+0pc{
\mathrm{sComm}_\mathrm{simplicial}{\mathrm{Module}_\circledcirc}_{\Gamma_k^{\times q}}(\mathcal{O}_{X_{R_k,A,q}})\ar[d]\ar[d]\ar[d]\ar[d]\ar[r]^{\mathrm{global}}\ar[r]\ar[r] &\mathrm{sComm}_\mathrm{simplicial}\varphi_q{\mathrm{Module}_\circledcirc}(\{\mathrm{Robba}^{\mathrm{extended},q}_{{R_k,A,I}}\}_I)\ar[d]\ar[d]\ar[d]\ar[d]\\
\mathrm{sComm}_\mathrm{simplicial}{\mathrm{Module}_\circledcirc}_{\Gamma_0^{\times q}}(\mathcal{O}_{X_{\mathbb{Q}_p(p^{1/p^\infty})^{\wedge\flat},A,q}})\ar[r]^{\mathrm{global}}\ar[r]\ar[r] &\mathrm{sComm}_\mathrm{simplicial}\varphi_q{\mathrm{Modules}_\circledcirc}_{\Gamma_0^{\times q}}(\{\mathrm{Robba}^{\mathrm{extended},q}_{{R_0,A,I}}\}_I).  
}
\]

\item Then as in \cite{LBV} we have a functor (global section) of the de Rham complex after \cite[Definition 5.9, Section 5.2.1]{KKM}\footnote{Here $\circledcirc=\text{solidquasicoherentsheaves}$.}:
\[\displayindent=-0.5in
\xymatrix@R+6pc@C+0pc{
\mathrm{deRham}_{\mathrm{sComm}_\mathrm{simplicial}{\mathrm{Modules}_\circledcirc}_{\Gamma_k^{\times q}}(\mathcal{O}_{X_{R_k,A,q}})\ar[r]^{\mathrm{global}}}(-)\ar[d]\ar[d]\ar[d]\ar[d]\ar[r]\ar[r] &\mathrm{deRham}_{\mathrm{sComm}_\mathrm{simplicial}\varphi_q{\mathrm{Modules}_\circledcirc}_{\Gamma_k^{\times q}}(\{\mathrm{Robba}^{\mathrm{extended},q}_{{R_k,A,I}}\}_I)}(-)\ar[d]\ar[d]\ar[d]\ar[d]\\
\mathrm{deRham}_{\mathrm{sComm}_\mathrm{simplicial}{\mathrm{Modules}_\circledcirc}_{\Gamma_0^{\times q}}(\mathcal{O}_{X_{\mathbb{Q}_p(p^{1/p^\infty})^{\wedge\flat},A,q}})\ar[r]^{\mathrm{global}}}(-)\ar[r]\ar[r] &\mathrm{deRham}_{\mathrm{sComm}_\mathrm{simplicial}\varphi_q{\mathrm{Modules}_\circledcirc}_{\Gamma_0^{\times q}}(\{\mathrm{Robba}^{\mathrm{extended},q}_{{R_0,A,I}}\}_I)}(-). 
}
\]

\item Then we have the following a functor (global section) of $K$-group $(\infty,1)$-spectrum from \cite{BGT}\footnote{Here $\circledcirc=\text{solidquasicoherentsheaves}$.}:
\[
\xymatrix@R+6pc@C+0pc{
\mathrm{K}^\mathrm{BGT}_{\mathrm{sComm}_\mathrm{simplicial}{\mathrm{Modules}_\circledcirc}_{\Gamma_k^{\times q}}(\mathcal{O}_{X_{R_k,A,q}})\ar[r]^{\mathrm{global}}}(-)\ar[d]\ar[d]\ar[d]\ar[d]\ar[r]\ar[r] &\mathrm{K}^\mathrm{BGT}_{\mathrm{sComm}_\mathrm{simplicial}\varphi_q{\mathrm{Modules}_\circledcirc}_{\Gamma_k^{\times q}}(\{\mathrm{Robba}^{\mathrm{extended},q}_{{R_k,A,I}}\}_I)}(-)\ar[d]\ar[d]\ar[d]\ar[d]\\
\mathrm{K}^\mathrm{BGT}_{\mathrm{sComm}_\mathrm{simplicial}{\mathrm{Modules}_\circledcirc}_{\Gamma_0^{\times q}}(\mathcal{O}_{X_{\mathbb{Q}_p(p^{1/p^\infty})^{\wedge\flat},A,q}})\ar[r]^{\mathrm{global}}}(-)\ar[r]\ar[r] &\mathrm{K}^\mathrm{BGT}_{\mathrm{sComm}_\mathrm{simplicial}\varphi_q{\mathrm{Modules}_\circledcirc}_{\Gamma_0^{\times q}}(\{\mathrm{Robba}^{\mathrm{extended},q}_{{R_0,A,I}}\}_I)}(-). 
}
\]
	
\end{itemize}

\begin{proposition}
All the global functors from \cite[Proposition 13.8, Theorem 14.9, Remark 14.10]{1CS2} achieve the equivalences on both sides.	
\end{proposition}

\newpage
\subsection{$\infty$-Categorical Analytic Stacks and Descents III}

\indent As before, we have the following. Let $\mathcal{A}$ vary in the category of all the Banach algebras over $\mathbb{Q}_p$ we have the following.

\begin{itemize}

\item (\text{Proposition}) There is a functor (global section) between the $\infty$-prestacks of inductive Banach quasicoherent presheaves:
\[
\xymatrix@R+0pc@C+0pc{
\mathrm{Ind}\mathrm{Banach}(\mathcal{O}_{X_{R,-,q}})\ar[r]^{\mathrm{global}}\ar[r]\ar[r] &\varphi_q\mathrm{Ind}\mathrm{Banach}(\{\mathrm{Robba}^\mathrm{extended}_{{R,-,I,q}}\}_I).  
}
\]
\item (\text{Proposition}) There is a functor (global section) between the $\infty$-prestacks of monomorphic inductive Banach quasicoherent presheaves:
\[
\xymatrix@R+0pc@C+0pc{
\mathrm{Ind}^m\mathrm{Banach}(\mathcal{O}_{X_{R,-,q}})\ar[r]^{\mathrm{global}}\ar[r]\ar[r] &\varphi_q\mathrm{Ind}^m\mathrm{Banach}(\{\mathrm{Robba}^\mathrm{extended}_{{R,-,I,q}}\}_I).  
}
\]

\item (\text{Proposition}) There is a functor (global section) between the $\infty$-prestacks of inductive Banach quasicoherent presheaves:
\[
\xymatrix@R+0pc@C+0pc{
\mathrm{Ind}\mathrm{Banach}(\mathcal{O}_{X_{R,-,q}})\ar[r]^{\mathrm{global}}\ar[r]\ar[r] &\varphi_q\mathrm{Ind}\mathrm{Banach}(\{\mathrm{Robba}^\mathrm{extended}_{{R,-,I,q}}\}_I).  
}
\]
\item (\text{Proposition}) There is a functor (global section) between the $\infty$-stacks of monomorphic inductive Banach quasicoherent presheaves:
\[
\xymatrix@R+0pc@C+0pc{
\mathrm{Ind}^m\mathrm{Banach}(\mathcal{O}_{X_{R,-,q}})\ar[r]^{\mathrm{global}}\ar[r]\ar[r] &\varphi_q\mathrm{Ind}^m\mathrm{Banach}(\{\mathrm{Robba}^\mathrm{extended}_{{R,-,I,q}}\}_I).  
}
\]
\item (\text{Proposition}) There is a functor (global section) between the $\infty$-prestacks of inductive Banach quasicoherent commutative algebra $E_\infty$ objects:
\[
\xymatrix@R+0pc@C+0pc{
\mathrm{sComm}_\mathrm{simplicial}\mathrm{Ind}\mathrm{Banach}(\mathcal{O}_{X_{R,-,q}})\ar[r]^{\mathrm{global}}\ar[r]\ar[r] &\mathrm{sComm}_\mathrm{simplicial}\varphi_q\mathrm{Ind}\mathrm{Banach}(\{\mathrm{Robba}^\mathrm{extended}_{{R,-,I,q}}\}_I).  
}
\]
\item (\text{Proposition}) There is a functor (global section) between the $\infty$-prestacks of monomorphic inductive Banach quasicoherent commutative algebra $E_\infty$ objects:
\[
\xymatrix@R+0pc@C+0pc{
\mathrm{sComm}_\mathrm{simplicial}\mathrm{Ind}^m\mathrm{Banach}(\mathcal{O}_{X_{R,-,q}})\ar[r]^{\mathrm{global}}\ar[r]\ar[r] &\mathrm{sComm}_\mathrm{simplicial}\varphi_q\mathrm{Ind}^m\mathrm{Banach}(\{\mathrm{Robba}^\mathrm{extended}_{{R,-,I,q}}\}_I).  
}
\]

\item Then parallel as in \cite{LBV} we have a functor (global section) of the de Rham complex after \cite[Definition 5.9, Section 5.2.1]{KKM}:
\[
\xymatrix@R+0pc@C+0pc{
\mathrm{deRham}_{\mathrm{sComm}_\mathrm{simplicial}\mathrm{Ind}\mathrm{Banach}(\mathcal{O}_{X_{R,-,q}})\ar[r]^{\mathrm{global}}}(-)\ar[r]\ar[r] &\mathrm{deRham}_{\mathrm{sComm}_\mathrm{simplicial}\varphi_q\mathrm{Ind}\mathrm{Banach}(\{\mathrm{Robba}^\mathrm{extended}_{{R,-,I,q}}\}_I)}(-), 
}
\]
\[
\xymatrix@R+0pc@C+0pc{
\mathrm{deRham}_{\mathrm{sComm}_\mathrm{simplicial}\mathrm{Ind}^m\mathrm{Banach}(\mathcal{O}_{X_{R,-,q}})\ar[r]^{\mathrm{global}}}(-)\ar[r]\ar[r] &\mathrm{deRham}_{\mathrm{sComm}_\mathrm{simplicial}\varphi_q\mathrm{Ind}^m\mathrm{Banach}(\{\mathrm{Robba}^\mathrm{extended}_{{R,-,I,q}}\}_I)}(-). 
}
\]

\item Then we have the following a functor (global section) of $K$-group $(\infty,1)$-spectrum from \cite{BGT}:
\[
\xymatrix@R+0pc@C+0pc{
\mathrm{K}^\mathrm{BGT}_{\mathrm{sComm}_\mathrm{simplicial}\mathrm{Ind}\mathrm{Banach}(\mathcal{O}_{X_{R,-,q}})\ar[r]^{\mathrm{global}}}(-)\ar[r]\ar[r] &\mathrm{K}^\mathrm{BGT}_{\mathrm{sComm}_\mathrm{simplicial}\varphi_q\mathrm{Ind}\mathrm{Banach}(\{\mathrm{Robba}^\mathrm{extended}_{{R,-,I,q}}\}_I)}(-), 
}
\]
\[
\xymatrix@R+0pc@C+0pc{
\mathrm{K}^\mathrm{BGT}_{\mathrm{sComm}_\mathrm{simplicial}\mathrm{Ind}^m\mathrm{Banach}(\mathcal{O}_{X_{R,-,q}})\ar[r]^{\mathrm{global}}}(-)\ar[r]\ar[r] &\mathrm{K}^\mathrm{BGT}_{\mathrm{sComm}_\mathrm{simplicial}\varphi_q\mathrm{Ind}^m\mathrm{Banach}(\{\mathrm{Robba}^\mathrm{extended}_{{R,-,I,q}}\}_I)}(-). 
}
\]
\end{itemize}

\noindent Now let $R=\mathbb{Q}_p(p^{1/p^\infty})^{\wedge\flat}$ and $R_k=\mathbb{Q}_p(p^{1/p^\infty})^{\wedge}\left<T_1^{\pm 1/p^{\infty}},...,T_k^{\pm 1/p^{\infty}}\right>^\flat$ we have the following Galois theoretic results with naturality along $f:\mathrm{Spa}(\mathbb{Q}_p(p^{1/p^\infty})^{\wedge}\left<T_1^{\pm 1/p^\infty},...,T_k^{\pm 1/p^\infty}\right>^\flat)\rightarrow \mathrm{Spa}(\mathbb{Q}_p(p^{1/p^\infty})^{\wedge\flat})$:

\begin{itemize}
\item (\text{Proposition}) There is a functor (global section) between the $\infty$-prestacks of inductive Banach quasicoherent presheaves:
\[
\xymatrix@R+6pc@C+0pc{
\mathrm{Ind}\mathrm{Banach}(\mathcal{O}_{X_{\mathbb{Q}_p(p^{1/p^\infty})^{\wedge}\left<T_1^{\pm 1/p^\infty},...,T_k^{\pm 1/p^\infty}\right>^\flat,-,q}})\ar[d]\ar[d]\ar[d]\ar[d] \ar[r]^{\mathrm{global}}\ar[r]\ar[r] &\varphi_q\mathrm{Ind}\mathrm{Banach}(\{\mathrm{Robba}^\mathrm{extended}_{{R_k,-,I,q}}\}_I) \ar[d]\ar[d]\ar[d]\ar[d].\\
\mathrm{Ind}\mathrm{Banach}(\mathcal{O}_{X_{\mathbb{Q}_p(p^{1/p^\infty})^{\wedge\flat},-,q}})\ar[r]^{\mathrm{global}}\ar[r]\ar[r] &\varphi_q\mathrm{Ind}\mathrm{Banach}(\{\mathrm{Robba}^\mathrm{extended}_{{R_0,-,I,q}}\}_I).\\ 
}
\]
\item (\text{Proposition}) There is a functor (global section) between the $\infty$-prestacks of monomorphic inductive Banach quasicoherent presheaves:
\[
\xymatrix@R+6pc@C+0pc{
\mathrm{Ind}^m\mathrm{Banach}(\mathcal{O}_{X_{R_k,-,q}})\ar[r]^{\mathrm{global}}\ar[d]\ar[d]\ar[d]\ar[d]\ar[r]\ar[r] &\varphi_q\mathrm{Ind}^m\mathrm{Banach}(\{\mathrm{Robba}^\mathrm{extended}_{{R_k,-,I,q}}\}_I)\ar[d]\ar[d]\ar[d]\ar[d]\\
\mathrm{Ind}^m\mathrm{Banach}(\mathcal{O}_{X_{\mathbb{Q}_p(p^{1/p^\infty})^{\wedge\flat},-,q}})\ar[r]^{\mathrm{global}}\ar[r]\ar[r] &\varphi_q\mathrm{Ind}^m\mathrm{Banach}(\{\mathrm{Robba}^\mathrm{extended}_{{R_0,-,I,q}}\}_I).\\  
}
\]
\item (\text{Proposition}) There is a functor (global section) between the $\infty$-prestacks of inductive Banach quasicoherent commutative algebra $E_\infty$ objects:
\[
\xymatrix@R+6pc@C+0pc{
\mathrm{sComm}_\mathrm{simplicial}\mathrm{Ind}\mathrm{Banach}(\mathcal{O}_{X_{R_k,-,q}})\ar[d]\ar[d]\ar[d]\ar[d]\ar[r]^{\mathrm{global}}\ar[r]\ar[r] &\mathrm{sComm}_\mathrm{simplicial}\varphi_q\mathrm{Ind}\mathrm{Banach}(\{\mathrm{Robba}^\mathrm{extended}_{{R_k,-,I,q}}\}_I)\ar[d]\ar[d]\ar[d]\ar[d]\\
\mathrm{sComm}_\mathrm{simplicial}\mathrm{Ind}\mathrm{Banach}(\mathcal{O}_{X_{\mathbb{Q}_p(p^{1/p^\infty})^{\wedge\flat},-,q}})\ar[r]^{\mathrm{global}}\ar[r]\ar[r] &\mathrm{sComm}_\mathrm{simplicial}\varphi_q\mathrm{Ind}\mathrm{Banach}(\{\mathrm{Robba}^\mathrm{extended}_{{R_0,-,I,q}}\}_I).  
}
\]
\item (\text{Proposition}) There is a functor (global section) between the $\infty$-prestacks of monomorphic inductive Banach quasicoherent commutative algebra $E_\infty$ objects:
\[\displayindent=-0.4in
\xymatrix@R+6pc@C+0pc{
\mathrm{sComm}_\mathrm{simplicial}\mathrm{Ind}^m\mathrm{Banach}(\mathcal{O}_{X_{R_k,-,q}})\ar[d]\ar[d]\ar[d]\ar[d]\ar[r]^{\mathrm{global}}\ar[r]\ar[r] &\mathrm{sComm}_\mathrm{simplicial}\varphi_q\mathrm{Ind}^m\mathrm{Banach}(\{\mathrm{Robba}^\mathrm{extended}_{{R_k,-,I,q}}\}_I)\ar[d]\ar[d]\ar[d]\ar[d]\\
 \mathrm{sComm}_\mathrm{simplicial}\mathrm{Ind}^m\mathrm{Banach}(\mathcal{O}_{X_{\mathbb{Q}_p(p^{1/p^\infty})^{\wedge\flat},-,q}})\ar[r]^{\mathrm{global}}\ar[r]\ar[r] &\mathrm{sComm}_\mathrm{simplicial}\varphi_q\mathrm{Ind}^m\mathrm{Banach}(\{\mathrm{Robba}^\mathrm{extended}_{{R_0,-,I,q}}\}_I).
}
\]

\item Then parallel as in \cite{LBV} we have a functor (global section) of the de Rham complex after \cite[Definition 5.9, Section 5.2.1]{KKM}:
\[\displayindent=-0.4in
\xymatrix@R+6pc@C+0pc{
\mathrm{deRham}_{\mathrm{sComm}_\mathrm{simplicial}\mathrm{Ind}\mathrm{Banach}(\mathcal{O}_{X_{R_k,-,q}})\ar[r]^{\mathrm{global}}}(-)\ar[d]\ar[d]\ar[d]\ar[d]\ar[r]\ar[r] &\mathrm{deRham}_{\mathrm{sComm}_\mathrm{simplicial}\varphi_q\mathrm{Ind}\mathrm{Banach}(\{\mathrm{Robba}^\mathrm{extended}_{{R_k,-,I,q}}\}_I)}(-)\ar[d]\ar[d]\ar[d]\ar[d]\\
\mathrm{deRham}_{\mathrm{sComm}_\mathrm{simplicial}\mathrm{Ind}\mathrm{Banach}(\mathcal{O}_{X_{\mathbb{Q}_p(p^{1/p^\infty})^{\wedge\flat},-,q}})\ar[r]^{\mathrm{global}}}(-)\ar[r]\ar[r] &\mathrm{deRham}_{\mathrm{sComm}_\mathrm{simplicial}\varphi_q\mathrm{Ind}\mathrm{Banach}(\{\mathrm{Robba}^\mathrm{extended}_{{R_0,-,I,q}}\}_I)}(-), 
}
\]
\[\displayindent=-0.4in
\xymatrix@R+6pc@C+0pc{
\mathrm{deRham}_{\mathrm{sComm}_\mathrm{simplicial}\mathrm{Ind}^m\mathrm{Banach}(\mathcal{O}_{X_{R_k,-,q}})\ar[r]^{\mathrm{global}}}(-)\ar[d]\ar[d]\ar[d]\ar[d]\ar[r]\ar[r] &\mathrm{deRham}_{\mathrm{sComm}_\mathrm{simplicial}\varphi_q\mathrm{Ind}^m\mathrm{Banach}(\{\mathrm{Robba}^\mathrm{extended}_{{R_k,-,I,q}}\}_I)}(-)\ar[d]\ar[d]\ar[d]\ar[d]\\
\mathrm{deRham}_{\mathrm{sComm}_\mathrm{simplicial}\mathrm{Ind}^m\mathrm{Banach}(\mathcal{O}_{X_{\mathbb{Q}_p(p^{1/p^\infty})^{\wedge\flat},-,q}})\ar[r]^{\mathrm{global}}}(-)\ar[r]\ar[r] &\mathrm{deRham}_{\mathrm{sComm}_\mathrm{simplicial}\varphi_q\mathrm{Ind}^m\mathrm{Banach}(\{\mathrm{Robba}^\mathrm{extended}_{{R_0,-,I,q}}\}_I)}(-). 
}
\]

\item Then we have the following a functor (global section) of $K$-group $(\infty,1)$-spectrum from \cite{BGT}:
\[
\xymatrix@R+6pc@C+0pc{
\mathrm{K}^\mathrm{BGT}_{\mathrm{sComm}_\mathrm{simplicial}\mathrm{Ind}\mathrm{Banach}(\mathcal{O}_{X_{R_k,-,q}})\ar[r]^{\mathrm{global}}}(-)\ar[d]\ar[d]\ar[d]\ar[d]\ar[r]\ar[r] &\mathrm{K}^\mathrm{BGT}_{\mathrm{sComm}_\mathrm{simplicial}\varphi_q\mathrm{Ind}\mathrm{Banach}(\{\mathrm{Robba}^\mathrm{extended}_{{R_k,-,I,q}}\}_I)}(-)\ar[d]\ar[d]\ar[d]\ar[d]\\
\mathrm{K}^\mathrm{BGT}_{\mathrm{sComm}_\mathrm{simplicial}\mathrm{Ind}\mathrm{Banach}(\mathcal{O}_{X_{\mathbb{Q}_p(p^{1/p^\infty})^{\wedge\flat},-,q}})\ar[r]^{\mathrm{global}}}(-)\ar[r]\ar[r] &\mathrm{K}^\mathrm{BGT}_{\mathrm{sComm}_\mathrm{simplicial}\varphi_q\mathrm{Ind}\mathrm{Banach}(\{\mathrm{Robba}^\mathrm{extended}_{{R_0,-,I,q}}\}_I)}(-), 
}
\]
\[
\xymatrix@R+6pc@C+0pc{
\mathrm{K}^\mathrm{BGT}_{\mathrm{sComm}_\mathrm{simplicial}\mathrm{Ind}^m\mathrm{Banach}(\mathcal{O}_{X_{R_k,-,q}})\ar[r]^{\mathrm{global}}}(-)\ar[d]\ar[d]\ar[d]\ar[d]\ar[r]\ar[r] &\mathrm{K}^\mathrm{BGT}_{\mathrm{sComm}_\mathrm{simplicial}\varphi_q\mathrm{Ind}^m\mathrm{Banach}(\{\mathrm{Robba}^\mathrm{extended}_{{R_k,-,I,q}}\}_I)}(-)\ar[d]\ar[d]\ar[d]\ar[d]\\
\mathrm{K}^\mathrm{BGT}_{\mathrm{sComm}_\mathrm{simplicial}\mathrm{Ind}^m\mathrm{Banach}(\mathcal{O}_{X_{\mathbb{Q}_p(p^{1/p^\infty})^{\wedge\flat},-,q}})\ar[r]^{\mathrm{global}}}(-)\ar[r]\ar[r] &\mathrm{K}^\mathrm{BGT}_{\mathrm{sComm}_\mathrm{simplicial}\varphi_q\mathrm{Ind}^m\mathrm{Banach}(\{\mathrm{Robba}^\mathrm{extended}_{{R_0,-,I,q}}\}_I)}(-). 
}
\]

\end{itemize}

\
\indent Then we consider further equivariance by considering the arithmetic profinite fundamental groups and actually its $q$-th power $\mathrm{Gal}(\overline{{Q}_p\left<T_1^{\pm 1},...,T_k^{\pm 1}\right>}/R_k)^{\times q}$ through the following diagram:\\

\[
\xymatrix@R+6pc@C+0pc{
\mathbb{Z}_p^k=\mathrm{Gal}(R_k/{\mathbb{Q}_p(p^{1/p^\infty})^\wedge\left<T_1^{\pm 1},...,T_k^{\pm 1}\right>})\ar[d]\ar[d]\ar[d]\ar[d] \ar[r]\ar[r] \ar[r]\ar[r] &\mathrm{Gal}(\overline{{Q}_p\left<T_1^{\pm 1},...,T_k^{\pm 1}\right>}/R_k) \ar[d]\ar[d]\ar[d] \ar[r]\ar[r] &\Gamma_{\mathbb{Q}_p} \ar[d]\ar[d]\ar[d]\ar[d]\\
(\mathbb{Z}_p^k=\mathrm{Gal}(R_k/{\mathbb{Q}_p(p^{1/p^\infty})^\wedge\left<T_1^{\pm 1},...,T_k^{\pm 1}\right>}))^{\times q} \ar[r]\ar[r] \ar[r]\ar[r] &\Gamma_k^{\times q}:=\mathrm{Gal}(R_k/{\mathbb{Q}_p\left<T_1^{\pm 1},...,T_k^{\pm 1}\right>})^{\times q} \ar[r] \ar[r]\ar[r] &\Gamma_{\mathbb{Q}_p}^{\times q}.
}
\]

\

We then have the correspond arithmetic profinite fundamental groups equivariant versions:
\begin{itemize}
\item (\text{Proposition}) There is a functor (global section) between the $\infty$-prestacks of inductive Banach quasicoherent presheaves:
\[
\xymatrix@R+6pc@C+0pc{
\mathrm{Ind}\mathrm{Banach}_{\Gamma_{k}^{\times q}}(\mathcal{O}_{X_{\mathbb{Q}_p(p^{1/p^\infty})^{\wedge}\left<T_1^{\pm 1/p^\infty},...,T_k^{\pm 1/p^\infty}\right>^\flat,-,q}})\ar[d]\ar[d]\ar[d]\ar[d] \ar[r]^{\mathrm{global}}\ar[r]\ar[r] &\varphi_q\mathrm{Ind}\mathrm{Banach}_{\Gamma_{k}^{\times q}}(\{\mathrm{Robba}^\mathrm{extended}_{{R_k,-,I,q}}\}_I) \ar[d]\ar[d]\ar[d]\ar[d].\\
\mathrm{Ind}\mathrm{Banach}(\mathcal{O}_{X_{\mathbb{Q}_p(p^{1/p^\infty})^{\wedge\flat},-,q}})\ar[r]^{\mathrm{global}}\ar[r]\ar[r] &\varphi_q\mathrm{Ind}\mathrm{Banach}(\{\mathrm{Robba}^\mathrm{extended}_{{R_0,-,I,q}}\}_I).\\ 
}
\]
\item (\text{Proposition}) There is a functor (global section) between the $\infty$-prestacks of monomorphic inductive Banach quasicoherent presheaves:
\[
\xymatrix@R+6pc@C+0pc{
\mathrm{Ind}^m\mathrm{Banach}_{\Gamma_{k}^{\times q}}(\mathcal{O}_{X_{R_k,-,q}})\ar[r]^{\mathrm{global}}\ar[d]\ar[d]\ar[d]\ar[d]\ar[r]\ar[r] &\varphi_q\mathrm{Ind}^m\mathrm{Banach}_{\Gamma_{k}^{\times q}}(\{\mathrm{Robba}^\mathrm{extended}_{{R_k,-,I,q}}\}_I)\ar[d]\ar[d]\ar[d]\ar[d]\\
\mathrm{Ind}^m\mathrm{Banach}_{\Gamma_{0}^{\times q}}(\mathcal{O}_{X_{\mathbb{Q}_p(p^{1/p^\infty})^{\wedge\flat},-,q}})\ar[r]^{\mathrm{global}}\ar[r]\ar[r] &\varphi_q\mathrm{Ind}^m\mathrm{Banach}_{\Gamma_{0}^{\times q}}(\{\mathrm{Robba}^\mathrm{extended}_{{R_0,-,I,q}}\}_I).\\  
}
\]
\item (\text{Proposition}) There is a functor (global section) between the $\infty$-stacks of inductive Banach quasicoherent commutative algebra $E_\infty$ objects:
\[\displayindent=-0.4in
\xymatrix@R+6pc@C+0pc{
\mathrm{sComm}_\mathrm{simplicial}\mathrm{Ind}\mathrm{Banach}_{\Gamma_{k}^{\times q}}(\mathcal{O}_{X_{R_k,-,q}})\ar[d]\ar[d]\ar[d]\ar[d]\ar[r]^{\mathrm{global}}\ar[r]\ar[r] &\mathrm{sComm}_\mathrm{simplicial}\varphi_q\mathrm{Ind}\mathrm{Banach}_{\Gamma_{k}^{\times q}}(\{\mathrm{Robba}^\mathrm{extended}_{{R_k,-,I,q}}\}_I)\ar[d]\ar[d]\ar[d]\ar[d]\\
\mathrm{sComm}_\mathrm{simplicial}\mathrm{Ind}\mathrm{Banach}_{\Gamma_{0}^{\times q}}(\mathcal{O}_{X_{\mathbb{Q}_p(p^{1/p^\infty})^{\wedge\flat},-,q}})\ar[r]^{\mathrm{global}}\ar[r]\ar[r] &\mathrm{sComm}_\mathrm{simplicial}\varphi_q\mathrm{Ind}\mathrm{Banach}_{\Gamma_{0}^{\times q}}(\{\mathrm{Robba}^\mathrm{extended}_{{R_0,-,I,q}}\}_I).  
}
\]
\item (\text{Proposition}) There is a functor (global section) between the $\infty$-prestacks of monomorphic inductive Banach quasicoherent commutative algebra $E_\infty$ objects:
\[\displayindent=-0.4in
\xymatrix@R+6pc@C+0pc{
\mathrm{sComm}_\mathrm{simplicial}\mathrm{Ind}^m\mathrm{Banach}_{\Gamma_{k}^{\times q}}(\mathcal{O}_{X_{R_k,-,q}})\ar[d]\ar[d]\ar[d]\ar[d]\ar[r]^{\mathrm{global}}\ar[r]\ar[r] &\mathrm{sComm}_\mathrm{simplicial}\varphi_q\mathrm{Ind}^m\mathrm{Banach}_{\Gamma_{k}^{\times q}}(\{\mathrm{Robba}^\mathrm{extended}_{{R_k,-,I,q}}\}_I)\ar[d]\ar[d]\ar[d]\ar[d]\\
 \mathrm{sComm}_\mathrm{simplicial}\mathrm{Ind}^m\mathrm{Banach}_{\Gamma_{0}^{\times q}}(\mathcal{O}_{X_{\mathbb{Q}_p(p^{1/p^\infty})^{\wedge\flat},-,q}})\ar[r]^{\mathrm{global}}\ar[r]\ar[r] &\mathrm{sComm}_\mathrm{simplicial}\varphi_q\mathrm{Ind}^m\mathrm{Banach}_{\Gamma_{0}^{\times q}}(\{\mathrm{Robba}^\mathrm{extended}_{{R_0,-,I,q}}\}_I). 
}
\]

\item Then parallel as in \cite{LBV} we have a functor (global section) of the de Rham complex after \cite[Definition 5.9, Section 5.2.1]{KKM}:
\[\displayindent=-0.4in
\xymatrix@R+6pc@C+0pc{
\mathrm{deRham}_{\mathrm{sComm}_\mathrm{simplicial}\mathrm{Ind}\mathrm{Banach}_{\Gamma_{k}^{\times q}}(\mathcal{O}_{X_{R_k,-,q}})\ar[r]^{\mathrm{global}}}(-)\ar[d]\ar[d]\ar[d]\ar[d]\ar[r]\ar[r] &\mathrm{deRham}_{\mathrm{sComm}_\mathrm{simplicial}\varphi_q\mathrm{Ind}\mathrm{Banach}_{\Gamma_{k}^{\times q}}(\{\mathrm{Robba}^\mathrm{extended}_{{R_k,-,I,q}}\}_I)}(-)\ar[d]\ar[d]\ar[d]\ar[d]\\
\mathrm{deRham}_{\mathrm{sComm}_\mathrm{simplicial}\mathrm{Ind}\mathrm{Banach}_{\Gamma_{0}^{\times q}}(\mathcal{O}_{X_{\mathbb{Q}_p(p^{1/p^\infty})^{\wedge\flat},-,q}})\ar[r]^{\mathrm{global}}}(-)\ar[r]\ar[r] &\mathrm{deRham}_{\mathrm{sComm}_\mathrm{simplicial}\varphi_q\mathrm{Ind}\mathrm{Banach}_{\Gamma_{0}^{\times q}}(\{\mathrm{Robba}^\mathrm{extended}_{{R_0,-,I,q}}\}_I)}(-), 
}
\]
\[\displayindent=-0.7in
\xymatrix@R+6pc@C+0pc{
\mathrm{deRham}_{\mathrm{sComm}_\mathrm{simplicial}\mathrm{Ind}^m\mathrm{Banach}_{\Gamma_{k}^{\times q}}(\mathcal{O}_{X_{R_k,-,q}})\ar[r]^{\mathrm{global}}}(-)\ar[d]\ar[d]\ar[d]\ar[d]\ar[r]\ar[r] &\mathrm{deRham}_{\mathrm{sComm}_\mathrm{simplicial}\varphi_q\mathrm{Ind}^m\mathrm{Banach}_{\Gamma_{k}^{\times q}}(\{\mathrm{Robba}^\mathrm{extended}_{{R_k,-,I,q}}\}_I)}(-)\ar[d]\ar[d]\ar[d]\ar[d]\\
\mathrm{deRham}_{\mathrm{sComm}_\mathrm{simplicial}\mathrm{Ind}^m\mathrm{Banach}_{\Gamma_{0}^{\times q}}(\mathcal{O}_{X_{\mathbb{Q}_p(p^{1/p^\infty})^{\wedge\flat},-,q}})\ar[r]^{\mathrm{global}}}(-)\ar[r]\ar[r] &\mathrm{deRham}_{\mathrm{sComm}_\mathrm{simplicial}\varphi_q\mathrm{Ind}^m\mathrm{Banach}_{\Gamma_{0}^{\times q}}(\{\mathrm{Robba}^\mathrm{extended}_{{R_0,-,I,q}}\}_I)}(-). 
}
\]

\item Then we have the following a functor (global section) of $K$-group $(\infty,1)$-spectrum from \cite{BGT}:
\[
\xymatrix@R+6pc@C+0pc{
\mathrm{K}^\mathrm{BGT}_{\mathrm{sComm}_\mathrm{simplicial}\mathrm{Ind}\mathrm{Banach}_{\Gamma_{k}^{\times q}}(\mathcal{O}_{X_{R_k,-,q}})\ar[r]^{\mathrm{global}}}(-)\ar[d]\ar[d]\ar[d]\ar[d]\ar[r]\ar[r] &\mathrm{K}^\mathrm{BGT}_{\mathrm{sComm}_\mathrm{simplicial}\varphi_q\mathrm{Ind}\mathrm{Banach}_{\Gamma_{k}^{\times q}}(\{\mathrm{Robba}^\mathrm{extended}_{{R_k,-,I,q}}\}_I)}(-)\ar[d]\ar[d]\ar[d]\ar[d]\\
\mathrm{K}^\mathrm{BGT}_{\mathrm{sComm}_\mathrm{simplicial}\mathrm{Ind}\mathrm{Banach}_{\Gamma_{0}^{\times q}}(\mathcal{O}_{X_{\mathbb{Q}_p(p^{1/p^\infty})^{\wedge\flat},-,q}})\ar[r]^{\mathrm{global}}}(-)\ar[r]\ar[r] &\mathrm{K}^\mathrm{BGT}_{\mathrm{sComm}_\mathrm{simplicial}\varphi_q\mathrm{Ind}\mathrm{Banach}_{\Gamma_{0}^{\times q}}(\{\mathrm{Robba}^\mathrm{extended}_{{R_0,-,I,q}}\}_I)}(-), 
}
\]
\[
\xymatrix@R+6pc@C+0pc{
\mathrm{K}^\mathrm{BGT}_{\mathrm{sComm}_\mathrm{simplicial}\mathrm{Ind}^m\mathrm{Banach}_{\Gamma_{k}^{\times q}}(\mathcal{O}_{X_{R_k,-,q}})\ar[r]^{\mathrm{global}}}(-)\ar[d]\ar[d]\ar[d]\ar[d]\ar[r]\ar[r] &\mathrm{K}^\mathrm{BGT}_{\mathrm{sComm}_\mathrm{simplicial}\varphi_q\mathrm{Ind}^m\mathrm{Banach}_{\Gamma_{k}^{\times q}}(\{\mathrm{Robba}^\mathrm{extended}_{{R_k,-,I,q}}\}_I)}(-)\ar[d]\ar[d]\ar[d]\ar[d]\\
\mathrm{K}^\mathrm{BGT}_{\mathrm{sComm}_\mathrm{simplicial}\mathrm{Ind}^m\mathrm{Banach}_{\Gamma_{0}^{\times q}}(\mathcal{O}_{X_{\mathbb{Q}_p(p^{1/p^\infty})^{\wedge\flat},-,q}})\ar[r]^{\mathrm{global}}}(-)\ar[r]\ar[r] &\mathrm{K}^\mathrm{BGT}_{\mathrm{sComm}_\mathrm{simplicial}\varphi_q\mathrm{Ind}^m\mathrm{Banach}_{\Gamma_{0}^{\times q}}(\{\mathrm{Robba}^\mathrm{extended}_{{R_0,-,I,q}}\}_I)}(-). 
}
\]

\end{itemize}

\

\begin{remark}
\noindent We can certainly consider the quasicoherent sheaves in \cite[Lemma 7.11, Remark 7.12]{1BBK}, therefore all the quasicoherent presheaves and modules will be those satisfying \cite[Lemma 7.11, Remark 7.12]{1BBK} if one would like to consider the the quasicoherent sheaves. That being all as this said, we would believe that the big quasicoherent presheaves are automatically quasicoherent sheaves (namely satisfying the corresponding \v{C}ech $\infty$-descent as in \cite[Section 9.3]{KKM} and \cite[Lemma 7.11, Remark 7.12]{1BBK}) and the corresponding global section functors are automatically equivalence of $\infty$-categories.\\ 
\end{remark}

\

\indent In Clausen-Scholze formalism we have the following:

\begin{itemize}

\item (\text{Proposition}) There is a functor (global section) between the $\infty$-prestacks of inductive Banach quasicoherent sheaves:
\[
\xymatrix@R+0pc@C+0pc{
{\mathrm{Modules}_\circledcirc}(\mathcal{O}_{X_{R,-,q}})\ar[r]^{\mathrm{global}}\ar[r]\ar[r] &\varphi_q{\mathrm{Modules}_\circledcirc}(\{\mathrm{Robba}^{\mathrm{extended},q}_{{R,-,I}}\}_I).  
}
\]

\item (\text{Proposition}) There is a functor (global section) between the $\infty$-prestacks of inductive Banach quasicoherent sheaves:
\[
\xymatrix@R+0pc@C+0pc{
{\mathrm{Modules}_\circledcirc}(\mathcal{O}_{X_{R,-,q}})\ar[r]^{\mathrm{global}}\ar[r]\ar[r] &\varphi_q{\mathrm{Modules}_\circledcirc}(\{\mathrm{Robba}^{\mathrm{extended},q}_{{R,-,I}}\}_I).  
}
\]

\item (\text{Proposition}) There is a functor (global section) between the $\infty$-prestacks of inductive Banach quasicoherent commutative algebra $E_\infty$ objects\footnote{Here $\circledcirc=\text{solidquasicoherentsheaves}$.}:
\[
\xymatrix@R+0pc@C+0pc{
\mathrm{sComm}_\mathrm{simplicial}{\mathrm{Modules}_\circledcirc}(\mathcal{O}_{X_{R,-,q}})\ar[r]^{\mathrm{global}}\ar[r]\ar[r] &\mathrm{sComm}_\mathrm{simplicial}\varphi_q{\mathrm{Modules}_\circledcirc}(\{\mathrm{Robba}^{\mathrm{extended},q}_{{R,-,I}}\}_I).  
}
\]

\item Then as in \cite{LBV} we have a functor (global section) of the de Rham complex after \cite[Definition 5.9, Section 5.2.1]{KKM}\footnote{Here $\circledcirc=\text{solidquasicoherentsheaves}$.}:
\[
\xymatrix@R+0pc@C+0pc{
\mathrm{deRham}_{\mathrm{sComm}_\mathrm{simplicial}{\mathrm{Modules}_\circledcirc}(\mathcal{O}_{X_{R,-,q}})\ar[r]^{\mathrm{global}}}(-)\ar[r]\ar[r] &\mathrm{deRham}_{\mathrm{sComm}_\mathrm{simplicial}\varphi_q{\mathrm{Modules}_\circledcirc}(\{\mathrm{Robba}^{\mathrm{extended},q}_{{R,-,I}}\}_I)}(-). 
}
\]

\item Then we have the following a functor (global section) of $K$-group $(\infty,1)$-spectrum from \cite{BGT}\footnote{Here $\circledcirc=\text{solidquasicoherentsheaves}$.}:
\[
\xymatrix@R+0pc@C+0pc{
\mathrm{K}^\mathrm{BGT}_{\mathrm{sComm}_\mathrm{simplicial}{\mathrm{Modules}_\circledcirc}(\mathcal{O}_{X_{R,-,q}})\ar[r]^{\mathrm{global}}}(-)\ar[r]\ar[r] &\mathrm{K}^\mathrm{BGT}_{\mathrm{sComm}_\mathrm{simplicial}\varphi_q{\mathrm{Modules}_\circledcirc}(\{\mathrm{Robba}^{\mathrm{extended},q}_{{R,-,I}}\}_I)}(-). 
}
\]
\end{itemize}

\noindent Now let $R=\mathbb{Q}_p(p^{1/p^\infty})^{\wedge\flat}$ and $R_k=\mathbb{Q}_p(p^{1/p^\infty})^{\wedge}\left<T_1^{\pm 1/p^{\infty}},...,T_k^{\pm 1/p^{\infty}}\right>^\flat$ we have the following Galois theoretic results with naturality along $f:\mathrm{Spa}(\mathbb{Q}_p(p^{1/p^\infty})^{\wedge}\left<T_1^{\pm 1/p^\infty},...,T_k^{\pm 1/p^\infty}\right>^\flat)\rightarrow \mathrm{Spa}(\mathbb{Q}_p(p^{1/p^\infty})^{\wedge\flat})$:

\begin{itemize}
\item (\text{Proposition}) There is a functor (global section) between the $\infty$-prestacks of inductive Banach quasicoherent sheaves\footnote{Here $\circledcirc=\text{solidquasicoherentsheaves}$.}:
\[
\xymatrix@R+6pc@C+0pc{
{\mathrm{Modules}_\circledcirc}(\mathcal{O}_{X_{\mathbb{Q}_p(p^{1/p^\infty})^{\wedge}\left<T_1^{\pm 1/p^\infty},...,T_k^{\pm 1/p^\infty}\right>^\flat,-,q}})\ar[d]\ar[d]\ar[d]\ar[d] \ar[r]^{\mathrm{global}}\ar[r]\ar[r] &\varphi_q{\mathrm{Modules}_\circledcirc}(\{\mathrm{Robba}^{\mathrm{extended},q}_{{R_k,-,I}}\}_I) \ar[d]\ar[d]\ar[d]\ar[d].\\
{\mathrm{Modules}_\circledcirc}(\mathcal{O}_{X_{\mathbb{Q}_p(p^{1/p^\infty})^{\wedge\flat},-,q}})\ar[r]^{\mathrm{global}}\ar[r]\ar[r] &\varphi_q{\mathrm{Modules}_\circledcirc}(\{\mathrm{Robba}^{\mathrm{extended},q}_{{R_0,-,I}}\}_I).\\ 
}
\]
\item (\text{Proposition}) There is a functor (global section) between the $\infty$-prestacks of inductive Banach quasicoherent commutative algebra $E_\infty$ objects\footnote{Here $\circledcirc=\text{solidquasicoherentsheaves}$.}:
\[
\xymatrix@R+6pc@C+0pc{
\mathrm{sComm}_\mathrm{simplicial}{\mathrm{Modules}_\circledcirc}(\mathcal{O}_{X_{R_k,-,q}})\ar[d]\ar[d]\ar[d]\ar[d]\ar[r]^{\mathrm{global}}\ar[r]\ar[r] &\mathrm{sComm}_\mathrm{simplicial}\varphi_q{\mathrm{Modules}_\circledcirc}(\{\mathrm{Robba}^{\mathrm{extended},q}_{{R_k,-,I}}\}_I)\ar[d]\ar[d]\ar[d]\ar[d]\\
\mathrm{sComm}_\mathrm{simplicial}{\mathrm{Modules}_\circledcirc}(\mathcal{O}_{X_{\mathbb{Q}_p(p^{1/p^\infty})^{\wedge\flat},-,q}})\ar[r]^{\mathrm{global}}\ar[r]\ar[r] &\mathrm{sComm}_\mathrm{simplicial}\varphi_q{\mathrm{Modules}_\circledcirc}(\{\mathrm{Robba}^{\mathrm{extended},q}_{{R_0,-,I}}\}_I).  
}
\]

\item Then as in \cite{LBV} we have a functor (global section) of the de Rham complex after \cite[Definition 5.9, Section 5.2.1]{KKM}\footnote{Here $\circledcirc=\text{solidquasicoherentsheaves}$.}:
\[\displayindent=-0.4in
\xymatrix@R+6pc@C+0pc{
\mathrm{deRham}_{\mathrm{sComm}_\mathrm{simplicial}{\mathrm{Modules}_\circledcirc}(\mathcal{O}_{X_{R_k,-,q}})\ar[r]^{\mathrm{global}}}(-)\ar[d]\ar[d]\ar[d]\ar[d]\ar[r]\ar[r] &\mathrm{deRham}_{\mathrm{sComm}_\mathrm{simplicial}\varphi_q{\mathrm{Modules}_\circledcirc}(\{\mathrm{Robba}^{\mathrm{extended},q}_{{R_k,-,I}}\}_I)}(-)\ar[d]\ar[d]\ar[d]\ar[d]\\
\mathrm{deRham}_{\mathrm{sComm}_\mathrm{simplicial}{\mathrm{Modules}_\circledcirc}(\mathcal{O}_{X_{\mathbb{Q}_p(p^{1/p^\infty})^{\wedge\flat},-,q}})\ar[r]^{\mathrm{global}}}(-)\ar[r]\ar[r] &\mathrm{deRham}_{\mathrm{sComm}_\mathrm{simplicial}\varphi_q{\mathrm{Modules}_\circledcirc}(\{\mathrm{Robba}^{\mathrm{extended},q}_{{R_0,-,I}}\}_I)}(-). 
}
\]

\item Then we have the following a functor (global section) of $K$-group $(\infty,1)$-spectrum from \cite{BGT}\footnote{Here $\circledcirc=\text{solidquasicoherentsheaves}$.}:
\[
\xymatrix@R+6pc@C+0pc{
\mathrm{K}^\mathrm{BGT}_{\mathrm{sComm}_\mathrm{simplicial}{\mathrm{Modules}_\circledcirc}(\mathcal{O}_{X_{R_k,-,q}})\ar[r]^{\mathrm{global}}}(-)\ar[d]\ar[d]\ar[d]\ar[d]\ar[r]\ar[r] &\mathrm{K}^\mathrm{BGT}_{\mathrm{sComm}_\mathrm{simplicial}\varphi_q{\mathrm{Modules}_\circledcirc}(\{\mathrm{Robba}^{\mathrm{extended},q}_{{R_k,-,I}}\}_I)}(-)\ar[d]\ar[d]\ar[d]\ar[d]\\
\mathrm{K}^\mathrm{BGT}_{\mathrm{sComm}_\mathrm{simplicial}{\mathrm{Modules}_\circledcirc}(\mathcal{O}_{X_{\mathbb{Q}_p(p^{1/p^\infty})^{\wedge\flat},-,q}})\ar[r]^{\mathrm{global}}}(-)\ar[r]\ar[r] &\mathrm{K}^\mathrm{BGT}_{\mathrm{sComm}_\mathrm{simplicial}\varphi_q{\mathrm{Modules}_\circledcirc}(\{\mathrm{Robba}^{\mathrm{extended},q}_{{R_0,-,I}}\}_I)}(-). 
}
\]

\end{itemize}

\
\indent Then we consider further equivariance by considering the arithmetic profinite fundamental groups $\Gamma_{\mathbb{Q}_p}$ and $\mathrm{Gal}(\overline{{Q}_p\left<T_1^{\pm 1},...,T_k^{\pm 1}\right>}/R_k)$ through the following diagram:

\[
\xymatrix@R+0pc@C+0pc{
\mathbb{Z}_p^k=\mathrm{Gal}(R_k/{\mathbb{Q}_p(p^{1/p^\infty})^\wedge\left<T_1^{\pm 1},...,T_k^{\pm 1}\right>}) \ar[r]\ar[r] \ar[r]\ar[r] &\Gamma_k:=\mathrm{Gal}(R_k/{\mathbb{Q}_p\left<T_1^{\pm 1},...,T_k^{\pm 1}\right>}) \ar[r] \ar[r]\ar[r] &\Gamma_{\mathbb{Q}_p}.
}
\]

\begin{itemize}
\item (\text{Proposition}) There is a functor (global section) between the $\infty$-prestacks of inductive Banach quasicoherent sheaves\footnote{Here $\circledcirc=\text{solidquasicoherentsheaves}$.}:
\[
\xymatrix@R+6pc@C+0pc{
{\mathrm{Modules}_\circledcirc}_{\Gamma_k^{\times q}}(\mathcal{O}_{X_{\mathbb{Q}_p(p^{1/p^\infty})^{\wedge}\left<T_1^{\pm 1/p^\infty},...,T_k^{\pm 1/p^\infty}\right>^\flat,-,q}})\ar[d]\ar[d]\ar[d]\ar[d] \ar[r]^{\mathrm{global}}\ar[r]\ar[r] &\varphi_q{\mathrm{Modules}_\circledcirc}_{\Gamma_k^{\times q}}(\{\mathrm{Robba}^{\mathrm{extended},q}_{{R_k,-,I}}\}_I) \ar[d]\ar[d]\ar[d]\ar[d].\\
{\mathrm{Modules}_\circledcirc}(\mathcal{O}_{X_{\mathbb{Q}_p(p^{1/p^\infty})^{\wedge\flat},-,q}})\ar[r]^{\mathrm{global}}\ar[r]\ar[r] &\varphi_q{\mathrm{Modules}_\circledcirc}(\{\mathrm{Robba}^{\mathrm{extended},q}_{{R_0,-,I}}\}_I).\\ 
}
\]

\item (\text{Proposition}) There is a functor (global section) between the $\infty$-stacks of inductive Banach quasicoherent commutative algebra $E_\infty$ objects\footnote{Here $\circledcirc=\text{solidquasicoherentsheaves}$.}:
\[\displayindent=-0.4in
\xymatrix@R+6pc@C+0pc{
\mathrm{sComm}_\mathrm{simplicial}{\mathrm{Modules}_\circledcirc}_{\Gamma_k^{\times q}}(\mathcal{O}_{X_{R_k,-,q}})\ar[d]\ar[d]\ar[d]\ar[d]\ar[r]^{\mathrm{global}}\ar[r]\ar[r] &\mathrm{sComm}_\mathrm{simplicial}\varphi_q{\mathrm{Modules}_\circledcirc}_{\Gamma_k^{\times q}}(\{\mathrm{Robba}^{\mathrm{extended},q}_{{R_k,-,I}}\}_I)\ar[d]\ar[d]\ar[d]\ar[d]\\
\mathrm{sComm}_\mathrm{simplicial}{\mathrm{Modules}_\circledcirc}_{\Gamma_0^{\times q}}(\mathcal{O}_{X_{\mathbb{Q}_p(p^{1/p^\infty})^{\wedge\flat},-,q}})\ar[r]^{\mathrm{global}}\ar[r]\ar[r] &\mathrm{sComm}_\mathrm{simplicial}\varphi_q{\mathrm{Modules}_\circledcirc}_{\Gamma_0^{\times q}}(\{\mathrm{Robba}^{\mathrm{extended},q}_{{R_0,-,I}}\}_I).  
}
\]

\item Then as in \cite{LBV} we have a functor (global section) of the de Rham complex after \cite[Definition 5.9, Section 5.2.1]{KKM}\footnote{Here $\circledcirc=\text{solidquasicoherentsheaves}$.}:
\[\displayindent=-0.7in
\xymatrix@R+6pc@C+0pc{
\mathrm{deRham}_{\mathrm{sComm}_\mathrm{simplicial}{\mathrm{Modules}_\circledcirc}_{\Gamma_k^{\times q}}(\mathcal{O}_{X_{R_k,-,q}})\ar[r]^{\mathrm{global}}}(-)\ar[d]\ar[d]\ar[d]\ar[d]\ar[r]\ar[r] &\mathrm{deRham}_{\mathrm{sComm}_\mathrm{simplicial}\varphi_q{\mathrm{Modules}_\circledcirc}_{\Gamma_k^{\times q}}(\{\mathrm{Robba}^{\mathrm{extended},q}_{{R_k,-,I}}\}_I)}(-)\ar[d]\ar[d]\ar[d]\ar[d]\\
\mathrm{deRham}_{\mathrm{sComm}_\mathrm{simplicial}{\mathrm{Modules}_\circledcirc}_{\Gamma_0^{\times q}}(\mathcal{O}_{X_{\mathbb{Q}_p(p^{1/p^\infty})^{\wedge\flat},-,q}})\ar[r]^{\mathrm{global}}}(-)\ar[r]\ar[r] &\mathrm{deRham}_{\mathrm{sComm}_\mathrm{simplicial}\varphi_q{\mathrm{Modules}_\circledcirc}_{\Gamma_0^{\times q}}(\{\mathrm{Robba}^{\mathrm{extended},q}_{{R_0,-,I}}\}_I)}(-). 
}
\]

\item Then we have the following a functor (global section) of $K$-group $(\infty,1)$-spectrum from \cite{BGT}\footnote{Here $\circledcirc=\text{solidquasicoherentsheaves}$.}:
\[
\xymatrix@R+6pc@C+0pc{
\mathrm{K}^\mathrm{BGT}_{\mathrm{sComm}_\mathrm{simplicial}{\mathrm{Modules}_\circledcirc}_{\Gamma_k^{\times q}}(\mathcal{O}_{X_{R_k,-,q}})\ar[r]^{\mathrm{global}}}(-)\ar[d]\ar[d]\ar[d]\ar[d]\ar[r]\ar[r] &\mathrm{K}^\mathrm{BGT}_{\mathrm{sComm}_\mathrm{simplicial}\varphi_q{\mathrm{Modules}_\circledcirc}_{\Gamma_k^{\times q}}(\{\mathrm{Robba}^{\mathrm{extended},q}_{{R_k,-,I}}\}_I)}(-)\ar[d]\ar[d]\ar[d]\ar[d]\\
\mathrm{K}^\mathrm{BGT}_{\mathrm{sComm}_\mathrm{simplicial}{\mathrm{Modules}_\circledcirc}_{\Gamma_0^{\times q}}(\mathcal{O}_{X_{\mathbb{Q}_p(p^{1/p^\infty})^{\wedge\flat},-,q}})\ar[r]^{\mathrm{global}}}(-)\ar[r]\ar[r] &\mathrm{K}^\mathrm{BGT}_{\mathrm{sComm}_\mathrm{simplicial}\varphi_q{\mathrm{Modules}_\circledcirc}_{\Gamma_0^{\times q}}(\{\mathrm{Robba}^{\mathrm{extended},q}_{{R_0,-,I}}\}_I)}(-). 
}
\]

\end{itemize}

\begin{proposition}
All the global functors from \cite[Proposition 13.8, Theorem 14.9, Remark 14.10]{1CS2} achieve the equivalences on both sides.	
\end{proposition}

\newpage
\subsection{$\infty$-Categorical Analytic Stacks and Descents IV}

\indent In the following the right had of each row in each diagram will be the corresponding quasicoherent Robba bundles over the Robba ring carrying the corresponding action from the Frobenius or the fundamental groups, defined by directly applying \cite[Section 9.3]{KKM} and \cite{BBM}. We now let $\mathcal{A}$ be any commutative algebra objects in the corresponding $\infty$-toposes over ind-Banach commutative algebra objects over $\mathbb{Q}_p$ or the corresponding born\'e commutative algebra objects over $\mathbb{Q}_p$ carrying the Grothendieck topology defined by essentially the corresponding monomorphism homotopy in the opposite category. Then we promote the construction to the corresponding $\infty$-stack over the same $\infty$-categories of affinoids. We now take the corresponding colimit through all the $(\infty,1)$-categories. Therefore all the corresponding $(\infty,1)$-functors into $(\infty,1)$-categories or $(\infty,1)$-groupoids are from the homotopy closure of $\mathbb{Q}_p\left<C_1,...,C_\ell\right>$ $\ell=1,q,...$ in $\mathrm{sComm}\mathrm{Ind}\mathrm{Banach}_{\mathbb{Q}_p}$ or $\mathbb{Q}_p\left<C_1,...,C_\ell\right>$ $\ell=1,q,...$ in $\mathrm{sComm}\mathrm{Ind}^m\mathrm{Banach}_{\mathbb{Q}_p}$ as in \cite[Section 4.2]{BBM}:
\begin{align}
\mathrm{Ind}^{\mathbb{Q}_p\left<C_1,...,C_\ell\right>,\ell=1,q,...}\mathrm{sComm}\mathrm{Ind}\mathrm{Banach}_{\mathbb{Q}_p},\\
\mathrm{Ind}^{\mathbb{Q}_p\left<C_1,...,C_\ell\right>,\ell=1,q,...}\mathrm{sComm}\mathrm{Ind}\mathrm{Banach}_{\mathbb{Q}_p}	.
\end{align}

\begin{itemize}

\item (\text{Proposition}) There is a functor (global section) between the $\infty$-prestacks of inductive Banach quasicoherent presheaves:
\[
\xymatrix@R+0pc@C+0pc{
\mathrm{Ind}\mathrm{Banach}(\mathcal{O}_{X_{R,-,q}})\ar[r]^{\mathrm{global}}\ar[r]\ar[r] &\varphi_q\mathrm{Ind}\mathrm{Banach}(\{\mathrm{Robba}^\mathrm{extended}_{{R,-,I,q}}\}_I).  
}
\]
The definition is given by the following:
\[
\xymatrix@R+0pc@C+0pc{
\mathrm{homotopycolimit}_i(\mathrm{Ind}\mathrm{Banach}(\mathcal{O}_{X_{R,-,q}})\ar[r]^{\mathrm{global}}\ar[r]\ar[r] &\varphi_q\mathrm{Ind}\mathrm{Banach}(\{\mathrm{Robba}^\mathrm{extended}_{{R,-,I,q}}\}_I))(\mathcal{O}_i),  
}
\]
each $\mathcal{O}_i$ is just as $\mathbb{Q}_p\left<C_1,...,C_\ell\right>,\ell=1,q,...$.
\item (\text{Proposition}) There is a functor (global section) between the $\infty$-prestacks of monomorphic inductive Banach quasicoherent presheaves:
\[
\xymatrix@R+0pc@C+0pc{
\mathrm{Ind}^m\mathrm{Banach}(\mathcal{O}_{X_{R,-,q}})\ar[r]^{\mathrm{global}}\ar[r]\ar[r] &\varphi_q\mathrm{Ind}^m\mathrm{Banach}(\{\mathrm{Robba}^\mathrm{extended}_{{R,-,I,q}}\}_I).  
}
\]
The definition is given by the following:
\[
\xymatrix@R+0pc@C+0pc{
\mathrm{homotopycolimit}_i(\mathrm{Ind}^m\mathrm{Banach}(\mathcal{O}_{X_{R,-,q}})\ar[r]^{\mathrm{global}}\ar[r]\ar[r] &\varphi_q\mathrm{Ind}^m\mathrm{Banach}(\{\mathrm{Robba}^\mathrm{extended}_{{R,-,I,q}}\}_I))(\mathcal{O}_i),  
}
\]
each $\mathcal{O}_i$ is just as $\mathbb{Q}_p\left<C_1,...,C_\ell\right>,\ell=1,q,...$.

\item (\text{Proposition}) There is a functor (global section) between the $\infty$-prestacks of inductive Banach quasicoherent presheaves:
\[
\xymatrix@R+0pc@C+0pc{
\mathrm{Ind}\mathrm{Banach}(\mathcal{O}_{X_{R,-,q}})\ar[r]^{\mathrm{global}}\ar[r]\ar[r] &\varphi_q\mathrm{Ind}\mathrm{Banach}(\{\mathrm{Robba}^\mathrm{extended}_{{R,-,I,q}}\}_I).  
}
\]
The definition is given by the following:
\[
\xymatrix@R+0pc@C+0pc{
\mathrm{homotopycolimit}_i(\mathrm{Ind}\mathrm{Banach}(\mathcal{O}_{X_{R,-,q}})\ar[r]^{\mathrm{global}}\ar[r]\ar[r] &\varphi_q\mathrm{Ind}\mathrm{Banach}(\{\mathrm{Robba}^\mathrm{extended}_{{R,-,I,q}}\}_I))(\mathcal{O}_i),  
}
\]
each $\mathcal{O}_i$ is just as $\mathbb{Q}_p\left<C_1,...,C_\ell\right>,\ell=1,q,...$.
\item (\text{Proposition}) There is a functor (global section) between the $\infty$-stacks of monomorphic inductive Banach quasicoherent presheaves:
\[
\xymatrix@R+0pc@C+0pc{
\mathrm{Ind}^m\mathrm{Banach}(\mathcal{O}_{X_{R,-,q}})\ar[r]^{\mathrm{global}}\ar[r]\ar[r] &\varphi_q\mathrm{Ind}^m\mathrm{Banach}(\{\mathrm{Robba}^\mathrm{extended}_{{R,-,I,q}}\}_I).  
}
\]
The definition is given by the following:
\[
\xymatrix@R+0pc@C+0pc{
\mathrm{homotopycolimit}_i(\mathrm{Ind}^m\mathrm{Banach}(\mathcal{O}_{X_{R,-,q}})\ar[r]^{\mathrm{global}}\ar[r]\ar[r] &\varphi_q\mathrm{Ind}^m\mathrm{Banach}(\{\mathrm{Robba}^\mathrm{extended}_{{R,-,I,q}}\}_I))(\mathcal{O}_i),  
}
\]
each $\mathcal{O}_i$ is just as $\mathbb{Q}_p\left<C_1,...,C_\ell\right>,\ell=1,q,...$.
\item (\text{Proposition}) There is a functor (global section) between the $\infty$-prestacks of inductive Banach quasicoherent commutative algebra $E_\infty$ objects:
\[
\xymatrix@R+0pc@C+0pc{
\mathrm{sComm}_\mathrm{simplicial}\mathrm{Ind}\mathrm{Banach}(\mathcal{O}_{X_{R,-,q}})\ar[r]^{\mathrm{global}}\ar[r]\ar[r] &\mathrm{sComm}_\mathrm{simplicial}\varphi_q\mathrm{Ind}\mathrm{Banach}(\{\mathrm{Robba}^\mathrm{extended}_{{R,-,I,q}}\}_I).  
}
\]
The definition is given by the following:
\[
\xymatrix@R+0pc@C+0pc{
\mathrm{homotopycolimit}_i(\mathrm{sComm}_\mathrm{simplicial}\mathrm{Ind}\mathrm{Banach}(\mathcal{O}_{X_{R,-,q}})\ar[r]^{\mathrm{global}}\ar[r]\ar[r] &\mathrm{sComm}_\mathrm{simplicial}\varphi_q\mathrm{Ind}\mathrm{Banach}(\{\mathrm{Robba}^\mathrm{extended}_{{R,-,I,q}}\}_I))(\mathcal{O}_i),  
}
\]
each $\mathcal{O}_i$ is just as $\mathbb{Q}_p\left<C_1,...,C_\ell\right>,\ell=1,q,...$.
\item (\text{Proposition}) There is a functor (global section) between the $\infty$-prestacks of monomorphic inductive Banach quasicoherent commutative algebra $E_\infty$ objects:
\[
\xymatrix@R+0pc@C+0pc{
\mathrm{sComm}_\mathrm{simplicial}\mathrm{Ind}^m\mathrm{Banach}(\mathcal{O}_{X_{R,-,q}})\ar[r]^{\mathrm{global}}\ar[r]\ar[r] &\mathrm{sComm}_\mathrm{simplicial}\varphi_q\mathrm{Ind}^m\mathrm{Banach}(\{\mathrm{Robba}^\mathrm{extended}_{{R,-,I,q}}\}_I).  
}
\]
The definition is given by the following:
\[
\xymatrix@R+0pc@C+0pc{
\mathrm{homotopycolimit}_i(\mathrm{sComm}_\mathrm{simplicial}\mathrm{Ind}^m\mathrm{Banach}(\mathcal{O}_{X_{R,-,q}})\ar[r]^{\mathrm{global}}\ar[r]\ar[r] &\mathrm{sComm}_\mathrm{simplicial}\varphi_q\mathrm{Ind}^m\mathrm{Banach}(\{\mathrm{Robba}^\mathrm{extended}_{{R,-,I,q}}\}_I))(\mathcal{O}_i),  
}
\]
each $\mathcal{O}_i$ is just as $\mathbb{Q}_p\left<C_1,...,C_\ell\right>,\ell=1,q,...$.

\item Then parallel as in \cite{LBV} we have a functor (global section ) of the de Rham complex after \cite[Definition 5.9, Section 5.2.1]{KKM}:
\[
\xymatrix@R+0pc@C+0pc{
\mathrm{deRham}_{\mathrm{sComm}_\mathrm{simplicial}\mathrm{Ind}\mathrm{Banach}(\mathcal{O}_{X_{R,-,q}})\ar[r]^{\mathrm{global}}}(-)\ar[r]\ar[r] &\mathrm{deRham}_{\mathrm{sComm}_\mathrm{simplicial}\varphi_q\mathrm{Ind}\mathrm{Banach}(\{\mathrm{Robba}^\mathrm{extended}_{{R,-,I,q}}\}_I)}(-), 
}
\]
\[
\xymatrix@R+0pc@C+0pc{
\mathrm{deRham}_{\mathrm{sComm}_\mathrm{simplicial}\mathrm{Ind}^m\mathrm{Banach}(\mathcal{O}_{X_{R,-,q}})\ar[r]^{\mathrm{global}}}(-)\ar[r]\ar[r] &\mathrm{deRham}_{\mathrm{sComm}_\mathrm{simplicial}\varphi_q\mathrm{Ind}^m\mathrm{Banach}(\{\mathrm{Robba}^\mathrm{extended}_{{R,-,I,q}}\}_I)}(-). 
}
\]
The definition is given by the following:
\[\displayindent=-0.1in
\xymatrix@R+0pc@C+0pc{
\mathrm{homotopycolimit}_i\\
(\mathrm{deRham}_{\mathrm{sComm}_\mathrm{simplicial}\mathrm{Ind}\mathrm{Banach}(\mathcal{O}_{X_{R,-,q}})\ar[r]^{\mathrm{global}}}(-)\ar[r]\ar[r] &\mathrm{deRham}_{\mathrm{sComm}_\mathrm{simplicial}\varphi_q\mathrm{Ind}\mathrm{Banach}(\{\mathrm{Robba}^\mathrm{extended}_{{R,-,I,q}}\}_I)}(-))(\mathcal{O}_i),  
}
\]
\[\displayindent=-0.2in
\xymatrix@R+0pc@C+0pc{
\mathrm{homotopycolimit}_i\\
(\mathrm{deRham}_{\mathrm{sComm}_\mathrm{simplicial}\mathrm{Ind}^m\mathrm{Banach}(\mathcal{O}_{X_{R,-,q}})\ar[r]^{\mathrm{global}}}(-)\ar[r]\ar[r] &\mathrm{deRham}_{\mathrm{sComm}_\mathrm{simplicial}\varphi_q\mathrm{Ind}^m\mathrm{Banach}(\{\mathrm{Robba}^\mathrm{extended}_{{R,-,I,q}}\}_I)}(-))(\mathcal{O}_i),  
}
\]
each $\mathcal{O}_i$ is just as $\mathbb{Q}_p\left<C_1,...,C_\ell\right>,\ell=1,q,...$.\item Then we have the following a functor (global section) of $K$-group $(\infty,1)$-spectrum from \cite{BGT}:
\[
\xymatrix@R+0pc@C+0pc{
\mathrm{K}^\mathrm{BGT}_{\mathrm{sComm}_\mathrm{simplicial}\mathrm{Ind}\mathrm{Banach}(\mathcal{O}_{X_{R,-,q}})\ar[r]^{\mathrm{global}}}(-)\ar[r]\ar[r] &\mathrm{K}^\mathrm{BGT}_{\mathrm{sComm}_\mathrm{simplicial}\varphi_q\mathrm{Ind}\mathrm{Banach}(\{\mathrm{Robba}^\mathrm{extended}_{{R,-,I,q}}\}_I)}(-), 
}
\]
\[
\xymatrix@R+0pc@C+0pc{
\mathrm{K}^\mathrm{BGT}_{\mathrm{sComm}_\mathrm{simplicial}\mathrm{Ind}^m\mathrm{Banach}(\mathcal{O}_{X_{R,-,q}})\ar[r]^{\mathrm{global}}}(-)\ar[r]\ar[r] &\mathrm{K}^\mathrm{BGT}_{\mathrm{sComm}_\mathrm{simplicial}\varphi_q\mathrm{Ind}^m\mathrm{Banach}(\{\mathrm{Robba}^\mathrm{extended}_{{R,-,I,q}}\}_I)}(-). 
}
\]
The definition is given by the following:
\[
\xymatrix@R+0pc@C+0pc{
\mathrm{homotopycolimit}_i(\mathrm{K}^\mathrm{BGT}_{\mathrm{sComm}_\mathrm{simplicial}\mathrm{Ind}\mathrm{Banach}(\mathcal{O}_{X_{R,-,q}})\ar[r]^{\mathrm{global}}}(-)\ar[r]\ar[r] &\mathrm{K}^\mathrm{BGT}_{\mathrm{sComm}_\mathrm{simplicial}\varphi_q\mathrm{Ind}\mathrm{Banach}(\{\mathrm{Robba}^\mathrm{extended}_{{R,-,I,q}}\}_I)}(-))(\mathcal{O}_i),  
}
\]
\[\displayindent=-0.5in
\xymatrix@R+0pc@C+0pc{
\mathrm{homotopycolimit}_i(\mathrm{K}^\mathrm{BGT}_{\mathrm{sComm}_\mathrm{simplicial}\mathrm{Ind}^m\mathrm{Banach}(\mathcal{O}_{X_{R,-,q}})\ar[r]^{\mathrm{global}}}(-)\ar[r]\ar[r] &\mathrm{K}^\mathrm{BGT}_{\mathrm{sComm}_\mathrm{simplicial}\varphi_q\mathrm{Ind}^m\mathrm{Banach}(\{\mathrm{Robba}^\mathrm{extended}_{{R,-,I,q}}\}_I)}(-))(\mathcal{O}_i),  
}
\]
each $\mathcal{O}_i$ is just as $\mathbb{Q}_p\left<C_1,...,C_\ell\right>,\ell=1,q,...$.
\end{itemize}

\noindent Now let $R=\mathbb{Q}_p(p^{1/p^\infty})^{\wedge\flat}$ and $R_k=\mathbb{Q}_p(p^{1/p^\infty})^{\wedge}\left<T_1^{\pm 1/p^{\infty}},...,T_k^{\pm 1/p^{\infty}}\right>^\flat$ we have the following Galois theoretic results with naturality along $f:\mathrm{Spa}(\mathbb{Q}_p(p^{1/p^\infty})^{\wedge}\left<T_1^{\pm 1/p^\infty},...,T_k^{\pm 1/p^\infty}\right>^\flat)\rightarrow \mathrm{Spa}(\mathbb{Q}_p(p^{1/p^\infty})^{\wedge\flat})$:

\begin{itemize}
\item (\text{Proposition}) There is a functor (global section) between the $\infty$-prestacks of inductive Banach quasicoherent presheaves:
\[
\xymatrix@R+6pc@C+0pc{
\mathrm{Ind}\mathrm{Banach}(\mathcal{O}_{X_{\mathbb{Q}_p(p^{1/p^\infty})^{\wedge}\left<T_1^{\pm 1/p^\infty},...,T_k^{\pm 1/p^\infty}\right>^\flat,-,q}})\ar[d]\ar[d]\ar[d]\ar[d] \ar[r]^{\mathrm{global}}\ar[r]\ar[r] &\varphi_q\mathrm{Ind}\mathrm{Banach}(\{\mathrm{Robba}^\mathrm{extended}_{{R_k,-,I,q}}\}_I) \ar[d]\ar[d]\ar[d]\ar[d].\\
\mathrm{Ind}\mathrm{Banach}(\mathcal{O}_{X_{\mathbb{Q}_p(p^{1/p^\infty})^{\wedge\flat},-,q}})\ar[r]^{\mathrm{global}}\ar[r]\ar[r] &\varphi_q\mathrm{Ind}\mathrm{Banach}(\{\mathrm{Robba}^\mathrm{extended}_{{R_0,-,I,q}}\}_I).\\ 
}
\]
The definition is given by the following:
\[
\xymatrix@R+0pc@C+0pc{
\mathrm{homotopycolimit}_i(\square)(\mathcal{O}_i),  
}
\]
each $\mathcal{O}_i$ is just as $\mathbb{Q}_p\left<C_1,...,C_\ell\right>,\ell=1,q,...$ and $\square$ is the relative diagram of $\infty$-functors.
\item (\text{Proposition}) There is a functor (global section) between the $\infty$-prestacks of monomorphic inductive Banach quasicoherent presheaves:
\[
\xymatrix@R+6pc@C+0pc{
\mathrm{Ind}^m\mathrm{Banach}(\mathcal{O}_{X_{R_k,-,q}})\ar[r]^{\mathrm{global}}\ar[d]\ar[d]\ar[d]\ar[d]\ar[r]\ar[r] &\varphi_q\mathrm{Ind}^m\mathrm{Banach}(\{\mathrm{Robba}^\mathrm{extended}_{{R_k,-,I,q}}\}_I)\ar[d]\ar[d]\ar[d]\ar[d]\\
\mathrm{Ind}^m\mathrm{Banach}(\mathcal{O}_{X_{\mathbb{Q}_p(p^{1/p^\infty})^{\wedge\flat},-,q}})\ar[r]^{\mathrm{global}}\ar[r]\ar[r] &\varphi_q\mathrm{Ind}^m\mathrm{Banach}(\{\mathrm{Robba}^\mathrm{extended}_{{R_0,-,I,q}}\}_I).\\  
}
\]
The definition is given by the following:
\[
\xymatrix@R+0pc@C+0pc{
\mathrm{homotopycolimit}_i(\square)(\mathcal{O}_i),  
}
\]
each $\mathcal{O}_i$ is just as $\mathbb{Q}_p\left<C_1,...,C_\ell\right>,\ell=1,q,...$ and $\square$ is the relative diagram of $\infty$-functors.

\item (\text{Proposition}) There is a functor (global section) between the $\infty$-prestacks of inductive Banach quasicoherent commutative algebra $E_\infty$ objects:
\[
\xymatrix@R+6pc@C+0pc{
\mathrm{sComm}_\mathrm{simplicial}\mathrm{Ind}\mathrm{Banach}(\mathcal{O}_{X_{R_k,-,q}})\ar[d]\ar[d]\ar[d]\ar[d]\ar[r]^{\mathrm{global}}\ar[r]\ar[r] &\mathrm{sComm}_\mathrm{simplicial}\varphi_q\mathrm{Ind}\mathrm{Banach}(\{\mathrm{Robba}^\mathrm{extended}_{{R_k,-,I,q}}\}_I)\ar[d]\ar[d]\ar[d]\ar[d]\\
\mathrm{sComm}_\mathrm{simplicial}\mathrm{Ind}\mathrm{Banach}(\mathcal{O}_{X_{\mathbb{Q}_p(p^{1/p^\infty})^{\wedge\flat},-,q}})\ar[r]^{\mathrm{global}}\ar[r]\ar[r] &\mathrm{sComm}_\mathrm{simplicial}\varphi_q\mathrm{Ind}\mathrm{Banach}(\{\mathrm{Robba}^\mathrm{extended}_{{R_0,-,I,q}}\}_I).  
}
\]
The definition is given by the following:
\[
\xymatrix@R+0pc@C+0pc{
\mathrm{homotopycolimit}_i(\square)(\mathcal{O}_i),  
}
\]
each $\mathcal{O}_i$ is just as $\mathbb{Q}_p\left<C_1,...,C_\ell\right>,\ell=1,q,...$ and $\square$ is the relative diagram of $\infty$-functors.

\item (\text{Proposition}) There is a functor (global section) between the $\infty$-prestacks of monomorphic inductive Banach quasicoherent commutative algebra $E_\infty$ objects:
\[\displayindent=-0.4in
\xymatrix@R+6pc@C+0pc{
\mathrm{sComm}_\mathrm{simplicial}\mathrm{Ind}^m\mathrm{Banach}(\mathcal{O}_{X_{R_k,-,q}})\ar[d]\ar[d]\ar[d]\ar[d]\ar[r]^{\mathrm{global}}\ar[r]\ar[r] &\mathrm{sComm}_\mathrm{simplicial}\varphi_q\mathrm{Ind}^m\mathrm{Banach}(\{\mathrm{Robba}^\mathrm{extended}_{{R_k,-,I,q}}\}_I)\ar[d]\ar[d]\ar[d]\ar[d]\\
 \mathrm{sComm}_\mathrm{simplicial}\mathrm{Ind}^m\mathrm{Banach}(\mathcal{O}_{X_{\mathbb{Q}_p(p^{1/p^\infty})^{\wedge\flat},-,q}})\ar[r]^{\mathrm{global}}\ar[r]\ar[r] &\mathrm{sComm}_\mathrm{simplicial}\varphi_q\mathrm{Ind}^m\mathrm{Banach}(\{\mathrm{Robba}^\mathrm{extended}_{{R_0,-,I,q}}\}_I).
}
\]
The definition is given by the following:
\[
\xymatrix@R+0pc@C+0pc{
\mathrm{homotopycolimit}_i(\square)(\mathcal{O}_i),  
}
\]
each $\mathcal{O}_i$ is just as $\mathbb{Q}_p\left<C_1,...,C_\ell\right>,\ell=1,q,...$ and $\square$ is the relative diagram of $\infty$-functors.

\item Then parallel as in \cite{LBV} we have a functor (global section) of the de Rham complex after \cite[Definition 5.9, Section 5.2.1]{KKM}:
\[\displayindent=-0.4in
\xymatrix@R+6pc@C+0pc{
\mathrm{deRham}_{\mathrm{sComm}_\mathrm{simplicial}\mathrm{Ind}\mathrm{Banach}(\mathcal{O}_{X_{R_k,-,q}})\ar[r]^{\mathrm{global}}}(-)\ar[d]\ar[d]\ar[d]\ar[d]\ar[r]\ar[r] &\mathrm{deRham}_{\mathrm{sComm}_\mathrm{simplicial}\varphi_q\mathrm{Ind}\mathrm{Banach}(\{\mathrm{Robba}^\mathrm{extended}_{{R_k,-,I,q}}\}_I)}(-)\ar[d]\ar[d]\ar[d]\ar[d]\\
\mathrm{deRham}_{\mathrm{sComm}_\mathrm{simplicial}\mathrm{Ind}\mathrm{Banach}(\mathcal{O}_{X_{\mathbb{Q}_p(p^{1/p^\infty})^{\wedge\flat},-,q}})\ar[r]^{\mathrm{global}}}(-)\ar[r]\ar[r] &\mathrm{deRham}_{\mathrm{sComm}_\mathrm{simplicial}\varphi_q\mathrm{Ind}\mathrm{Banach}(\{\mathrm{Robba}^\mathrm{extended}_{{R_0,-,I,q}}\}_I)}(-), 
}
\]
\[\displayindent=-0.4in
\xymatrix@R+6pc@C+0pc{
\mathrm{deRham}_{\mathrm{sComm}_\mathrm{simplicial}\mathrm{Ind}^m\mathrm{Banach}(\mathcal{O}_{X_{R_k,-,q}})\ar[r]^{\mathrm{global}}}(-)\ar[d]\ar[d]\ar[d]\ar[d]\ar[r]\ar[r] &\mathrm{deRham}_{\mathrm{sComm}_\mathrm{simplicial}\varphi_q\mathrm{Ind}^m\mathrm{Banach}(\{\mathrm{Robba}^\mathrm{extended}_{{R_k,-,I,q}}\}_I)}(-)\ar[d]\ar[d]\ar[d]\ar[d]\\
\mathrm{deRham}_{\mathrm{sComm}_\mathrm{simplicial}\mathrm{Ind}^m\mathrm{Banach}(\mathcal{O}_{X_{\mathbb{Q}_p(p^{1/p^\infty})^{\wedge\flat},-,q}})\ar[r]^{\mathrm{global}}}(-)\ar[r]\ar[r] &\mathrm{deRham}_{\mathrm{sComm}_\mathrm{simplicial}\varphi_q\mathrm{Ind}^m\mathrm{Banach}(\{\mathrm{Robba}^\mathrm{extended}_{{R_0,-,I,q}}\}_I)}(-). 
}
\]

\item Then we have the following a functor (global section) of $K$-group $(\infty,1)$-spectrum from \cite{BGT}:
\[
\xymatrix@R+6pc@C+0pc{
\mathrm{K}^\mathrm{BGT}_{\mathrm{sComm}_\mathrm{simplicial}\mathrm{Ind}\mathrm{Banach}(\mathcal{O}_{X_{R_k,-,q}})\ar[r]^{\mathrm{global}}}(-)\ar[d]\ar[d]\ar[d]\ar[d]\ar[r]\ar[r] &\mathrm{K}^\mathrm{BGT}_{\mathrm{sComm}_\mathrm{simplicial}\varphi_q\mathrm{Ind}\mathrm{Banach}(\{\mathrm{Robba}^\mathrm{extended}_{{R_k,-,I,q}}\}_I)}(-)\ar[d]\ar[d]\ar[d]\ar[d]\\
\mathrm{K}^\mathrm{BGT}_{\mathrm{sComm}_\mathrm{simplicial}\mathrm{Ind}\mathrm{Banach}(\mathcal{O}_{X_{\mathbb{Q}_p(p^{1/p^\infty})^{\wedge\flat},-,q}})\ar[r]^{\mathrm{global}}}(-)\ar[r]\ar[r] &\mathrm{K}^\mathrm{BGT}_{\mathrm{sComm}_\mathrm{simplicial}\varphi_q\mathrm{Ind}\mathrm{Banach}(\{\mathrm{Robba}^\mathrm{extended}_{{R_0,-,I,q}}\}_I)}(-), 
}
\]
\[
\xymatrix@R+6pc@C+0pc{
\mathrm{K}^\mathrm{BGT}_{\mathrm{sComm}_\mathrm{simplicial}\mathrm{Ind}^m\mathrm{Banach}(\mathcal{O}_{X_{R_k,-,q}})\ar[r]^{\mathrm{global}}}(-)\ar[d]\ar[d]\ar[d]\ar[d]\ar[r]\ar[r] &\mathrm{K}^\mathrm{BGT}_{\mathrm{sComm}_\mathrm{simplicial}\varphi_q\mathrm{Ind}^m\mathrm{Banach}(\{\mathrm{Robba}^\mathrm{extended}_{{R_k,-,I,q}}\}_I)}(-)\ar[d]\ar[d]\ar[d]\ar[d]\\
\mathrm{K}^\mathrm{BGT}_{\mathrm{sComm}_\mathrm{simplicial}\mathrm{Ind}^m\mathrm{Banach}(\mathcal{O}_{X_{\mathbb{Q}_p(p^{1/p^\infty})^{\wedge\flat},-,q}})\ar[r]^{\mathrm{global}}}(-)\ar[r]\ar[r] &\mathrm{K}^\mathrm{BGT}_{\mathrm{sComm}_\mathrm{simplicial}\varphi_q\mathrm{Ind}^m\mathrm{Banach}(\{\mathrm{Robba}^\mathrm{extended}_{{R_0,-,I,q}}\}_I)}(-). 
}
\]
The definition is given by the following:
\[
\xymatrix@R+0pc@C+0pc{
\mathrm{homotopycolimit}_i(\square)(\mathcal{O}_i),  
}
\]
each $\mathcal{O}_i$ is just as $\mathbb{Q}_p\left<C_1,...,C_\ell\right>,\ell=1,q,...$ and $\square$ is the relative diagram of $\infty$-functors.

\end{itemize}

\
\indent Then we consider further equivariance by considering the arithmetic profinite fundamental groups and actually its $q$-th power $\mathrm{Gal}(\overline{{Q}_p\left<T_1^{\pm 1},...,T_k^{\pm 1}\right>}/R_k)^{\times q}$ through the following diagram:\\

\[
\xymatrix@R+6pc@C+0pc{
\mathbb{Z}_p^k=\mathrm{Gal}(R_k/{\mathbb{Q}_p(p^{1/p^\infty})^\wedge\left<T_1^{\pm 1},...,T_k^{\pm 1}\right>})\ar[d]\ar[d]\ar[d]\ar[d] \ar[r]\ar[r] \ar[r]\ar[r] &\mathrm{Gal}(\overline{{Q}_p\left<T_1^{\pm 1},...,T_k^{\pm 1}\right>}/R_k) \ar[d]\ar[d]\ar[d] \ar[r]\ar[r] &\Gamma_{\mathbb{Q}_p} \ar[d]\ar[d]\ar[d]\ar[d]\\
(\mathbb{Z}_p^k=\mathrm{Gal}(R_k/{\mathbb{Q}_p(p^{1/p^\infty})^\wedge\left<T_1^{\pm 1},...,T_k^{\pm 1}\right>}))^{\times q} \ar[r]\ar[r] \ar[r]\ar[r] &\Gamma_k^{\times q}:=\mathrm{Gal}(R_k/{\mathbb{Q}_p\left<T_1^{\pm 1},...,T_k^{\pm 1}\right>})^{\times q} \ar[r] \ar[r]\ar[r] &\Gamma_{\mathbb{Q}_p}^{\times q}.
}
\]

\

We then have the correspond arithmetic profinite fundamental groups equivariant versions:

\begin{itemize}
\item (\text{Proposition}) There is a functor (global section) between the $\infty$-prestacks of inductive Banach quasicoherent presheaves:
\[
\xymatrix@R+6pc@C+0pc{
\mathrm{Ind}\mathrm{Banach}_{\Gamma_{k}^{\times q}}(\mathcal{O}_{X_{\mathbb{Q}_p(p^{1/p^\infty})^{\wedge}\left<T_1^{\pm 1/p^\infty},...,T_k^{\pm 1/p^\infty}\right>^\flat,-,q}})\ar[d]\ar[d]\ar[d]\ar[d] \ar[r]^{\mathrm{global}}\ar[r]\ar[r] &\varphi_q\mathrm{Ind}\mathrm{Banach}_{\Gamma_{k}^{\times q}}(\{\mathrm{Robba}^\mathrm{extended}_{{R_k,-,I,q}}\}_I) \ar[d]\ar[d]\ar[d]\ar[d].\\
\mathrm{Ind}\mathrm{Banach}_{\Gamma_{0}^{\times q}}(\mathcal{O}_{X_{\mathbb{Q}_p(p^{1/p^\infty})^{\wedge\flat},-,q}})\ar[r]^{\mathrm{global}}\ar[r]\ar[r] &\varphi_q\mathrm{Ind}\mathrm{Banach}_{\Gamma_{0}^{\times q}}(\{\mathrm{Robba}^\mathrm{extended}_{{R_0,-,I,q}}\}_I).\\ 
}
\]
The definition is given by the following:
\[
\xymatrix@R+0pc@C+0pc{
\mathrm{homotopycolimit}_i(\square)(\mathcal{O}_i),  
}
\]
each $\mathcal{O}_i$ is just as $\mathbb{Q}_p\left<C_1,...,C_\ell\right>,\ell=1,q,...$ and $\square$ is the relative diagram of $\infty$-functors.

\item (\text{Proposition}) There is a functor (global section) between the $\infty$-prestacks of monomorphic inductive Banach quasicoherent presheaves:
\[
\xymatrix@R+6pc@C+0pc{
\mathrm{Ind}^m\mathrm{Banach}_{\Gamma_{k}^{\times q}}(\mathcal{O}_{X_{R_k,-,q}})\ar[r]^{\mathrm{global}}\ar[d]\ar[d]\ar[d]\ar[d]\ar[r]\ar[r] &\varphi_q\mathrm{Ind}^m\mathrm{Banach}_{\Gamma_{k}^{\times q}}(\{\mathrm{Robba}^\mathrm{extended}_{{R_k,-,I,q}}\}_I)\ar[d]\ar[d]\ar[d]\ar[d]\\
\mathrm{Ind}^m\mathrm{Banach}_{\Gamma_{0}^{\times q}}(\mathcal{O}_{X_{\mathbb{Q}_p(p^{1/p^\infty})^{\wedge\flat},-,q}})\ar[r]^{\mathrm{global}}\ar[r]\ar[r] &\varphi_q\mathrm{Ind}^m\mathrm{Banach}_{\Gamma_{0}^{\times q}}(\{\mathrm{Robba}^\mathrm{extended}_{{R_0,-,I,q}}\}_I).\\  
}
\]
The definition is given by the following:
\[
\xymatrix@R+0pc@C+0pc{
\mathrm{homotopycolimit}_i(\square)(\mathcal{O}_i),  
}
\]
each $\mathcal{O}_i$ is just as $\mathbb{Q}_p\left<C_1,...,C_\ell\right>,\ell=1,q,...$ and $\square$ is the relative diagram of $\infty$-functors.

\item (\text{Proposition}) There is a functor (global section) between the $\infty$-stacks of inductive Banach quasicoherent commutative algebra $E_\infty$ objects:
\[\displayindent=-0.4in
\xymatrix@R+6pc@C+0pc{
\mathrm{sComm}_\mathrm{simplicial}\mathrm{Ind}\mathrm{Banach}_{\Gamma_{k}^{\times q}}(\mathcal{O}_{X_{R_k,-,q}})\ar[d]\ar[d]\ar[d]\ar[d]\ar[r]^{\mathrm{global}}\ar[r]\ar[r] &\mathrm{sComm}_\mathrm{simplicial}\varphi_q\mathrm{Ind}\mathrm{Banach}_{\Gamma_{k}^{\times q}}(\{\mathrm{Robba}^\mathrm{extended}_{{R_k,-,I,q}}\}_I)\ar[d]\ar[d]\ar[d]\ar[d]\\
\mathrm{sComm}_\mathrm{simplicial}\mathrm{Ind}\mathrm{Banach}_{\Gamma_{0}^{\times q}}(\mathcal{O}_{X_{\mathbb{Q}_p(p^{1/p^\infty})^{\wedge\flat},-,q}})\ar[r]^{\mathrm{global}}\ar[r]\ar[r] &\mathrm{sComm}_\mathrm{simplicial}\varphi_q\mathrm{Ind}\mathrm{Banach}_{\Gamma_{0}^{\times q}}(\{\mathrm{Robba}^\mathrm{extended}_{{R_0,-,I,q}}\}_I).  
}
\]
The definition is given by the following:
\[
\xymatrix@R+0pc@C+0pc{
\mathrm{homotopycolimit}_i(\square)(\mathcal{O}_i),  
}
\]
each $\mathcal{O}_i$ is just as $\mathbb{Q}_p\left<C_1,...,C_\ell\right>,\ell=1,q,...$ and $\square$ is the relative diagram of $\infty$-functors.

\item (\text{Proposition}) There is a functor (global section) between the $\infty$-prestacks of monomorphic inductive Banach quasicoherent commutative algebra $E_\infty$ objects:
\[\displayindent=-0.4in
\xymatrix@R+6pc@C+0pc{
\mathrm{sComm}_\mathrm{simplicial}\mathrm{Ind}^m\mathrm{Banach}_{\Gamma_{k}^{\times q}}(\mathcal{O}_{X_{R_k,-,q}})\ar[d]\ar[d]\ar[d]\ar[d]\ar[r]^{\mathrm{global}}\ar[r]\ar[r] &\mathrm{sComm}_\mathrm{simplicial}\varphi_q\mathrm{Ind}^m\mathrm{Banach}_{\Gamma_{k}^{\times q}}(\{\mathrm{Robba}^\mathrm{extended}_{{R_k,-,I,q}}\}_I)\ar[d]\ar[d]\ar[d]\ar[d]\\
 \mathrm{sComm}_\mathrm{simplicial}\mathrm{Ind}^m\mathrm{Banach}_{\Gamma_{0}^{\times q}}(\mathcal{O}_{X_{\mathbb{Q}_p(p^{1/p^\infty})^{\wedge\flat},-,q}})\ar[r]^{\mathrm{global}}\ar[r]\ar[r] &\mathrm{sComm}_\mathrm{simplicial}\varphi_q\mathrm{Ind}^m\mathrm{Banach}_{\Gamma_{0}^{\times q}}(\{\mathrm{Robba}^\mathrm{extended}_{{R_0,-,I,q}}\}_I). 
}
\]
The definition is given by the following:
\[
\xymatrix@R+0pc@C+0pc{
\mathrm{homotopycolimit}_i(\square)(\mathcal{O}_i),  
}
\]
each $\mathcal{O}_i$ is just as $\mathbb{Q}_p\left<C_1,...,C_\ell\right>,\ell=1,q,...$ and $\square$ is the relative diagram of $\infty$-functors.

\item Then parallel as in \cite{LBV} we have a functor (global section) of the de Rham complex after \cite[Definition 5.9, Section 5.2.1]{KKM}:
\[\displayindent=-0.6in
\xymatrix@R+6pc@C+0pc{
\mathrm{deRham}_{\mathrm{sComm}_\mathrm{simplicial}\mathrm{Ind}\mathrm{Banach}_{\Gamma_{k}^{\times q}}(\mathcal{O}_{X_{R_k,-,q}})\ar[r]^{\mathrm{global}}}(-)\ar[d]\ar[d]\ar[d]\ar[d]\ar[r]\ar[r] &\mathrm{deRham}_{\mathrm{sComm}_\mathrm{simplicial}\varphi_q\mathrm{Ind}\mathrm{Banach}_{\Gamma_{k}^{\times q}}(\{\mathrm{Robba}^\mathrm{extended}_{{R_k,-,I,q}}\}_I)}(-)\ar[d]\ar[d]\ar[d]\ar[d]\\
\mathrm{deRham}_{\mathrm{sComm}_\mathrm{simplicial}\mathrm{Ind}\mathrm{Banach}_{\Gamma_{0}^{\times q}}(\mathcal{O}_{X_{\mathbb{Q}_p(p^{1/p^\infty})^{\wedge\flat},-,q}})\ar[r]^{\mathrm{global}}}(-)\ar[r]\ar[r] &\mathrm{deRham}_{\mathrm{sComm}_\mathrm{simplicial}\varphi_q\mathrm{Ind}\mathrm{Banach}_{\Gamma_{0}^{\times q}}(\{\mathrm{Robba}^\mathrm{extended}_{{R_0,-,I,q}}\}_I)}(-), 
}
\]
\[\displayindent=-0.7in
\xymatrix@R+6pc@C+0pc{
\mathrm{deRham}_{\mathrm{sComm}_\mathrm{simplicial}\mathrm{Ind}^m\mathrm{Banach}_{\Gamma_{k}^{\times q}}(\mathcal{O}_{X_{R_k,-,q}})\ar[r]^{\mathrm{global}}}(-)\ar[d]\ar[d]\ar[d]\ar[d]\ar[r]\ar[r] &\mathrm{deRham}_{\mathrm{sComm}_\mathrm{simplicial}\varphi_q\mathrm{Ind}^m\mathrm{Banach}_{\Gamma_{k}^{\times q}}(\{\mathrm{Robba}^\mathrm{extended}_{{R_k,-,I,q}}\}_I)}(-)\ar[d]\ar[d]\ar[d]\ar[d]\\
\mathrm{deRham}_{\mathrm{sComm}_\mathrm{simplicial}\mathrm{Ind}^m\mathrm{Banach}_{\Gamma_{0}^{\times q}}(\mathcal{O}_{X_{\mathbb{Q}_p(p^{1/p^\infty})^{\wedge\flat},-,q}})\ar[r]^{\mathrm{global}}}(-)\ar[r]\ar[r] &\mathrm{deRham}_{\mathrm{sComm}_\mathrm{simplicial}\varphi_q\mathrm{Ind}^m\mathrm{Banach}_{\Gamma_{0}^{\times q}}(\{\mathrm{Robba}^\mathrm{extended}_{{R_0,-,I,q}}\}_I)}(-). 
}
\]
The definition is given by the following:
\[
\xymatrix@R+0pc@C+0pc{
\mathrm{homotopycolimit}_i(\square)(\mathcal{O}_i),  
}
\]
each $\mathcal{O}_i$ is just as $\mathbb{Q}_p\left<C_1,...,C_\ell\right>,\ell=1,q,...$ and $\square$ is the relative diagram of $\infty$-functors.

\item Then we have the following a functor (global section) of $K$-group $(\infty,1)$-spectrum from \cite{BGT}:
\[
\xymatrix@R+6pc@C+0pc{
\mathrm{K}^\mathrm{BGT}_{\mathrm{sComm}_\mathrm{simplicial}\mathrm{Ind}\mathrm{Banach}_{\Gamma_{k}^{\times q}}(\mathcal{O}_{X_{R_k,-,q}})\ar[r]^{\mathrm{global}}}(-)\ar[d]\ar[d]\ar[d]\ar[d]\ar[r]\ar[r] &\mathrm{K}^\mathrm{BGT}_{\mathrm{sComm}_\mathrm{simplicial}\varphi_q\mathrm{Ind}\mathrm{Banach}_{\Gamma_{k}^{\times q}}(\{\mathrm{Robba}^\mathrm{extended}_{{R_k,-,I,q}}\}_I)}(-)\ar[d]\ar[d]\ar[d]\ar[d]\\
\mathrm{K}^\mathrm{BGT}_{\mathrm{sComm}_\mathrm{simplicial}\mathrm{Ind}\mathrm{Banach}_{\Gamma_{0}^{\times q}}(\mathcal{O}_{X_{\mathbb{Q}_p(p^{1/p^\infty})^{\wedge\flat},-,q}})\ar[r]^{\mathrm{global}}}(-)\ar[r]\ar[r] &\mathrm{K}^\mathrm{BGT}_{\mathrm{sComm}_\mathrm{simplicial}\varphi_q\mathrm{Ind}\mathrm{Banach}_{\Gamma_{0}^{\times q}}(\{\mathrm{Robba}^\mathrm{extended}_{{R_0,-,I,q}}\}_I)}(-), 
}
\]
\[
\xymatrix@R+6pc@C+0pc{
\mathrm{K}^\mathrm{BGT}_{\mathrm{sComm}_\mathrm{simplicial}\mathrm{Ind}^m\mathrm{Banach}_{\Gamma_{k}^{\times q}}(\mathcal{O}_{X_{R_k,-,q}})\ar[r]^{\mathrm{global}}}(-)\ar[d]\ar[d]\ar[d]\ar[d]\ar[r]\ar[r] &\mathrm{K}^\mathrm{BGT}_{\mathrm{sComm}_\mathrm{simplicial}\varphi_q\mathrm{Ind}^m\mathrm{Banach}_{\Gamma_{k}^{\times q}}(\{\mathrm{Robba}^\mathrm{extended}_{{R_k,-,I,q}}\}_I)}(-)\ar[d]\ar[d]\ar[d]\ar[d]\\
\mathrm{K}^\mathrm{BGT}_{\mathrm{sComm}_\mathrm{simplicial}\mathrm{Ind}^m\mathrm{Banach}_{\Gamma_{0}^{\times q}}(\mathcal{O}_{X_{\mathbb{Q}_p(p^{1/p^\infty})^{\wedge\flat},-,q}})\ar[r]^{\mathrm{global}}}(-)\ar[r]\ar[r] &\mathrm{K}^\mathrm{BGT}_{\mathrm{sComm}_\mathrm{simplicial}\varphi_q\mathrm{Ind}^m\mathrm{Banach}_{\Gamma_{0}^{\times q}}(\{\mathrm{Robba}^\mathrm{extended}_{{R_0,-,I,q}}\}_I)}(-). 
}
\]
The definition is given by the following:
\[
\xymatrix@R+0pc@C+0pc{
\mathrm{homotopycolimit}_i(\square)(\mathcal{O}_i),  
}
\]
each $\mathcal{O}_i$ is just as $\mathbb{Q}_p\left<C_1,...,C_\ell\right>,\ell=1,q,...$ and $\square$ is the relative diagram of $\infty$-functors.

\end{itemize}

\

\begin{remark}
\noindent We can certainly consider the quasicoherent sheaves in \cite[Lemma 7.11, Remark 7.12]{1BBK}, therefore all the quasicoherent presheaves and modules will be those satisfying \cite[Lemma 7.11, Remark 7.12]{1BBK} if one would like to consider the the quasicoherent sheaves. That being all as this said, we would believe that the big quasicoherent presheaves are automatically quasicoherent sheaves (namely satisfying the corresponding \v{C}ech $\infty$-descent as in \cite[Section 9.3]{KKM} and \cite[Lemma 7.11, Remark 7.12]{1BBK}) and the corresponding global section functors are automatically equivalence of $\infty$-categories. \\
\end{remark}

\

\indent In Clausen-Scholze formalism we have the following\footnote{Certainly the homotopy colimit in the rings side will be within the condensed solid animated analytic rings from \cite{1CS2}.}:

\begin{itemize}
\item (\text{Proposition}) There is a functor (global section) between the $\infty$-prestacks of inductive Banach quasicoherent sheaves:
\[
\xymatrix@R+0pc@C+0pc{
{\mathrm{Modules}_\circledcirc}(\mathcal{O}_{X_{R,-,q}})\ar[r]^{\mathrm{global}}\ar[r]\ar[r] &\varphi_q{\mathrm{Modules}_\circledcirc}(\{\mathrm{Robba}^{\mathrm{extended},q}_{{R,-,I}}\}_I).  
}
\]
The definition is given by the following:
\[
\xymatrix@R+0pc@C+0pc{
\mathrm{homotopycolimit}_i({\mathrm{Modules}_\circledcirc}(\mathcal{O}_{X_{R,-,q}})\ar[r]^{\mathrm{global}}\ar[r]\ar[r] &\varphi_q{\mathrm{Modules}_\circledcirc}(\{\mathrm{Robba}^{\mathrm{extended},q}_{{R,-,I}}\}_I))(\mathcal{O}_i),  
}
\]
each $\mathcal{O}_i$ is just as $\mathbb{Q}_p\left<C_1,...,C_\ell\right>,\ell=1,2,...$.

\item (\text{Proposition}) There is a functor (global section) between the $\infty$-prestacks of inductive Banach quasicoherent sheaves:
\[
\xymatrix@R+0pc@C+0pc{
{\mathrm{Modules}_\circledcirc}(\mathcal{O}_{X_{R,-,q}})\ar[r]^{\mathrm{global}}\ar[r]\ar[r] &\varphi_q{\mathrm{Modules}_\circledcirc}(\{\mathrm{Robba}^{\mathrm{extended},q}_{{R,-,I}}\}_I).  
}
\]
The definition is given by the following:
\[
\xymatrix@R+0pc@C+0pc{
\mathrm{homotopycolimit}_i({\mathrm{Modules}_\circledcirc}(\mathcal{O}_{X_{R,-,q}})\ar[r]^{\mathrm{global}}\ar[r]\ar[r] &\varphi_q{\mathrm{Modules}_\circledcirc}(\{\mathrm{Robba}^{\mathrm{extended},q}_{{R,-,I}}\}_I))(\mathcal{O}_i),  
}
\]
each $\mathcal{O}_i$ is just as $\mathbb{Q}_p\left<C_1,...,C_\ell\right>,\ell=1,2,...$.

\item (\text{Proposition}) There is a functor (global section) between the $\infty$-prestacks of inductive Banach quasicoherent commutative algebra $E_\infty$ objects\footnote{Here $\circledcirc=\text{solidquasicoherentsheaves}$.}:
\[
\xymatrix@R+0pc@C+0pc{
\mathrm{sComm}_\mathrm{simplicial}{\mathrm{Modules}_\circledcirc}(\mathcal{O}_{X_{R,-,q}})\ar[r]^{\mathrm{global}}\ar[r]\ar[r] &\mathrm{sComm}_\mathrm{simplicial}\varphi_q{\mathrm{Modules}_\circledcirc}(\{\mathrm{Robba}^{\mathrm{extended},q}_{{R,-,I}}\}_I).  
}
\]
The definition is given by the following:
\[
\xymatrix@R+0pc@C+0pc{
\mathrm{homotopycolimit}_i(\mathrm{sComm}_\mathrm{simplicial}{\mathrm{Modules}_\circledcirc}(\mathcal{O}_{X_{R,-,q}})\ar[r]^{\mathrm{global}}\ar[r]\ar[r] &\mathrm{sComm}_\mathrm{simplicial}\varphi_q{\mathrm{Modules}_\circledcirc}(\{\mathrm{Robba}^{\mathrm{extended},q}_{{R,-,I}}\}_I))(\mathcal{O}_i),  
}
\]
each $\mathcal{O}_i$ is just as $\mathbb{Q}_p\left<C_1,...,C_\ell\right>,\ell=1,2,...$.
\item Then as in \cite{LBV} we have a functor (global section ) of the de Rham complex after \cite[Definition 5.9, Section 5.2.1]{KKM}\footnote{Here $\circledcirc=\text{solidquasicoherentsheaves}$.}:
\[
\xymatrix@R+0pc@C+0pc{
\mathrm{deRham}_{\mathrm{sComm}_\mathrm{simplicial}{\mathrm{Modules}_\circledcirc}(\mathcal{O}_{X_{R,-,q}})\ar[r]^{\mathrm{global}}}(-)\ar[r]\ar[r] &\mathrm{deRham}_{\mathrm{sComm}_\mathrm{simplicial}\varphi_q{\mathrm{Modules}_\circledcirc}(\{\mathrm{Robba}^{\mathrm{extended},q}_{{R,-,I}}\}_I)}(-), 
}
\]
The definition is given by the following:
\[
\xymatrix@R+0pc@C+0pc{
\mathrm{homotopycolimit}_i\\
(\mathrm{deRham}_{\mathrm{sComm}_\mathrm{simplicial}{\mathrm{Modules}_\circledcirc}(\mathcal{O}_{X_{R,-,q}})\ar[r]^{\mathrm{global}}}(-)\ar[r]\ar[r] &\mathrm{deRham}_{\mathrm{sComm}_\mathrm{simplicial}\varphi_q{\mathrm{Modules}_\circledcirc}(\{\mathrm{Robba}^{\mathrm{extended},q}_{{R,-,I}}\}_I)}(-))(\mathcal{O}_i),  
}
\]
each $\mathcal{O}_i$ is just as $\mathbb{Q}_p\left<C_1,...,C_\ell\right>,\ell=1,2,...$.\item Then we have the following a functor (global section) of $K$-group $(\infty,1)$-spectrum from \cite{BGT}\footnote{Here $\circledcirc=\text{solidquasicoherentsheaves}$.}:
\[
\xymatrix@R+0pc@C+0pc{
\mathrm{K}^\mathrm{BGT}_{\mathrm{sComm}_\mathrm{simplicial}{\mathrm{Modules}_\circledcirc}(\mathcal{O}_{X_{R,-,q}})\ar[r]^{\mathrm{global}}}(-)\ar[r]\ar[r] &\mathrm{K}^\mathrm{BGT}_{\mathrm{sComm}_\mathrm{simplicial}\varphi_q{\mathrm{Modules}_\circledcirc}(\{\mathrm{Robba}^{\mathrm{extended},q}_{{R,-,I}}\}_I)}(-). 
}
\]
The definition is given by the following:
\[
\xymatrix@R+0pc@C+0pc{
\mathrm{homotopycolimit}_i(\mathrm{K}^\mathrm{BGT}_{\mathrm{sComm}_\mathrm{simplicial}{\mathrm{Modules}_\circledcirc}(\mathcal{O}_{X_{R,-,q}})\ar[r]^{\mathrm{global}}}(-)\ar[r]\ar[r] &\mathrm{K}^\mathrm{BGT}_{\mathrm{sComm}_\mathrm{simplicial}\varphi_q{\mathrm{Modules}_\circledcirc}(\{\mathrm{Robba}^{\mathrm{extended},q}_{{R,-,I}}\}_I)}(-))(\mathcal{O}_i),  
}
\]
each $\mathcal{O}_i$ is just as $\mathbb{Q}_p\left<C_1,...,C_\ell\right>,\ell=1,2,...$.
\end{itemize}

\noindent Now let $R=\mathbb{Q}_p(p^{1/p^\infty})^{\wedge\flat}$ and $R_k=\mathbb{Q}_p(p^{1/p^\infty})^{\wedge}\left<T_1^{\pm 1/p^{\infty}},...,T_k^{\pm 1/p^{\infty}}\right>^\flat$ we have the following Galois theoretic results with naturality along $f:\mathrm{Spa}(\mathbb{Q}_p(p^{1/p^\infty})^{\wedge}\left<T_1^{\pm 1/p^\infty},...,T_k^{\pm 1/p^\infty}\right>^\flat)\rightarrow \mathrm{Spa}(\mathbb{Q}_p(p^{1/p^\infty})^{\wedge\flat})$:

\begin{itemize}
\item (\text{Proposition}) There is a functor (global section) between the $\infty$-prestacks of inductive Banach quasicoherent sheaves\footnote{Here $\circledcirc=\text{solidquasicoherentsheaves}$.}:
\[
\xymatrix@R+6pc@C+0pc{
{\mathrm{Modules}_\circledcirc}(\mathcal{O}_{X_{\mathbb{Q}_p(p^{1/p^\infty})^{\wedge}\left<T_1^{\pm 1/p^\infty},...,T_k^{\pm 1/p^\infty}\right>^\flat,-,q}})\ar[d]\ar[d]\ar[d]\ar[d] \ar[r]^{\mathrm{global}}\ar[r]\ar[r] &\varphi_q{\mathrm{Modules}_\circledcirc}(\{\mathrm{Robba}^{\mathrm{extended},q}_{{R_k,-,I}}\}_I) \ar[d]\ar[d]\ar[d]\ar[d].\\
{\mathrm{Modules}_\circledcirc}(\mathcal{O}_{X_{\mathbb{Q}_p(p^{1/p^\infty})^{\wedge\flat},-,q}})\ar[r]^{\mathrm{global}}\ar[r]\ar[r] &\varphi_q{\mathrm{Modules}_\circledcirc}(\{\mathrm{Robba}^{\mathrm{extended},q}_{{R_0,-,I}}\}_I).\\ 
}
\]
The definition is given by the following:
\[
\xymatrix@R+0pc@C+0pc{
\mathrm{homotopycolimit}_i(\square)(\mathcal{O}_i),  
}
\]
each $\mathcal{O}_i$ is just as $\mathbb{Q}_p\left<C_1,...,C_\ell\right>,\ell=1,2,...$ and $\square$ is the relative diagram of $\infty$-functors.

\item (\text{Proposition}) There is a functor (global section) between the $\infty$-prestacks of inductive Banach quasicoherent commutative algebra $E_\infty$ objects\footnote{Here $\circledcirc=\text{solidquasicoherentsheaves}$.}:
\[
\xymatrix@R+6pc@C+0pc{
\mathrm{sComm}_\mathrm{simplicial}{\mathrm{Modules}_\circledcirc}(\mathcal{O}_{X_{R_k,-,q}})\ar[d]\ar[d]\ar[d]\ar[d]\ar[r]^{\mathrm{global}}\ar[r]\ar[r] &\mathrm{sComm}_\mathrm{simplicial}\varphi_q{\mathrm{Modules}_\circledcirc}(\{\mathrm{Robba}^{\mathrm{extended},q}_{{R_k,-,I}}\}_I)\ar[d]\ar[d]\ar[d]\ar[d]\\
\mathrm{sComm}_\mathrm{simplicial}{\mathrm{Modules}_\circledcirc}(\mathcal{O}_{X_{\mathbb{Q}_p(p^{1/p^\infty})^{\wedge\flat},-,q}})\ar[r]^{\mathrm{global}}\ar[r]\ar[r] &\mathrm{sComm}_\mathrm{simplicial}\varphi_q{\mathrm{Modules}_\circledcirc}(\{\mathrm{Robba}^{\mathrm{extended},q}_{{R_0,-,I}}\}_I).  
}
\]
The definition is given by the following:
\[
\xymatrix@R+0pc@C+0pc{
\mathrm{homotopycolimit}_i(\square)(\mathcal{O}_i),  
}
\]
each $\mathcal{O}_i$ is just as $\mathbb{Q}_p\left<C_1,...,C_\ell\right>,\ell=1,2,...$ and $\square$ is the relative diagram of $\infty$-functors.

\item Then as in \cite{LBV} we have a functor (global section) of the de Rham complex after \cite[Definition 5.9, Section 5.2.1]{KKM}\footnote{Here $\circledcirc=\text{solidquasicoherentsheaves}$.}:
\[\displayindent=-0.4in
\xymatrix@R+6pc@C+0pc{
\mathrm{deRham}_{\mathrm{sComm}_\mathrm{simplicial}{\mathrm{Modules}_\circledcirc}(\mathcal{O}_{X_{R_k,-,q}})\ar[r]^{\mathrm{global}}}(-)\ar[d]\ar[d]\ar[d]\ar[d]\ar[r]\ar[r] &\mathrm{deRham}_{\mathrm{sComm}_\mathrm{simplicial}\varphi_q{\mathrm{Modules}_\circledcirc}(\{\mathrm{Robba}^{\mathrm{extended},q}_{{R_k,-,I}}\}_I)}(-)\ar[d]\ar[d]\ar[d]\ar[d]\\
\mathrm{deRham}_{\mathrm{sComm}_\mathrm{simplicial}{\mathrm{Modules}_\circledcirc}(\mathcal{O}_{X_{\mathbb{Q}_p(p^{1/p^\infty})^{\wedge\flat},-,q}})\ar[r]^{\mathrm{global}}}(-)\ar[r]\ar[r] &\mathrm{deRham}_{\mathrm{sComm}_\mathrm{simplicial}\varphi_q{\mathrm{Modules}_\circledcirc}(\{\mathrm{Robba}^{\mathrm{extended},q}_{{R_0,-,I}}\}_I)}(-), 
}
\]

\item Then we have the following a functor (global section) of $K$-group $(\infty,1)$-spectrum from \cite{BGT}\footnote{Here $\circledcirc=\text{solidquasicoherentsheaves}$.}:
\[
\xymatrix@R+6pc@C+0pc{
\mathrm{K}^\mathrm{BGT}_{\mathrm{sComm}_\mathrm{simplicial}{\mathrm{Modules}_\circledcirc}(\mathcal{O}_{X_{R_k,-,q}})\ar[r]^{\mathrm{global}}}(-)\ar[d]\ar[d]\ar[d]\ar[d]\ar[r]\ar[r] &\mathrm{K}^\mathrm{BGT}_{\mathrm{sComm}_\mathrm{simplicial}\varphi_q{\mathrm{Modules}_\circledcirc}(\{\mathrm{Robba}^{\mathrm{extended},q}_{{R_k,-,I}}\}_I)}(-)\ar[d]\ar[d]\ar[d]\ar[d]\\
\mathrm{K}^\mathrm{BGT}_{\mathrm{sComm}_\mathrm{simplicial}{\mathrm{Modules}_\circledcirc}(\mathcal{O}_{X_{\mathbb{Q}_p(p^{1/p^\infty})^{\wedge\flat},-,q}})\ar[r]^{\mathrm{global}}}(-)\ar[r]\ar[r] &\mathrm{K}^\mathrm{BGT}_{\mathrm{sComm}_\mathrm{simplicial}\varphi_q{\mathrm{Modules}_\circledcirc}(\{\mathrm{Robba}^{\mathrm{extended},q}_{{R_0,-,I}}\}_I)}(-), 
}
\]

The definition is given by the following:
\[
\xymatrix@R+0pc@C+0pc{
\mathrm{homotopycolimit}_i(\square)(\mathcal{O}_i),  
}
\]
each $\mathcal{O}_i$ is just as $\mathbb{Q}_p\left<C_1,...,C_\ell\right>,\ell=1,2,...$ and $\square$ is the relative diagram of $\infty$-functors.

\end{itemize}

\
\indent Then we consider further equivariance by considering the arithmetic profinite fundamental groups $\Gamma_{\mathbb{Q}_p}$ and $\mathrm{Gal}(\overline{{Q}_p\left<T_1^{\pm 1},...,T_k^{\pm 1}\right>}/R_k)$ through the following diagram:

\[
\xymatrix@R+0pc@C+0pc{
\mathbb{Z}_p^k=\mathrm{Gal}(R_k/{\mathbb{Q}_p(p^{1/p^\infty})^\wedge\left<T_1^{\pm 1},...,T_k^{\pm 1}\right>}) \ar[r]\ar[r] \ar[r]\ar[r] &\Gamma_k^{\times q}:=\mathrm{Gal}(R_k/{\mathbb{Q}_p\left<T_1^{\pm 1},...,T_k^{\pm 1}\right>}) \ar[r] \ar[r]\ar[r] &\Gamma_{\mathbb{Q}_p}.
}
\]

\begin{itemize}
\item (\text{Proposition}) There is a functor (global section) between the $\infty$-prestacks of inductive Banach quasicoherent sheaves\footnote{Here $\circledcirc=\text{solidquasicoherentsheaves}$.}:
\[
\xymatrix@R+6pc@C+0pc{
{\mathrm{Modules}_\circledcirc}_{\Gamma_k^{\times q}}(\mathcal{O}_{X_{\mathbb{Q}_p(p^{1/p^\infty})^{\wedge}\left<T_1^{\pm 1/p^\infty},...,T_k^{\pm 1/p^\infty}\right>^\flat,-,q}})\ar[d]\ar[d]\ar[d]\ar[d] \ar[r]^{\mathrm{global}}\ar[r]\ar[r] &\varphi_q{\mathrm{Modules}_\circledcirc}_{\Gamma_k^{\times q}}(\{\mathrm{Robba}^{\mathrm{extended},q}_{{R_k,-,I}}\}_I) \ar[d]\ar[d]\ar[d]\ar[d].\\
{\mathrm{Modules}_\circledcirc}(\mathcal{O}_{X_{\mathbb{Q}_p(p^{1/p^\infty})^{\wedge\flat},-,q}})\ar[r]^{\mathrm{global}}\ar[r]\ar[r] &\varphi_q{\mathrm{Modules}_\circledcirc}(\{\mathrm{Robba}^{\mathrm{extended},q}_{{R_0,-,I}}\}_I).\\ 
}
\]
The definition is given by the following:
\[
\xymatrix@R+0pc@C+0pc{
\mathrm{homotopycolimit}_i(\square)(\mathcal{O}_i),  
}
\]
each $\mathcal{O}_i$ is just as $\mathbb{Q}_p\left<C_1,...,C_\ell\right>,\ell=1,2,...$ and $\square$ is the relative diagram of $\infty$-functors.

\item (\text{Proposition}) There is a functor (global section) between the $\infty$-prestacks of inductive Banach quasicoherent commutative algebra $E_\infty$ objects\footnote{Here $\circledcirc=\text{solidquasicoherentsheaves}$.}:
\[\displayindent=-0.4in
\xymatrix@R+6pc@C+0pc{
\mathrm{sComm}_\mathrm{simplicial}{\mathrm{Modules}_\circledcirc}_{\Gamma_k^{\times q}}(\mathcal{O}_{X_{R_k,-,q}})\ar[d]\ar[d]\ar[d]\ar[d]\ar[r]^{\mathrm{global}}\ar[r]\ar[r] &\mathrm{sComm}_\mathrm{simplicial}\varphi_q{\mathrm{Modules}_\circledcirc}_{\Gamma_k^{\times q}}(\{\mathrm{Robba}^{\mathrm{extended},q}_{{R_k,-,I}}\}_I)\ar[d]\ar[d]\ar[d]\ar[d]\\
\mathrm{sComm}_\mathrm{simplicial}{\mathrm{Modules}_\circledcirc}_{\Gamma_0^{\times q}}(\mathcal{O}_{X_{\mathbb{Q}_p(p^{1/p^\infty})^{\wedge\flat},-,q}})\ar[r]^{\mathrm{global}}\ar[r]\ar[r] &\mathrm{sComm}_\mathrm{simplicial}\varphi_q{\mathrm{Modules}_\circledcirc}_{\Gamma_0^{\times q}}(\{\mathrm{Robba}^{\mathrm{extended},q}_{{R_0,-,I}}\}_I).  
}
\]
The definition is given by the following:
\[
\xymatrix@R+0pc@C+0pc{
\mathrm{homotopycolimit}_i(\square)(\mathcal{O}_i),  
}
\]
each $\mathcal{O}_i$ is just as $\mathbb{Q}_p\left<C_1,...,C_\ell\right>,\ell=1,2,...$ and $\square$ is the relative diagram of $\infty$-functors.

\item Then as in \cite{LBV} we have a functor (global section) of the de Rham complex after \cite[Definition 5.9, Section 5.2.1]{KKM}\footnote{Here $\circledcirc=\text{solidquasicoherentsheaves}$.}:
\[\displayindent=-0.6in
\xymatrix@R+6pc@C+0pc{
\mathrm{deRham}_{\mathrm{sComm}_\mathrm{simplicial}{\mathrm{Modules}_\circledcirc}_{\Gamma_k^{\times q}}(\mathcal{O}_{X_{R_k,-,q}})\ar[r]^{\mathrm{global}}}(-)\ar[d]\ar[d]\ar[d]\ar[d]\ar[r]\ar[r] &\mathrm{deRham}_{\mathrm{sComm}_\mathrm{simplicial}\varphi_q{\mathrm{Modules}_\circledcirc}_{\Gamma_k^{\times q}}(\{\mathrm{Robba}^{\mathrm{extended},q}_{{R_k,-,I}}\}_I)}(-)\ar[d]\ar[d]\ar[d]\ar[d]\\
\mathrm{deRham}_{\mathrm{sComm}_\mathrm{simplicial}{\mathrm{Modules}_\circledcirc}_{\Gamma_0^{\times q}}(\mathcal{O}_{X_{\mathbb{Q}_p(p^{1/p^\infty})^{\wedge\flat},-,q}})\ar[r]^{\mathrm{global}}}(-)\ar[r]\ar[r] &\mathrm{deRham}_{\mathrm{sComm}_\mathrm{simplicial}\varphi_q{\mathrm{Modules}_\circledcirc}_{\Gamma_0^{\times q}}(\{\mathrm{Robba}^{\mathrm{extended},q}_{{R_0,-,I}}\}_I)}(-), 
}
\]

The definition is given by the following:
\[
\xymatrix@R+0pc@C+0pc{
\mathrm{homotopycolimit}_i(\square)(\mathcal{O}_i),  
}
\]
each $\mathcal{O}_i$ is just as $\mathbb{Q}_p\left<C_1,...,C_\ell\right>,\ell=1,2,...$ and $\square$ is the relative diagram of $\infty$-functors.

\item Then we have the following a functor (global section) of $K$-group $(\infty,1)$-spectrum from \cite{BGT}\footnote{Here $\circledcirc=\text{solidquasicoherentsheaves}$.}:
\[
\xymatrix@R+6pc@C+0pc{
\mathrm{K}^\mathrm{BGT}_{\mathrm{sComm}_\mathrm{simplicial}{\mathrm{Modules}_\circledcirc}_{\Gamma_k^{\times q}}(\mathcal{O}_{X_{R_k,-,q}})\ar[r]^{\mathrm{global}}}(-)\ar[d]\ar[d]\ar[d]\ar[d]\ar[r]\ar[r] &\mathrm{K}^\mathrm{BGT}_{\mathrm{sComm}_\mathrm{simplicial}\varphi_q{\mathrm{Modules}_\circledcirc}_{\Gamma_k^{\times q}}(\{\mathrm{Robba}^{\mathrm{extended},q}_{{R_k,-,I}}\}_I)}(-)\ar[d]\ar[d]\ar[d]\ar[d]\\
\mathrm{K}^\mathrm{BGT}_{\mathrm{sComm}_\mathrm{simplicial}{\mathrm{Modules}_\circledcirc}_{\Gamma_0^{\times q}}(\mathcal{O}_{X_{\mathbb{Q}_p(p^{1/p^\infty})^{\wedge\flat},-,q}})\ar[r]^{\mathrm{global}}}(-)\ar[r]\ar[r] &\mathrm{K}^\mathrm{BGT}_{\mathrm{sComm}_\mathrm{simplicial}\varphi_q{\mathrm{Modules}_\circledcirc}_{\Gamma_0^{\times q}}(\{\mathrm{Robba}^{\mathrm{extended},q}_{{R_0,-,I}}\}_I)}(-), 
}
\]

The definition is given by the following:
\[
\xymatrix@R+0pc@C+0pc{
\mathrm{homotopycolimit}_i(\square)(\mathcal{O}_i),  
}
\]
each $\mathcal{O}_i$ is just as $\mathbb{Q}_p\left<C_1,...,C_\ell\right>,\ell=1,2,...$ and $\square$ is the relative diagram of $\infty$-functors.\\

\end{itemize}

\

\begin{proposition}
All the global functors from \cite[Proposition 13.8, Theorem 14.9, Remark 14.10]{1CS2} achieve the equivalences on both sides.	
\end{proposition}

\newpage
\subsection{$\infty$-Categorical Analytic Stacks and Descents V}

Here we consider the corresponding archimedean picture, after \cite[Problem A.4, Kedlaya's Lecture]{1CBCKSW}. Recall for any algebraic variety $R$ over $\mathbb{R}$ this $X_R(\mathbb{C})$ is defined to be the corresponding quotient:
\begin{align}
X_{R}(\mathbb{C}):=R(\mathbb{C})\times P^1(\mathbb{C})/\varphi,\\
Y_R(\mathbb{C}):=R(\mathbb{C})\times P^1(\mathbb{C}).	
\end{align}
The Hodge structure is given by $\varphi$. We define the relative version by considering a further algebraic variety over $\mathbb{C}$, say $A$ as in the following:
\begin{align}
X_{R,A}(\mathbb{C}):=R(\mathbb{C})\times P^1(\mathbb{C})\times A(\mathbb{C})/\varphi,\\
Y_{R,A}(\mathbb{C}):=R(\mathbb{C})\times P^1(\mathbb{C})\times A(\mathbb{C}).	
\end{align}

\indent We then take $q$-th self power to achieve $X_{R,q}(\mathbb{C})$ as 
\begin{align}
X_{R,q}(\mathbb{C}):=(R(\mathbb{C})\times P^1(\mathbb{C}))^q/\varphi_q,\\
Y_{R,q}(\mathbb{C}):=(R(\mathbb{C})\times P^1(\mathbb{C}))^q.	
\end{align}
The multi hyperk\"ahler Hodge structure is given by $\varphi_q$. We define the relative version by considering a further algebraic variety over $\mathbb{C}$, say $A$ as in the following:
\begin{align}
X_{R,A}(\mathbb{C}):=(R(\mathbb{C})\times P^1(\mathbb{C}))^q\times A(\mathbb{C})/\varphi_q,\\
Y_{R,A}(\mathbb{C}):=(R(\mathbb{C})\times P^1(\mathbb{C}))^q\times A(\mathbb{C}).	
\end{align}

Then by \cite{1BBK} and \cite{1CS2} we have the corresponding $\infty$-category of $\infty$-sheaves of simplicial ind-Banach quasicoherent modules which in our situation will be assumed to the modules in \cite{1BBK}, as well as the corresponding associated Clausen-Scholze spaces:
\begin{align}
X_{R}(\mathbb{C}):=R(\mathbb{C})\times P^1(\mathbb{C})^\blacksquare/\varphi,\\
Y_R(\mathbb{C}):=R(\mathbb{C})\times P^1(\mathbb{C})^\blacksquare.	
\end{align}
\begin{align}
X_{R,A}(\mathbb{C}):=R(\mathbb{C})\times P^1(\mathbb{C})\times A(\mathbb{C})^\blacksquare/\varphi,\\
Y_{R,A}(\mathbb{C}):=R(\mathbb{C})\times P^1(\mathbb{C})\times A(\mathbb{C})^\blacksquare,	
\end{align}
with the $\infty$-category of $\infty$-sheaves of simplicial liquid quasicoherent modules, liquid vector bundles and liquid perfect complexes, with further descent \cite[Proposition 13.8, Theorem 14.9, Remark 14.10]{1CS2}.\\
 
\indent We then take $q$-th self power to achieve $X_{R,q}(\mathbb{C})$ as 
\begin{align}
X_{R,q}(\mathbb{C}):=(R(\mathbb{C})\times P^1(\mathbb{C}))^{q,\blacksquare}/\varphi_q,\\
Y_{R,q}(\mathbb{C}):=(R(\mathbb{C})\times P^1(\mathbb{C}))^{q,\blacksquare}.	
\end{align}
The multi hyperk\"ahler Hodge structure is given by $\varphi_q$. We define the relative version by considering a further algebraic variety over $\mathbb{C}$, say $A$ as in the following:
\begin{align}
X_{R,A}(\mathbb{C}):=(R(\mathbb{C})\times P^1(\mathbb{C}))^{q,\blacksquare}\times A(\mathbb{C})/\varphi_q,\\
Y_{R,A}(\mathbb{C}):=(R(\mathbb{C})\times P^1(\mathbb{C}))^{q,\blacksquare}\times A(\mathbb{C}).	
\end{align}

We call the resulting global sections are the corresponding $c$-equivariant Hodge Modules. Then we have the following direct analogy:

\begin{itemize}
\item (\text{Proposition}) There is an equivalence between the $\infty$-categories of inductive Banach quasicoherent presheaves:
\[
\xymatrix@R+0pc@C+0pc{
\mathrm{Ind}\mathrm{Banach}(\mathcal{O}_{X_{R,A,q}})\ar[r]^{\mathrm{equi}}\ar[r]\ar[r] &\varphi_q\mathrm{Ind}\mathrm{Banach}(\mathcal{O}_{Y_{R,A,q}}).  
}
\]
\item (\text{Proposition}) There is an equivalence between the $\infty$-categories of monomorphic inductive Banach quasicoherent presheaves:
\[
\xymatrix@R+0pc@C+0pc{
\mathrm{Ind}^m\mathrm{Banach}(\mathcal{O}_{X_{R,A,q}})\ar[r]^{\mathrm{equi}}\ar[r]\ar[r] &\varphi_q\mathrm{Ind}^m\mathrm{Banach}(\mathcal{O}_{Y_{R,A,q}}).  
}
\]
\end{itemize}

\begin{itemize}

\item (\text{Proposition}) There is an equivalence between the $\infty$-categories of inductive Banach quasicoherent presheaves:
\[
\xymatrix@R+0pc@C+0pc{
\mathrm{Ind}\mathrm{Banach}(\mathcal{O}_{X_{R,A,q}})\ar[r]^{\mathrm{equi}}\ar[r]\ar[r] &\varphi_q\mathrm{Ind}\mathrm{Banach}(\mathcal{O}_{Y_{R,A,q}}).  
}
\]
\item (\text{Proposition}) There is an equivalence between the $\infty$-categories of monomorphic inductive Banach quasicoherent presheaves:
\[
\xymatrix@R+0pc@C+0pc{
\mathrm{Ind}^m\mathrm{Banach}(\mathcal{O}_{X_{R,A,q}})\ar[r]^{\mathrm{equi}}\ar[r]\ar[r] &\varphi_q\mathrm{Ind}^m\mathrm{Banach}(\mathcal{O}_{Y_{R,A,q}}).  
}
\]
\item (\text{Proposition}) There is an equivalence between the $\infty$-categories of inductive Banach quasicoherent commutative algebra $E_\infty$ objects:
\[
\xymatrix@R+0pc@C+0pc{
\mathrm{sComm}_\mathrm{simplicial}\mathrm{Ind}\mathrm{Banach}(\mathcal{O}_{X_{R,A,q}})\ar[r]^{\mathrm{equi}}\ar[r]\ar[r] &\mathrm{sComm}_\mathrm{simplicial}\varphi_q\mathrm{Ind}\mathrm{Banach}(\mathcal{O}_{Y_{R,A,q}}).  
}
\]
\item (\text{Proposition}) There is an equivalence between the $\infty$-categories of monomorphic inductive Banach quasicoherent commutative algebra $E_\infty$ objects:
\[
\xymatrix@R+0pc@C+0pc{
\mathrm{sComm}_\mathrm{simplicial}\mathrm{Ind}^m\mathrm{Banach}(\mathcal{O}_{X_{R,A,q}})\ar[r]^{\mathrm{equi}}\ar[r]\ar[r] &\mathrm{sComm}_\mathrm{simplicial}\varphi_q\mathrm{Ind}^m\mathrm{Banach}(\mathcal{O}_{Y_{R,A,q}}).  
}
\]

\item Then parallel as in \cite{LBV} we have the equivalence of the de Rham complex after \cite[Definition 5.9, Section 5.2.1]{KKM}:
\[
\xymatrix@R+0pc@C+0pc{
\mathrm{deRham}_{\mathrm{sComm}_\mathrm{simplicial}\mathrm{Ind}\mathrm{Banach}(\mathcal{O}_{X_{R,A,q}})\ar[r]^{\mathrm{equi}}}(-)\ar[r]\ar[r] &\mathrm{deRham}_{\mathrm{sComm}_\mathrm{simplicial}\varphi_q\mathrm{Ind}\mathrm{Banach}(\mathcal{O}_{Y_{R,A,q}})}(-), 
}
\]
\[
\xymatrix@R+0pc@C+0pc{
\mathrm{deRham}_{\mathrm{sComm}_\mathrm{simplicial}\mathrm{Ind}^m\mathrm{Banach}(\mathcal{O}_{X_{R,A,q}})\ar[r]^{\mathrm{equi}}}(-)\ar[r]\ar[r] &\mathrm{deRham}_{\mathrm{sComm}_\mathrm{simplicial}\varphi_q\mathrm{Ind}^m\mathrm{Banach}(\mathcal{O}_{Y_{R,A,q}})}(-). 
}
\]

\item Then we have the following equivalence of $K$-group $(\infty,1)$-spectrum from \cite{BGT}:
\[
\xymatrix@R+0pc@C+0pc{
\mathrm{K}^\mathrm{BGT}_{\mathrm{sComm}_\mathrm{simplicial}\mathrm{Ind}\mathrm{Banach}(\mathcal{O}_{X_{R,A,q}})\ar[r]^{\mathrm{equi}}}(-)\ar[r]\ar[r] &\mathrm{K}^\mathrm{BGT}_{\mathrm{sComm}_\mathrm{simplicial}\varphi_q\mathrm{Ind}\mathrm{Banach}(\mathcal{O}_{Y_{R,A,q}})}(-), 
}
\]
\[
\xymatrix@R+0pc@C+0pc{
\mathrm{K}^\mathrm{BGT}_{\mathrm{sComm}_\mathrm{simplicial}\mathrm{Ind}^m\mathrm{Banach}(\mathcal{O}_{X_{R,A,q}})\ar[r]^{\mathrm{equi}}}(-)\ar[r]\ar[r] &\mathrm{K}^\mathrm{BGT}_{\mathrm{sComm}_\mathrm{simplicial}\varphi_q\mathrm{Ind}^m\mathrm{Banach}(\mathcal{O}_{Y_{R,A,q}})}(-). 
}
\]
\end{itemize}

\begin{assumption}\label{assumtionpresheaves}
All the functors of modules or algebras below are Clausen-Scholze sheaves \cite[Proposition 13.8, Theorem 14.9, Remark 14.10]{1CS2}. 	
\end{assumption}

\begin{itemize}
\item (\text{Proposition}) There is an equivalence between the $\infty$-categories of inductive liquid sheaves:
\[
\xymatrix@R+0pc@C+0pc{
\mathrm{Module}_\circledcirc(\mathcal{O}_{X_{R,A,q}})\ar[r]^{\mathrm{equi}}\ar[r]\ar[r] &\varphi_q\mathrm{Module}_\circledcirc(\mathcal{O}_{Y_{R,A,q}}).  
}
\]
\end{itemize}

\begin{itemize}

\item (\text{Proposition}) There is an equivalence between the $\infty$-categories of inductive Banach quasicoherent commutative algebra $E_\infty$ objects:
\[
\xymatrix@R+0pc@C+0pc{
\mathrm{sComm}_\mathrm{simplicial}\mathrm{Module}_{\text{liquidquasicoherentsheaves}}(\mathcal{O}_{X_{R,A,q}})\ar[r]^{\mathrm{equi}}\ar[r]\ar[r] &\mathrm{sComm}_\mathrm{simplicial}\varphi_q\mathrm{Module}_{\text{liquidquasicoherentsheaves}}(\mathcal{O}_{Y_{R,A,q}}).  
}
\]

\item Then as in \cite{LBV} we have the equivalence of the de Rham complex after \cite[Definition 5.9, Section 5.2.1]{KKM}\footnote{Here $\circledcirc=\text{liquidquasicoherentsheaves}$.}:
\[
\xymatrix@R+0pc@C+0pc{
\mathrm{deRham}_{\mathrm{sComm}_\mathrm{simplicial}\mathrm{Module}_\circledcirc(\mathcal{O}_{X_{R,A,q}})\ar[r]^{\mathrm{equi}}}(-)\ar[r]\ar[r] &\mathrm{deRham}_{\mathrm{sComm}_\mathrm{simplicial}\varphi_q\mathrm{Module}_\circledcirc(\mathcal{O}_{Y_{R,A,q}})}(-). 
}
\]

\item Then we have the following equivalence of $K$-group $(\infty,1)$-spectrum from \cite{BGT}\footnote{Here $\circledcirc=\text{liquidquasicoherentsheaves}$.}:
\[
\xymatrix@R+0pc@C+0pc{
\mathrm{K}^\mathrm{BGT}_{\mathrm{sComm}_\mathrm{simplicial}\mathrm{Module}_\circledcirc(\mathcal{O}_{X_{R,A,q}})\ar[r]^{\mathrm{equi}}}(-)\ar[r]\ar[r] &\mathrm{K}^\mathrm{BGT}_{\mathrm{sComm}_\mathrm{simplicial}\varphi_q\mathrm{Module}_\circledcirc(\mathcal{O}_{Y_{R,A,q}})}(-). 
}
\]
\end{itemize}

\newpage

\subsection*{Acknowledgements}
The discussion we made in this paper is based on the author's talk at UCSD number theory seminar. And the presentation and the choice of the topics here were completed and finished under the suggestions from Professor Kedlaya. We would like to thank Professor Kedlaya for helpful discussion before giving the talk in 2020.

\bibliographystyle{splncs}

\end{document}